\documentclass[11pt,reqno]{amsart}
\usepackage[left=1in,top=1in,right=1in,bottom=1in]{geometry}
\usepackage{amsfonts}
\usepackage{enumerate,subfigure}%,xcolor}
\usepackage{amssymb,latexsym, amsmath,graphicx}%,epsf}
\usepackage{epstopdf}
\usepackage{framed}
\usepackage{hyperref}
\usepackage{cite}
\usepackage{booktabs}
\usepackage[dvipsnames]{xcolor}

\numberwithin{equation}{section}
\usepackage[pagewise]{lineno}%\linenumbers

\usepackage{multirow}%,booktabs}
\DeclareMathOperator{\R}{{\mathbb R}}
\DeclareMathOperator{\C}{{\mathbb C}}
\DeclareMathOperator{\by}{\times}
\DeclareMathOperator{\bndry}{\partial\Omega}
\DeclareMathOperator{\texp}{\mathbf{t}^{\mbox{\tiny{\textbf{exp}}}}}
\DeclareMathOperator{\tzero}{\mathbf{t}^{0}}
\DeclareMathOperator{\qexp}{\mathbf{q}^{\mbox{\tiny{\textbf{exp}}}}}
\DeclareMathOperator{\qzero}{\mathbf{q}^0}
\DeclareMathOperator{\sigcal}{\sigma^{\mbox{\tiny{\textbf{CAL}}}}}

\DeclareMathOperator{\sigexp}{\sigma^{\mbox{\tiny{\textbf{exp}}}}}
\DeclareMathOperator{\sigmaexp}{\sigma^{\mbox{\tiny{\textbf{exp}}}}}

\DeclareMathOperator{\Fhat}{\hat{F}}

\DeclareMathOperator{\sigbest}{\sigma_{\mbox{\tiny{best}}}}

\makeatletter
\def\thickhline{%
  \noalign{\ifnum0=`}\fi\hrule \@height \thickarrayrulewidth \futurelet
   \reserved@a\@xthickhline}
\def\@xthickhline{\ifx\reserved@a\thickhline
               \vskip\doublerulesep
               \vskip-\thickarrayrulewidth
             \fi
      \ifnum0=`{\fi}}
\makeatother

\newlength{\thickarrayrulewidth}
\usepackage{boldline}
\definecolor{brown}{rgb}{0.7,0.3,0}

\definecolor{darkgreen}{rgb}{0,0.6,0}
\newcommand{\trev}[1]{{#1}}
\begin{document}
%--------------------------------------------------------------------
%\bibliographystyle{plainnat} % for Phys Meas.

\title[3D CGO-Based Experimental EIT]{Fast Absolute 3D CGO-Based Electrical Impedance Tomography on Experimental Tank Data}

\author[Hamilton et al.]{S.~J. Hamilton, P.~A. Muller, D. Isaacson, ~V. Kolehmainen, ~J. Newell, ~O. Rajabi~Shishvan, ~G. Saulnier, and ~J. Toivanen}

%\date{April 22, 2022}

\thanks{S.~J. Hamilton is with the Department of Mathematical and Statistical Sciences; Marquette University, Milwaukee, WI 53233 USA,  email: \texttt{sarah.hamilton@marquette.edu}}
\thanks{D. Isaacson is with the Department of Mathematics, Rensselaer Polytechnic Institute, Troy, NY 12180, USA}%, email: \texttt{isaacd@rpi.edu}}
\thanks{V. Kolehmainen and J. Toivanen are with the Department of Applied Physics, University of Eastern Finland, FI-70210 Kuopio, Finland}%, email: \texttt{ville.kolehmainen@uef.fi}}
\thanks{P.~A. Muller is with the Department of Mathematics \& Statistics; Villanova University, Villanova, PA 19085 USA}%,  email: \texttt{peter.muller@villanova.edu}}
\thanks{J. Newell is with the Department of Biomedical Engineering, Rensselaer Polytechnic Institute, Troy, NY 12180, USA}
\thanks{O. Rajabi~Shishvan and G. Saulnier are with the Department of Electrical and Computer Engineering, University at Albany - SUNY, Albany, NY 12222, USA}
\thanks{J. Toivanen is with the Department of Applied Physics, University of Eastern Finland, FI-70210 Kuopio, Finland}

%--------------------------------------------------------------------
\begin{abstract}
 {\it Objective:} To present the first 3D CGO-based absolute EIT reconstructions from experimental tank data. {\it Approach:} CGO-based methods for absolute EIT imaging are compared to traditional TV regularized non-linear least squares reconstruction methods.  Additional robustness testing is performed by considering incorrect model\textbf{}ing of domain shape.  {\it Main Results:} The CGO-based methods are fast, and show strong robustness to incorrect domain modeling comparable to classic difference EIT imaging and fewer boundary artefacts than the TV regularized non-linear least squares reference reconstructions. {\it Significance:}  This work is the first to demonstrate fully 3D CGO-based absolute EIT reconstruction on experimental data and also compares to TV-regularized absolute reconstruction. The speed (1-5 seconds) and quality of the reconstructions is encouraging for future work in absolute EIT.
\end{abstract}
%--------------------------------------------------------------------
%\pacs{PUT NUMBERS IN}
%\submitto{\PM}
%--------------------------------------------------------------------
\maketitle % uncomment if separate title page is required

%\tableofcontents

%--------------------------------------------------------------------
\section{Introduction}\label{sec:Intro}
%--------------------------------------------------------------------
The main objective of this paper is to demonstrate the feasibility, speed, and robustness of \trev{producing} 3D absolute (static) images of the electrical conductivity inside a tank, from experimental \trev{Electrical Impedance Tomography (EIT)} data measured on an array of electrodes on the tank’s surface by the ACT5~\cite{rajabishishvan2021} adaptive current tomography system, using \trev{CGO-based} reconstruction algorithms.  \trev{The primary contributions of this work are that it presents the first use of Complex Geometrical Optics (CGO) based methods to produce absolute (static) images of the conductivity from experimentally measured EIT data in 3D, and studies their robustness under modeling errors.}  

\trev{``Absolute", or ``static,'' conductivity images are images of the conductivity inside a body made from \emph{one} set of experimental measurements made on the surface of the body at \emph{one} time \cite{Brown_absoluteEIT}. Alternatively, ``dynamic," or ``time-difference," imaging uses \emph{two} sets of data measured at \emph{two} different times to make an image of the {\em change} in the conductivity that took place between the two times the measurements were made. We present both types of images made from experimental data by CGO-based algorithms in 3D in this paper.}

\trev{Applications of EIT began with geophysical exploration over a century ago}~\cite{Allaud1977,Keller1966} and have since expanded into several fields including the transport of fluids and gases \cite{Eggleston,Kaipio2004a,Wang2015}, and biomedical imaging \cite{Swanson1976,Henderson1978,Barber1984,Holder2005,Kao2020,Adler2021,Frerichs2017,deCastroMartins2019}.  Systems for monitoring lung function in real-time are now commercially available and clinical trials are in progress to determine the extent to which some of these systems might be used to guide mechanical ventilation \trev{\cite{Timpel,RecruitTrial,PlugGroup}}. These systems typically make and display two dimensional images of the differences in conductivity between two states, \trev{e.g.}, lungs filled with air and lungs depleted of air, in order to take advantage of \trev{dynamic (difference)} imaging methods and algorithms~\trev{\cite{Adler2021,Frerichs2017}}.  The desire to improve the diagnosis of cancer and stroke has motivated the development of systems and methods capable of imaging the absolute or static internal conductivity and permittivity \trev{in 3D} \cite{Cherepenin2001,Kerner2002,Choi2007,goren2018,Kao2007a}.  \trev{The interested reader is referred to \cite{Cheney1999} for further review of applications of EIT.}

\trev{Most absolute (static) EIT reconstruction focuses on solving a simplified linearized problem or iteratively solving an optimization-based method which requires repeated solutions to the forward problem which can become costly for highly dense meshes.  CGO-based methods are {\em direct methods} in that they do not require iteration. They have the capability to solve the full nonlinear mathematical inverse problem, and do not require repeated solutions to the forward problem, e.g. via the finite element method (FEM).  The D-bar inversion algorithm, which is a specific type of CGO-based inversion algorithm in 2D, has been used to make both absolute and difference images in 2D in real-time from experimental data measured on tanks, and for difference imaging on human subjects in laboratory and clinical settings.  See \cite{Mueller2020} for a recent review of the 2D D-bar method and its applications, and \cite{Nachman1996} for the theoretical foundation.}  

\trev{In 3D, existence and uniqueness of solutions \cite{Nachman1988a,Novikov1988} can be shown for a 3D D-bar type equation but the constructive proof, upon which the reconstruction algorithms are built, bypasses the 3D D-bar equation instead using a high (non-physical) frequency limit connecting the nonlinear scattering data and the linear Fourier data.  Advances in the numerical implementation of 3D CGO-based methods are more recent. The first numerical implementation of 3D CGO-based methods on simulated electrode data using current injection on the surface of the sphere was presented in \cite{Hamiltonetal_2021}.  A 3D CGO-based inversion algorithm, regularization scheme, rigorous proof of stability under certain hypotheses, and examples reconstructing the conductivity inside a sphere from numerically simulated Dirichlet data on the entire surface of the sphere (without electrodes) were given in \cite{Knudsen2022}. Until now it has been an open question if 3D CGO-based algorithms could be used with experimental data. This paper answers this question affirmatively by showing that the conductivity can be recovered rapidly and robustly from experimental data measured on 32~electrodes on the surface of a rectangular prism building on the work of \cite{Hamiltonetal_2021}.}
 
It was demonstrated in \cite{Goble1990,Blue1997,Holder2005,Kao2020} that it is possible to make electrical impedance tomography (EIT) Systems that produce both 3D static and dynamic images of the interior of the chest showing heart and lung regions, as well as changes in those regions due to ventilation and perfusion. \trev{These systems used linearized and iterative methods. The former are fast but less accurate and the latter are slower. Recent 3D CGO inversion algorithms and their analysis, when applied to synthetic or simulated data, suggested that they have the potential to be fast and more accurate  \cite{Hamiltonetal_2021, Delbary2014,Knudsen2022}.  Here we test their capabilities on experimental EIT data.} 

In this paper we compare the following methods for making static 3D images from experimental tank data, collected by the ACT5 imaging system:
\vspace{-0.5em}
\begin{itemize}
    \item A CGO linearization method based on Calder\'on’s original proposal introducing CGO ideas into the subject \cite{Calderon1980}. We call this algorithm Calder\'on’s method.
    \item Two CGO\trev{-based} methods for solving the full non-linear inversion problem based on Sylvester and Uhlmann’s original uniqueness proof and Nachman’s constructive uniqueness proof in \cite{Sylvester1987,Nachman1988a}.  \trev{The first is the $\texp$ algorithm, and the second is called $\tzero$, \cite{Delbary2012,Hamiltonetal_2021} both of which are simplifications of the original constructive proof.}
    \item A more traditional iterative inversion method based on optimization using a Total Variation (TV) regularization.%\tblue{[?? references include Wang from refs report??]}\tpeter{I can't find a reference to Wang in the ref reports.}
\end{itemize}

%\vspace{-0.5em}
The EIT data were measured on 32 electrodes attached to the six surfaces of a rectangular tank using the ACT5 system. The system \trev{can adaptively determine patterns of currents to apply to all 32 electrodes} that result in voltages on those electrodes that are proportional to the applied currents. These ``eigen currents" form a discrete orthogonal set and provide improved voltage data for reconstructing the conductivity inside the tank from a limited amount of current or power that can be applied to the tank. The theory of adaptive current tomography systems is given in \cite{Isaacson1986,Gisser1990,Gisser1988,Cheney1992,Li2013,Li2014,Newell2005,Saulnier2005}.

The remainder of the paper is organized as follows.  Section~\ref{sec:Methods} \trev{provides a brief review of the mathematics of EIT and} encompasses the methods used in the work: the CGO and reference reconstruction algorithms, experimental setup for the ACT5 data collection, robustness tests that will be explored, and metrics used to evaluate the results.  The results section, Section~\ref{sec:results}, presents slice, 3D, and isosurface renderings of the recovered conductivities for correct and incorrect domain modeling.  Section~\ref{sec:discussion} contains a discussion of the results and conclusions are drawn in Section~\ref{sec:conclusion}.

%--------------------------------------------------------------------
\section{Methods}\label{sec:Methods}
%--------------------------------------------------------------------
In this work we compare three CGO-based methods (Calder\'on's Method, the $\texp$ Method, and the $\tzero$ Method) to a more common iterative Total Variation (TV) regularized non-linear least squares method for absolute EIT image reconstruction.  For comparison, we also include time difference EIT images for the CGO-based methods and compare them to a linearized \trev{difference imaging} method.  \trev{We begin with a brief review of the mathematical problem.}
%--------------------------------------------------------------------
\subsection{\trev{Mathematical Background}}\label{sec-background}
%--------------------------------------------------------------------

The mathematical problem of reconstructing the internal conductivity, when measurements can be made with infinite precision everywhere on the boundary of a body, is currently called ``Calder\'on’s Problem” by much of the mathematical community since A. Calder\'on formulated this inverse problem as follows \cite{Calderon1980}: In two or more dimensions, can one find the conductivity, $\sigma(x)$, inside a body, $\Omega$, from all possible electrical measurements made on the surface, $\bndry$, of the body? Here the voltage or potential,  $u(x)$, inside the body, due to an applied surface current density,  $j(x)$, is assumed to satisfy the following low frequency approximation to Maxwell’s Equations, which we will refer to as the {\em conductivity equation}, with a Neumann boundary condition:
\begin{eqnarray}
    \nabla\cdot\sigma(x)\nabla u(x) &=& 0, \qquad x\in\Omega\subset \R^3 \label{eq:cond}\\
\sigma(x)\frac{\partial u(x)}{\partial \nu(x)} &=& j(x), \qquad x\in\bndry.\nonumber
\end{eqnarray}
Here $\nu =\nu(x)$ denotes the outward pointing unit normal to the surface at $x\in\bndry$, and, $\frac{\partial u(x)}{\partial \nu(x)}=\nu(x)\cdot\nabla u(x)$.  We denote the mapping from applied current density to resulting voltage on the surface, called the ``Neumann-to-Dirchlet'' (ND) map, by $\mathcal{R}_\sigma$, where $\mathcal{R}_\sigma j\equiv u(x)$ for $x\in\bndry$, and $u(x)$ solves the conductivity equation~\eqref{eq:cond} with Neumann data $j(x)$.

The Calder\'on problem can also be formulated as the mathematical problem.  If one is given the ND map or operator, $\mathcal{R}_\sigma$, or equivalently its inverse, the Dirichlet-to-Neumann (DN) operator,  $\Lambda_\sigma:=\mathcal{R}_\sigma^{-1}$,  can one find $\sigma(x)$? Calder\'on showed that if $\sigma(x)$ does not differ too much from a constant then one can recover an approximation to it from the ND map. His short paper showed this could be done \trev{using} Fourier Transforms in a very clever way, by introducing special solutions, $e^{\zeta\cdot x}$, to the conductivity equation with constant conductivity, i.e. the Laplace equation, where $\zeta$ is a complex valued vector with $\zeta\cdot\zeta=0$\trev{.} This is possible in two or more dimensions \trev{when the conductivity is close to a constant} and gave birth to Calder\'on's method as described in this paper, as well as more powerful methods for solving the full non-linear problem by generalizing Calder\'on’s solutions to what we now call Complex Geometrical Optics (CGO) solutions, introduced by \cite{Sylvester1986,Sylvester1987} in their landmark paper proving that the Calder\'on problem in three or more dimensions has a unique solution\trev{, where the inversion problem is reduced to a Fourier transform in the limit $|\zeta|\to\infty$. The $\texp$ and $\tzero$ methods described below in \S~\ref{sec:texpG0}  follow this strategy, as opposed to the 2D D-bar methods described in \cite{Mueller2020}. This strategy is not possible in 2D where the CGO solutions are found by solving a first order linear PDE involving D-bar operators in the complex vector $\zeta$ and taking it to 0.} Constructive proofs using these CGO solutions and D-bar ideas from inverse scattering theory, along with applications to scattering and acoustics were given in \cite{Nachman1988a, Novikov1988}. 

CGO methods were used to prove uniqueness in the more difficult case of 2D where constructive methods for reconstructing the conductivity were given in detail in \cite{Nachman1996} \trev{for $\sigma\in\mathcal{C}^2(\Omega)$} and later in \cite{Astala2006a} \trev{for $\sigma\in{L}^\infty(\Omega)$}. Other pioneering work proving uniqueness under a variety of hypotheses on the conductivity include \cite{Langer1933,Kohn1984,Novikov1988,Parker1984,Tikhonov1949,Druskin1982}. Extensive references to more recent progress in the analysis and numerical analysis of the Calder\'on problem can be found in \cite{Borcea2002a,Delbary2014,Mueller2012,Uhlmann2014,Feldman2019}; \trev{\cite{Mueller2020,Knudsen2022}}.

In what follows we will be interested in the problem of reconstructing an approximation to the internal conductivity from finitely many experimental measurements made with finite precision. Unfortunately, this is an ill-posed problem and, unlike the purely mathematical problem, it does not have a unique solution. Nevertheless, it is sometimes possible to reconstruct useful approximations to the internal conductivity with a finite number of degrees of freedom, or voxels, which we will illustrate by making images from experimental data and comparing them to the actual interior conductivity within the tanks.

The EIT problem for a body with $L$ electrodes, $e_\ell$, $\ell=1,2,\ldots,L$, on its surface, is to find an approximation to the internal conductivity from all the possible electrical measurements made on these $L$ electrodes. In particular we will assume that we apply $L-1$ linearly independent patterns of currents, $\vec{I}^{(k)}$, $k=1,2,\ldots, L-1$, to the $L$ electrodes, and measure the resulting $L-1$ voltage patterns, $\vec{V}^{(k)}$, where $\vec{I}^{(k)}_\ell$ and $\vec{V}^{(k)}_\ell$ denote the applied current, and measured voltage from the $k^{\mbox{\tiny th}}$ pattern, on the $\ell^{\mbox{\tiny th}}$ electrode, for $\ell = 1,2,\ldots,L$.  From conservation of charge, and our choice of ground, we assume
\[\sum_{\ell=1}^L \vec{I}^{(k)}_\ell = \sum_{\ell=1}^L \vec{V}^{(k)}_\ell = 0.\]
The $k^{\mbox{\tiny th}}$ voltage or potential, $u^{(k)}(x)$, resulting inside the body is determined by the conductivity equation, $\nabla\cdot\sigma \nabla u^{(k)}=0$, and the {\it ``Complete Electrode Model''}, \cite{Cheng1989,Somersalo1992} where:
\begin{align}
\int_{e_\ell} \sigma(x)\frac{\partial u^{(k)}(x)}{\partial \nu(x)} dS(x) &=I^{(k)}_\ell, \qquad x\in e_\ell\nonumber\\
\sigma(x)\frac{\partial u^{(k)}(x)}{\partial \nu(x)}  & = 0, \qquad x\notin \bigcup_{\ell=1}^L e_\ell\label{eq:CEM}\\
u^{(k)}(x) + z_\ell \sigma(x)\frac{\partial u^{(k)}(x)}{\partial \nu(x)} & = V^{(k)}_\ell, \qquad x\in e_\ell.\nonumber
\end{align}
Here $z_\ell$ is the effective contact, or surface, impedance on the $\ell^{\mbox{\tiny th}}$ electrode.  The current patterns used will be eigenvectors of the \trev{Current to Voltage map, which is a matrix approximation to the ND map, for the homogeneous saline tank. They are found numerically by simulating a homogeneous tank for static/absolute imaging, and experimentally by adaptive methods for difference/dynamic imaging. The matrix approximations to the ND maps from the conductivity distribution $\sigma$ are denoted by the $L\by L$ matrices}, $\mathbf{R}_\sigma$, where $\mathbf{R}_\sigma\vec{I}^{(k)}=\vec{V}^{(k)}$, and we define $\mathbf{R}_\sigma\vec{1}=\vec{0}$.  Here the vectors $\vec{1}$, $\vec{0}$, denote vectors all of whose components are $1$, or $0$, respectively. The discrete analog of the DN map used in the CGO-based methods is given by $\mathbf{L}_\sigma:=\mathbf{R}_\sigma^{-1}$.

%--------------------------------------------------------------------
\subsection{Calder\'on's Method}\label{sec-cald-orig}
%--------------------------------------------------------------------
Following Calder\'on's original paper \cite{Calderon1980}, Calder\'on's method approximates the conductivity, $\sigma(x)$, from its Fourier transform.  Here we present a brief description of the method and refer the reader to \cite{Calderon1980, Bikowski2008,Muller2017,Hamiltonetal_2021} for \trev{further} details. Calder\'on's method assumes the conductivity is a small perturbation, $\delta\sigma(x)$, from a constant background, $\sigma_b$, i.e. $\sigma(x)=\sigma_b+\delta\sigma(x)$. In this paper, we assume that the background conductivity $\sigma_b=1$.  If the background constant is not one, then the problem can be scaled and unscaled as in \cite{Isaacson2004, Hamiltonetal_2021}.

The three steps of Calder\'on's method in 3-D, as described in \cite{Hamiltonetal_2021}, are:

{\bf Step 1:} Use the DN maps $\Lambda_\sigma$ and $\Lambda_1$ to approximate the Fourier transform of the small perturbation in conductivity, $\widehat{\delta\sigma}(z)$, by
\begin{equation}\label{eq:CaldFhat}
\widehat{\delta\sigma}(z)\approx\trev{\Fhat}(z):=-\frac{1}{2\pi^2|z|^2}\int_{\partial\Omega}e^{\pi i(z\cdot x)+\pi(a\cdot x)}\left(\Lambda_\sigma-\Lambda_1\right) e^{\pi i(z\cdot x)-\pi(a\cdot x)}dS(x),
\end{equation}
where $z$ and $a$ satisfy 
\begin{equation}\label{eq:za_conditions}
    z,a\in\R^3, |z|=|a|,\text{ and }z\cdot a=0,
\end{equation}
and $\Lambda_1$ is the DN map for a constant conductivity of $1$.

{\bf Step 2:} Take the inverse Fourier transform of $\trev{\Fhat}(z)$: 
\begin{equation}\label{eq:cal-FhattoGammaCart}
\delta\sigma^{\mbox{\tiny{\textbf{CAL}}}}(x) \approx\mathcal{F}^{-1}\{\trev{\Fhat}(z)\}(x)=\int_{\mathbb{R}^3}{\trev{\Fhat}(z)e^{-2\pi i(x\cdot z)}dz}.
\end{equation}

{\bf Step 3:} Add the background to the perturbation to recover the approximate conductivity, $\sigma^{\mbox{\tiny{\textbf{CAL}}}}(x)$:
\begin{equation}\label{eq:sigCal_pert}
\sigma^{\mbox{\tiny{\textbf{CAL}}}}(x)=\sigma_b+\delta\sigma^{\mbox{\tiny{\textbf{CAL}}}}(x).
\end{equation}

\normalsize
\vspace{1em}

The definition of the Fourier transform in Calder\'on's method is different from that used in the $\mathbf{t}^{\mbox{\tiny{\textbf{exp}}}}$ and $\tzero$ methods described below.  However, each method is consistent with its definition and is consistent with the respective literature on that method.

In this paper, we compute $\Fhat(z)$ in spherical Fourier coordinates. As such, we choose \[z=|z|(\cos\tilde\phi\sin\tilde\theta,\sin\tilde\phi\sin\tilde\theta,\cos\tilde\theta)\text{ and } a=|z|(\cos\tilde\phi\cos\tilde\theta,\sin\tilde\phi\cos\tilde\theta,-\sin\tilde\theta),\]
for $|z|\geq0$, $0\leq\tilde\phi\leq2\pi$ and $0\leq\tilde\theta\leq\pi$ so that $z$ and $a$ satisfy \eqref{eq:za_conditions}. Then, the inverse Fourier transform in \eqref{eq:cal-FhattoGammaCart} becomes
\begin{equation}\label{eq:cal-FhattoGamma}
\delta\sigma^{\mbox{\tiny{\textbf{CAL}}}}(x) =\int_{0}^{\infty}{\int_0^{2\pi}{\int_0^\pi{|z|^2\sin\tilde\theta\trev{\Fhat}(|z|,\tilde\phi,\tilde\theta)e^{-2\pi i(x\cdot z)}d\tilde\theta d\tilde\phi d|z|}}}.
\end{equation}

Additionally, we implement the use of a mollifier, $\hat{\eta}\left(\frac{z}{y}\right)$, as introduced in \cite{Calderon1980} for some \trev{parameter} $y\in\mathbb{R}$ to reduce Gibbs phenomenon \trev{caused by jump discontinuities in $\sigma(x)$ while} recovering $\delta\sigma^{\mbox{\tiny{\textbf{CAL}}}}$  
\begin{equation}\label{eq-sigcal-mollified}
\delta\sigma^{\mbox{\tiny{\textbf{CAL}}}}(x) =\int_{0}^{\infty}{\int_0^{2\pi}{\int_0^\pi{|z|^2\sin\tilde\theta\trev{\Fhat}(|z|,\tilde\phi,\tilde\theta)\hat{\eta}\left(\frac{z}{y}\right)e^{-2\pi i(x\cdot z)}d\tilde\theta d\tilde\phi d|z|}}}.
\end{equation}
We implement the same mollifier as used in \cite{Bikowski2008, Hamiltonetal_2021},
\begin{equation}\label{eq-Cald-moll}
\hat{\eta}\left(\frac{z}{y}\right)=e^{-\pi t|z|^2},
\end{equation}
where $y=1/\sqrt{t}$ and $t$ acts as a smoothing parameter with larger $t$ values producing smoother reconstructions with smaller jumps at points of discontinuity in $\sigma(x)$.

Since noise causes \eqref{eq:CaldFhat} to blow up at large $|z|$, we use a non-uniform truncation regularization strategy. A similar regularization strategy was proved stable for the 2-D D-bar method in \cite{Knudsen2009}, which also noted non-uniform truncation also produces reliable reconstructions.  In our case, we will first compute $\Fhat(z)$ for $|z|$ within an outer radius of $T_{z_2}$. We keep values of $\Fhat(z)$ whose real and imaginary amplitudes are below a threshold determined by the amplitudes of $\Fhat$ within a smaller radius $|z|\leq T_{z_1}$; $\Fhat$ is set to 0 everywhere else. Both radii are chosen empirically. The inner radius $T_{z_1}$ is chosen as a region in $z$-space where noise does not cause $\Fhat$ to blow up and $T_{z_2}$ is chosen to keep as much reasonable information from $\Fhat$ without introducing holes in the non-zero region of $\Fhat$. As such, our truncated $\Fhat$ is computed by

\begin{equation}\label{eq:truncFhat}
\small \Fhat_R(z)=\left\{\begin{array}{ll}
\Fhat(z),& \text{if }|z|\leq T_{z_2} \text{, and }
|\text{Re}(\Fhat(z))|\leq \underset{|z|\leq T_{z_1}}{\max}|\text{Re}(\Fhat(z))|,\,\, |\text{Im}(\Fhat(z))|\leq \underset{|z|\leq T_{z_1}}{\max}|\text{Im}(\Fhat(z))|\\
0,& \mathrm{else}.
\end{array}\right.
\end{equation}\normalsize

With our truncated $\Fhat$, equation \eqref{eq-sigcal-mollified} is truncated in the radial variable, leading to the approximation
\begin{equation}\label{eq-cal-FhattoGamma-reg}
\delta\sigma_R^{\mbox{\tiny{\textbf{CAL}}}}(x) =\int_{0}^{T_{z_2}}{\int_0^{2\pi}{\int_0^\pi{|z|^2\sin\tilde\theta\Fhat_R(|z|,\tilde\phi,\tilde\theta)e^{-\pi t|z|^2}e^{-2\pi i(x\cdot z)}d\tilde\theta d\tilde\phi d|z|}}}.
\end{equation}

For difference images shown in this paper, we only perform Steps 1 and 2 of the method and replace $\Lambda_1$ in \eqref{eq:CaldFhat} with a reference DN map, $\Lambda_{\sigma_{\mbox{\tiny ref}}}$ before computing \eqref{eq-cal-FhattoGamma-reg}. Thus, the flow for absolute images is 
\[\left(\Lambda_\sigma,\Lambda_1\right)\overset{1}{\longrightarrow}\trev{\Fhat}(z) \overset{2}{\longrightarrow}\delta\sigma^{\mbox{\tiny{\textbf{CAL}}}}\overset{3}{\longrightarrow}\sigma^{\mbox{\tiny{\textbf{CAL}}}}\]
and the flow for difference images in this paper is
\[\left(\Lambda_\sigma,\Lambda_{\sigma_{\mbox{\tiny ref}}}\right)\overset{1}{\longrightarrow}\trev{\Fhat}(z) \overset{2}{\longrightarrow}\delta\sigma^{\mbox{\tiny{\textbf{CAL}}}}.\]

%--------------------------------------------------------------------
\subsubsection{Numerical Implementation}\label{sec-cal-Fhat-implement}
%--------------------------------------------------------------------

In this section, we review the implementation details of Calder\'on's method introduced in \cite{Hamiltonetal_2021}.

For Step~1, we compute \eqref{eq:CaldFhat} for $|z|\leq T_{z_2}$ by discretizing the boundary integral as follows,
\begin{eqnarray}
\Fhat(z)&=&-\frac{1}{2\pi^2|z|^2}\int_{\partial\Omega}e^{\pi i(z\cdot x)+\pi(a\cdot x)}\left(\Lambda_\sigma-\Lambda_1\right) e^{\pi i(z\cdot x)-\pi(a\cdot x)}dS(x)\nonumber\\
&\approx& -\frac{1}{2\pi^2|z|^2}\left(\frac{|\partial\Omega|}{L}\right)(e^{\pi i(z\cdot \mathbf{x})+\pi(a\cdot \mathbf{x})})^{\mbox{\tiny \bf T}}\mathbf{Q}\left(\mathbf{L}_\sigma-\mathbf{L}_1\right)\mathbf{Q}^{\mbox{\tiny \bf T}}\left[e^{\pi i(z\cdot \mathbf{x})-\pi(a\cdot \mathbf{x})}\right],\label{eq-discrete-Fhat_bdry}
\end{eqnarray}
where $|\partial\Omega|$ is the surface area of the domain; $L$ is the number of  electrodes; $\mathbf{x}\in\R^{(L\by3)}$ denotes the vector of the Cartesian centers of the $L$ electrodes; $^{\mbox{\tiny \bf T}}$ denotes the traditional, non-conjugate, matrix transpose; $\mathbf{L}_\sigma$ and $\mathbf{L}_1$ denote the discrete matrix approximations to the DN maps $\Lambda_\sigma$ and $\Lambda_1$ respectively; and $\mathbf{Q}\in\R^{L\by \mbox{\tiny $N_{li}$}}$ denotes an orthonormal basis created using  $N_{li}$ linearly independent applied currents over $L$ electrodes as was done in \cite{Hamiltonetal_2021}. The matrix $\mathbf{L}_1$ is based on the FEM solution of the CEM \eqref{eq:CEM}.  Problem-specific mesh details are given below in section~\ref{sec:robustness}.

We then compute $\Fhat_R$ according to equations \eqref{eq:truncFhat} and \eqref{eq-discrete-Fhat_bdry}. For Step 2, on an equally-spaced $16\times16\times16$ rectangular grid in $x$, the conductivity difference $\delta\sigma_R^{\mbox{\tiny{\textbf{CAL}}}}(x)$ is computed via \eqref{eq-cal-FhattoGamma-reg} using a 3D Simpson's rule using $N_{|z|}=10, N_{\tilde{\theta}}=10$, and $N_{\tilde{\phi}}=30$ uniformly-spaced nodes on the $|z|, \tilde{\theta}$, and $\tilde{\phi}$ grids, respectively. As the number of nodes in the Fourier domain increases, so does computation time, but some artefact reduction can be achieved. Empirically, the artefact reduction did not seem significant enough to warrant an increased computational time beyond these parameter choices.

Difference images are computed using \eqref{eq-cal-FhattoGamma-reg} and \eqref{eq-discrete-Fhat_bdry}, replacing $\mathbf{L}_1$ with a discrete reference DN map, $\mathbf{L}_{\sigma_{\mbox{\tiny ref}}}$ in \eqref{eq-discrete-Fhat_bdry}. For the phantom tank experiments in this paper, this reference map is from data collected with a tank filled only with saline \trev{matching the experiment} and no other inclusions.

The absolute reconstructions of $\sigma_R^{\mbox{\tiny{\textbf{CAL}}}}(x)$ in this paper are produced using \eqref{eq:sigCal_pert} replacing $\sigma_b$ with $\sigma_{\mbox{\tiny{best}}}$
\begin{equation*}
\sigma_R^{\mbox{\tiny{\textbf{CAL}}}}(x)=\sigma_{\mbox{\tiny{best}}}+\delta\sigma_R^{\mbox{\tiny{\textbf{CAL}}}}(x).
\end{equation*}
The solution is then interpolated to a $64\times64\times64$ rectangular grid. 

Following \cite{Isaacson2004}, $\sigma_{\mbox{\tiny{best}}}$ is given by
\begin{equation}\label{eq-gam-best}
\sigma_{\mbox{\tiny{best}}}=\frac{\sum_{k=1}^K\sum_{\ell=1}^LU_\ell^kU_\ell^k}{\sum_{k=1}^K\sum_{\ell=1}^LU_\ell^kV_\ell^k},
\end{equation}
where $U_\ell^k$ is the $k^{th}$ {\it simulated} voltage pattern measured on electrode $\ell$ with a homogeneous conductivity of $1$ S/m and $V_\ell^k$ is the $k^{th}$ voltage pattern {\it measured} on electrode $\ell$ for the inhomogeneous conductivity $\sigma$. The simulated voltages are the same voltages used to compute $\mathbf{L}_1$, as described in \S\ref{sec-fwd_model}.

%--------------------------------------------------------------------
\subsection{The $\texp$ and $\tzero$ methods}\label{sec:texpG0}
%--------------------------------------------------------------------
Both the $\texp$ and $\tzero$ methods are derived from the constructive proofs of \cite{Nachman1988a,Novikov1987,Novikov1988} and involve special solutions called {\it Complex Geometrical Optics} (CGO) solutions \cite{Sylvester1987}.  A brief summary is included here for the reader's convenience.  For further details see \cite{Delbary2012,Hamiltonetal_2021,Delbary2014,Nachman1988a}.

%\subsubsection{$\texp$ method} \tville{OK to add subsecs for texp and t0?}
Assuming that the conductivity is a constant $\sigma_c=1$ in a neighborhood of the boundary $\bndry$, the real-valued conductivity equation~\eqref{eq:cond} can be transformed to the Schr\"odinger equation
\begin{equation}\label{eq:schro}
    \left(-\Delta + q(x)\right) \widetilde{u}(x)=0, \qquad x\in\R^3,
\end{equation}
via the change of variables $q(x)=\frac{\Delta\sqrt{\sigma(x)}}{\sqrt{\sigma(x)}}$ and $\widetilde{u}(x)=u(x)\sqrt{\sigma(x)}$, by extending $\sigma(x)\equiv1$ for all $x\in\R^3\setminus\bar{\Omega}$.  For $\zeta(\xi)\in \mathcal{V}_\xi$, unique CGO solutions exist to the transformed problem 
\[\left(-\Delta + q(x)\right) \psi(x,\zeta)=0, \qquad x\in\R^3,\]
where $\psi(x,\zeta)\sim e^{ix\cdot\zeta}$ for large $|x|$ or $|\zeta|$, \trev{and}
\begin{equation}\label{eq:Vxi}
\trev{\mathcal{V}_\xi} = \left\{\zeta\in\mathbb{C}^3\middle|\zeta^2=0, \;\; \left(\xi+\zeta\right)^2=0\right\}, \quad \text{for each $\xi\in\mathbb{R}^3$},
\end{equation}
where $\zeta^2 = \zeta\cdot\zeta$ and $\zeta$ is a purely auxiliary parameter.  The conductivity $\sigma(x)$ can then be recovered from the DN map $\Lambda_\sigma$ as follows.

For each $x\in\bndry$ and $\zeta\in \trev{\mathcal{V}_\xi}$, solve the Fredholm integral equation of the Second Kind, 
\begin{equation}\label{eq:psi-bndry}
    \psi(x,\zeta) = e^{ix\cdot\zeta} - \int_{\bndry} G_\zeta(x-y)\left(\Lambda_\sigma - \Lambda_1\right)\psi(y,\zeta)\;dS(y),
\end{equation}
where 
\[G_\zeta(x)= \frac{e^{ix\cdot\zeta}}{(2\pi)^3}\int_{\R^3} \frac{e^{ix\cdot k}}{|k|^2 + 2k\cdot\zeta}dk, \quad x\in\R^3\setminus\{0\},\]
denotes the Faddeev Green's function \cite{Faddeev1966}.  Then, evaluate the scattering data
\begin{equation}\label{eq:t-scat-tBIE}
\mathbf{t}(\xi,\zeta)=\int_{\partial\Omega} e^{-ix\cdot(\xi+\zeta)}\left(\Lambda_{\sigma}-\Lambda_1\right)\psi(x,\zeta)\;dS(x).
\end{equation}
For $|\zeta|$ large, the Schr\"odinger potential $q(x)$ can be recovered via the inverse Fourier transform
\begin{equation}\label{eq:t-q}
q(x)\approx \mathcal{F}^{-1}\left\{\mathbf{t}(\xi,\zeta)\right\}(x)=\frac{1}{(2\pi)^3}\int_{\mathbb{R}^3} e^{i x\cdot\xi} \mathbf{t}(\xi,\zeta)\;d\xi, \quad x\in\mathbb{R}^3.
\end{equation}
The conductivity is then recovered by solving the boundary value problem
\begin{equation}\label{eq:QtoSigma}
\left\{\begin{array}{rclcl}
(-\Delta + \trev{q}(x)) \tilde{u}(x) & =& 0 & \qquad & x\in\Omega\subset\mathbb{R}^3\\
\tilde{u}(x) & =& 1 && x\in\partial\Omega,
\end{array}
\right.
\end{equation}
for $\widetilde{u}(x)$ and evaluating $\sigma(x)=\left(\tilde{u}(x)\right)^2$.  This is the full nonlinear reconstruction method.  

Replacing the CGO solutions $\psi(x,\zeta)$ by their asymptotic behavior $e^{ix\cdot\zeta}$ in the scattering data \eqref{eq:t-scat-tBIE} via 
\begin{equation}\label{eq:texp-scat}
\mathbf{t}^{\mbox{\tiny{\textbf{exp}}}}(\xi,\zeta)=\int_{\partial\Omega} e^{-ix\cdot(\xi+\zeta)}\left(\Lambda_{\sigma}-\Lambda_1\right)e^{ix\cdot\zeta}\;dS(x).
\end{equation}
yields a {\it `Born approximation'} typically called the $\texp$ approximation for EIT.  Using this approximate scattering data in place of the fully nonlinear $\mathbf{t}(\xi,\zeta)$, one proceeds with the recovery of an approximate potential $\qexp(x)$ via \eqref{eq:t-q} and conductivity $\sigexp(x)$ via \eqref{eq:QtoSigma}.  The flow is:
\[\left(\Lambda_{\sigma},\Lambda_1\right)\overset{1}{\longrightarrow}\texp(\xi,\zeta) \overset{2}{\longrightarrow}\qexp(x) \overset{3}{\longrightarrow}\sigexp(x).\]

An intermediate approximation can be computed by replacing the Faddeev Green's function $G_\zeta(x)$ in the single layer potential, in \eqref{eq:psi-bndry}, for the traces of the CGOs, with the standard Green's function $G_0(x) = \frac{1}{4\pi|x|}$ for the Laplacian operator.  Thus, one solves
\begin{equation}\label{eq:psi0-bndry}
    \psi^0(x,\zeta) = e^{ix\cdot\zeta} - \int_{\bndry} G_0(x-y)\left(\Lambda_\sigma - \Lambda_1\right)\psi^0(y,\zeta)\;dS(y),
\end{equation}
for the CGOs $\psi^0(x,\zeta)$, avoiding the exponentially growing Faddeev Green's function.  A corresponding approximation to the scattering data is then computed by using $\psi^0(x,\zeta)$ in place of $\psi(x,\zeta)$ in \eqref{eq:t-scat-tBIE} and then continuing to recover $\qzero(x)$ and $\sigma^0(x)$.  The flow is then
\[\left(\Lambda_{\sigma},\Lambda_1\right)\overset{1}{\longrightarrow}\psi^0(x,\zeta)\overset{2}{\longrightarrow}\tzero(\xi,\zeta) \overset{3}{\longrightarrow}\qzero(x) \overset{4}{\longrightarrow}\sigma^0(x).\]

We point out that the methods, as outlined above, assumed that the conductivity was a constant $\sigma_c=1$ near the boundary of the domain.  As mentioned in section \ref{sec-cald-orig}, the problem can be scaled and unscaled as in \cite{Isaacson2004, Hamiltonetal_2021} when the constant  $\sigma_c$ is not one.  In practice, we estimate the best-fit constant conductivity fit to the data $\sigbest$ as given by \eqref{eq-gam-best}. \trev{Explicitly, as in \cite{Hamiltonetal_2021}, we scale the DN map by using $\frac{1}{\sigbest}\Lambda_\sigma$ in place of $\Lambda_\sigma$, and re-scale at the end using $\sigexp(x)=\sigbest\left(\tilde{u}^{\mbox{\tiny exp}}(x)\right)$ and $\sigma^0(x)=\sigbest\left(\tilde{u}^0(x)\right)$.}

In this work we consider both the $\texp$ and $\tzero$ methods.  We note that this is the first time that $\tzero$ has been implemented on non-continuum DN data and the first time that either $\texp$ or $\tzero$ have been demonstrated on experimental 3D EIT data, for absolute or time-difference EIT imaging.

%--------------------------------------------------------------------
\subsubsection{Numerical Implementation}\label{sec-texp-tzero-implement}
%--------------------------------------------------------------------
Here we provide the numerical details pertinent to the implementation of the $\texp$ and $\tzero$ algorithms outlined above.  As with the Calder\'on method above, the main idea is to expand functions in the same orthonormal basis $\mathbf{Q}$ as described in section \ref{sec-cald-orig}.
%for $\R^{L\by \mbox{\tiny $N_{li}$}}$ derived from the applied current patterns, where $N_{li}$ denotes the number of linearly independent current patterns applied.

Following \cite{DeAngelo2010}, for each electrode center $x_\ell$, expand $\psi(x,\zeta)$ and $e^{ix\cdot\zeta}$ as 
\begin{equation}\label{eq:basis_expansions}
\psi(x_\ell,\zeta)\approx \sum_{j=1}^{N_{li}} b_j(\zeta)\mathbf{Q}^j_\ell, \qquad% \ell=1, \ldots L \label{eq:psi_expansion}\\
e^{ix_\ell\cdot\zeta}\approx \sum_{j=1}^{N_{li}} c_j(\zeta)\mathbf{Q}^j_\ell, \qquad \ell=1, \ldots L, %\label{eq:exp_expansion}
\end{equation}
where $\mathbf{Q}^j_\ell$ denotes the $(\ell,j)$ entry of $\mathbf{Q}$.  Then, the boundary integral equation \eqref{eq:psi0-bndry} can be approximated as follows
\begin{eqnarray}
    \psi^0(x_\ell,\zeta) 
    &=&= e^{ix_\ell\cdot\zeta} - \int_{\bndry} G_0(x_\ell-y)\left(\Lambda_\sigma - \Lambda_1\right)\psi^0(y,\zeta)\;dS(y)\nonumber\\
    &\approx& e^{ix_\ell\cdot\zeta} - \sum_{\ell'=1}^L\int_{E_{\ell'}} G_0(x_\ell-y)\left(\Lambda_\sigma - \Lambda_1\right)\psi^0(y,\zeta)\;dS(y)\nonumber\\
    &\approx& e^{ix_\ell\cdot\zeta} - \left[\sum_{\ell'=1}^L\int_{E_{\ell'}} G_0(x_\ell-y)dS(y)\right]\left[\left(\mathbf{L}_\sigma - \mathbf{L}_1\right)\psi^0(y_{\ell'},\zeta)\right],\label{eq:psi_bie_expan1}
\end{eqnarray}
where $E_{\ell'}$ denotes the $\ell'^{\mbox{\tiny th}}$ {\bf \trev{extended}} electrode where $\bigcup_{\ell'=1}^L E_{\ell'}=\bndry$ and the $E_{\ell'}$ are mutually disjoint \cite{Hyvnen2009ApproximatingIB}.  \trev{Note the true electrodes need not cover the surface $\bndry$, only the extended (mathematical) electrodes that we will use to discretize the integral.}  Then, using the expansions from \eqref{eq:basis_expansions} in \eqref{eq:psi_bie_expan1}, we have
\begin{eqnarray}
    \sum_{j=1}^{N_{li}} b_j(\zeta)\mathbf{Q}^j_\ell 
    &\approx& \sum_{j=1}^{N_{li}} c_j(\zeta)\mathbf{Q}^j_\ell - \left[\sum_{\ell'=1}^L\int_{E_{\ell'}} G_0(x_\ell-y)dS(y)\right]\left[\left(\mathbf{L}_\sigma - \mathbf{L}_1\right)\sum_{j=1}^{N_{li}} b_j(\zeta)\mathbf{Q}^j_{\ell'}\right] \nonumber\\
    &=& \sum_{j=1}^{N_{li}} c_j(\zeta)\mathbf{Q}^j_\ell - \left[\sum_{\ell'=1}^L\int_{E_{\ell'}} G_0(x_\ell-y)dS(y)\right]\left[\sum_{j=1}^{N_{li}} b_j(\zeta)f_j\left(y_{\ell'}\right)\right], \label{eq:eq:psi_bie_expan2}
\end{eqnarray}
where $f_j\left(y_{\ell'}\right)$ represents the action of $\left(\mathbf{L}_\sigma - \mathbf{L}_1\right)$ on $\mathbf{Q}^j$ evaluated at $y_{\ell'}$, which can be computed as the $(\ell',j)$ entry of $\mathbf{Q}\left(\mathbf{L}_\sigma - \mathbf{L}_1\right)$.  Define
\begin{equation}\label{eq:G0-mat}
    \widetilde{\mathbf{G}}_0(\ell,\ell')=\begin{cases} \frac{1}{4\pi\left|x_\ell - y_{\ell'}\right|} & \ell\neq\ell'\\
    0 & \ell=\ell',
    \end{cases}
\end{equation}
where we have removed the singularities at $x_\ell=y_\ell$.  Assuming $|E_\ell|=\frac{|\bndry|}{L}$ for each $\ell=1,\ldots,L$, and using \eqref{eq:G0-mat} in \eqref{eq:eq:psi_bie_expan2} we find
\begin{eqnarray}
    \sum_{j=1}^{N_{li}} b_j(\zeta)\mathbf{Q}^j_\ell 
    &\approx& \sum_{j=1}^{N_{li}} c_j(\zeta)\mathbf{Q}^j_\ell - \frac{|\bndry|}{L}\sum_{j=1}^{N_{li}} b_j(\zeta) \sum_{\ell'=1}^L \widetilde{\mathbf{G}}_0(\ell,\ell')f_j\left(y_{\ell'}\right) \nonumber\\
    &\approx& \sum_{j=1}^{N_{li}} c_j(\zeta)\mathbf{Q}^j_\ell - \frac{|\bndry|}{L}\sum_{j=1}^{N_{li}} b_j(\zeta) \left[\widetilde{\mathbf{G}}_0\mathbf{Q}\left(\mathbf{L}_\sigma - \mathbf{L}_1\right)\right](\ell,j) \nonumber
\end{eqnarray}
or in matrix form,
\[\mathbf{Q}\vec{b} = \mathbf{Q}\vec{c} - \frac{|\bndry|}{L}\widetilde{\mathbf{G}}_0\mathbf{Q}\left(\mathbf{L}_\sigma - \mathbf{L}_1\right)\vec{b}.\]
The solution to this equation can be found by solving the following system for the unknowns $\vec{b}$,
\begin{equation}\label{eq:psi0-bie-matrixform}
\left(I+A\right)\vec{b} = \vec{c},    
\end{equation}
where $A=\frac{|\bndry|}{L}\mathbf{Q}^{\mbox{\tiny \bf T}}\widetilde{\mathbf{G}}_0\mathbf{Q}\left(\mathbf{L}_\sigma - \mathbf{L}_1\right)$, and $I$ is the identity matrix.  \trev{If the extended electrodes $E_\ell$ are not uniform in size, then one could compute as weighted sum replacing the uniform weight $|E_\ell|=\frac{|\bndry|}{L}$ as appropriate.}

Next, following \cite{Hamiltonetal_2021}, the scattering data $\tzero(\xi,\zeta(\xi))$ can be computed for all $|\xi|$ less than a chosen truncation radius $T_{\xi}$ via 
\begin{equation}\label{eq:tzero-matrix}
\tzero(\xi,\zeta)\approx
\begin{cases}
\frac{|\bndry|}{L} \left[e^{-i\mathbf{x}\cdot(\xi + \zeta)}\right]^{\mbox{\tiny \bf T}}\mathbf{Q}\left(\mathbf{L}_\sigma - \mathbf{L}_1\right)\vec{b} & |\xi|<T_{\xi}\\
0 & else,
\end{cases}
\end{equation}
since $\vec{b}=\mathbf{Q}^{\mbox{\tiny \bf T}}\psi^0(\mathbf{x},\zeta)$, where $\mathbf{x}\in\R^{(L\by3)}$ is the same as in Section \ref{sec-cal-Fhat-implement}, i.e. a vector storing the centers of the electrodes.

Next, the approximate potential $\qzero$ is recovered by computing the inverse Fourier transform of the truncated scattering data $\tzero$ using a Simpson's rule in 3D,
\[\qzero(x)=\frac{1}{(2\pi)^3}\int_{[-T_\xi,T_\xi]^3} e^{i x\cdot\xi} \tzero(\xi,\zeta(\xi))\;d\xi.\]
Alternatively, one could use an IFFT to achieve additional speedup, taking care with quadrature points and the particular form of the kernel.  Following \cite{Hamiltonetal_2021}, the conductivity $\sigma^0(x)$ was recovered by first solving the boundary value problem \eqref{eq:QtoSigma}, using the {\sc PDE toolbox} in {\sc Matlab} using a mesh with approximately $21,000$ 3D elements, then computing $\sigma^0(x)=\sigma{\mbox{\tiny{best}}}\left(\tilde{u}(x)\right)^2$.  The 3D visualizations were obtained by interpolation to a $64\times64\times64$ rectangular grid using {\sc Matlab}'s {\tt scatteredInterpolant} function.
The $\zeta$ values were computed following the minimal-zeta approach outlined in \cite{Delbary2014}.  We remark that while the \trev{non-}existence of exceptional points for the solution of the boundary integral equation \eqref{eq:psi-bndry} is proven for large $|\zeta|$ \cite{Nachman1988a, Sylvester1987}, as well as small $|\zeta|$  \cite{Cornean2006}; it is still an open question for the intermediate values required here to perform the computation on a computer. 

The $\texp$ reconstructions were obtained in an analogous fashion to those of $\tzero$, this time bypassing the boundary integral equation \eqref{eq:t-scat-tBIE} and \trev{directly} computing 
\[
\texp(\xi,\zeta)\approx
\begin{cases}
\frac{|\bndry|}{L} \left[e^{-i\mathbf{x}\cdot(\xi + \zeta)}\right]^{\mbox{\tiny \bf T}}\mathbf{Q}\left(\mathbf{L}_\sigma - \mathbf{L}_1\right)\mathbf{Q}^{\mbox{\tiny \bf T}}\left[e^{i\mathbf{x}\cdot\zeta}\right] & |\xi|<T_{\xi}\\
0 & else,
\end{cases}\]
where we have replaced the vector of coefficients $\vec{b}$ with the expansion of the asymptotic behavior of $e^{ix\cdot\zeta}$ given by $\mathbf{Q}^{\mbox{\tiny \bf T}}\left[e^{i\mathbf{x}\cdot\zeta}\right]$.

Difference imaging can be performed with the $\tzero$ and $\texp$ methods by replacing the matrix $\mathbf{L}_1$ with $\mathbf{L}_{\sigma_{\mbox{\tiny ref}}}$ and computing $\sigma_{\mbox{\tiny diff}(x)} = \sigma(x) - \sigbest$ in the final step.

\subsection{Reference methods}
Total Variation regularized non-linear least squares reconstructions will serve as the reference reconstructions for the CGO-based absolute imaging cases considered here.  A classic linear difference imaging scheme, also reviewed below, will serve as the reference for the CGO-based difference images.

\subsubsection{Absolute imaging with TV regularization} 
A widely used numerical approach for absolute EIT is the total variation (TV) regularized non-linear least squares minimization 
\begin{equation}
  \hat\sigma = {\rm arg} \min_{\sigma>0}\{\Vert L_e \left( V-U(\sigma,z_\ast) \right)\Vert^2
+ \alpha {\rm TV}_\beta (\sigma)\},
\label{lssol}
\end{equation}
where $U(\sigma,z)$ is the finite dimensional forward map, $z_\ast \in \mathbb{R}^L$ are the fixed electrode contact impedances obtained from an initialization step of the minimization, $L_e$ is the Cholesky factor of the noise precision matrix of $e$ so that $L_e^{\rm T}L_e = \Gamma_e^{-1}$, scalar valued $\alpha$ is the regularization parameter and ${\rm TV}_\beta (\sigma)$ is the (smooth) TV regularization functional \cite{Rudin1992}
\begin{equation}\label{stvpsi}
{\rm TV}_\beta (\sigma) = \int_\Omega \sqrt{\| \nabla \sigma \|^2 + 
\beta }\ {\rm d} r,
\end{equation}
where $\beta$ is the (fixed) smoothing parameter. The forward model $U(\sigma,z)$ in \eqref{lssol} is based on the finite element (FEM) discretization of the complete electrode model \cite{Somersalo1992}. For details of the FEM model, see
\cite{Kaipio2000,Vauhkonen1998a,Vauhkonen1997}. In the FEM model, the electric conductivity is approximated as a linear combination of the piecewise linear nodal basis functions in a uniform tetrahedral mesh of $N$ nodes, leading to vector of unknowns $\sigma \in \mathbb{R}^N$, and the electric potential is approximated similarly in a significantly more dense tetrahedral mesh with refinements near the electrodes. 
The non-linear optimization in \eqref{lssol} is solved by a lagged Gauss-Newton method equipped with
a line search algorithm. The line search is implemented using bounded minimization such that the non-negativity $\sigma > 0$ is enforced. For more details of the method, see \cite{Toivanen2021}.

The fixed contact impedances $z_\ast$ and initial (constant) conductivity estimate for \eqref{lssol} are obtained from the solution of the non-linear least squares problem
\begin{equation}
  (\sigma_0, z_\ast) = {\rm arg} \min_{\sigma_c, z >0} \{\Vert L_e \left( V-U(\sigma_c,z) \right)\Vert^2\},
\label{lsinit}
\end{equation}
where the scalar $\sigma_c \in \mathbb{R}$ is the coefficient of a spatially constant conductivity image $\sigma_c \mathbf{1}$ and $z\in \mathbb{R}^L$. The non-linear least squares problem \eqref{lsinit} is solved by a Gauss Newton optimization.

\subsubsection{Linear difference imaging}
\label{sec:LD}
\trev{In linear difference imaging, see e.g. \cite{Bagshaw2003, Barber1987fast},} the objective is to reconstruct the change in conductivity between two measurements (V1; V2). \trev{The reference linear difference imaging approach of this paper uses linearized approximations of the observation models}
\begin{eqnarray}
    V_1 \approx U(\sigma_0,z_\ast) + J_\sigma (\sigma_1-\sigma_0) + e_1\\
    V_2 \approx U(\sigma_0,z_\ast) + J_\sigma (\sigma_2-\sigma_0) + e_2,
\end{eqnarray}
where the Jacobian matrix $J_\sigma$ of $U(\sigma,z_\ast)$ is evaluated at $\sigma = \sigma_0$. With these linearized models, the difference in the measurements becomes
\begin{eqnarray}
    \nonumber &\delta V &= V_2 - V_1\\
    \nonumber &&= \left(U(\sigma_0) + J_\sigma (\sigma_2-\sigma_0) + e_2\right) - \left(U(\sigma_0) + J_\sigma (\sigma_1-\sigma_0) + e_1\right) \\
    \label{equ:dV}
     &&= J_\sigma \delta\sigma + \delta e,
\end{eqnarray}
where $\delta \sigma = \sigma_2 - \sigma_1$ and $\delta e = e_2 - e_1$.
Now, the inverse problem is to reconstruct $\delta \sigma$ based on the difference data $\delta V$ and the model \eqref{equ:dV}. A widely used formulation for the linear problem is to use the (generalized) Tikhonov regularization with a smoothness promoting regularization functional
\begin{equation}
\label{equ:minLD}
  \hat{\delta\sigma} = \arg \min_{\delta \sigma} \left\{ || L_{\delta e}(\delta V - J_\sigma \delta\sigma) ||^2 + % p_{LD}(\delta\sigma) 
  ||L_p \delta \sigma||^2
  \right\},
\end{equation}
where $L_{\delta e}$ is the Cholesky factor of the noise precision matrix of $\delta e$ so that $L_{\delta e}^{\rm T}L_{\delta e} = \Gamma_{\delta e}^{-1} = \left(\Gamma_{e_1} + \Gamma_{e_2}\right)^{-1}$.
In this paper, the regularization matrix $L_p$ is constructed by utilization of a distance based covariance function. More specifically, we set $L_p^{\mathrm{T}}L_p = \Gamma_p^{-1}$, where the (prior covariance) matrix $\Gamma_p$ is constructed using the distance based correlation function 
\cite{Lieberman10}
\begin{align}
\label{equ:GammaSmooth}
  \Gamma_{p}(i,j) = \mathrm{std(}\sigma)^2\exp\left(-\frac{\|x_i - x_j \|^2}{2a^2}\right), \ \ \ i, j = 1,\ldots, N,
\end{align}
where the parameter $a$ controls the correlation length and can be solved by setting the distance $\|x_i - x_j \|$ to a selected value $d$ (e.g. half the radius of the target) and setting $\Gamma_{p}(i,j)$ to the desired covariance for that distance (e.g. 1\% of variance).

%--------------------------------------------------------------------
\subsection{Experimental Setup}\label{sec:experiments}
%--------------------------------------------------------------------
ACT5 \cite{rajabishishvan2021} is a 32 electrode parallel EIT instrument that uses 32 current sources to apply patterns of current to the target and 32 voltmeters to measure the resulting voltages. The current sources in ACT5 adjust the delivered current to compensate for current lost through shunt impedance, including that introduced by the capacitance of cables that connect the instrument to the electrodes, enabling the desired current patterns to be applied with high precision \cite{saulnier2020}.  ACT5 can apply sinusoidal currents at frequencies in the range of 11~kHz to 1~MHz with a peak amplitude of 0.25~mA. The voltmeters measure voltages  up to 0.5~V peak with a maximum signal-to-noise ratio (SNR) of 96~dB. Because the current source compensates for shunt capacitance of the cables, the system operates with grounded-shield cables \cite{abdelwahab2020}. 

A test tank was built from 3/8-inch-thick Plexiglas with interior dimensions of 17.0cm~x~25.5cm x~17.0cm.  Thirty-two electrodes were fabricated of 16-gauge 316 stainless steel, each 80.0~mm square.  The inside surfaces of the tank were milled out in the shape and depth of the electrodes so that the interior surfaces were flat.  The nominal gaps between electrodes were 5.0~mm. The 32 electrodes were placed four on each end and six on each side.  Bolts through the center of each electrode passed through the Plexiglas, allowing connections to the cables from the ACT5 electronics.  Five sides of the tank were glued together, and provision was made to secure the top with threaded rods at all four corners.  Access holes of 0.6~mm at various sites in the top to allow threads to suspend targets, and 25-Ga hypodermic needles to complete the filling of the tank after the lid was in place.  In use, a mixture of saline at room temperature at the desired background conductivity is made and measured using a Oakton CON~6+ handheld conductivity meter.  

Targets were made the day before experiments by adding NaCl and a few drops of food dye to distilled water until the desired conductivity was obtained, as measured by an immersible conductivity meter.  The solution was then heated slowly with stirring as 4\% agar (Fisher Scientific) was added until the temperature reached $85^{\circ}$C$-90^{\circ}$C.  The mixture was poured into 5.27~cm inside-diameter spherical molds and allowed to harden overnight.  Test cells were also filled at this time so the final conductivities could be verified using the ACT5 instrument.

Data were first collected of the tank filled with saline of the desired conductivity (24 mS/m) with no targets present.  Then, targets were added to the tank supported on fine toothpicks from below.  Photos were captured of the tank, with each target in place, to verify its/their position.  The conductivity values of the spheres used were approximately 290~mS/m. Averaged data, over 100 frames, were used in the subsequent conductivity reconstructions.  The experimental data will be made publicly available on github.\footnote{See \url{https://github.com/sarahjhamilton/open3D_EIT_data}}

\trev{Note that even though we are testing on a rectangular prism tank, the mathematical reconstruction algorithms are not limited to this geometry.}

\begin{figure}
    \centering
    \includegraphics[width=0.25\textwidth]{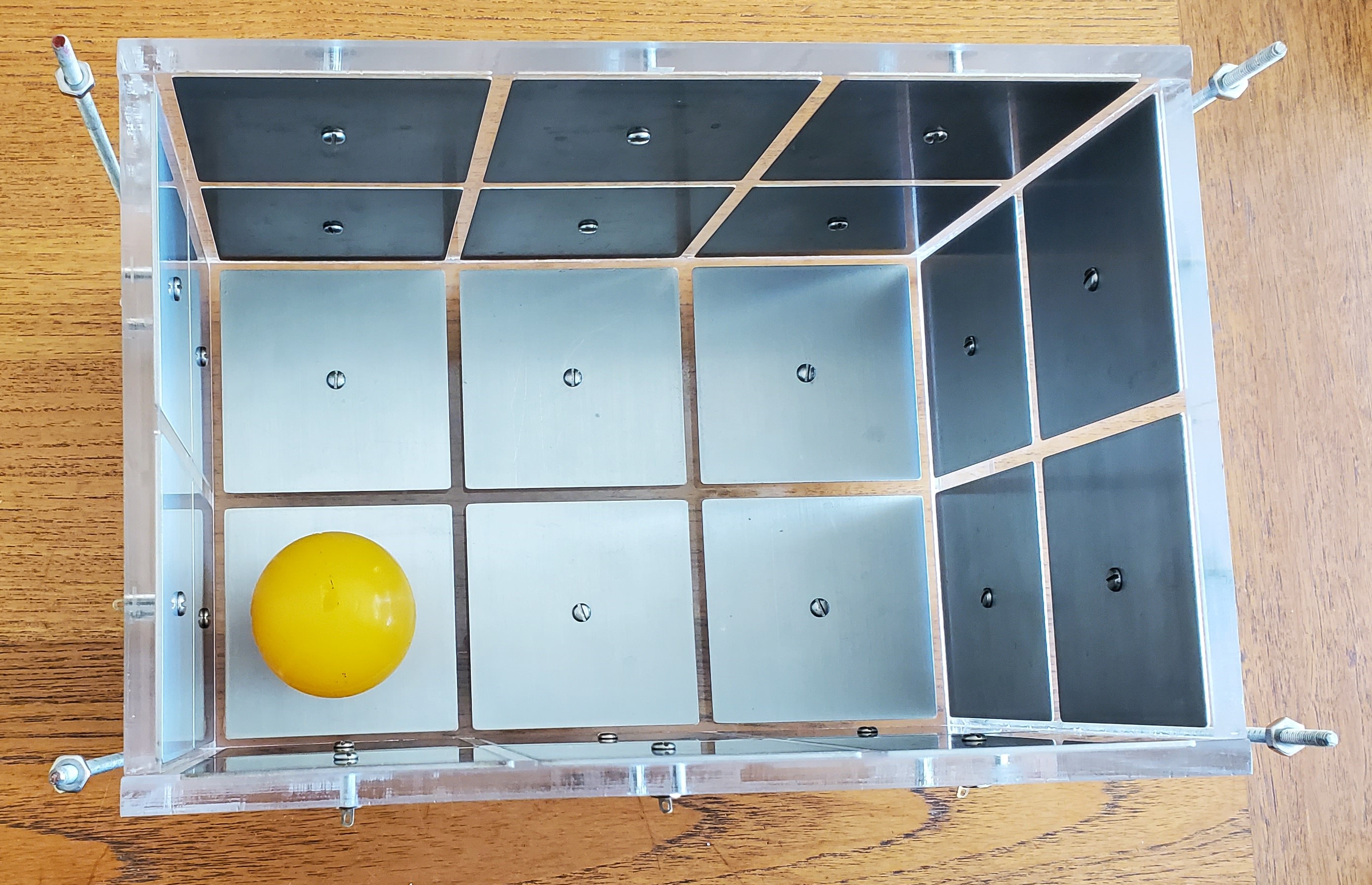}
    \hspace{3 em}
    \includegraphics[width=0.25\textwidth]{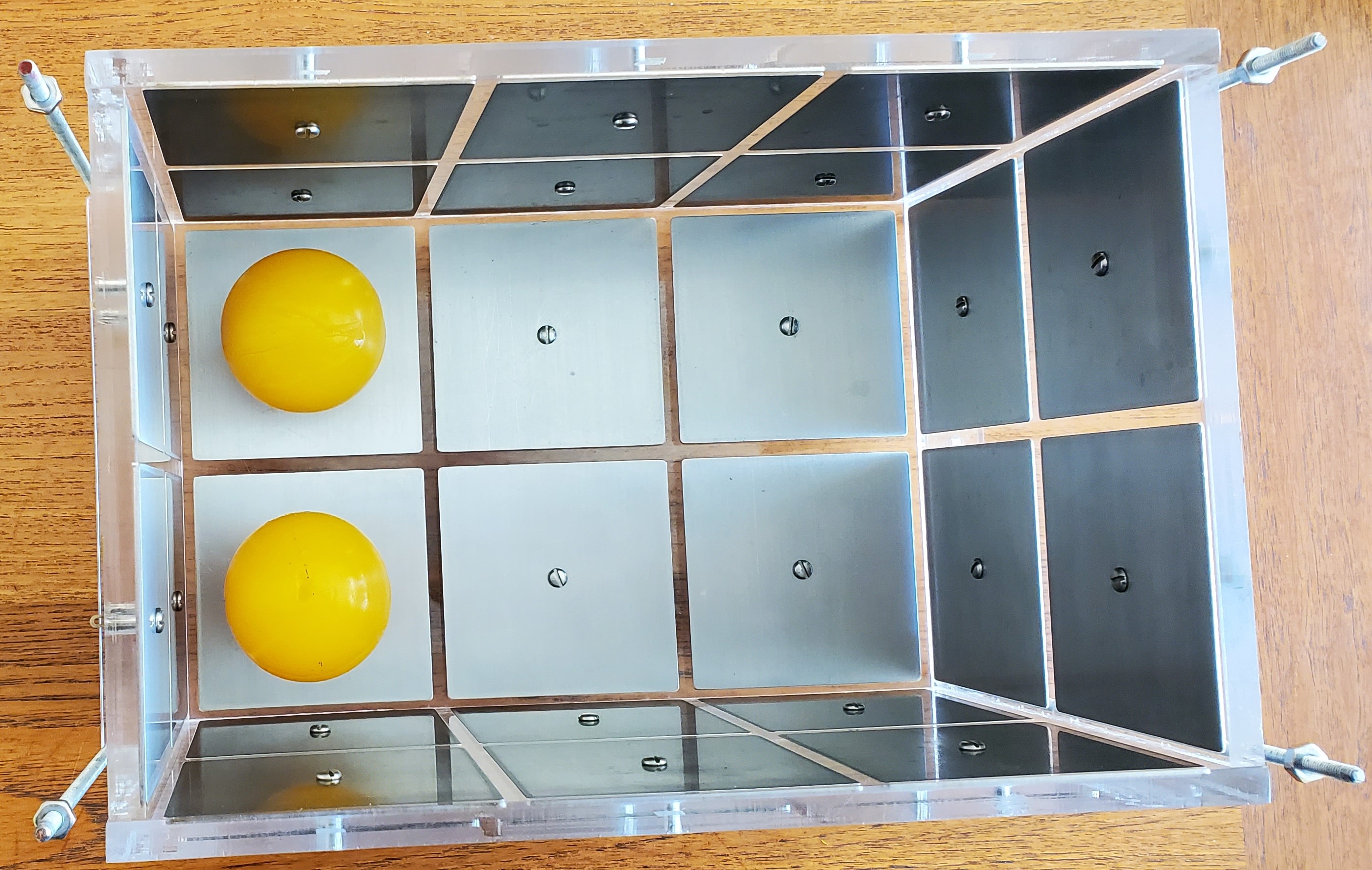}
    \hspace{3em}
    \includegraphics[width=0.25\textwidth]{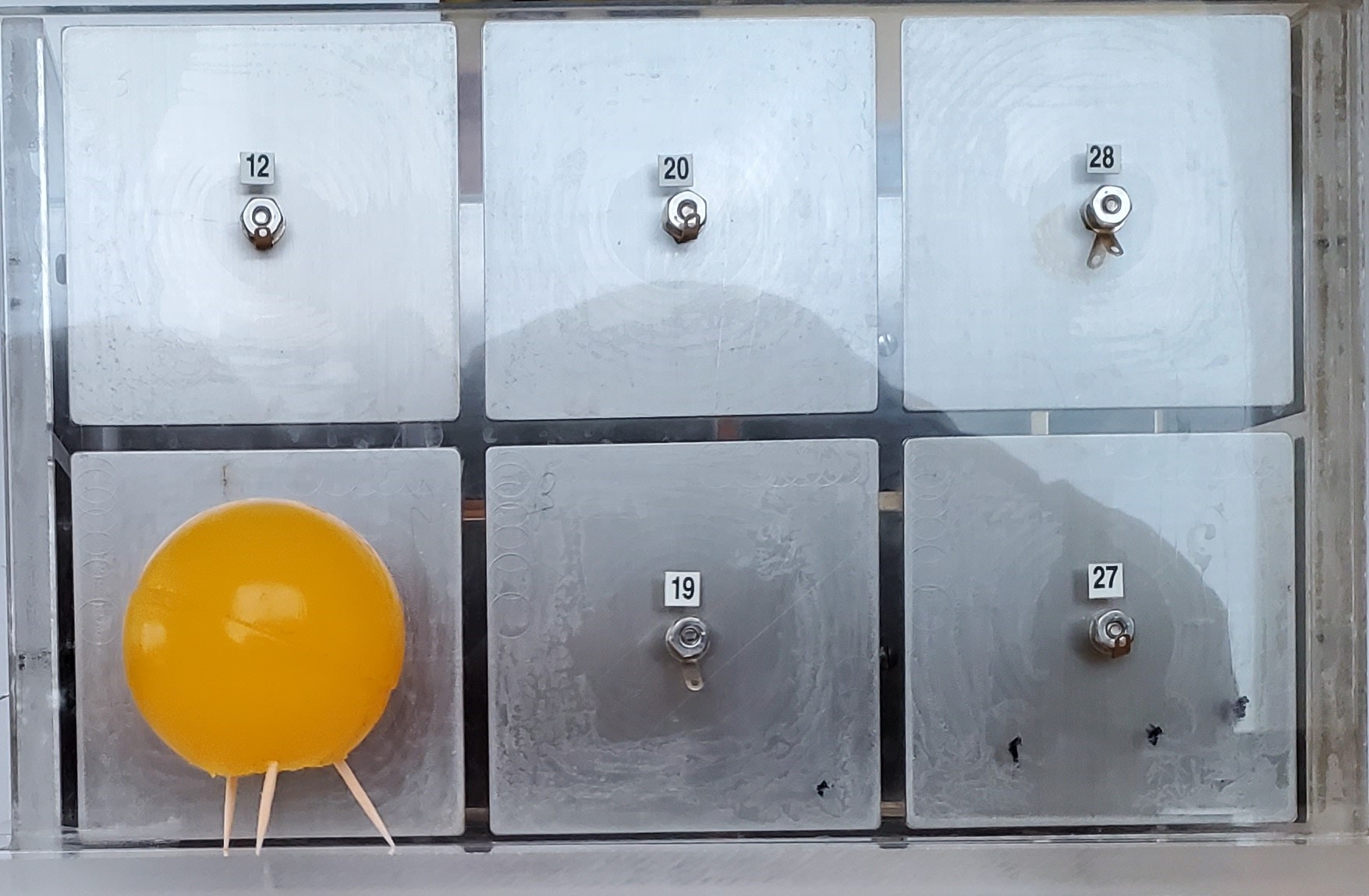}
    \caption{Experimental setups. Left: Top view of the one target experiment. Middle: Top view of two target experiment. Right: side view of height of targets above floor of tank. Note that the targets were measured at 290~mS/m in test cells.}
    \label{fig:truths}
\end{figure}

%--------------------------------------------------------------------
\subsection{Robustness Tests}\label{sec:robustness}
%--------------------------------------------------------------------
In addition to testing the reconstruction methods with the correct domain modeling, we consider a moderate, as well as strong, mis-modeling of the domain.  The incorrect domain modeling comes into play in the following places for the CGO-based methods.  First in the simulation of the $\mathbf{L}_1$ DN data that requires solving the forward EIT problem with a known conductivity $\sigma\equiv 1$\trev{~S/m} and the applied current patterns.  Next, the incorrect domain modeling presents as incorrect information about the centers of the electrodes $\mathbf{x}$ and surface area for the domain $|\bndry|$; the subsequent reconstructions of $\sigcal$, $\sigmaexp$, and $\sigma^0$ are recovered on the incorrect domain. For the TV regularized method, the incorrect domain modeling is present in the FEM based forward map $U(\sigma,z)$ used in the minimization problem \eqref{lssol} and computational mesh. For the linear difference imaging method, the incorrect domain modeling is present in the Jacobian matrix of the forward map $J_{\delta \sigma}$ in the minimization problem \eqref{equ:minLD} and computational mesh.

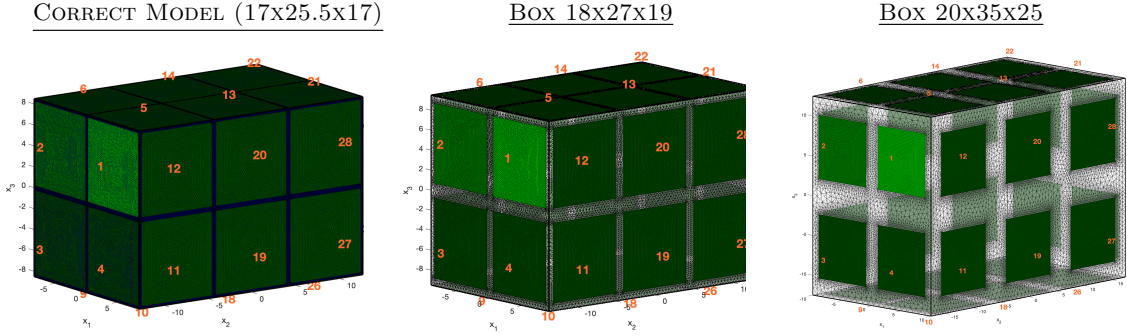
\begin{figure}[h!]
    \centering
    \begin{picture}(420,125)
    
    \put(-20,-10){\includegraphics[height=120pt]{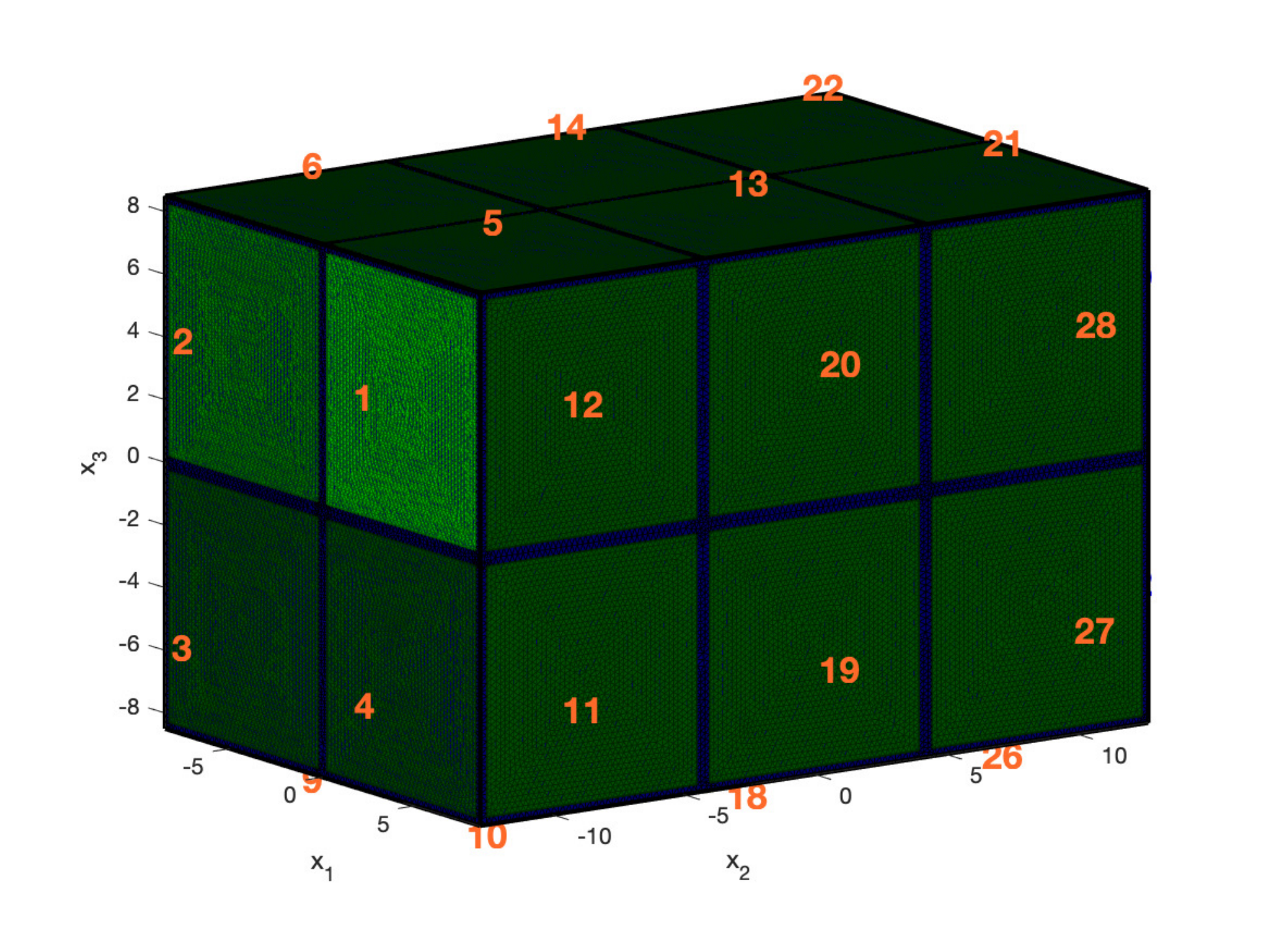}}
    
    \put(130,-10){\includegraphics[height=120pt]{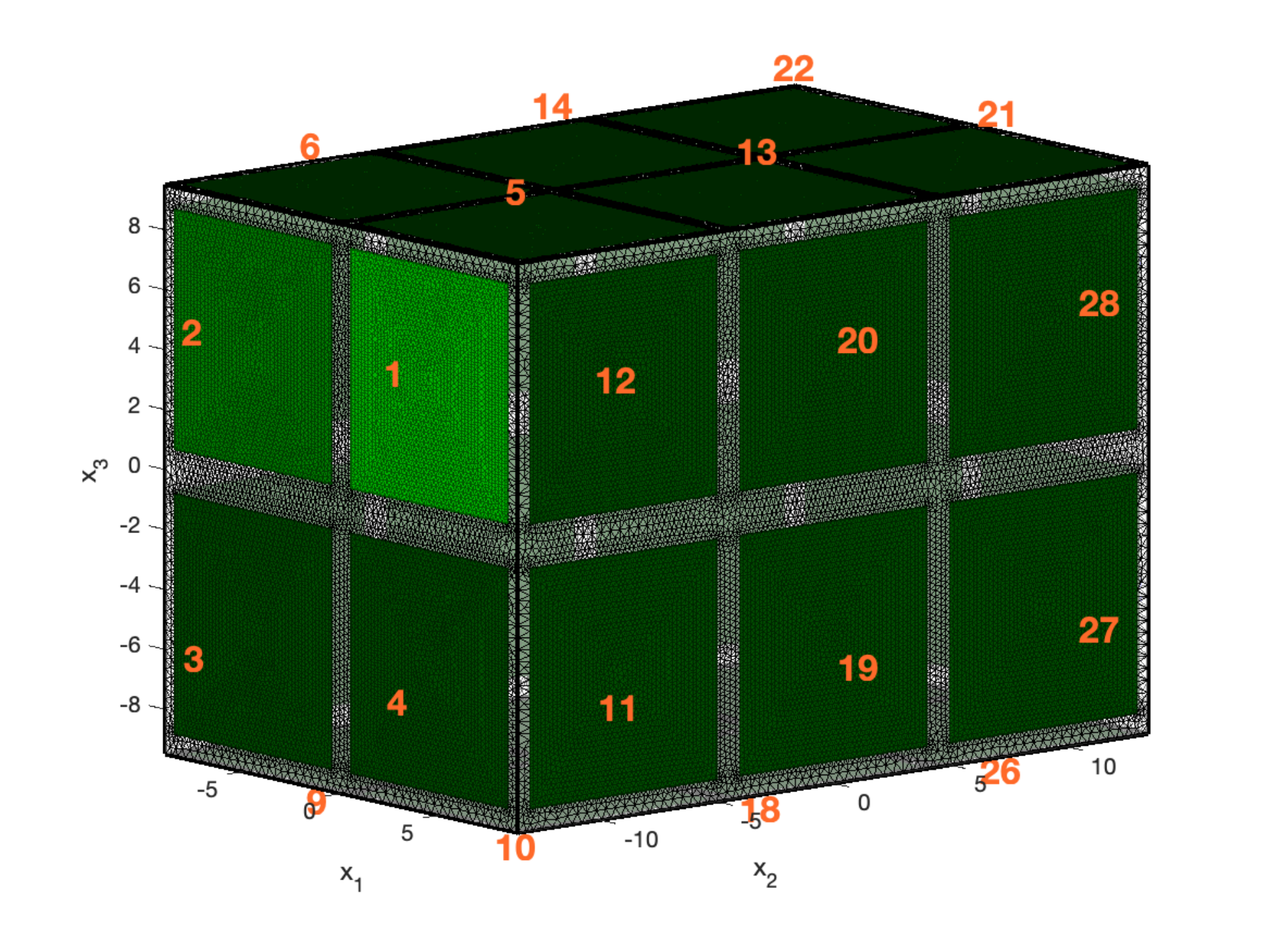}}
    
    \put(270,-10){\includegraphics[height=120pt]{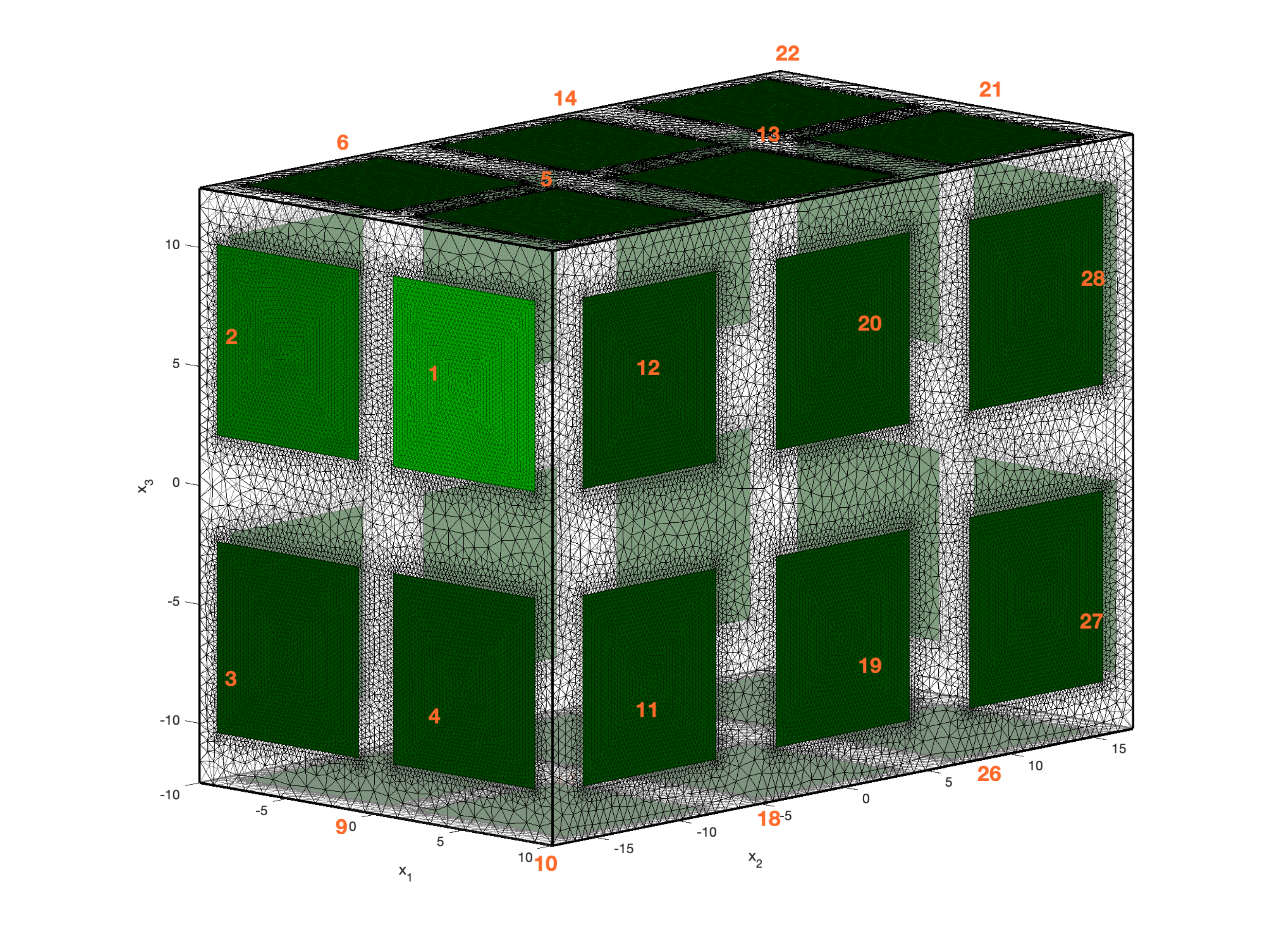}}
    
    \put(0,115){\underline{\footnotesize\sc Correct Model (17x25.5x17)}}
    \put(180,115){\underline{\footnotesize\sc Box 18x27x19}}
    \put(320,115){\underline{\footnotesize\sc Box 20x35x25}}
    
    \end{picture}
\caption{3D renderings of simulated boxes used in robustness testing.}
    \label{fig:robusttests}
    
\end{figure}

Accurate domain modeling can be challenging in practice.  To address this we explored a domain with a moderately incorrect domain shape, $18$cm~x~$27$cm~x~$19$cm, as well as a domain with a significantly incorrect domain shape, $20$cm~x~$35$cm~x~$25$cm\footnote{Due to the smaller ratio of longest side to shortest side, the $16\times16\times16$ $X$-grid used for Calder\'on's method as mentioned in \S\ref{sec-cal-Fhat-implement} does not capture the whole domain, so a $32\times32\times32$ $x$-grid is used in the significantly incorrect domain case}.  Electrodes were uniformly distributed as in the original box design, but clearly had different locations and spacing than the truth due to the new box dimensions.  Figure~\ref{fig:robusttests} (center and right) depict the larger box domains. Each reconstruction method was then tested assuming that the actual experimental voltage and current measurements were coming from the boxes shown in Figure~\ref{fig:robusttests}.

%--------------------------------------------------------------------
%--------------------------------------------------------------------

Time-difference reconstructions were performed using basal measurements with only \trev{24~mS/m} saline in the experimental tank.  This data was then used in the reconstruction methods discussed in \S\ref{sec-texp-tzero-implement} and \S\ref{sec:LD}.

%--------------------------------------------------------------------
\subsubsection{Forward modeling and Regularization Parameters}\label{sec-fwd_model}
%--------------------------------------------------------------------
For the CGO-based methods, the matrix approximation $\mathbf{L}_1$ to the DN map $\Lambda_1$ was formed using simulated voltage data, using EIDORS \cite{AdlerLionheart_2006,Vauhkonen1997}.  The EIDORS data is produced by solving the forward conductivity problem with the finite element approximation of the Complete Electrode Model~\eqref{eq:CEM} with $\sigma\equiv 1$\trev{~S/m}.  In each domain modeling case, see Figure~\ref{fig:robusttests}, the forward problem was solved on a box of the corresponding dimensions, with 8cm by 8cm electrodes, using approximately $250,000$ nodes and 1.2 million elements with high refinement used near the electrodes, and the default contact impedance in EIDORS. For the reference methods, the FEM based forward problem was solved using the same meshes for the electric potential as with the CGO-based methods and the conductivities were approximated in sparser, uniform meshes with $8,918$; $9,410$; $14,678$ nodes and $43,690$; $45,106$; $72,438$ elements, respectively, for the different domain modeling cases.

The regularization parameters for each of the reconstruction methods were chosen as follows. The regularization parameters $T_{z_1}$ and $T_{z_2}$ for Calder\'on's method were empirically chosen to be 1.4 and 1.7, respectively, and the mollifying parameter was chosen to be $t=0.1$ for all reconstructions. These values were chosen empirically to provide the best visual reconstructions.  
The truncation radius $T_\xi$ for the nonlinear scattering data used in the $\texp$ and $\tzero$ methods was chosen empirically in the range [10.5, 12] for each case shown here.  The localization of the target(s) did not appear to change significantly with $T_\xi$, however the contrast does appear sensitive to this parameter choice.  Further discussion is given below in Section~\ref{sec:discussion}. For the total variation regularized method, the regularization parameter $\alpha$ and the smoothing parameter $\beta$ were \trev{selected by computing a series of reconstructions with varying $\alpha$ and $\beta$ values, in the range [1e-6, 100] for $\alpha$ and [1e-4, 0.1] for $\beta$, and by choosing the parameters that gave the best visual quality of reconstructions in the correct domain model case. The chosen} values $\alpha = 0.01$ and $\beta = 0.001$ were used in all the test cases with different domain models.  For the linear difference imaging method, the visually best reconstructions were obtained using a standard deviation of conductivity $\text{std}(\sigma) = 16\sigma_0$ and parameter $a$ calculated using a covariance of 1 $\%$ of the variance at a correlation distance $d = l_{\text{box}}/8$, where $l_{\text{box}}$ is the length of the longest side of the box that is used as the domain model.

%--------------------------------------------------------------------
\subsection{Evaluation Metrics}\label{sec:metrics}
%--------------------------------------------------------------------
As the main focus of this work is experimental data reconstruction, comparing to a known `truth' with zero error is not possible. The conductivity values for the targets were measured at approximately 290~mS/m and the targets were placed to be roughly centered in a sub-cube of the overall box.  These data allow us to approximate the localization error (LE) of our reconstructions and maximum conductivity in each reconstructed target. 

LE, as in \cite{Hamiltonetal_2021}, measures the  distance between the centroids of the reconstructed targets and those of the true targets by 
\begin{equation}\label{eq:LE}
        \text{LE}=\sqrt{(x_1^{\mbox{\tiny recon}}-x_1^{\mbox{\tiny truth}})^2+(x_2^{\mbox{\tiny recon}}-x_2^{\mbox{\tiny truth}})^2+(x_3^{\mbox{\tiny recon}}-x_3^{\mbox{\tiny truth}})^2}.
\end{equation}
A localization error of $0$ is ideal, and signifies the targets are reconstructed in the correct location. In our experiments, the centers of the spherical targets were placed at the centers of the nearest electrode in each direction and in \eqref{eq:LE},  $(x_1^{\mbox{\tiny truth}},x_2^{\mbox{\tiny truth}},x_3^{\mbox{\tiny truth}})$ is our best estimate of the true location of the target center, acknowledging possible errors within millimeters of the truth. Following \cite{Hamiltonetal_2021}, to obtain the centers of the reconstructed targets, we use a thresholded segmentation of the reconstructions to identify targets with {\sc MATLAB}'s \texttt{regionprops3}. The thresholds for segmentation were used to empirically identify regions with more than a certain percentage of the maximum reconstructed conductivity. We note that the volume of the segmented targets changes with the choice of threshold, but as we increase the percentage (i.e. threshold) the locations of the centroids stabilize (up to millimeters), from which we computed the LE.  

\trev{As a major focus of this work is to study how the algorithms tolerate modeling errors, we report a {\bf scaled} version of the LE, 
\begin{equation}\label{eq:LE_scaled}
        \text{Scaled LE}=\frac{\sqrt{(x_1^{\mbox{\tiny recon}}-x_1^{\mbox{\tiny truth}})^2+(x_2^{\mbox{\tiny recon}}-x_2^{\mbox{\tiny truth}})^2+(x_3^{\mbox{\tiny recon}}-x_3^{\mbox{\tiny truth}})^2}}{25.5},
\end{equation}
where we have divided the LE from \eqref{eq:LE} by 25.5cm, the length of the longest side of the true tank.  This scaling value for the computing the LE metric is kept fixed across incorrect modeling scenarios.} 

%--------------------------------------------------------------------
\section{Results}\label{sec:results}
\begin{figure}[ht]
\centering
\begin{picture}(450,410)
\linethickness{.3mm}
%----------------- 
%----------------- 
% Truth
\put(300,320){\includegraphics[width=60pt]{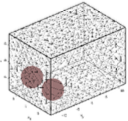}}
%\put(350,320){\includegraphics[width=100pt]{cald_true_1_targ_ABS_3D-eps-converted-to.pdf}}
\put(240,320){\includegraphics[angle=-90,origin=c,width=42pt]{TwoTargs_TopView.jpeg}}

\put(70,320){\includegraphics[width=60pt]{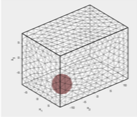}}
%\put(120,240){\includegraphics[width=100pt]{cald_true_2_targ_ABS_3D-eps-converted-to.pdf}}
\put(10,320){\includegraphics[angle=-90,origin=c,width=42pt]{OneTarg_TopView.jpeg}}
% %----------------- 
% cald
\put(290,240){\includegraphics[width=90pt]{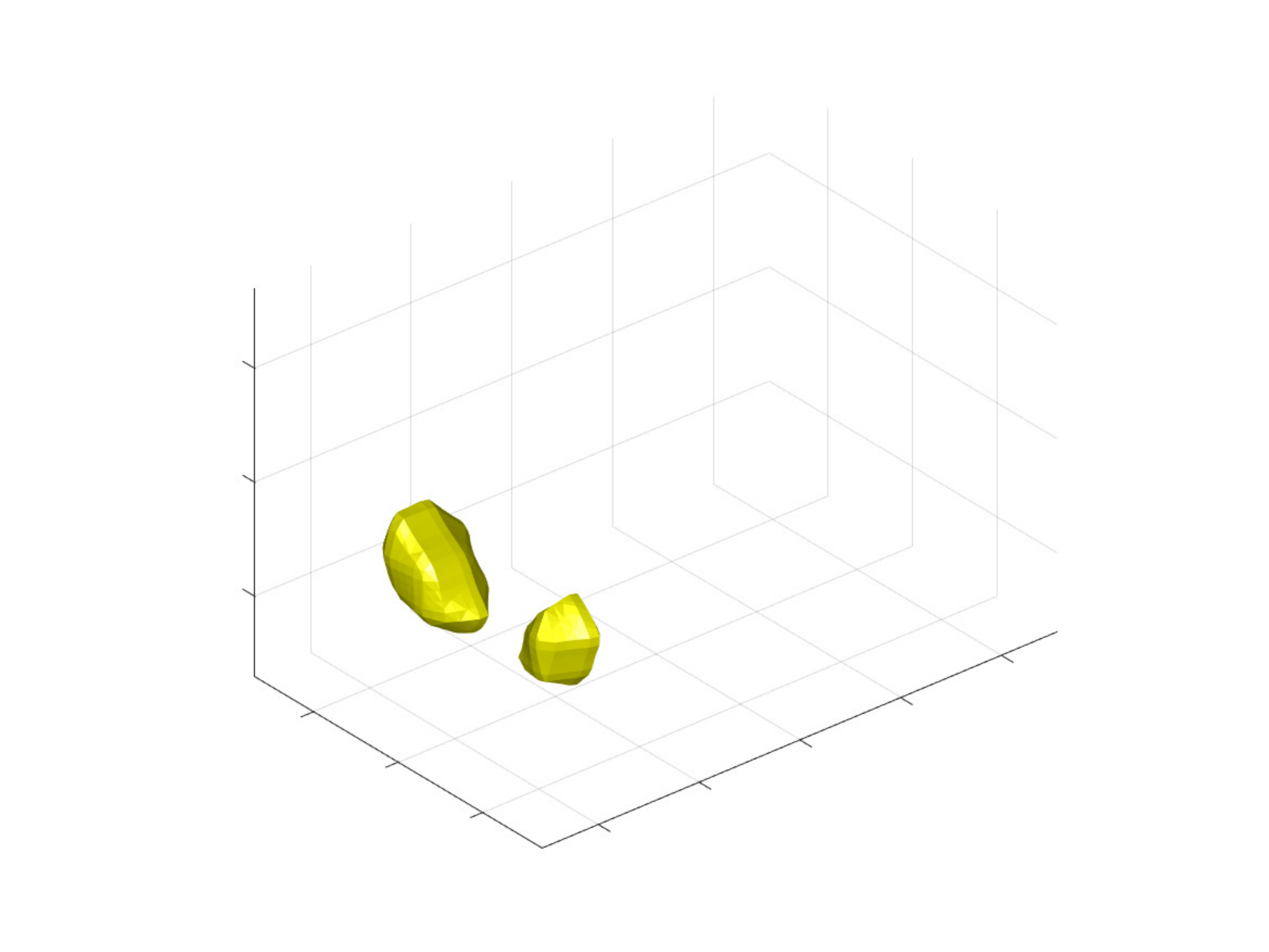}}
\put(350,240){\includegraphics[width=100pt]{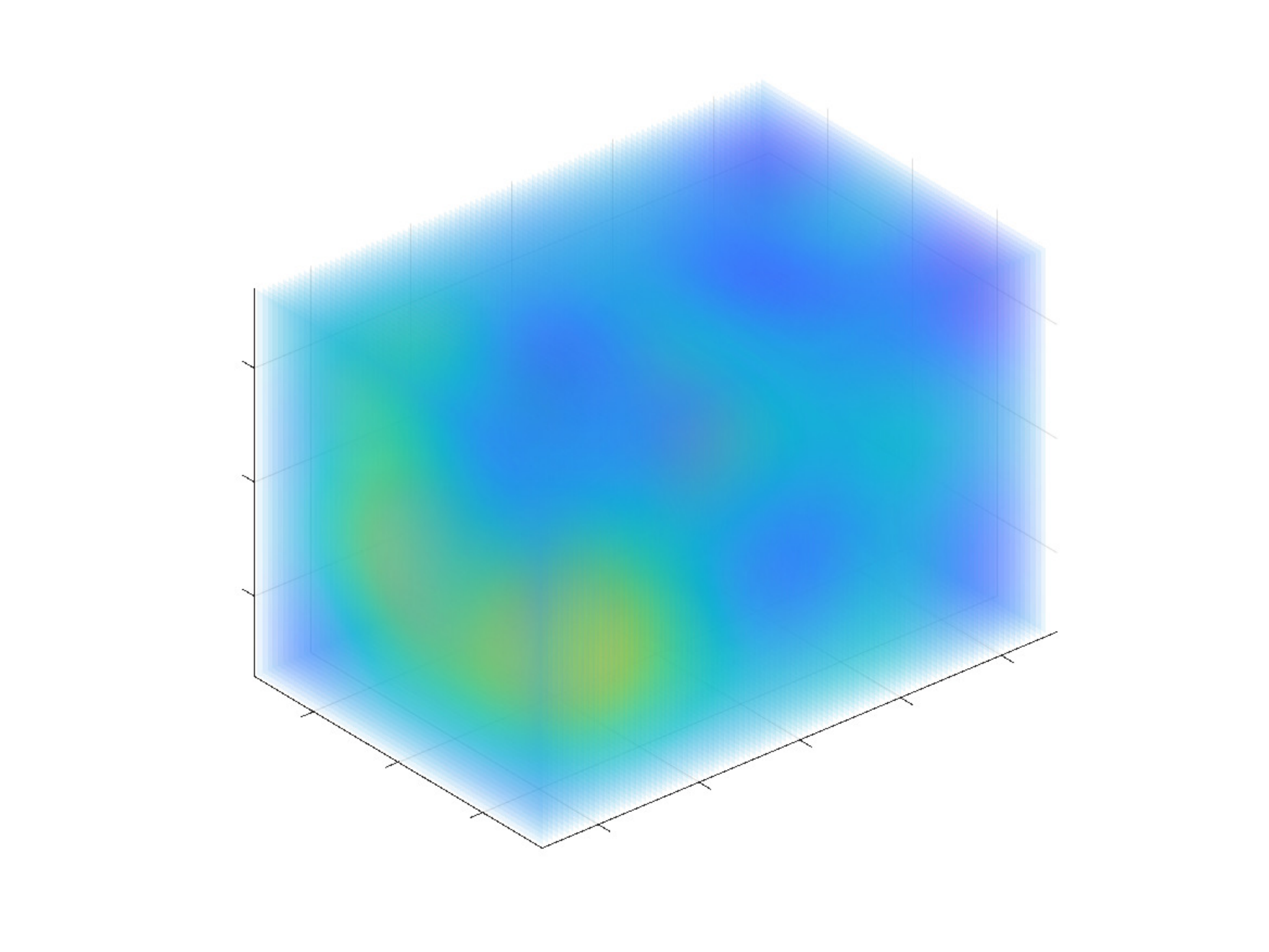}}
\put(230,240){\includegraphics[width=70pt]{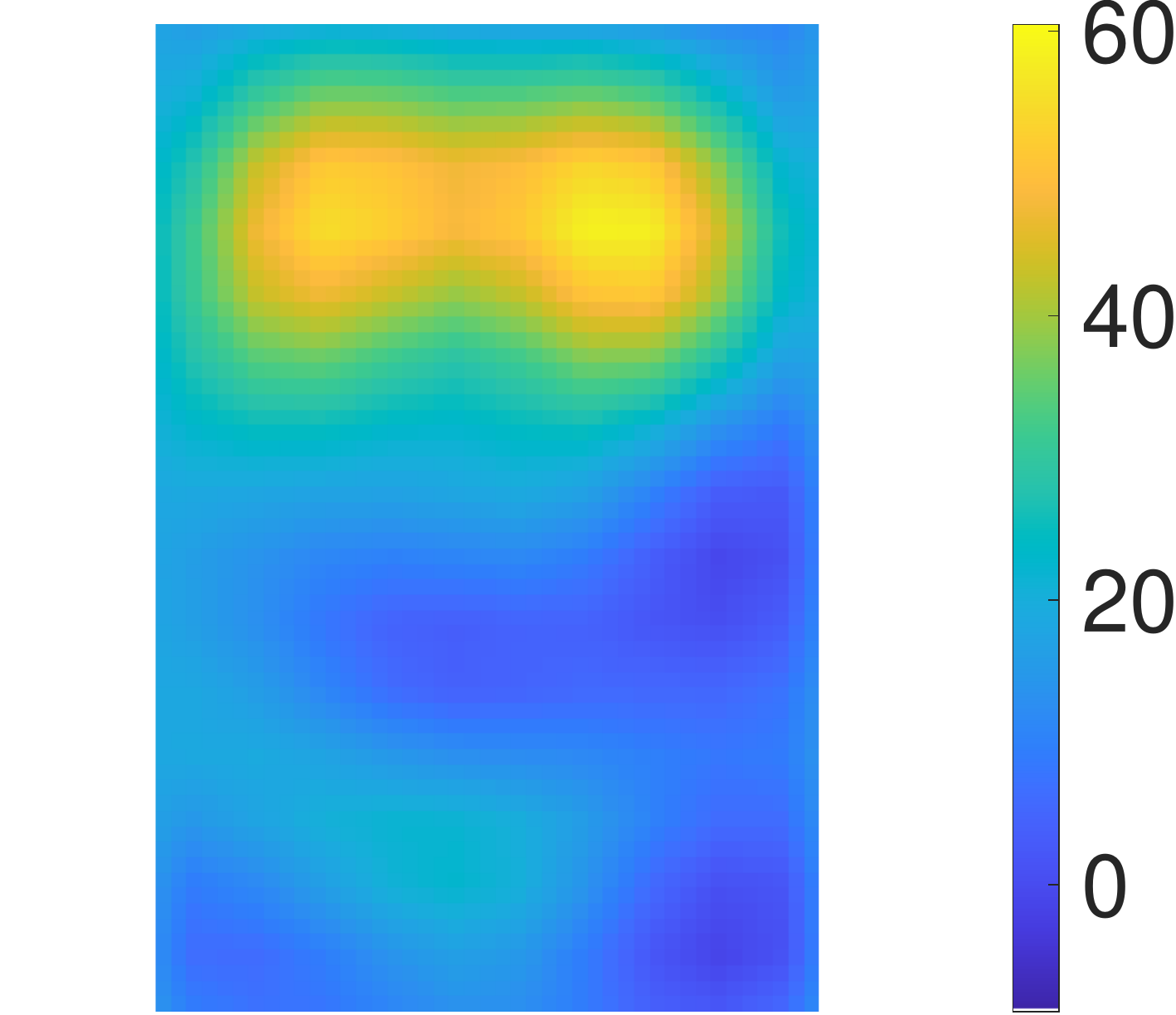}}

\put(60,240){\includegraphics[width=90pt]{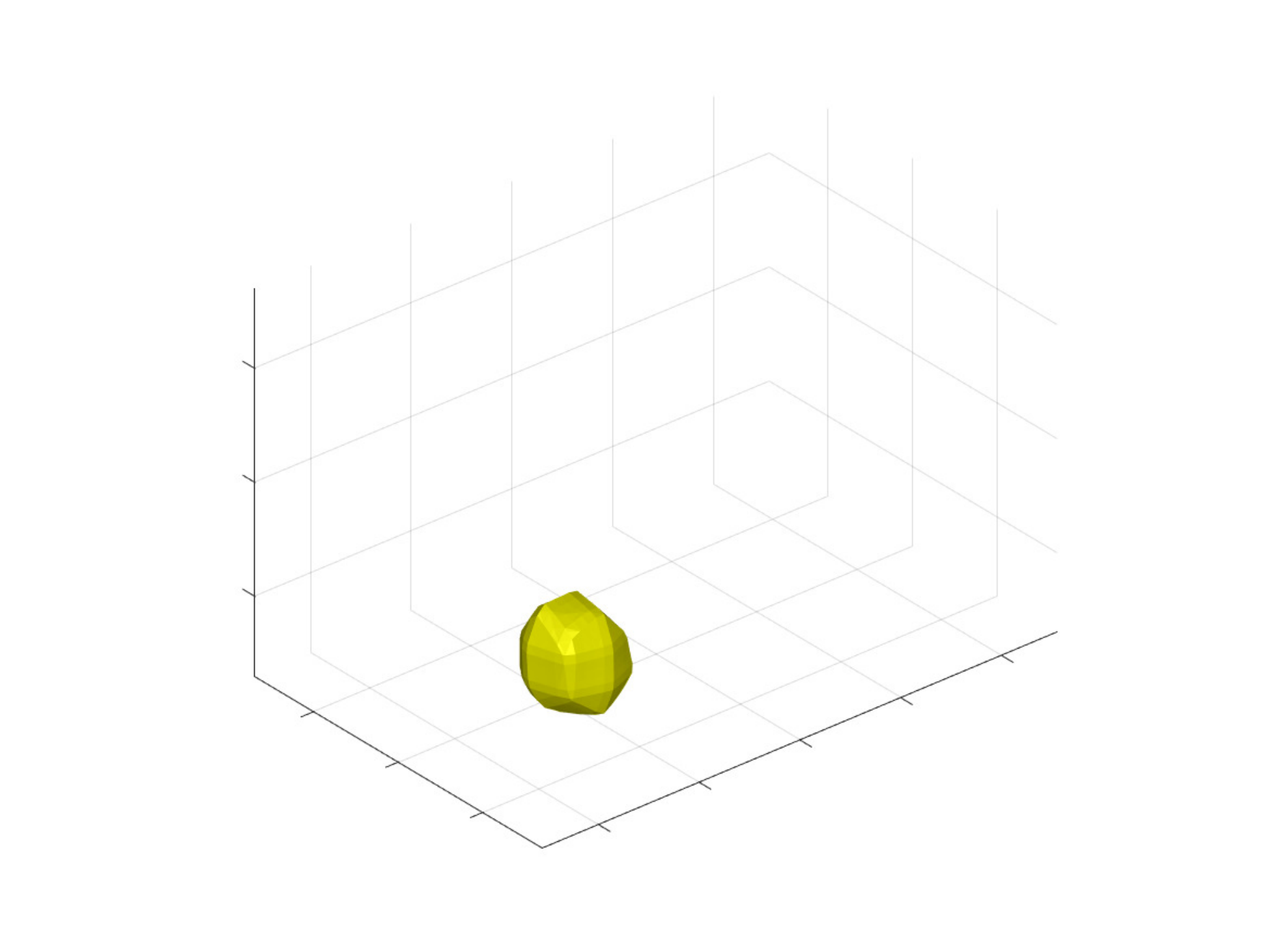}}
\put(120,240){\includegraphics[width=100pt]{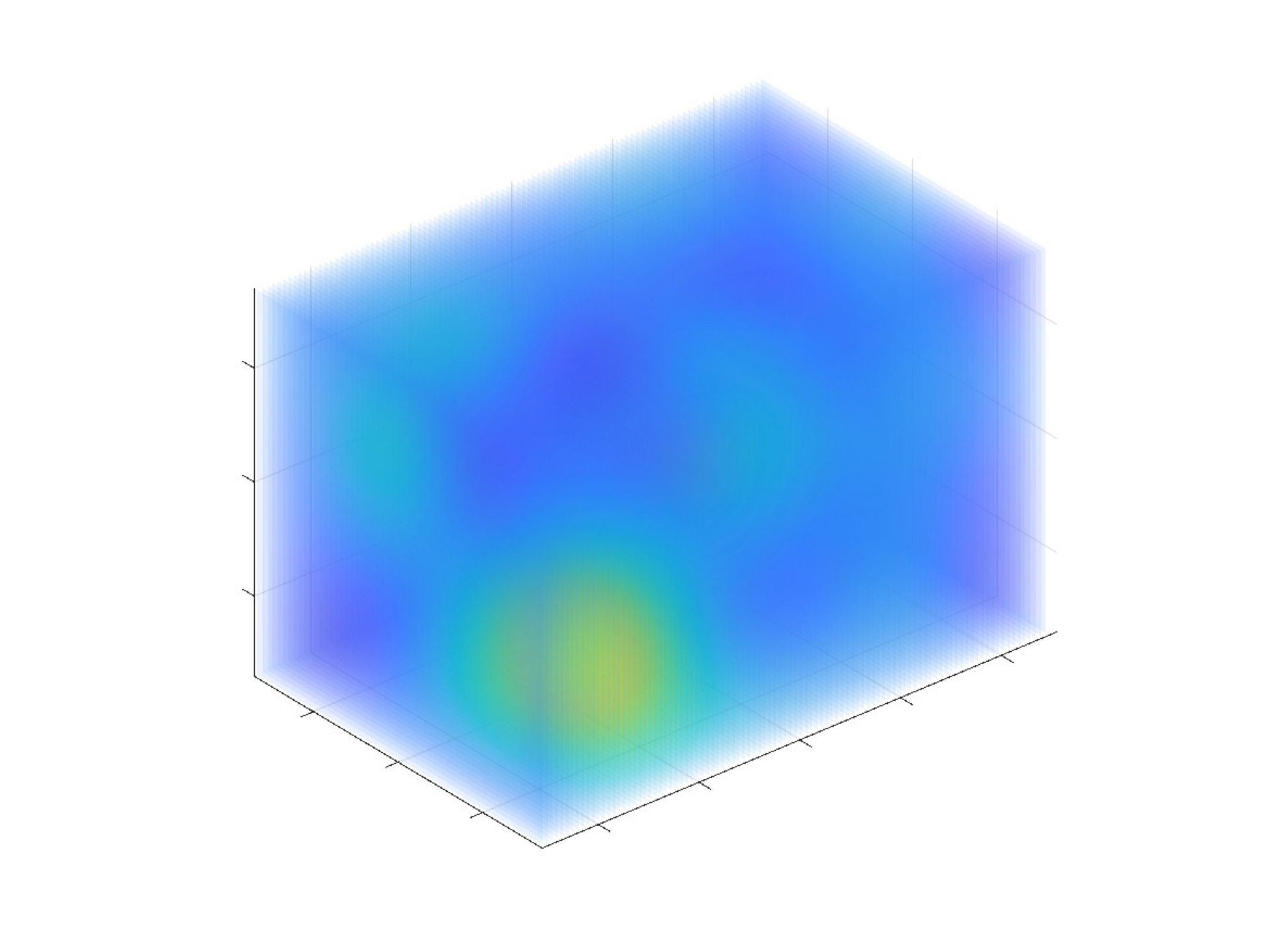}}
\put(0,240){\includegraphics[width=70pt]{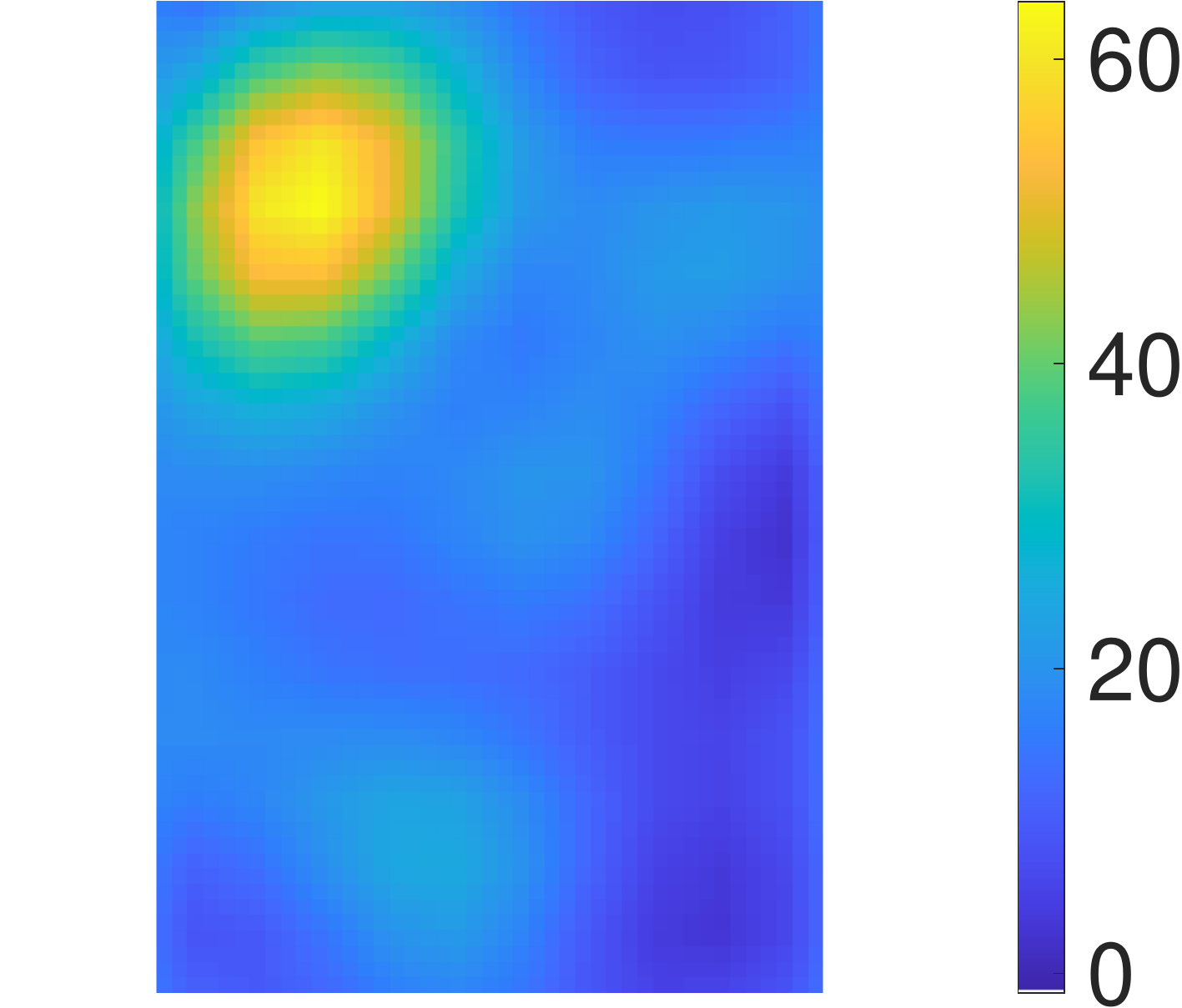}}
% %----------------- 
% texp
\put(290,160){\includegraphics[width=90pt]{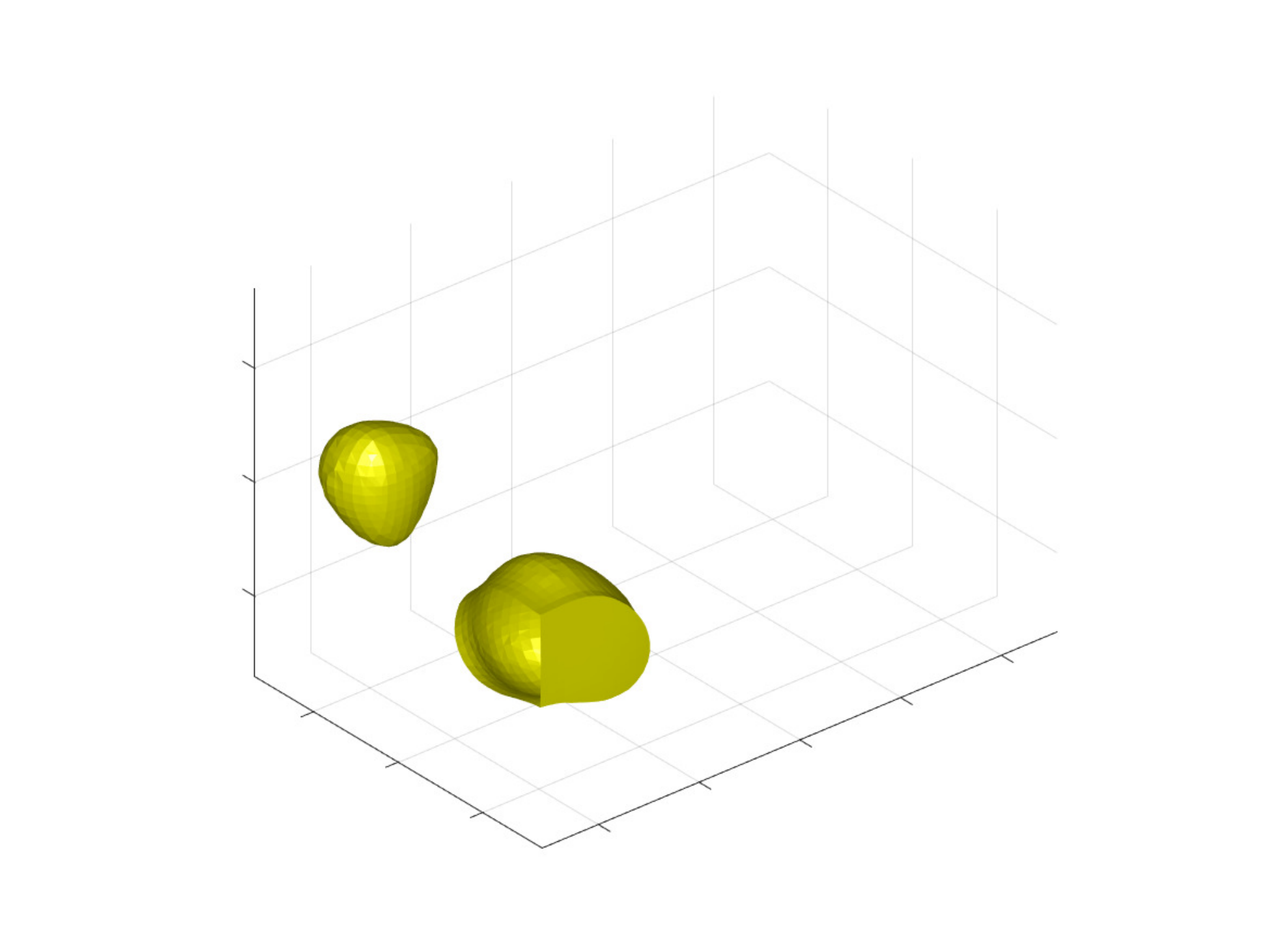}}
\put(350,160){\includegraphics[width=100pt]{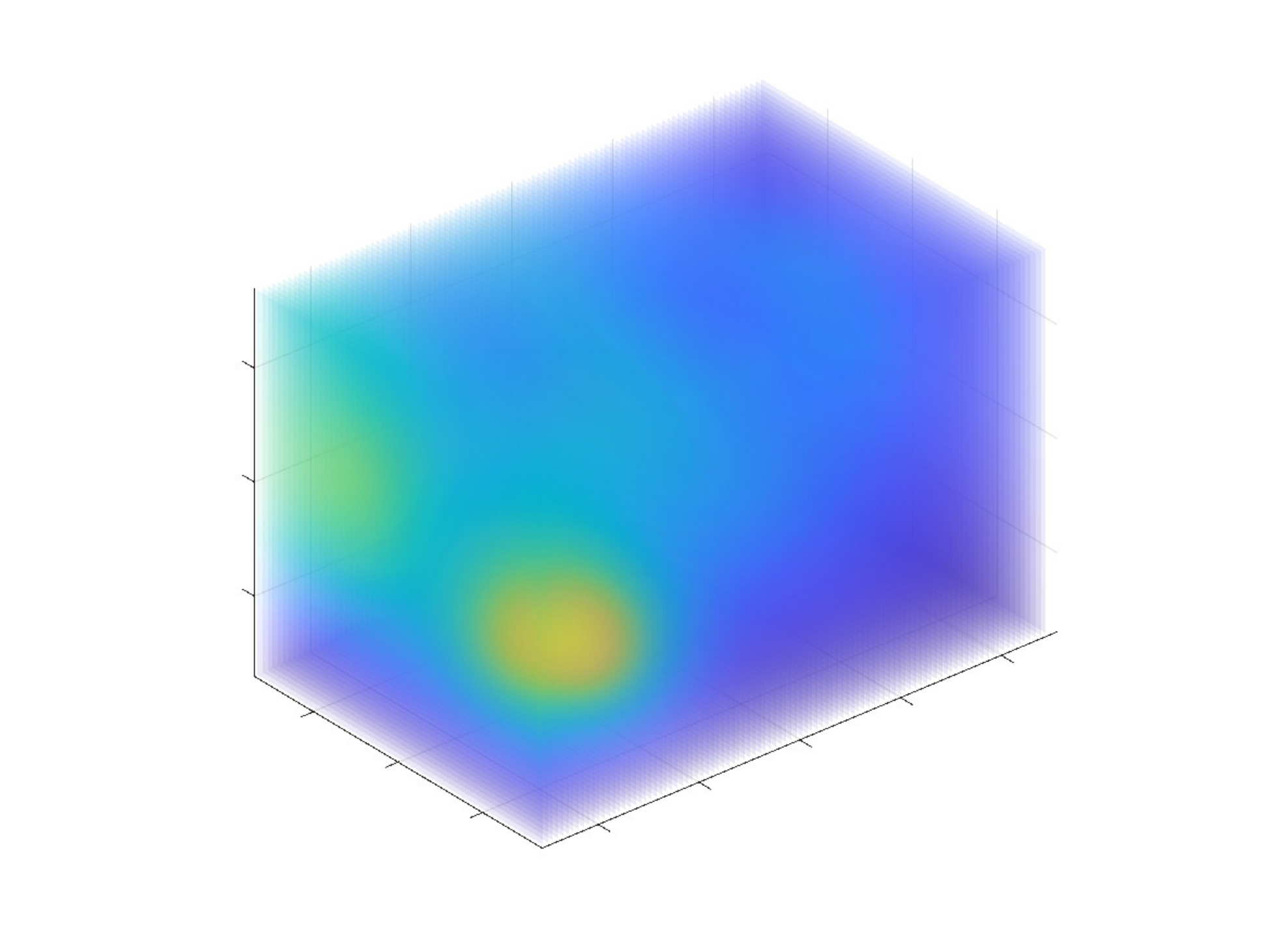}}
\put(230,160){\includegraphics[width=70pt]{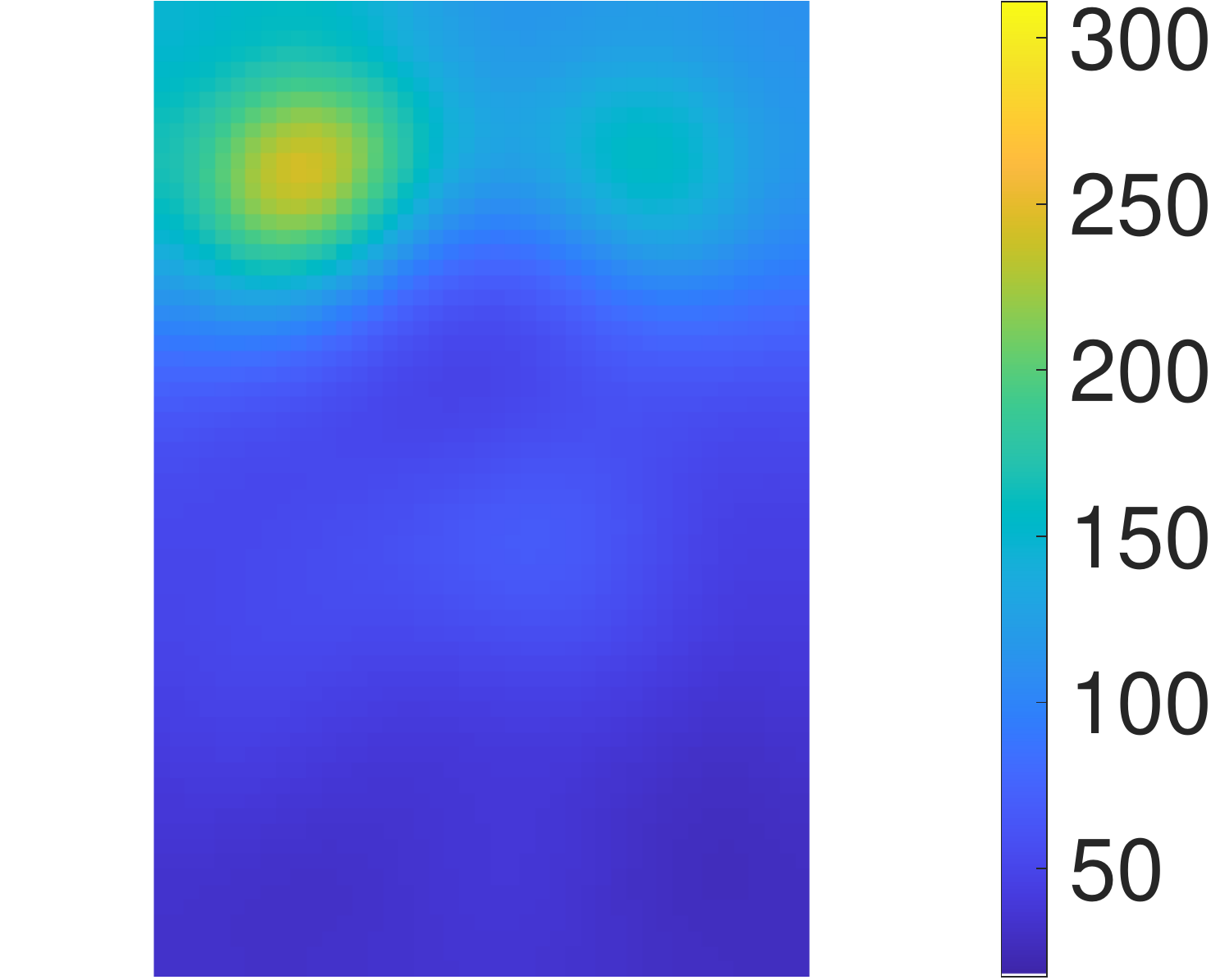}}

\put(60,160){\includegraphics[width=90pt]{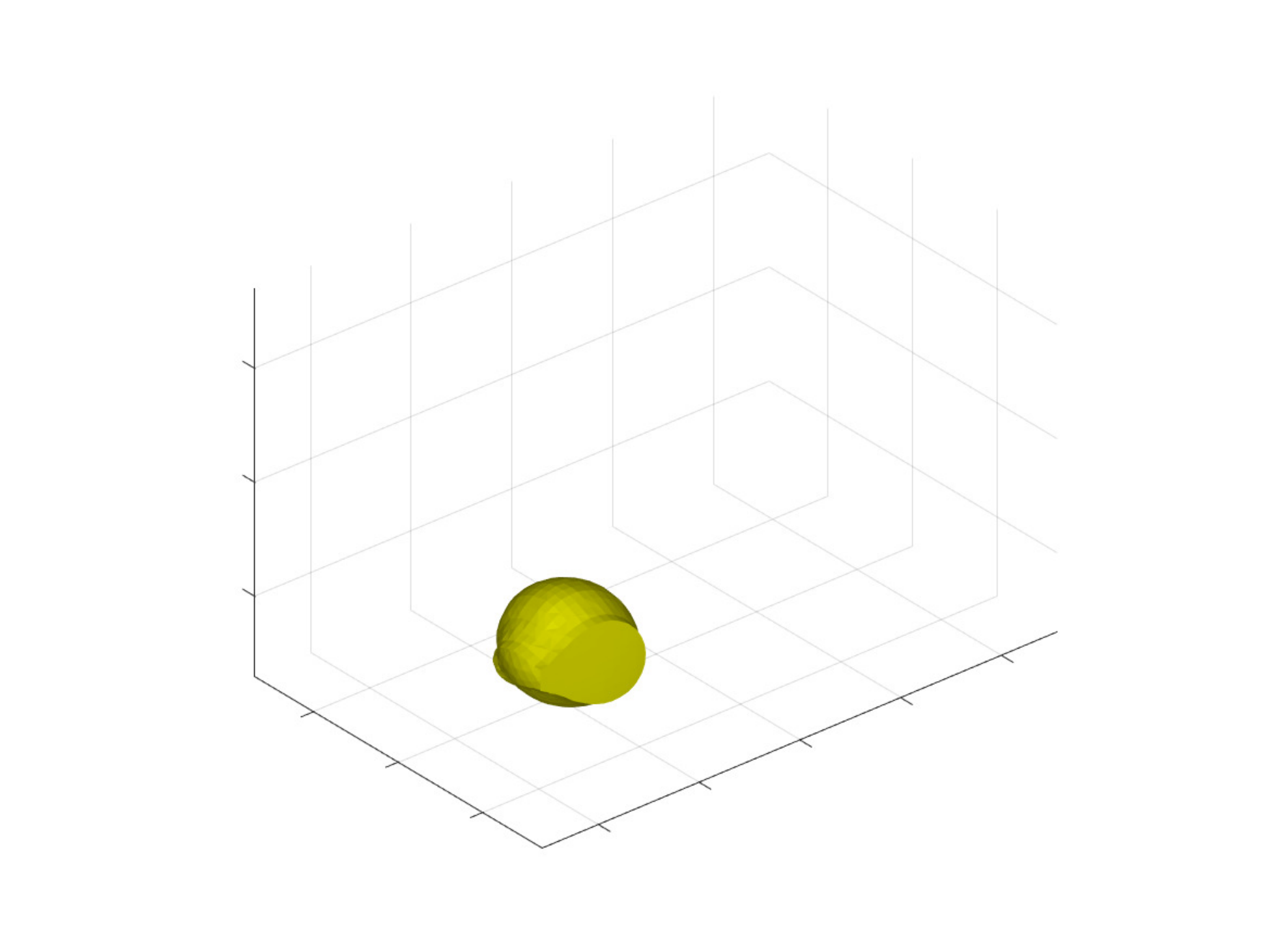}}
\put(120,160){\includegraphics[width=100pt]{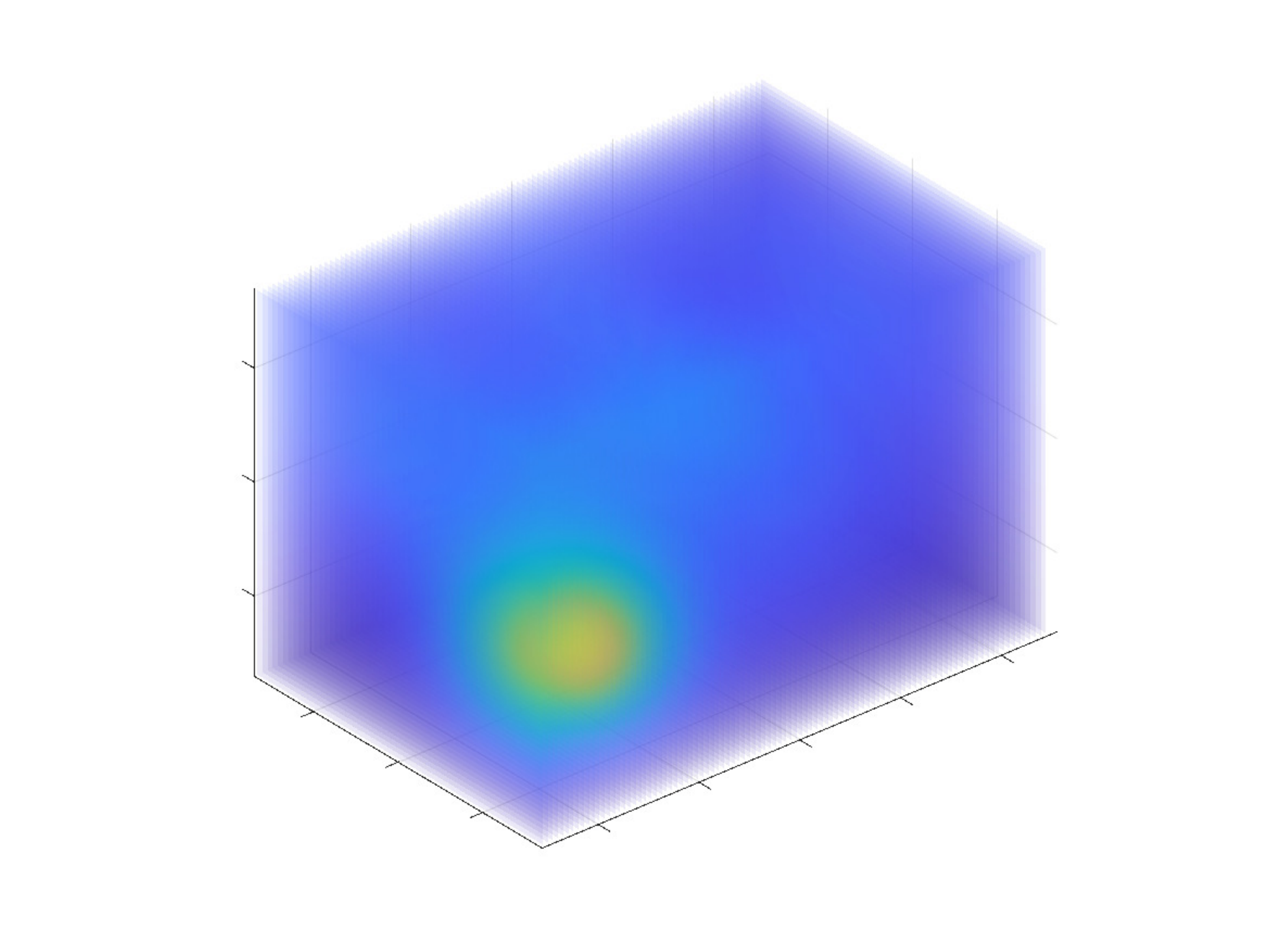}}
\put(0,160){\includegraphics[width=70pt]{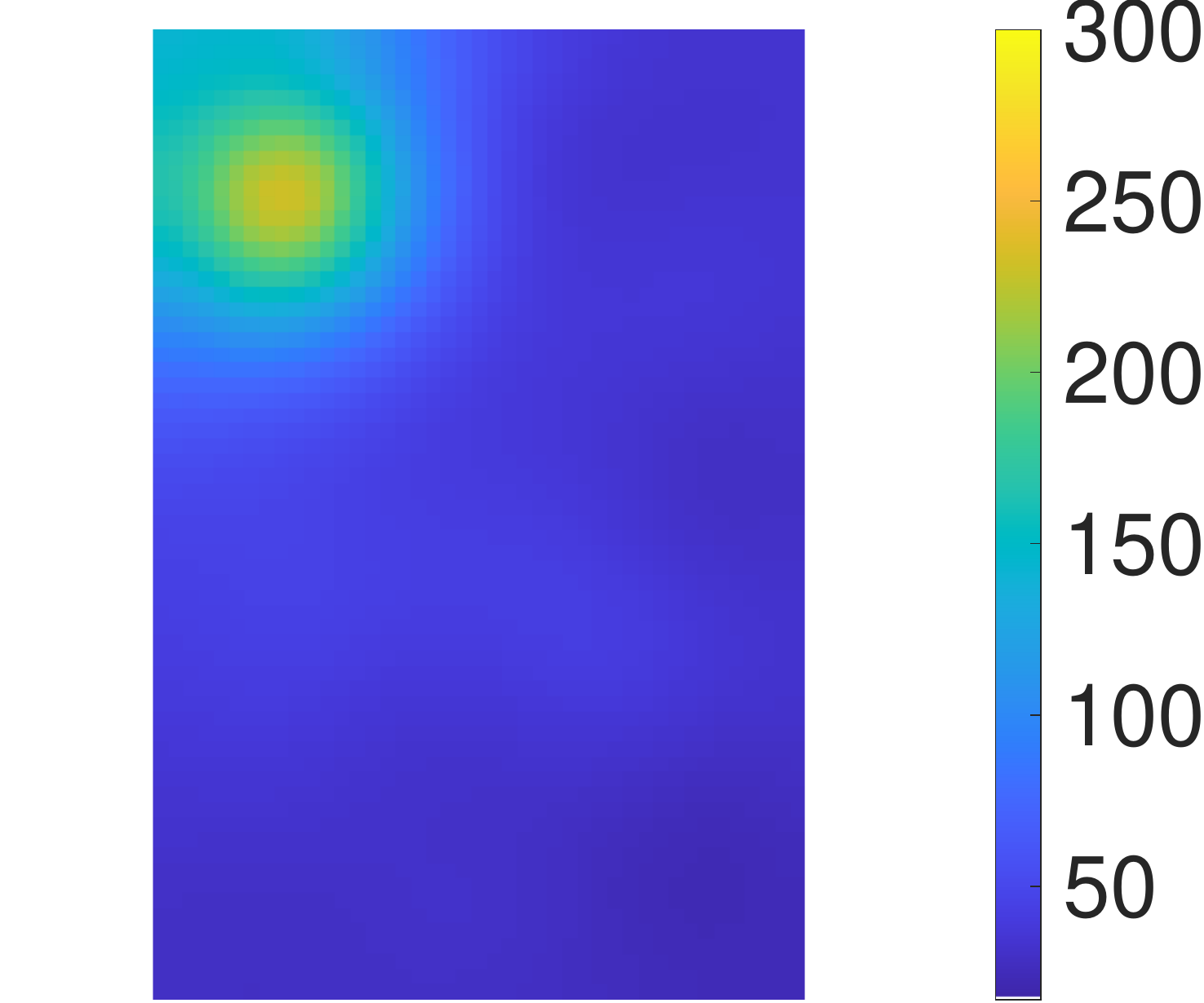}}

% %----------------- 
% t0
\put(290,80){\includegraphics[width=90pt]{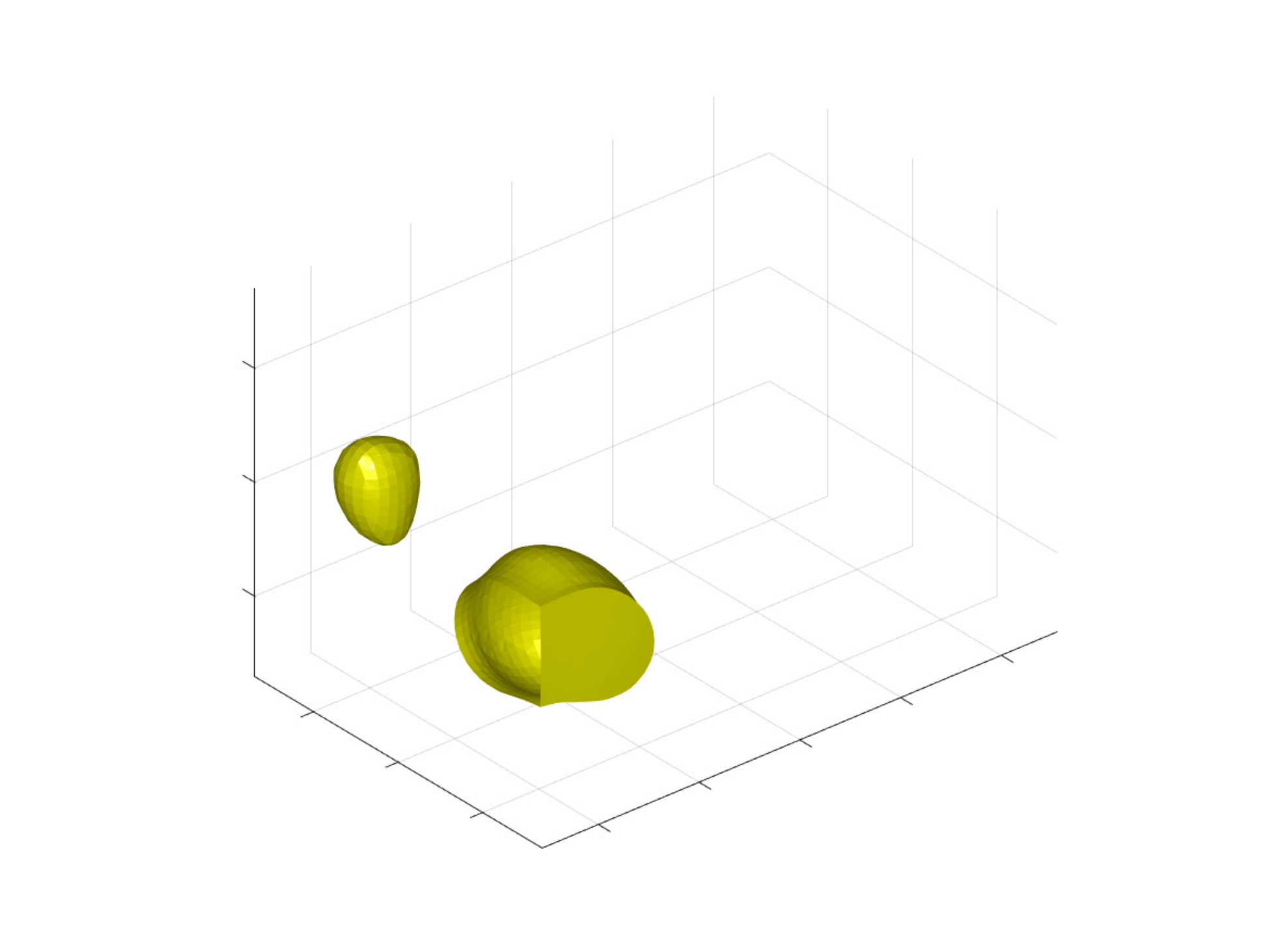}}
\put(350,80){\includegraphics[width=100pt]{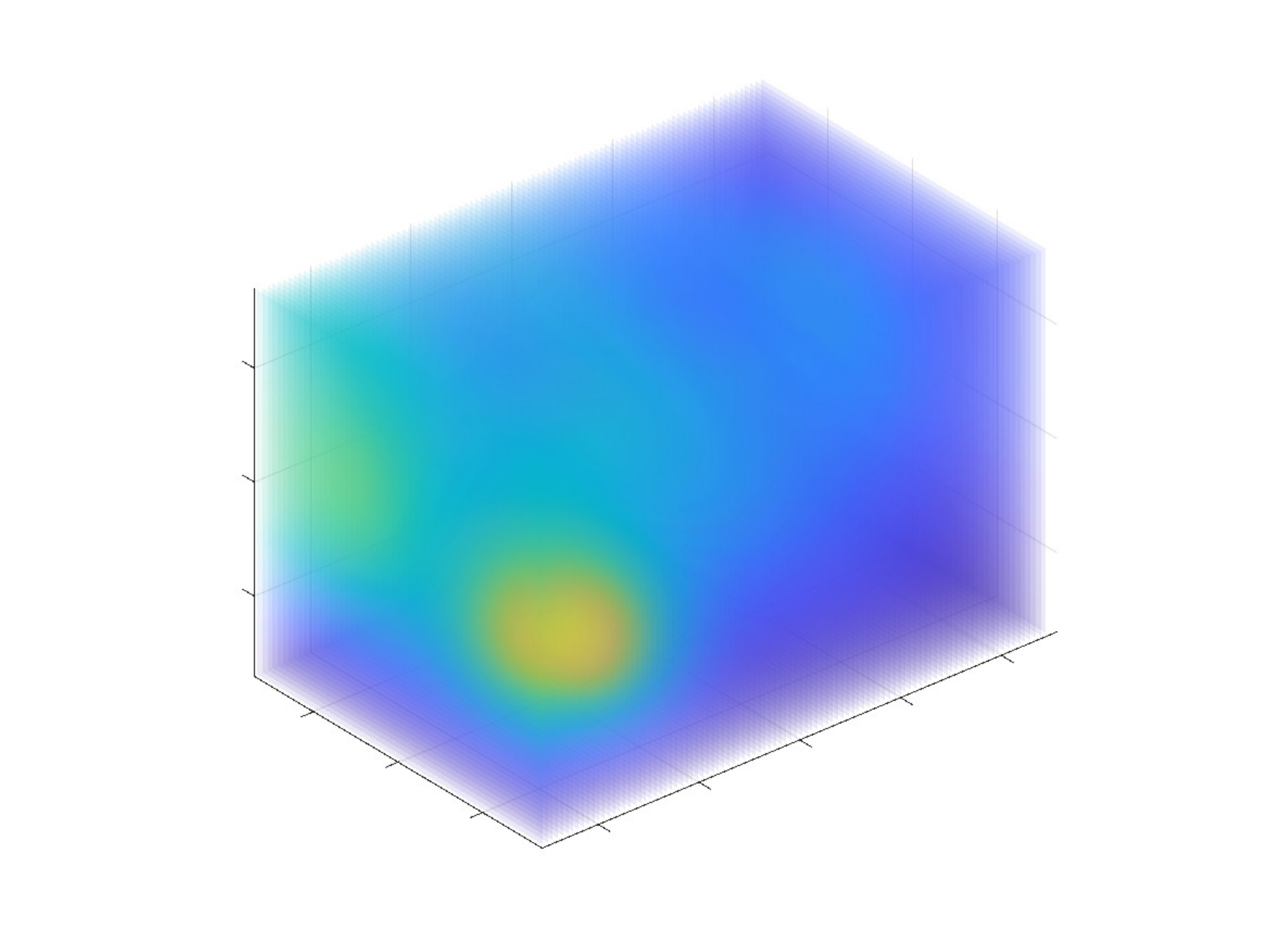}}
\put(230,80){\includegraphics[width=70pt]{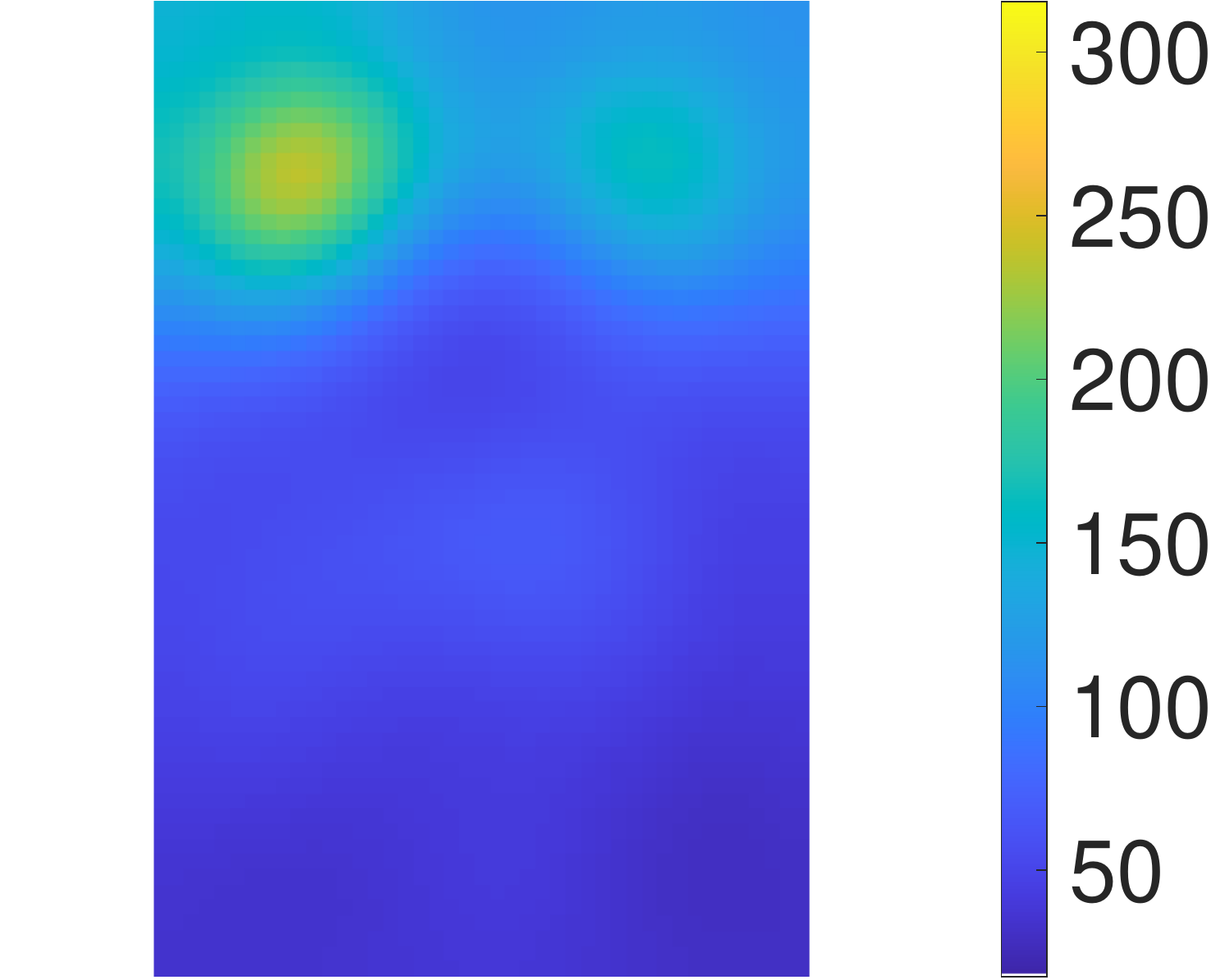}}

\put(60,80){\includegraphics[width=90pt]{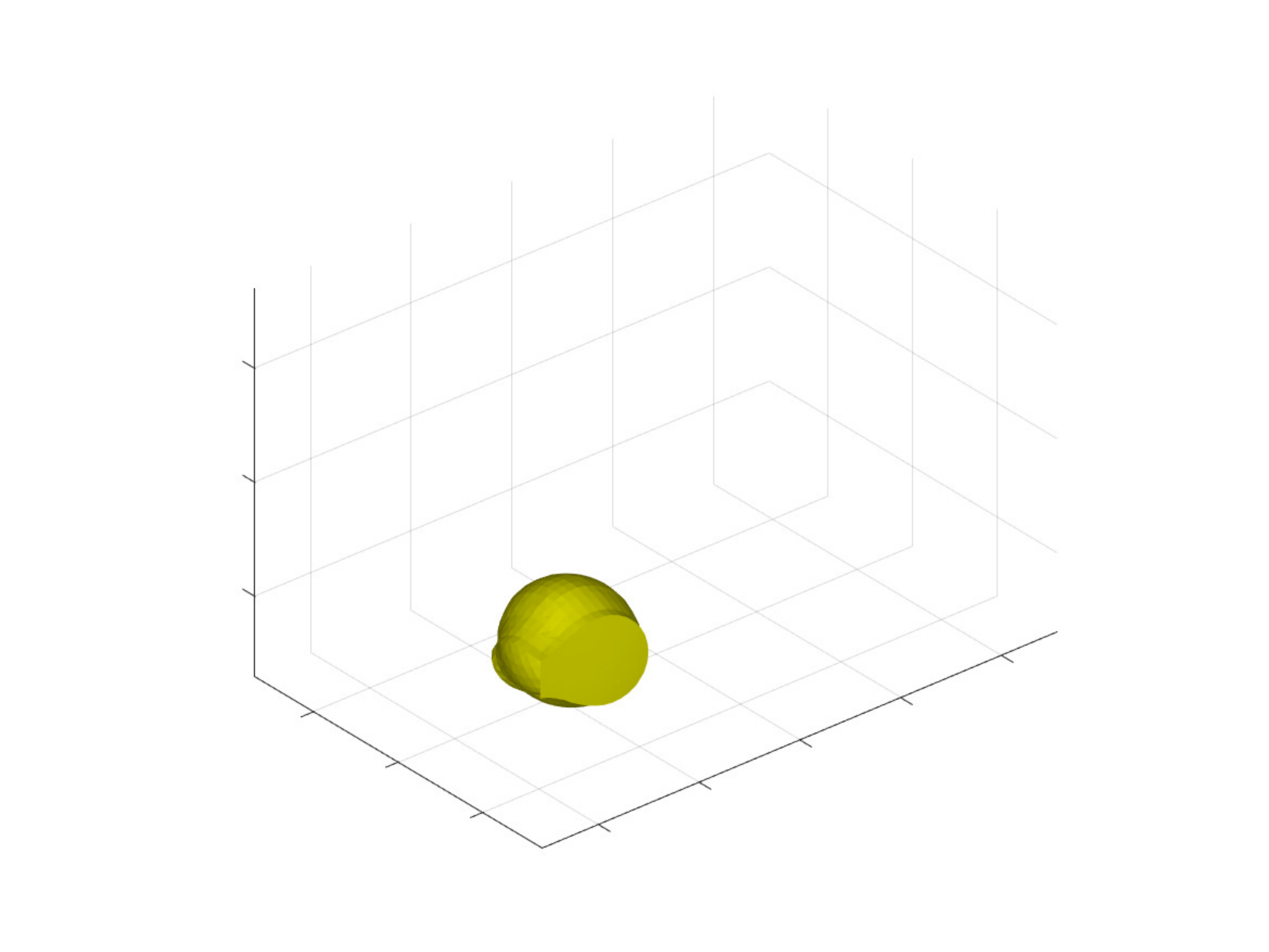}}
\put(120,80){\includegraphics[width=100pt]{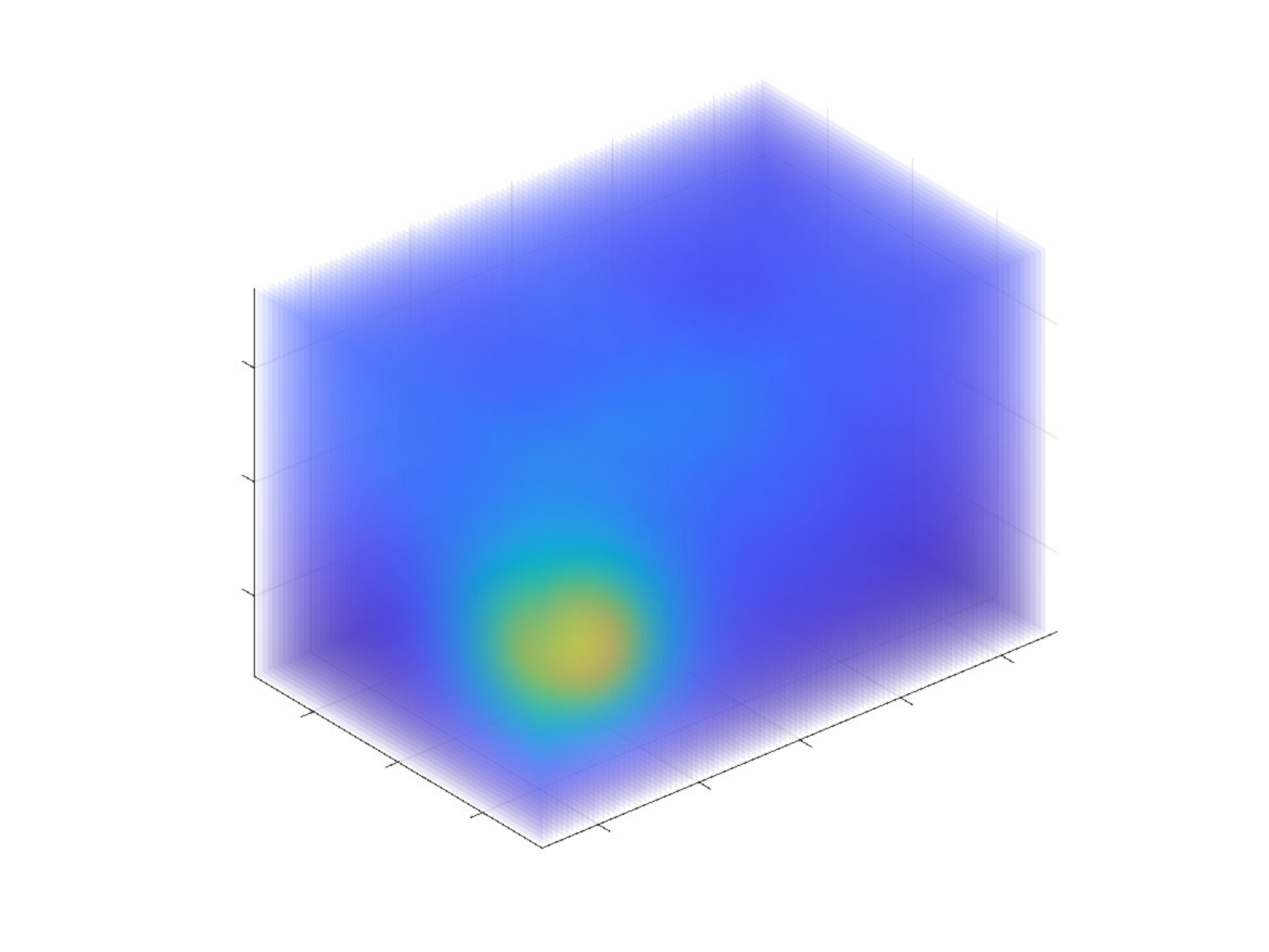}}
\put(0,80){\includegraphics[width=70pt]{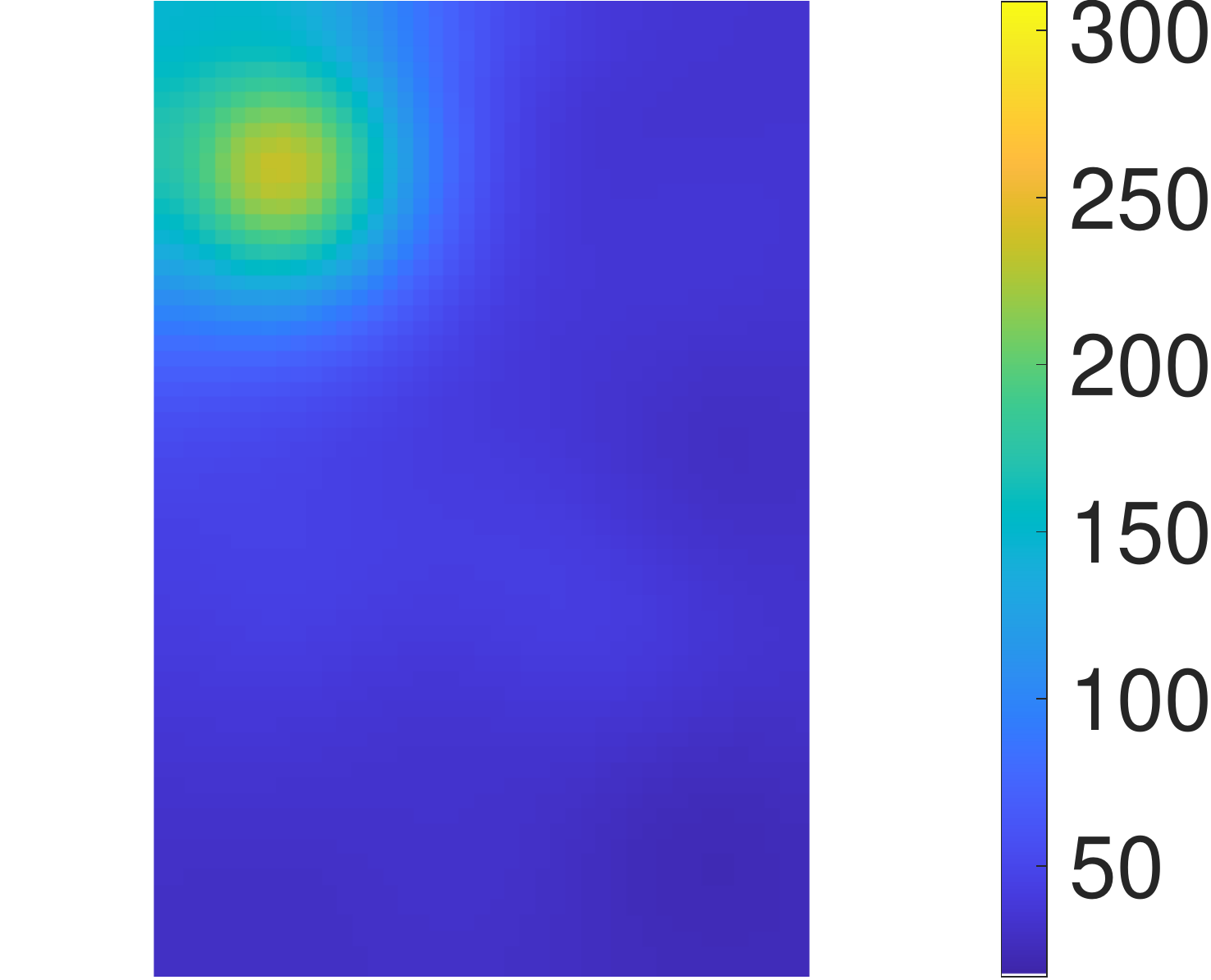}}

%----------------- 
% TV
\put(290,0){\includegraphics[width=90pt]{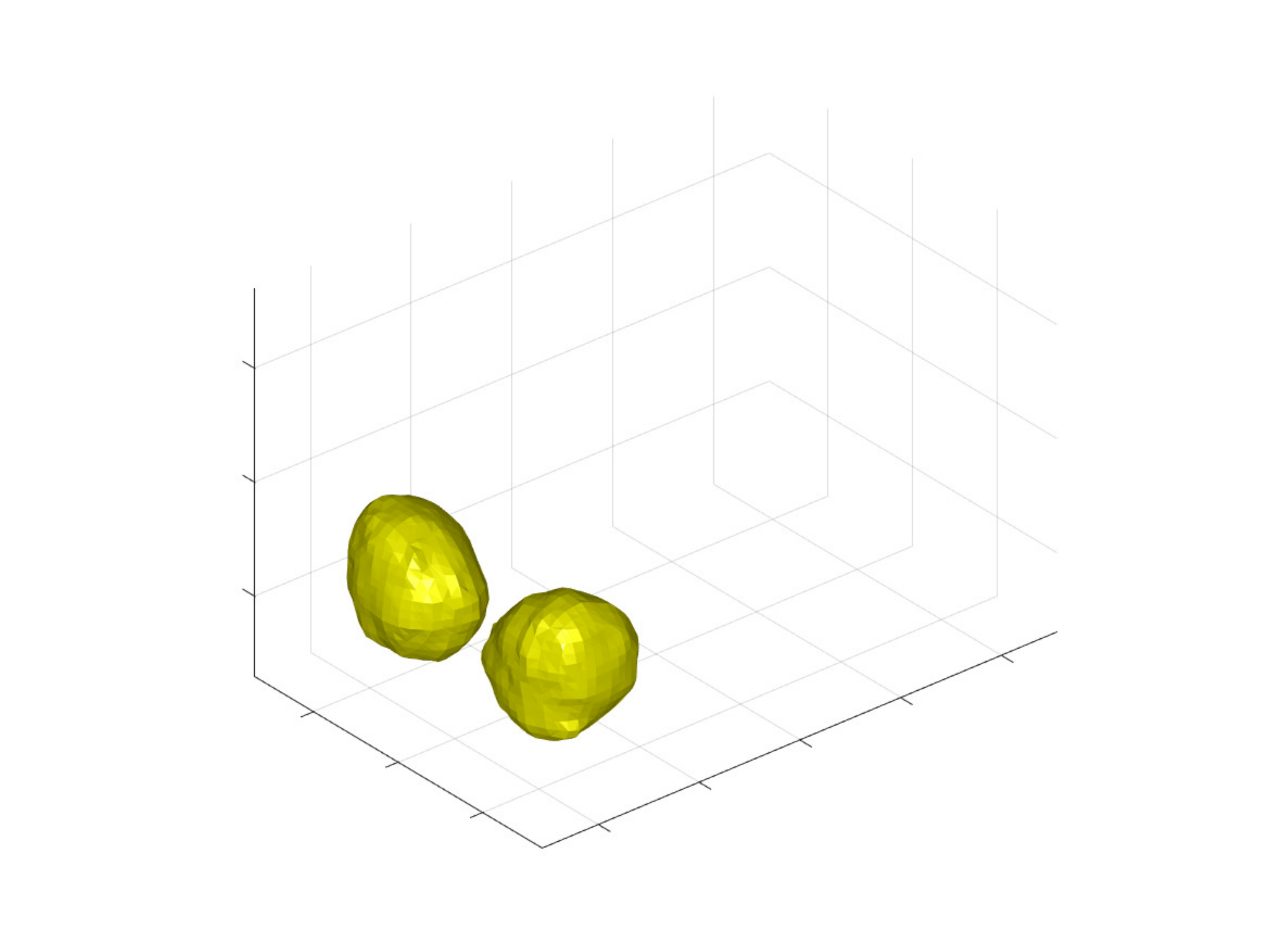}}
\put(350,0){\includegraphics[width=100pt]{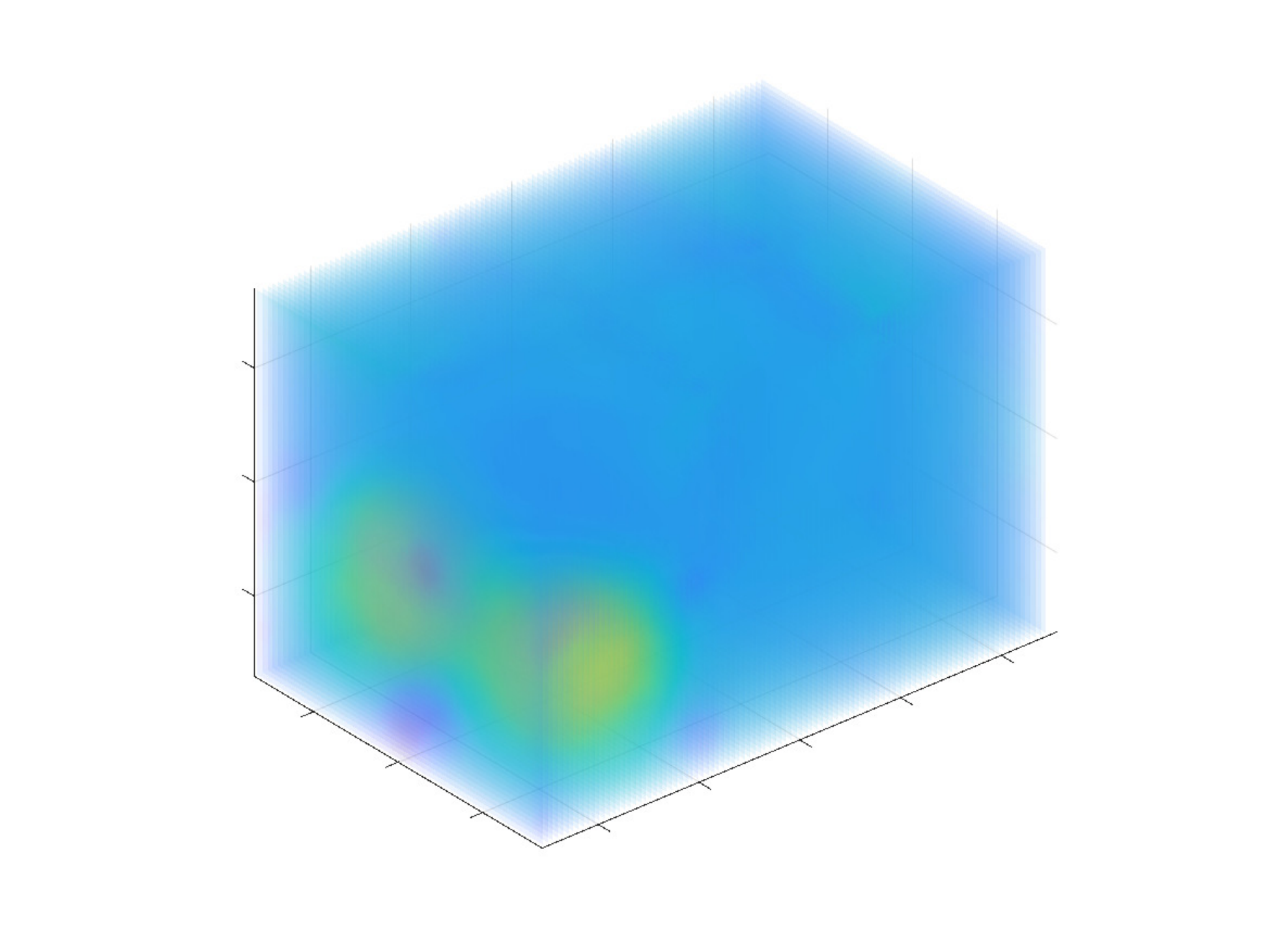}}
\put(230,0){\includegraphics[width=70pt]{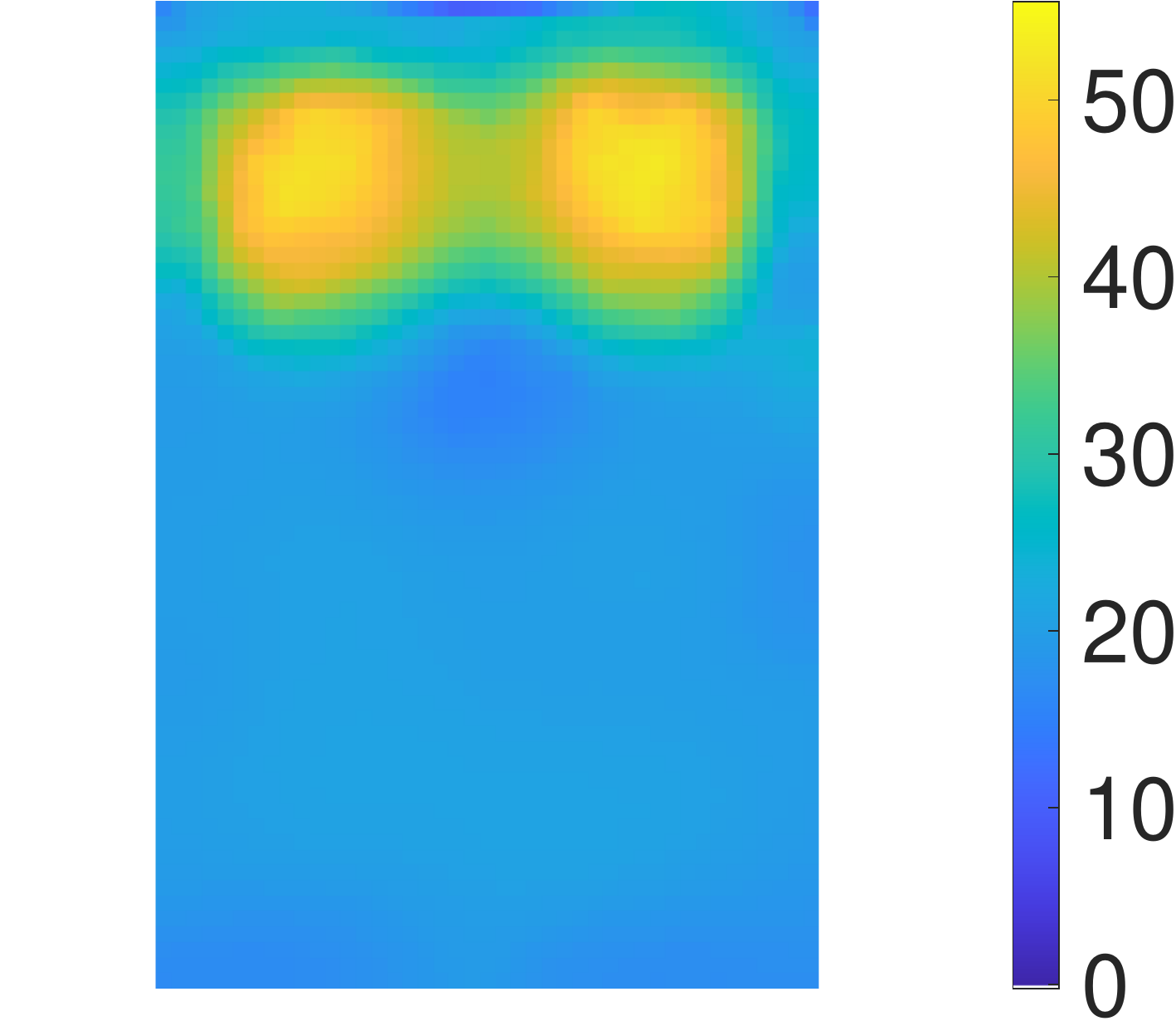}}

\put(60,0){\includegraphics[width=90pt]{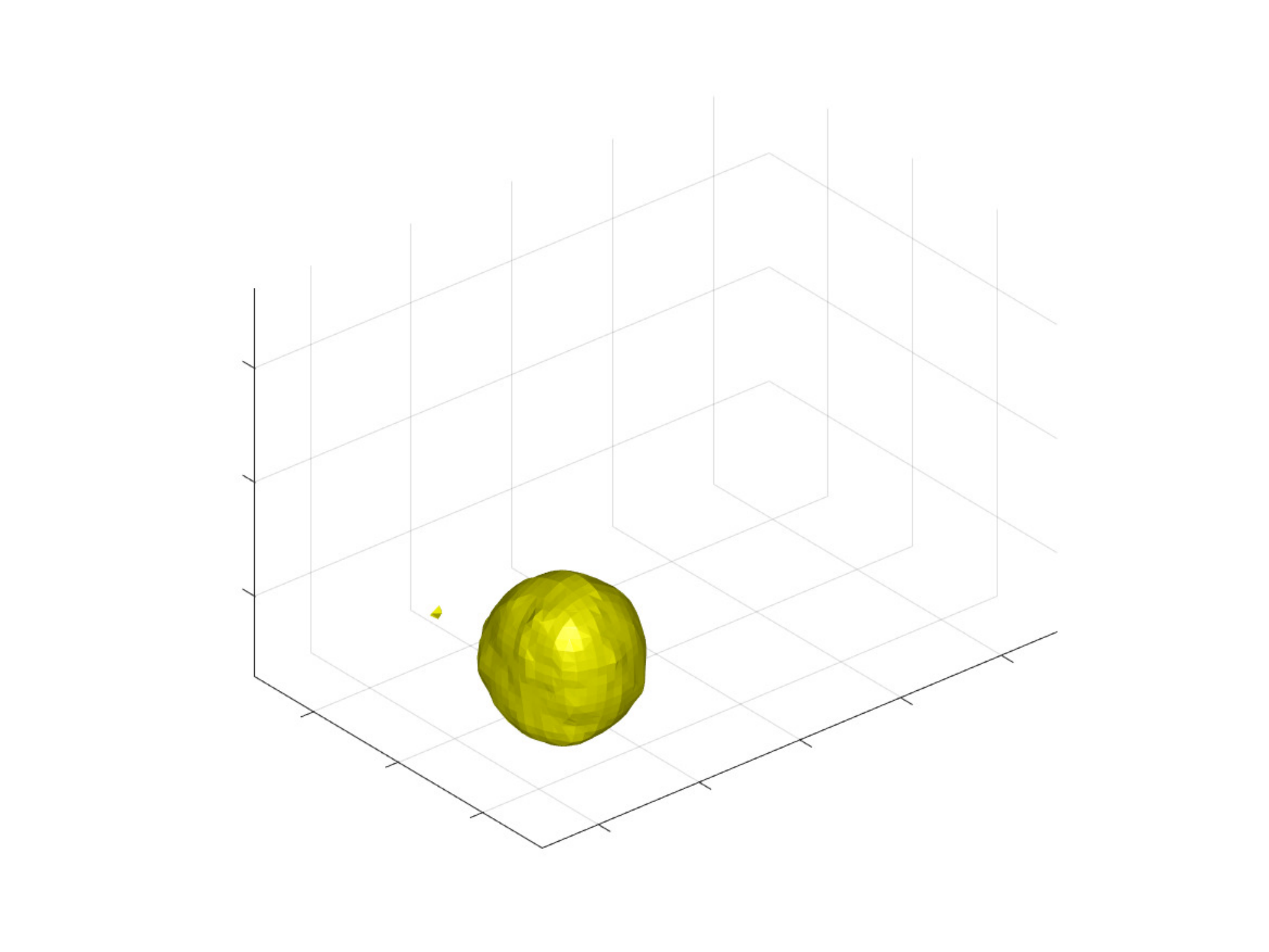}}
\put(120,0){\includegraphics[width=100pt]{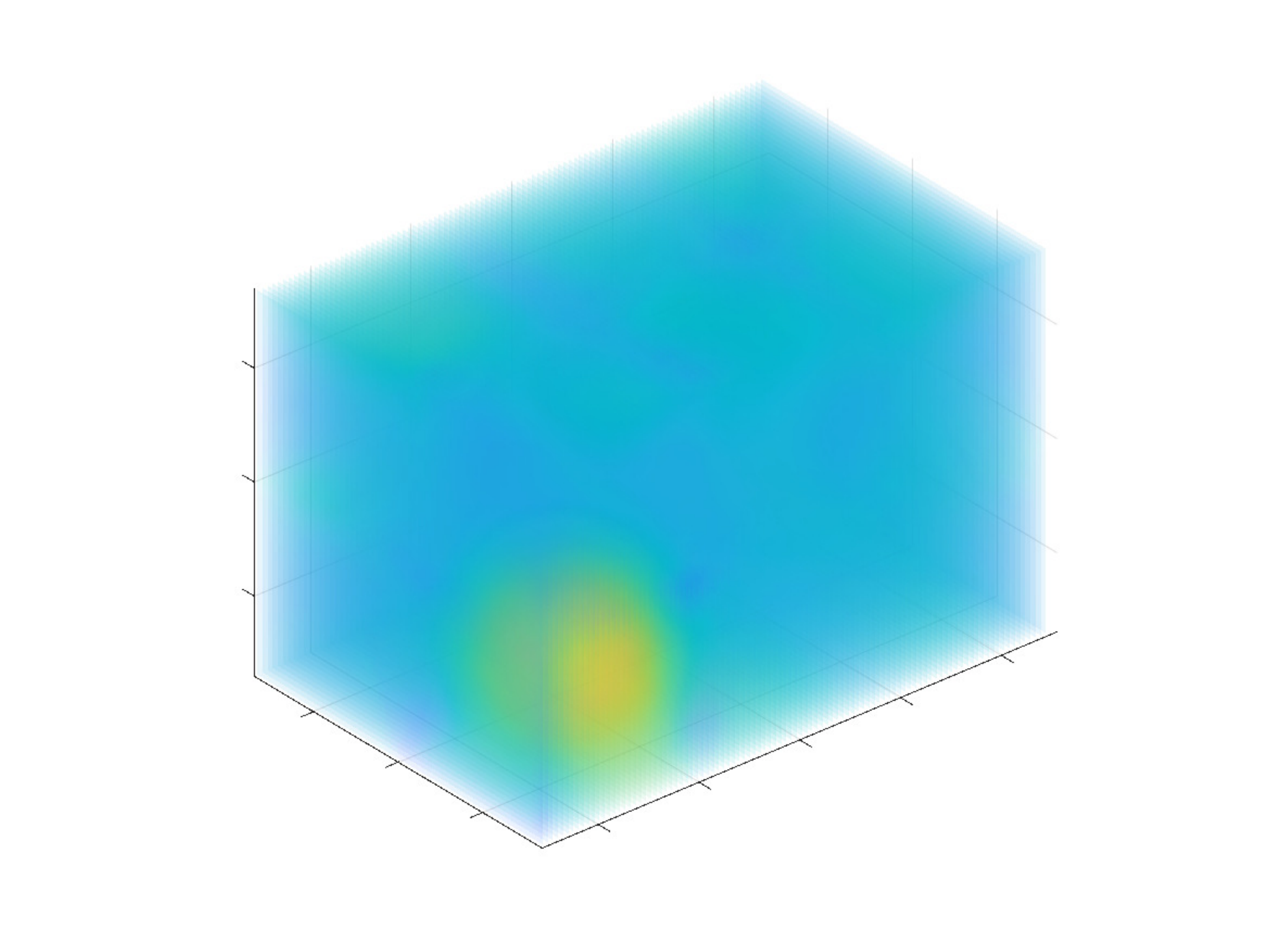}}
\put(0,0){\includegraphics[width=70pt]{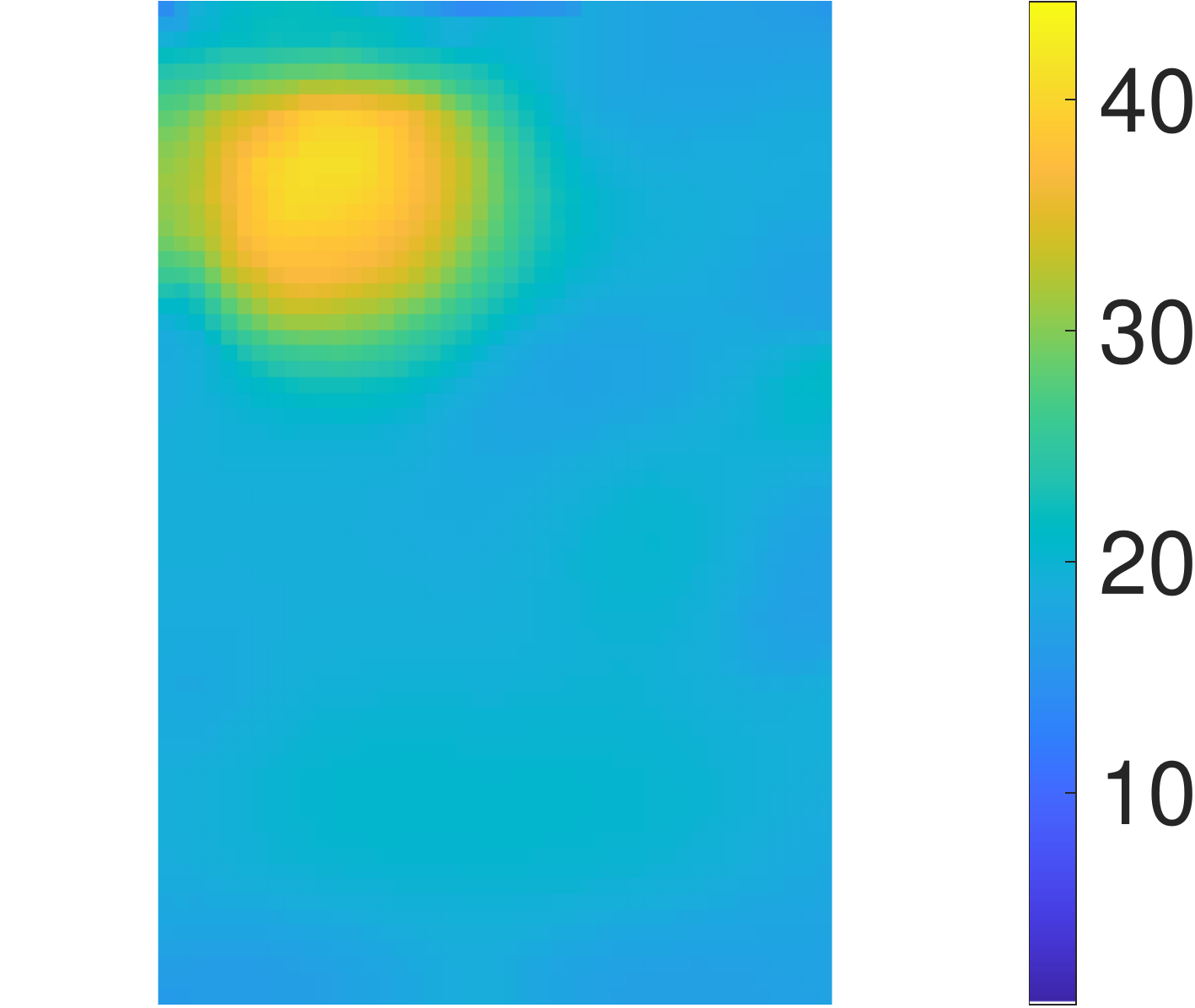}}

\put(3,0){\line(0,1){400}}
%\put(145,0){\line(0,1){400}}
\put(220,0){\line(0,1){400}}
%\put(445,0){\line(0,1){400}}
%----------------- 
% Top labels:
\put(20,385){{{\sc\footnotesize Slice}}}
\put(75,385){{{\sc\footnotesize Isosurface}}}
\put(160,385){{{\sc\footnotesize 3D}}}
\put(250,385){{{\sc\footnotesize Slice}}}
\put(310,385){{{\sc\footnotesize Isosurface}}}
\put(390,385){{{\sc\footnotesize 3D}}}

\put(68,400){{\underline{\sc One Target}}}
\put(298,400){{\underline{\sc Two Targets}}}
%\put(320,395){{\underline{\sc 20x35x25 modeling}}}

%----------------- 
% Side labels:
\put(-10,330){\rotatebox{90}{\sc Truth}}
\put(-10,250){\rotatebox{90}{\sc Cald}}
\put(-10,180){\rotatebox{90}{\sc $\texp$}}
\put(-10,105){\rotatebox{90}{\sc $\tzero$}}
\put(-10,25){\rotatebox{90}{\sc TV}}
%----------------- 

\put(445,0){\line(0,1){400}}

\end{picture}
\caption{\label{fig:abs_correct} {\bf Absolute image} reconstructions comparing the CGO methods to the regularized method with {\it correct domain modeling}.  Slices, isosurfaces, and 3D renderings of the conductivity are shown.  Note the truth targets had a measured conductivity of approx 290 mS/m.}
\end{figure}

%--------------------------------------------------------------------
We first present absolute EIT reconstructions with correct domain modeling, for the one and two targets, and for the four reconstruction methods (Figure~\ref{fig:abs_correct}).  Cross-sections in the center of the bottom row of electrodes of the box are shown on the left; this slice cuts through the center(s) of the physical target(s).  Isosurfaces, shown in the middle column, are extracted from the 3D reconstructions to indicate the locations, size, and number of targets.  The isosurfaces were extracted using {\sc MATLAB}'s \texttt{regionprops3} with threshold segmentation in the range of 60\% to 85\% of the maximum recovered conductivity in the tank.  A uniform threshold across all examples and reconstruction methods was not used as the blurring varied across methods.  Note that since the isosurfaces are dependent on this threshold, they can omit information from the full reconstruction. Therefore, the third column shows a full 3D rendering of conductivity in the domain. This 3D rendering is produced in {\sc MATLAB} by displaying 100 equally-spaced $x_1$-$x_3$ slices set to $0.1$ transparency of each slice.  

% ---------------------------------------------------------------------------------
% TABLE FOR ONE TARGET RECONSTRUCTIONS
% ---------------------------------------------------------------------------------
\begin{table}[ht]
\centering 
     \footnotesize
     \caption{One target {\bf absolute} imaging evaluation metrics across all domain modeling scenarios.}
\begin{tabular}{|c|c||c|c|c|c|}
\hline
                                                &       & Calder\'on & $\texp$  & $\tzero$    & TV \\ \hline
    \multirow{3}[0]{*}{\trev{Scaled LE}}       & correct domain   & \trev{{\bf 0.042}}     & \trev{0.094}   & \trev{0.098}   & \trev{0.043}  \\ \cline{2-6} 
                                               % & 0.2\% electrodes & 1.44     & 2.09   & 2.11   & {\bf 0.82}  \\ \cline{2-6} 
                                                & box 18x27x19     & \trev{0.075}     & \trev{0.131}   & \trev{0.119}   & \trev{{\bf 0.021}}  \\ \cline{2-6} 
                                                & box 20x35x25     & \trev{0.241}     & \trev{0.250}  & \trev{{\bf 0.231}}   & n/a   \\ \hline
    \hline
\multirow{3}[0]{*}{$\begin{array}{c}\text{Max. Target}\\ \text{Conductivity (mS/m)}\end{array}$}& correct domain   & 63.77    & {\bf 300.00} & 308.67 & 44.26 \\ \cline{2-6} 
                                                %& 0.2\% electrodes & 65.62    & {\bf 297.36} & 276.38 & 59.98 \\ \cline{2-6} 
                                                & box 18x27x19     & 76.94    & 303.32 & {\bf 284.67} & 70.14 \\ \cline{2-6} 
                                                & box 20x35x25     & 85.51    & {\bf 305.63} & 260.88 & n/a   \\ \hline
\end{tabular}
  \label{tab:1Targ_errs}
\end{table}
% ---------------------------------------------------------------------------------

% ---------------------------------------------------------------------------------
% TABLE FOR TWO TARGET RECONSTRUCTIONS
% ---------------------------------------------------------------------------------
\begin{table}[ht]
  \centering
  \footnotesize
  \caption{Two target {\bf absolute} imaging evaluation metrics across all domain modeling scenarios. Note: Target~1 is the same target as in the one-target case.}
    \begin{tabular}{|c|c|c||c|c|c|c|}
    \hline
                                                \multicolumn{1}{|r|}{} &       &       & Calder\'on &$\mathbf{t}^{\mbox{\tiny{\textbf{exp}}}}$ &  \multicolumn{1}{c|}{$\mathbf{t}^{\mbox{\tiny{\textbf{0}}}}$} & \multicolumn{1}{c|}{TV} \\
    \hline
    \hline
    \multirow{6}[0]{*}{\trev{Scaled LE}}         & \multirow{3}{*}{Target 1} & correct domain   & \trev{0.056}     & \trev{0.098}   & \trev{0.105}   & \trev{{\bf 0.034}}  \\ \cline{3-7} 
                                                %&                           & 0.2\% electrodes & 2.13     & 2.38   & 2.59   & {\bf 0.78}  \\ \cline{3-7} 
                                                &                           & box 18x27x19     & \trev{0.103}     & \trev{0.135}   & \trev{0.133}   & \trev{{\bf 0.021}}  \\ \cline{3-7} 
                                                &                           & box 20x35x25     & \trev{0.268}     & \trev{0.267}   & \trev{{\bf0.257}}   & n/a   \\ \cline{2-7} 
   \clineB{2-7}{2.5}                                             & \multirow{3}{*}{Target 2} & correct domain   & \trev{0.062}     &\trev{0.195}   & \trev{0.180}   & \trev{{\bf 0.031}}  \\ \cline{3-7} 
                                            %    &                           & 0.2\% electrodes & 2.40     & 5.46   & 5.67   & {\bf 0.64}  \\ \cline{3-7} 
                                                &                           & box 18x27x19     & \trev{0.120}     & \trev{0.191}   & \trev{0.214}   & \trev{{\bf 0.057}}  \\ \cline{3-7} 
                                                &                           & box 20x35x25     & \trev{{\bf 0.277}}     & \trev{0.311}   & \trev{0.295 }  & n/a   \\ \hline
\hline
    \multirow{7}[0]{*}{$\begin{array}{c}\text{Max. Target}\\\text{Conductivity (mS/m)}\end{array}$} & \multirow{3}{*}{Target 1} & correct domain   & 55.61    & {\bf 310.95} & 315.50 & 55.57 \\ \cline{3-7} 
                                            %    &                           & 0.2\% electrodes & 58.22    & 284.78 & {\bf 285.95} & 65.17 \\ \cline{3-7} 
                                                &                           & box 18x27x19     & 73.53    & 216.29 & {\bf 322.39} & 68.29 \\ \cline{3-7} 
                                                &                           & box 20x35x25     & 80.00    & 334.58 & {\bf 264.72} & n/a   \\ \cline{2-7} 
       \clineB{2-7}{2.5}                                         & \multirow{3}{*}{Target 2} & correct domain   & 60.47    & {\bf 234.85} & 226.35 & 53.94 \\ \cline{3-7} 
                                            %    &                           & 0.2\% electrodes & 59.35    & 263.72 & {\bf 283.11} & 66.61 \\ \cline{3-7} 
                                                &                           & box 18x27x19     & 67.94    & 153.82 & {\bf 203.84} & 64.22 \\ \cline{3-7} 
                                                &                           & box 20x35x25     & 77.81    & {\bf 285.06} & 244.80 & n/a   \\ \hline
\end{tabular}
  \label{tab:2Targ_errs}%
\end{table}
% ---------------------------------------------------------------------------------

% ---------------------------------------------------------------------------------
% TABLE FOR TWO TARGET RECONSTRUCTIONS - DIFFERENCE IMAGES
% ---------------------------------------------------------------------------------
\begin{table}[ht]
  \centering
  \footnotesize
  \caption{Two target {\bf difference} imaging evaluation metrics across all domain modeling scenarios. Note: Target~1 is the same target as in the one-target case.}
    \begin{tabular}{|c|c|c||c|c|c|c|}
    \hline
                                                \multicolumn{1}{|r|}{} &       &       & Calder\'on &$\mathbf{t}^{\mbox{\tiny{\textbf{exp}}}}$ &  \multicolumn{1}{c|}{$\mathbf{t}^{\mbox{\tiny{\textbf{0}}}}$} & \multicolumn{1}{c|}{Linear} \\
    \hline
    \hline
    \multirow{6}[0]{*}{\trev{Scaled LE}}         & \multirow{3}{*}{Target 1} & correct domain   & \trev{{\bf 0.034}}     & \trev{0.095}   & \trev{0.103}   & \trev{0.049}   \\ \cline{3-7} 
                                             %   &                           & 0.2\% electrodes & {\bf 0.83}     & 2.62   & 2.74   & 1.57   \\ \cline{3-7} 
                                                &                           & box 18x27x19     & \trev{{\bf 0.022}}     & \trev{0.116}   & \trev{0.124 }  & \trev{0.065}   \\ \cline{3-7} 
                                                &                           & box 20x35x25     & \trev{{\bf 0.113}}     & \trev{0.242}   & \trev{0.243 }  & n/a    \\ \cline{2-7} 
   \clineB{2-7}{2.5}                                             & \multirow{3}{*}{Target 2} & correct domain   & \trev{0.053}     & \trev{0.100}   & \trev{0.094}   & \trev{{\bf 0.050}}   \\ \cline{3-7} 
                                           %     &                           & 0.2\% electrodes & 1.46     & 2.53   & 2.41   & {\bf 1.21}   \\ \cline{3-7} 
                                                &                           & box 18x27x19     & \trev{0.043}     & \trev{0.116}   & \trev{0.115}   & \trev{{\bf 0.029}}   \\ \cline{3-7} 
                                                &                           & box 20x35x25     & \trev{{\bf 0.112}}     & \trev{0.230}  & \trev{0.233}   & n/a    \\ \hline
\hline
    \multirow{6}[0]{*}{$\begin{array}{c}\text{Max. Target}\\\text{Conductivity (mS/m)}\end{array}$} & \multirow{3}{*}{Target 1} & correct domain   & 57.18    & {\bf 277.83} & 295.63 & 484.83 \\ \cline{3-7} 
                                          %      &                           & 0.2\% electrodes & 57.07    & 295.99 & {\bf 290.70} & 423.95 \\ \cline{3-7} 
                                                &                           & box 18x27x19     & 61.72    & {\bf 272.80} & 282.05 & 227.95 \\ \cline{3-7} 
                                                &                           & box 20x35x25     & 56.30    & {\bf 253.46} & 243.88 & 227.95 \\ \cline{2-7}
       \clineB{2-7}{2.5}                                         & \multirow{3}{*}{Target 2} & correct domain   & 53.84    & {\bf 266.91} & 284.25 & 490.88 \\ \cline{3-7} 
                                         %       &                           & 0.2\% electrodes & 54.43    & 298.07 & {\bf 287.78} & 420.63 \\ \cline{3-7} 
                                                &                           & box 18x27x19     & 59.91    & {\bf 294.17} & 306.30 & 208.05 \\ \cline{3-7} 
                                                &                           & box 20x35x25     & 55.07    & 293.61 & {\bf 282.38} & n/a    \\ \hline
\end{tabular}
  \label{tab:2Targ_errs_diff}%
\end{table}
% ---------------------------------------------------------------------------------

\begin{figure}[ht]
\centering
\begin{picture}(450,410)
\linethickness{.3mm}
%----------------- 
%----------------- 
% Truth
\put(300,320){\includegraphics[width=60pt]{truth_2targ_iso.png}}
%\put(350,320){\includegraphics[width=100pt]{cald_true_1_targ_ABS_3D-eps-converted-to.pdf}}
\put(240,320){\includegraphics[angle=-90,origin=c,width=42pt]{TwoTargs_TopView.jpeg}}

\put(70,320){\includegraphics[width=60pt]{truth_1targ_iso.png}}
%\put(120,240){\includegraphics[width=100pt]{cald_true_2_targ_ABS_3D-eps-converted-to.pdf}}
\put(10,320){\includegraphics[angle=-90,origin=c,width=42pt]{OneTarg_TopView.jpeg}}
% %----------------- 
% cald
\put(290,240){\includegraphics[width=90pt]{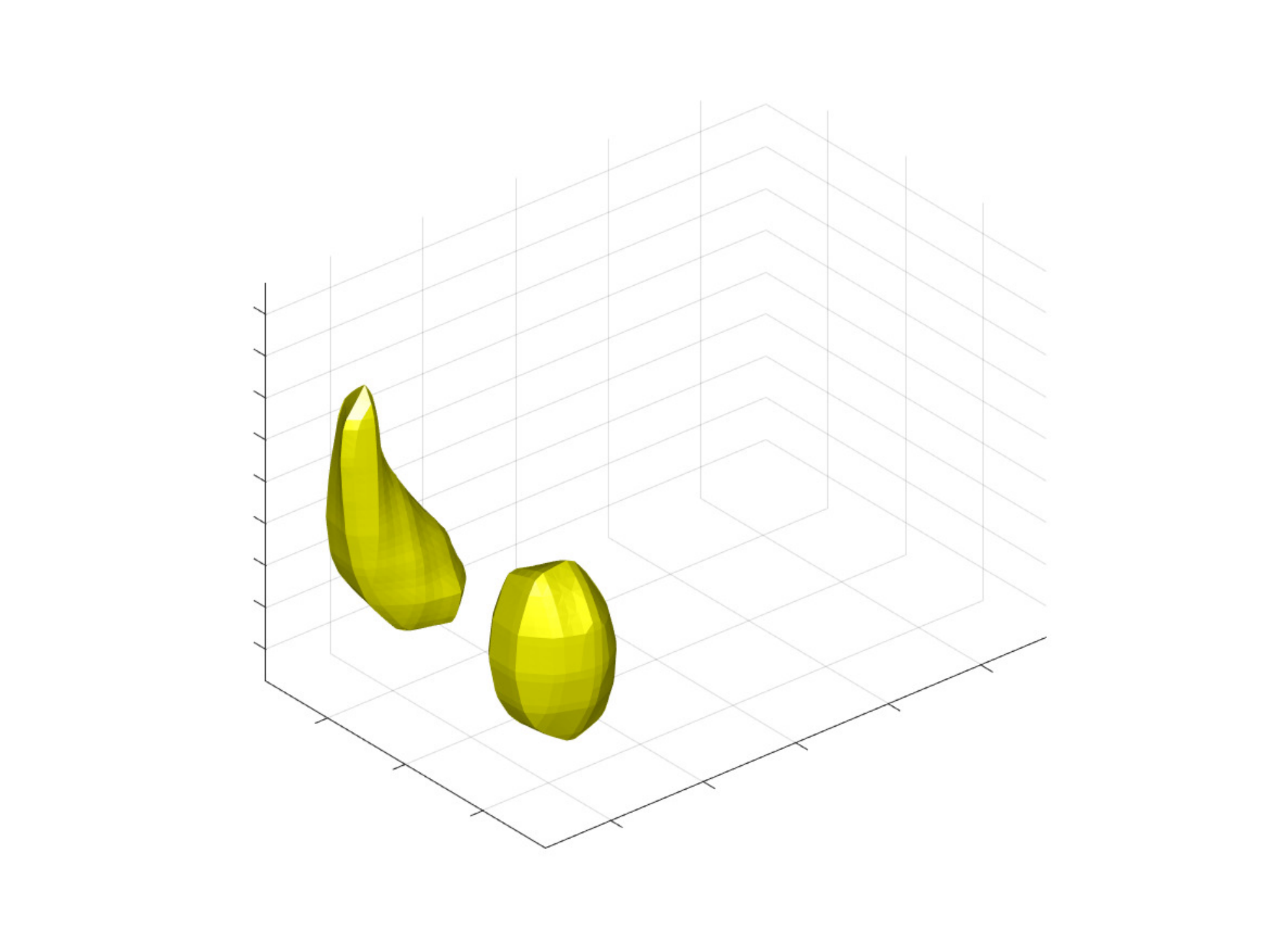}}
\put(350,240){\includegraphics[width=100pt]{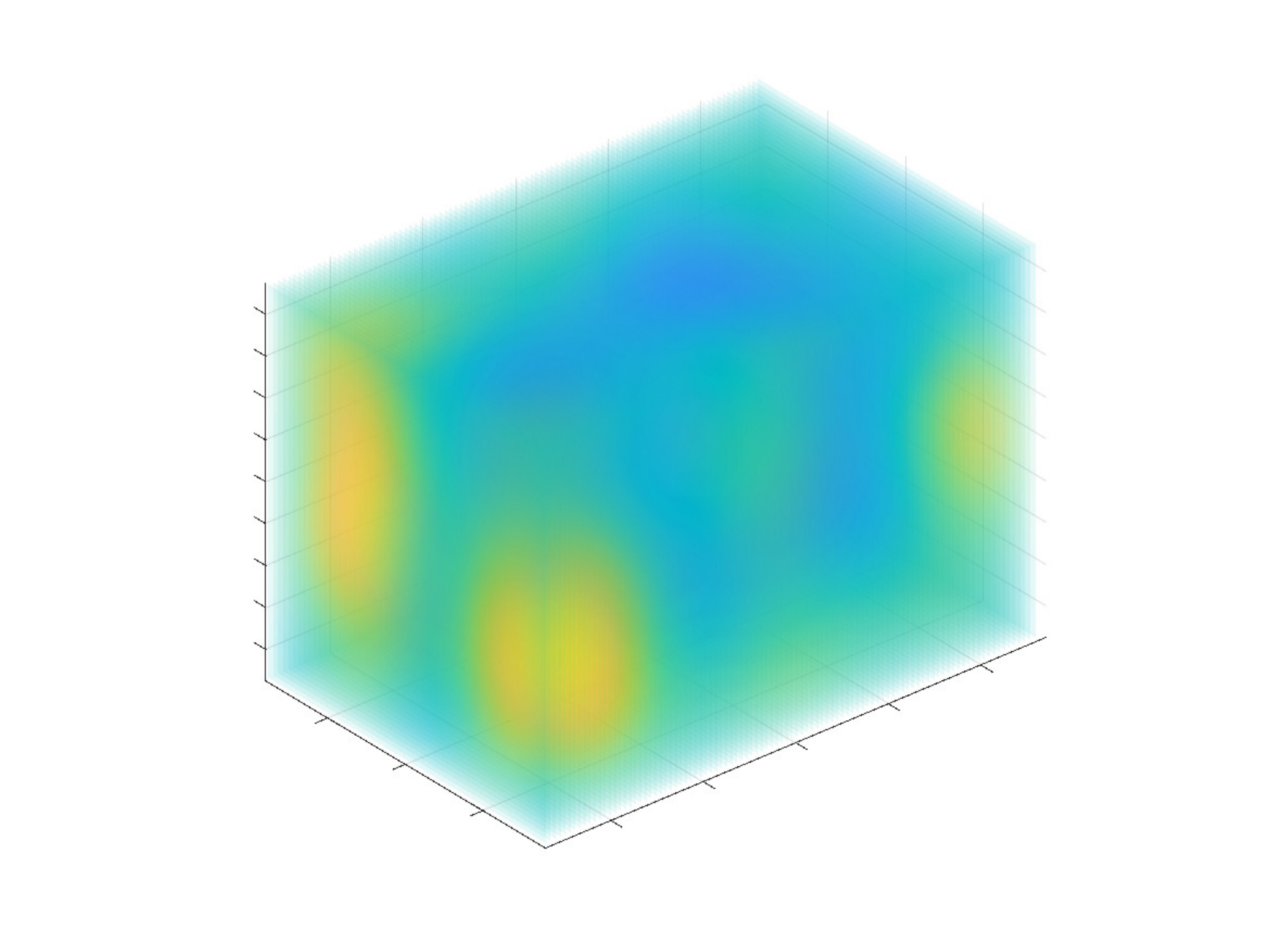}}
\put(230,240){\includegraphics[width=70pt]{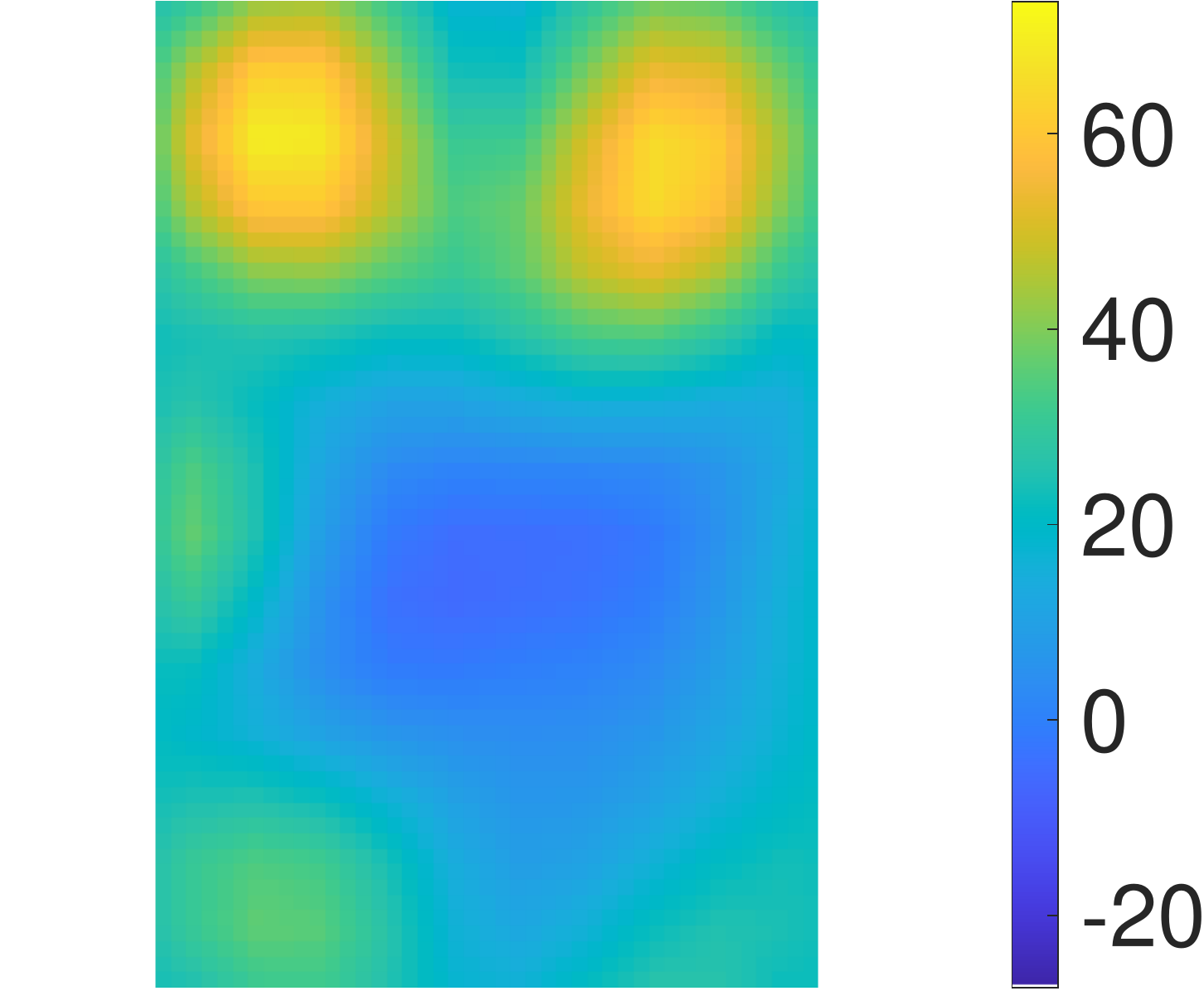}}

\put(60,240){\includegraphics[width=90pt]{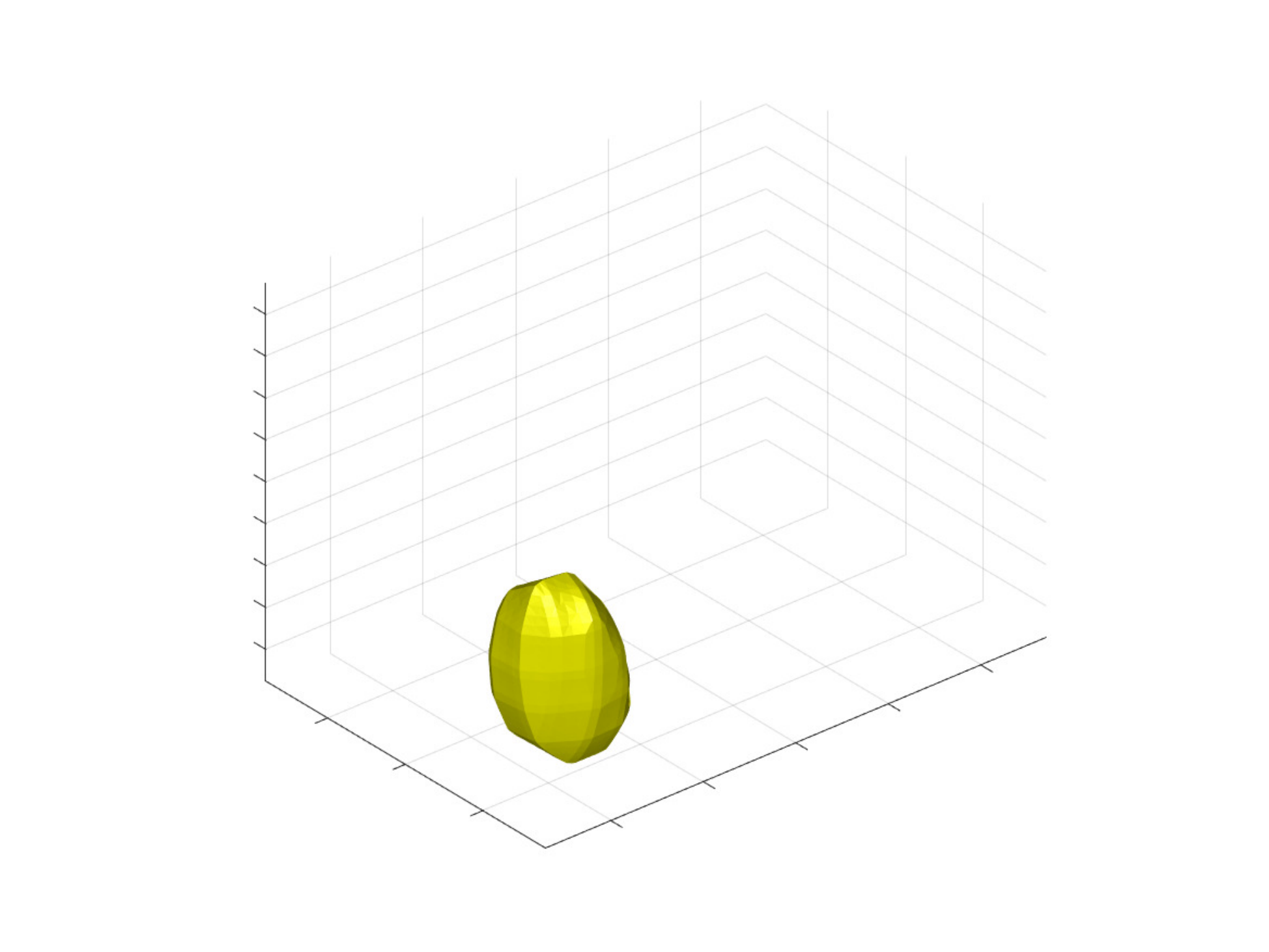}}
\put(120,240){\includegraphics[width=100pt]{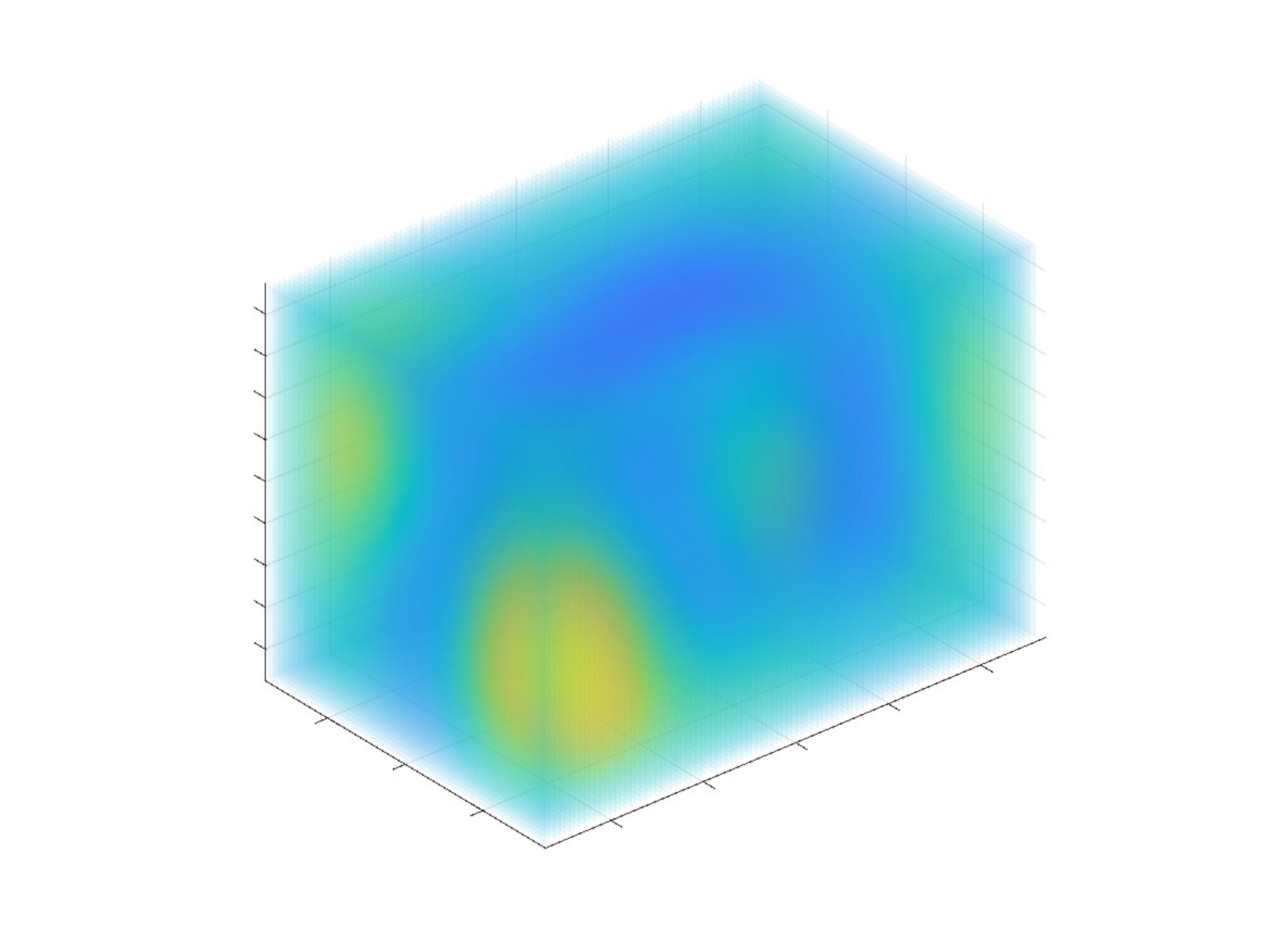}}
\put(0,240){\includegraphics[width=70pt]{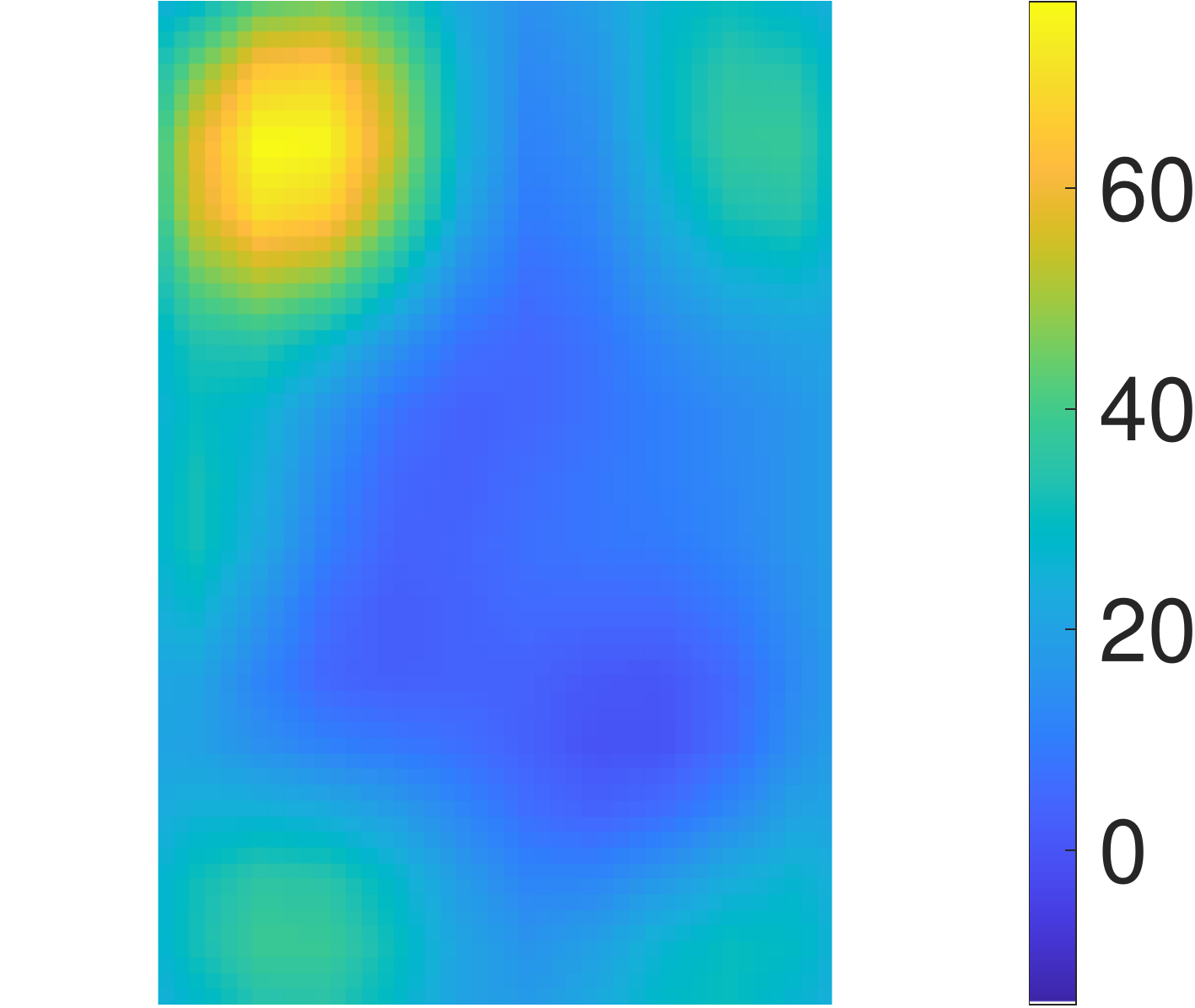}}
% %----------------- 
% texp
\put(290,160){\includegraphics[width=90pt]{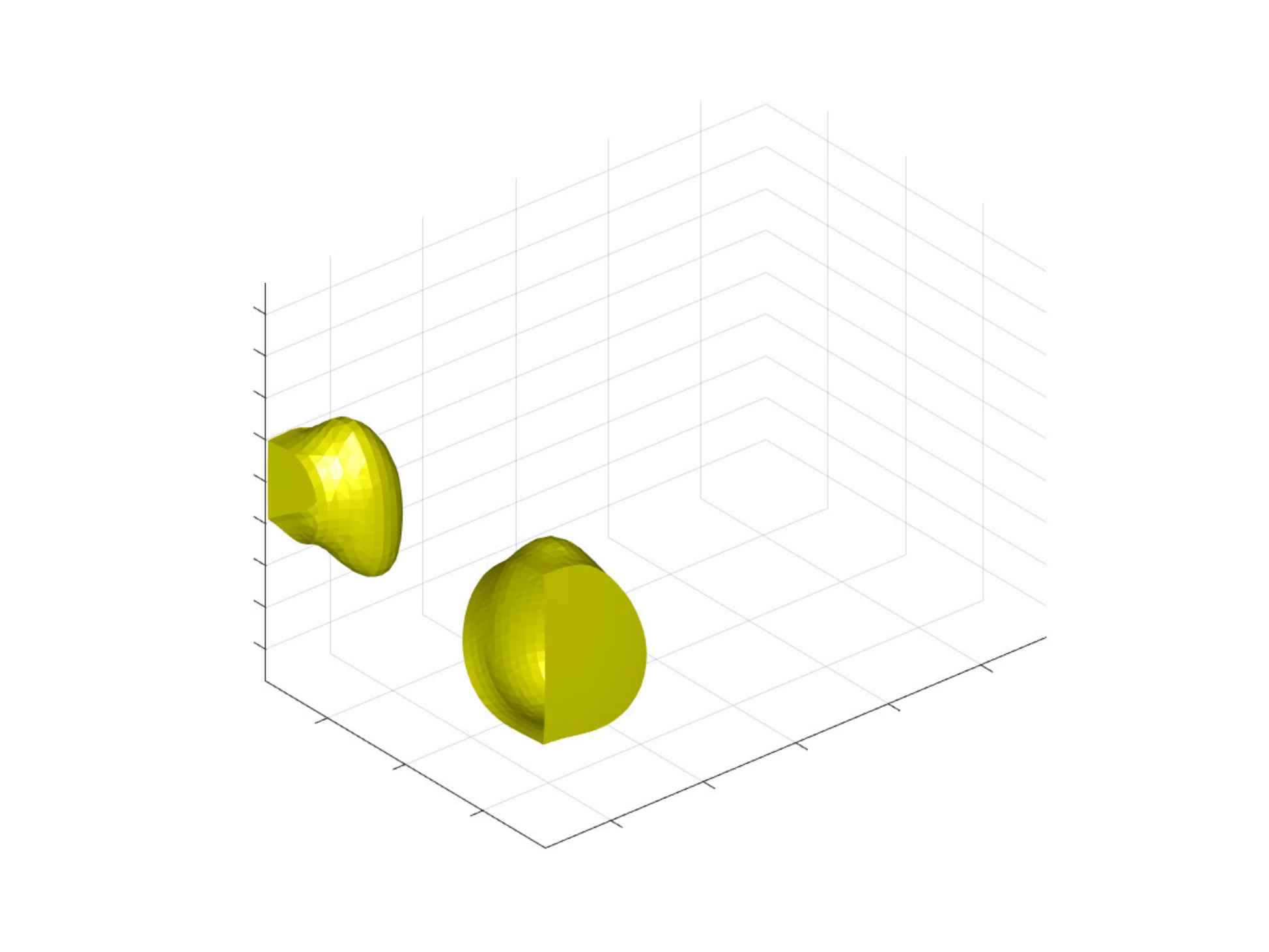}}
\put(350,160){\includegraphics[width=100pt]{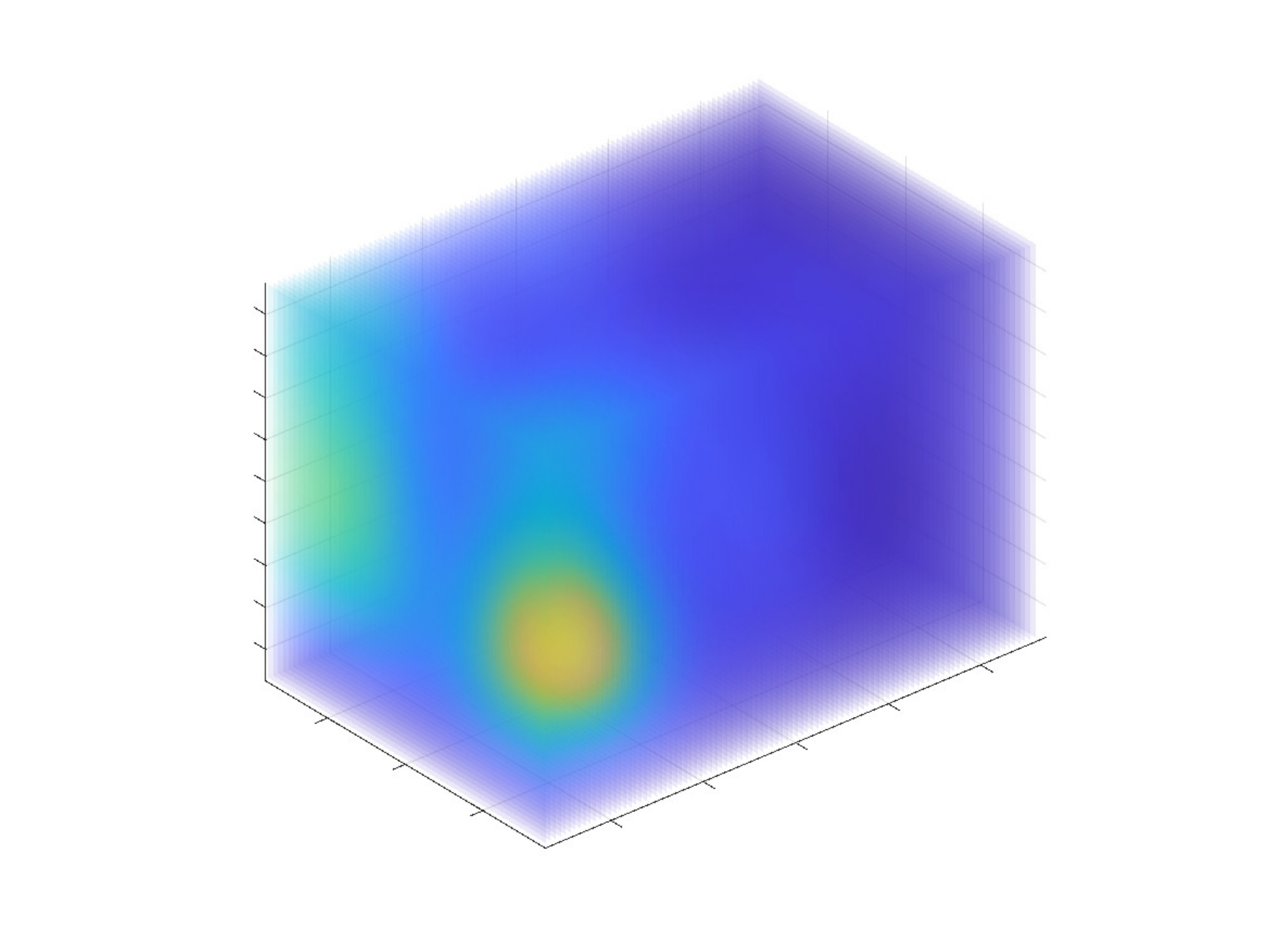}}
\put(230,160){\includegraphics[width=70pt]{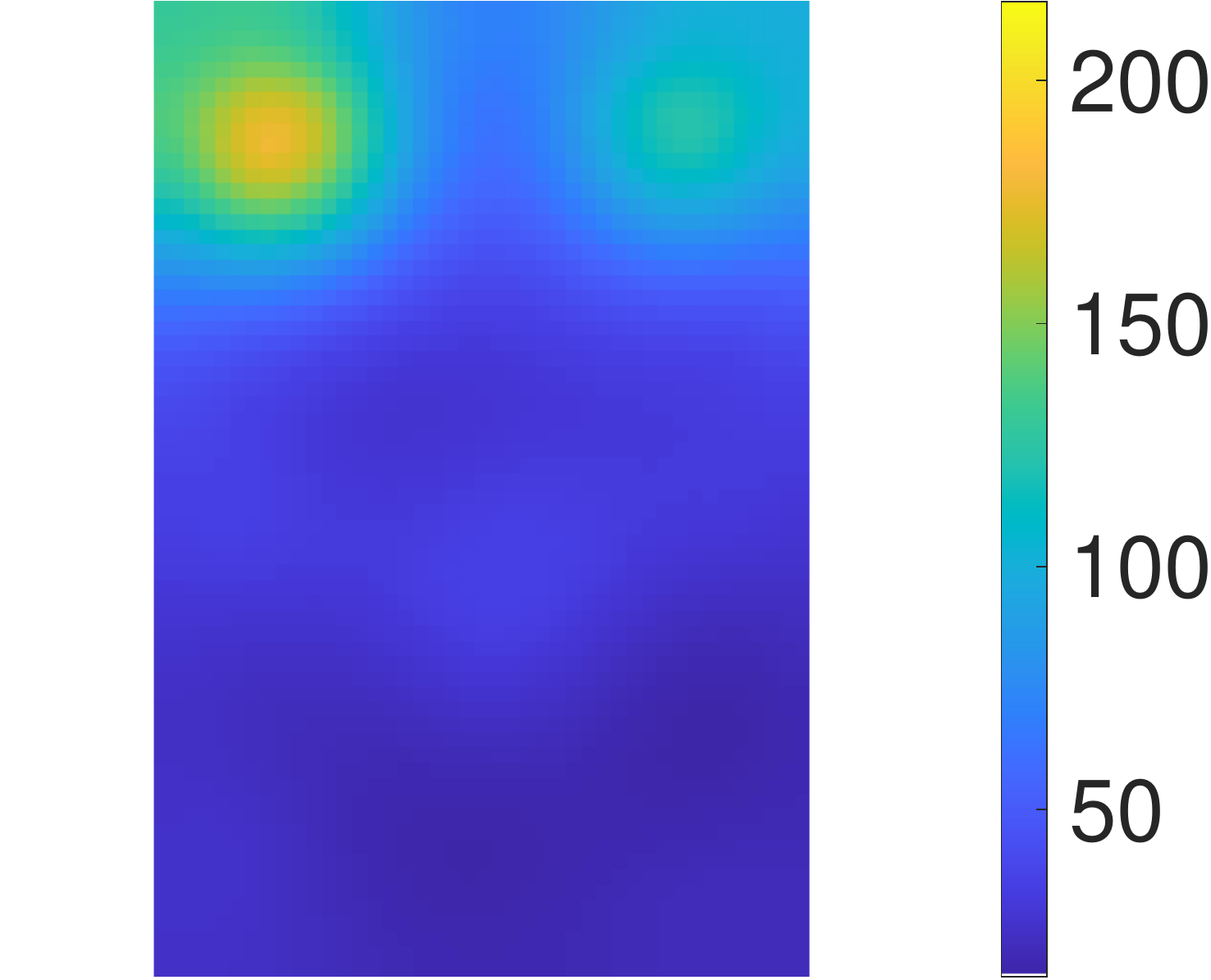}}

\put(60,160){\includegraphics[width=90pt]{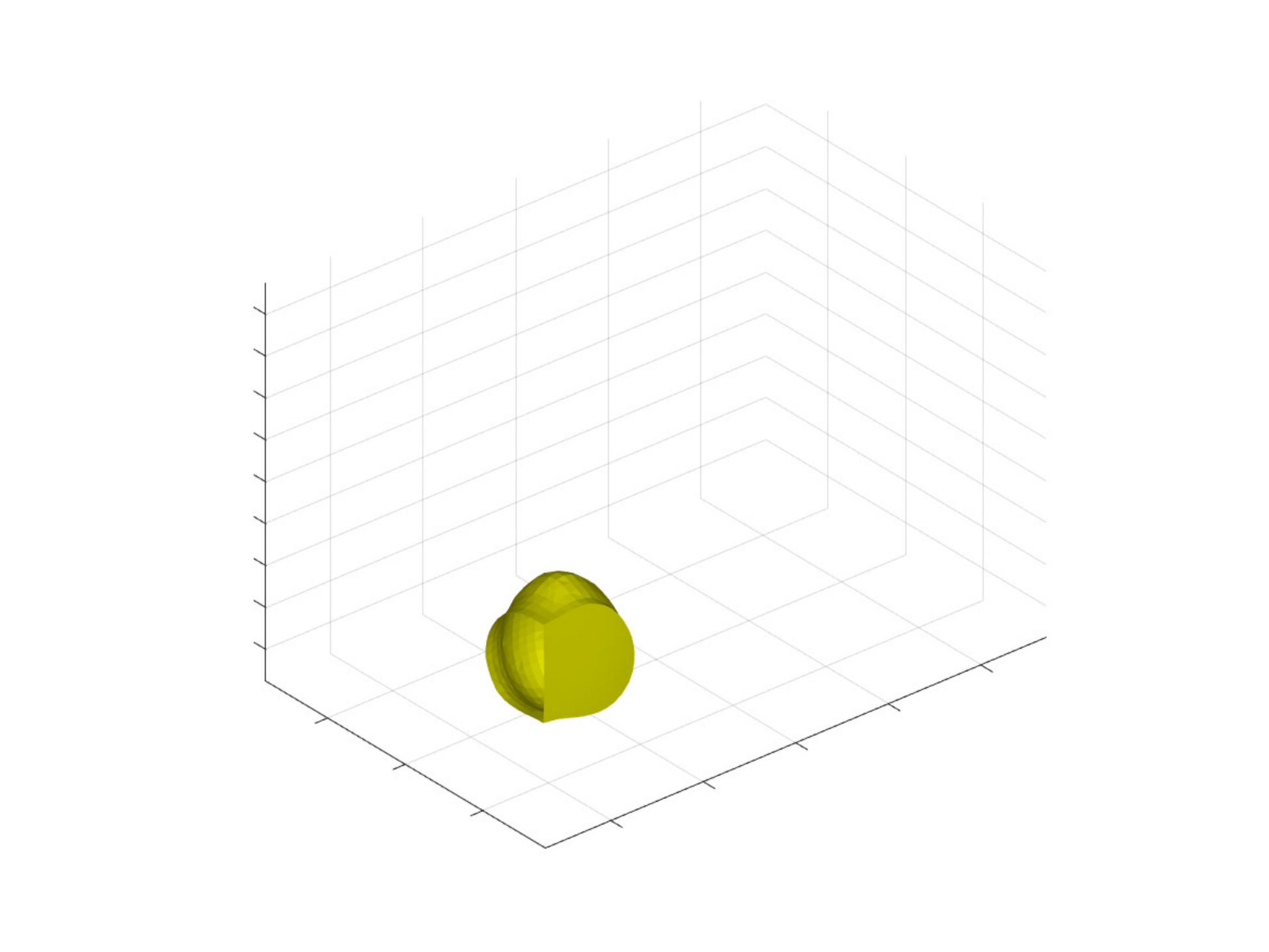}}
\put(120,160){\includegraphics[width=100pt]{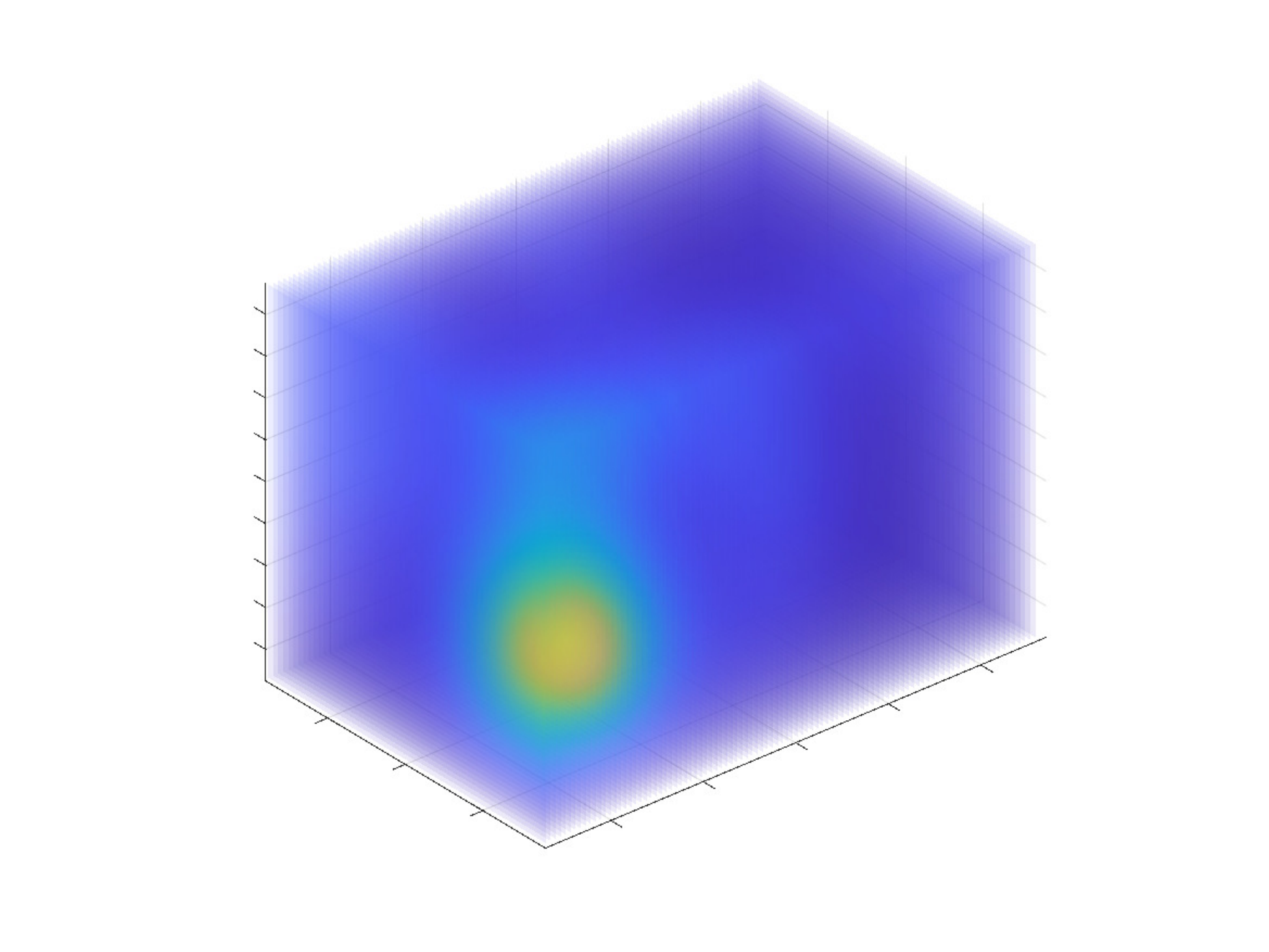}}
\put(0,160){\includegraphics[width=70pt]{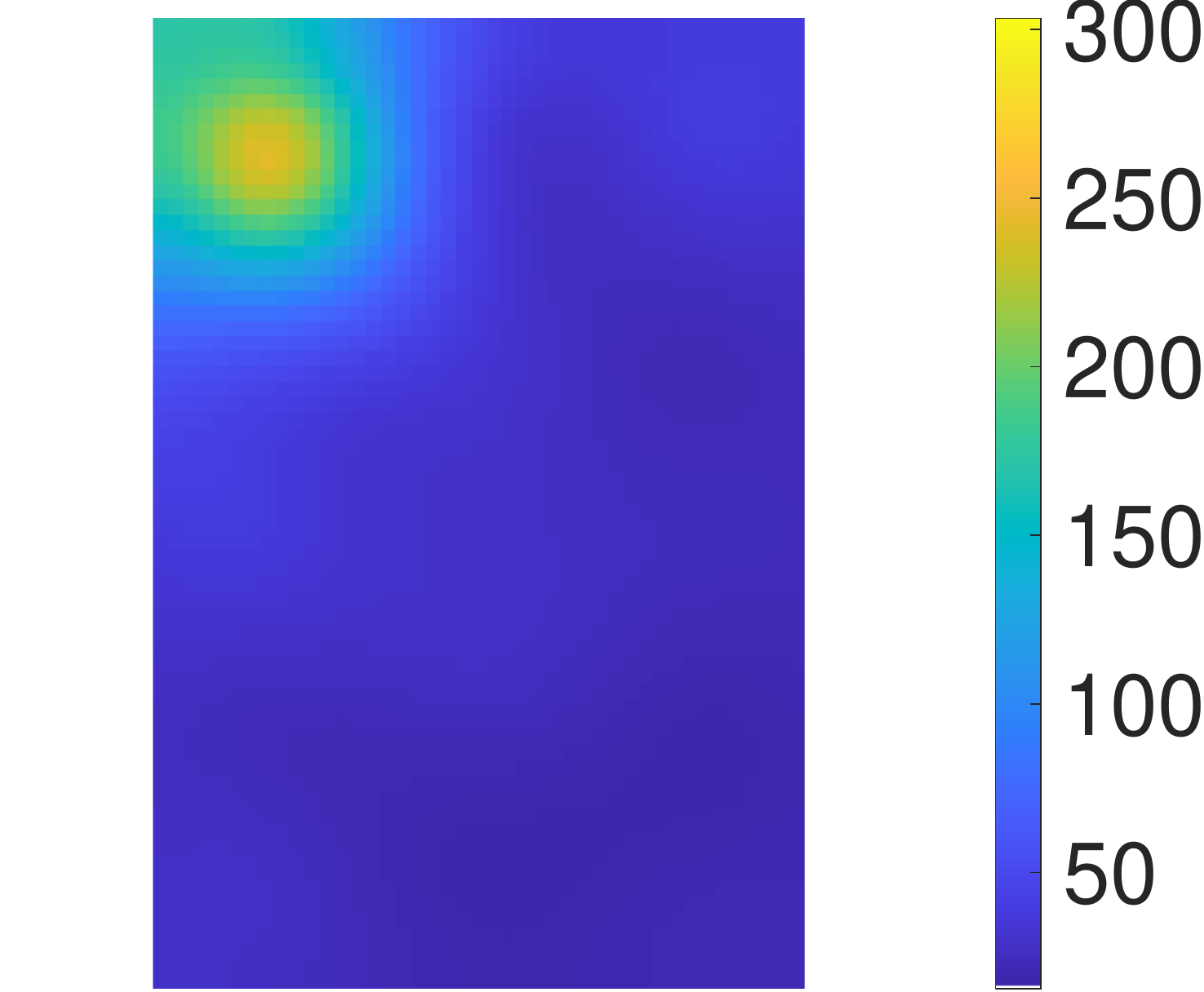}}

% %----------------- 
% t0
\put(290,80){\includegraphics[width=90pt]{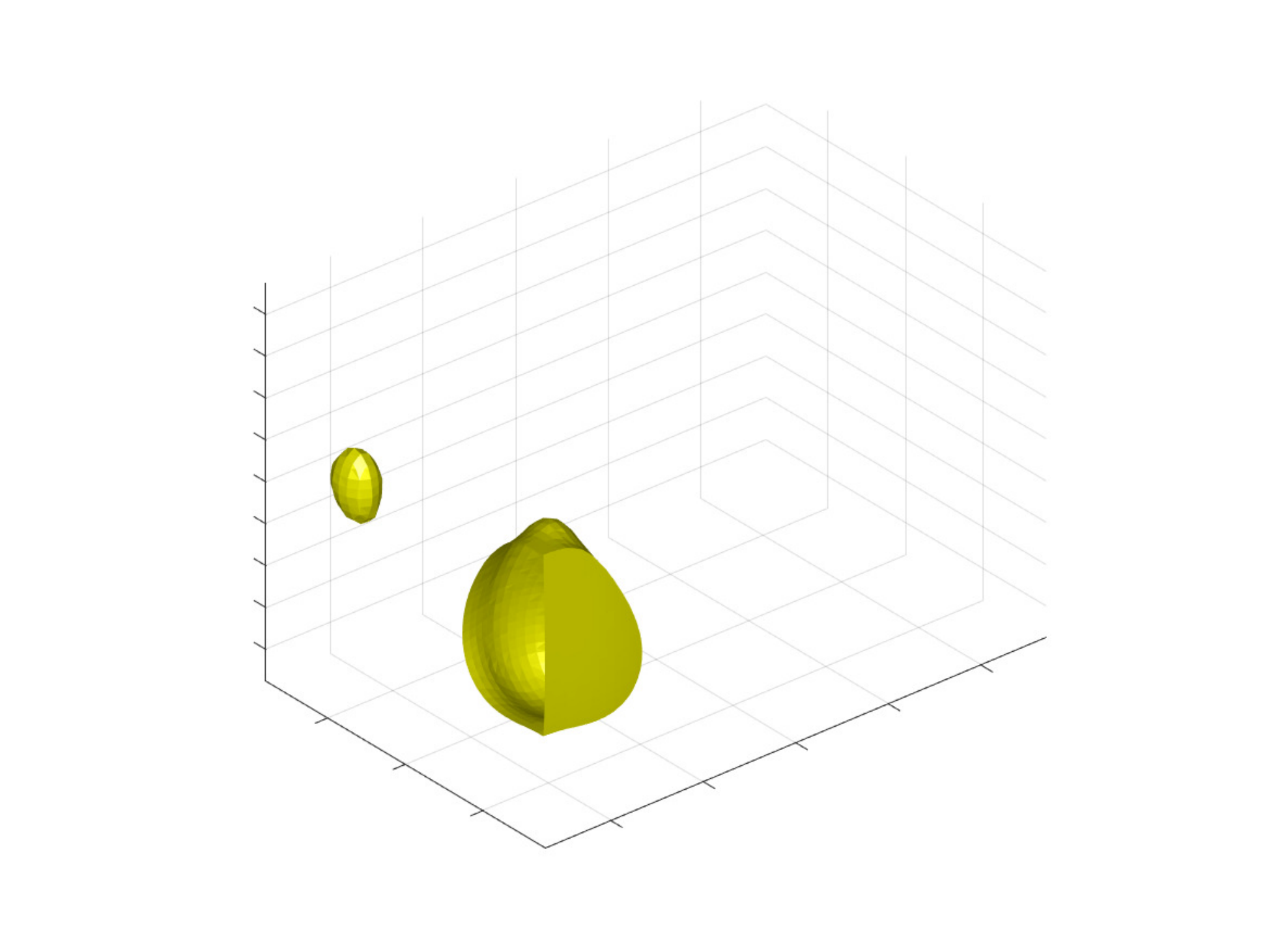}}
\put(350,80){\includegraphics[width=100pt]{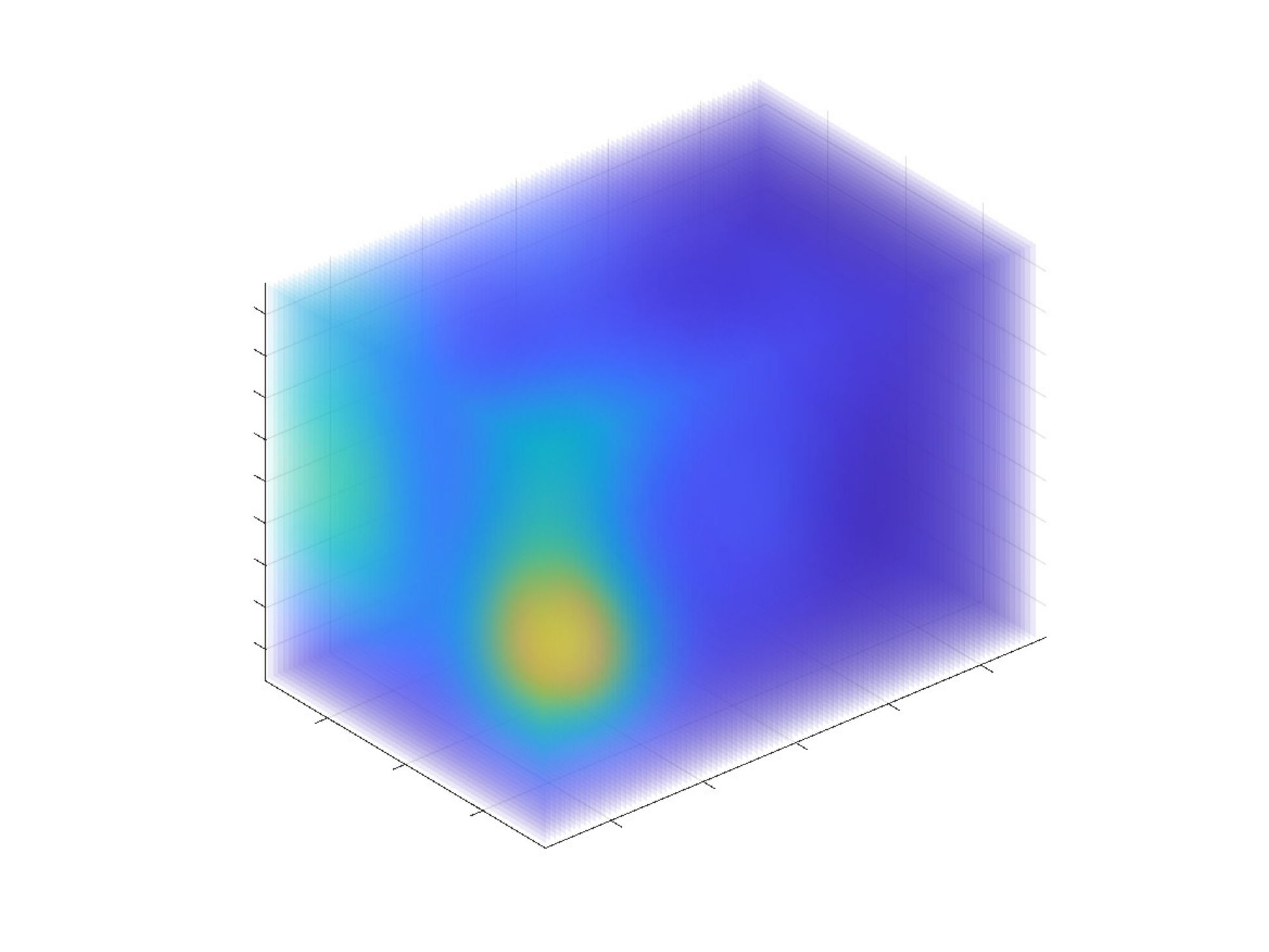}}
\put(230,80){\includegraphics[width=70pt]{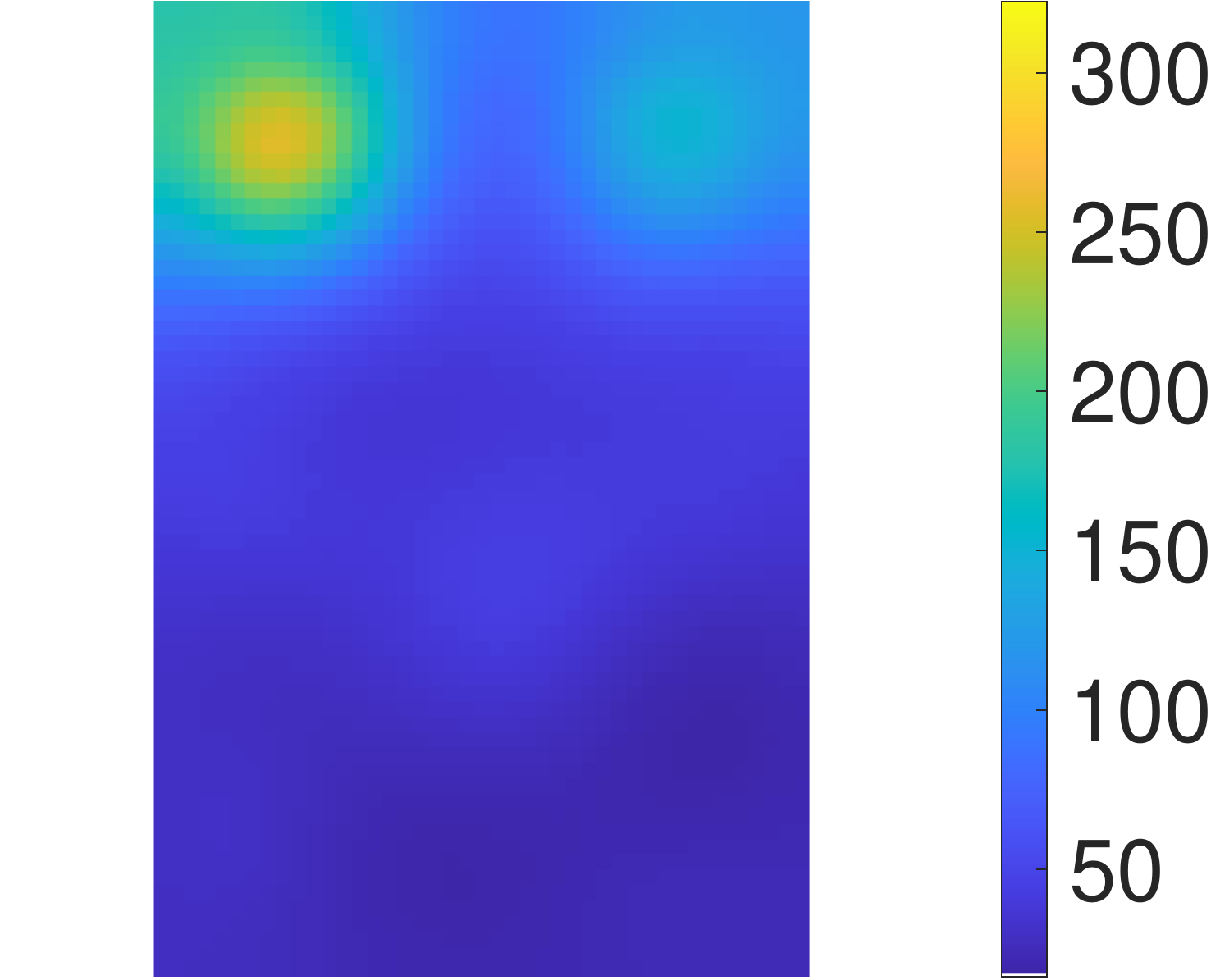}}

\put(60,80){\includegraphics[width=90pt]{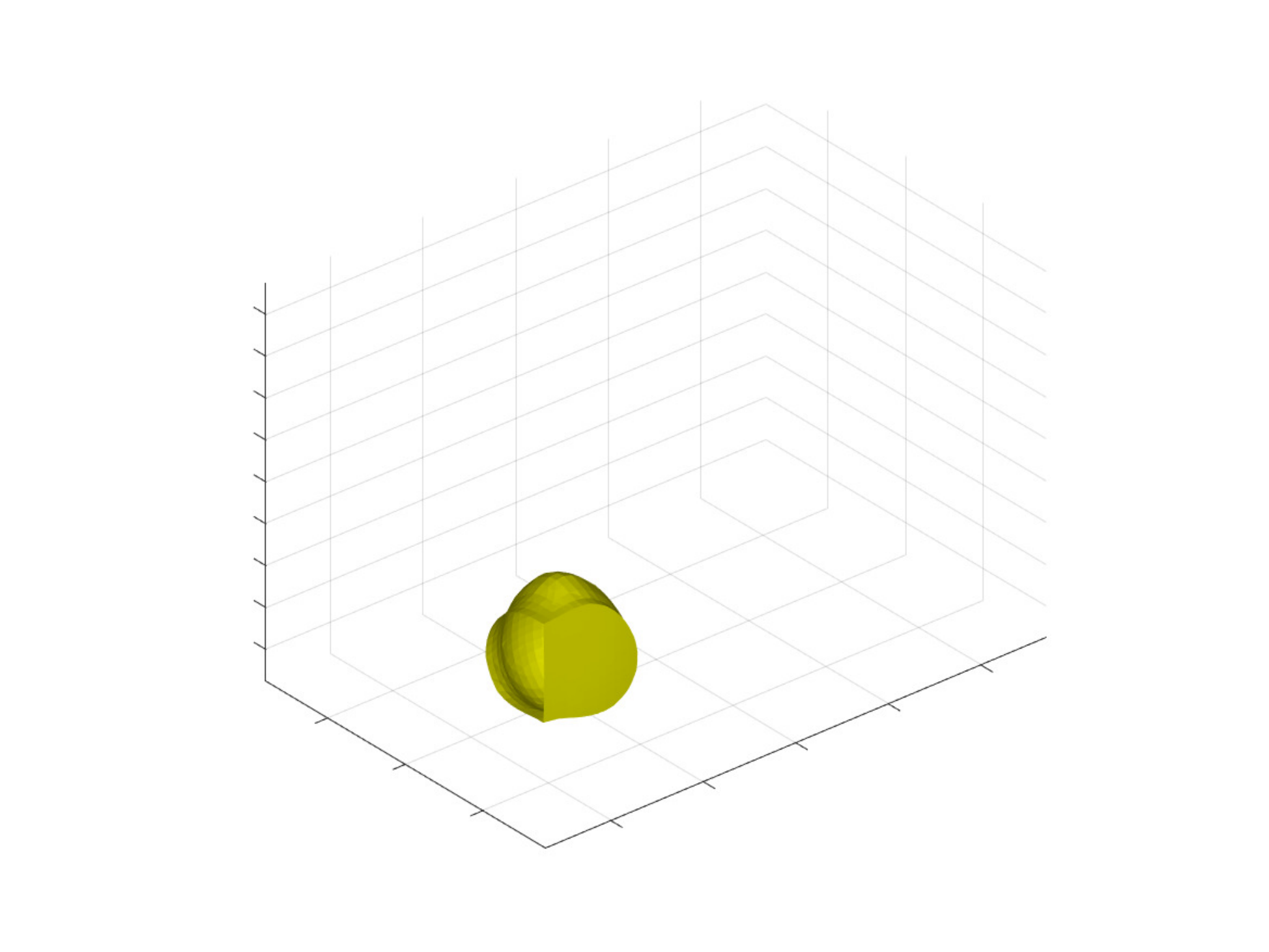}}
\put(120,80){\includegraphics[width=100pt]{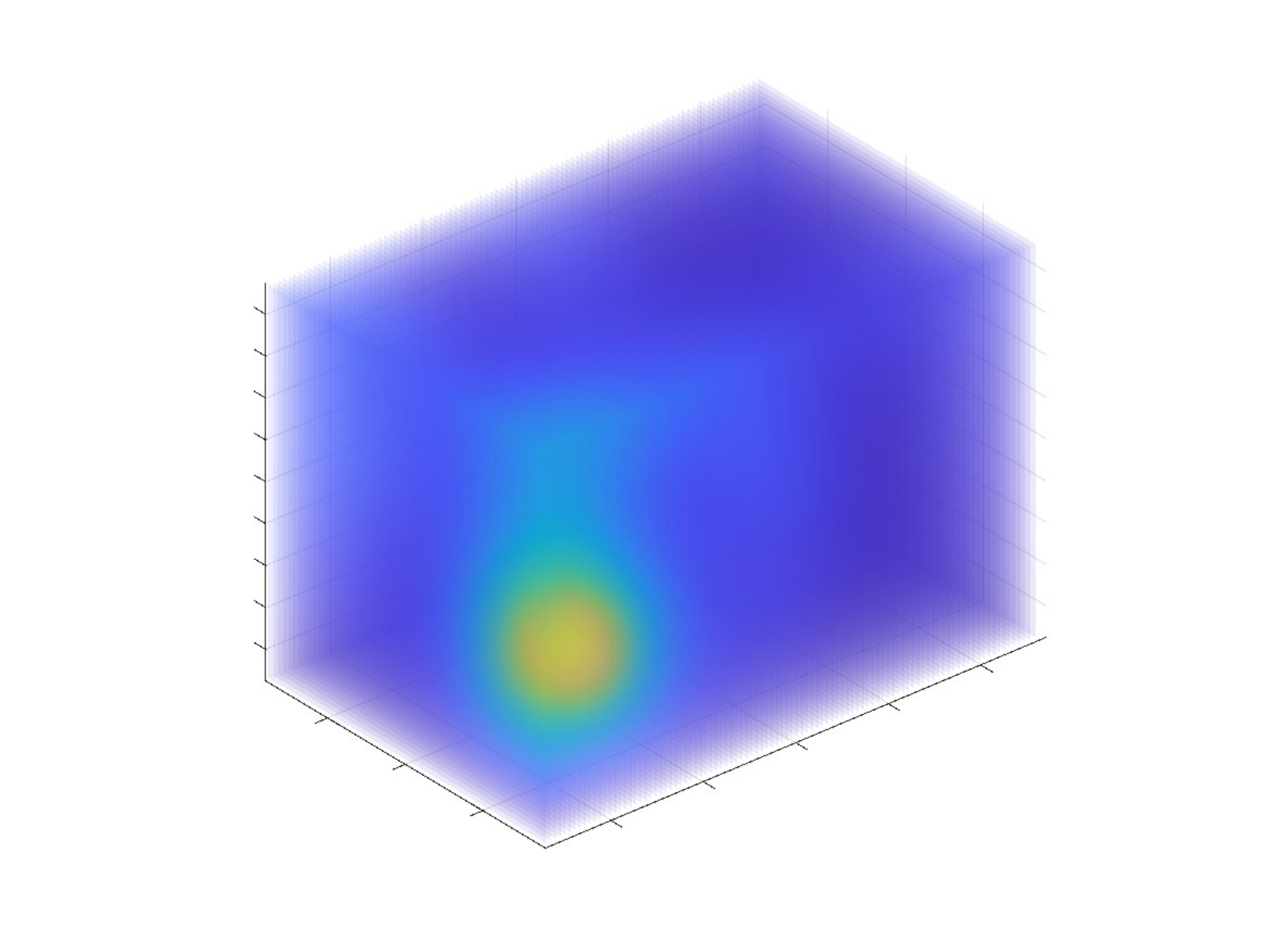}}
\put(0,80){\includegraphics[width=70pt]{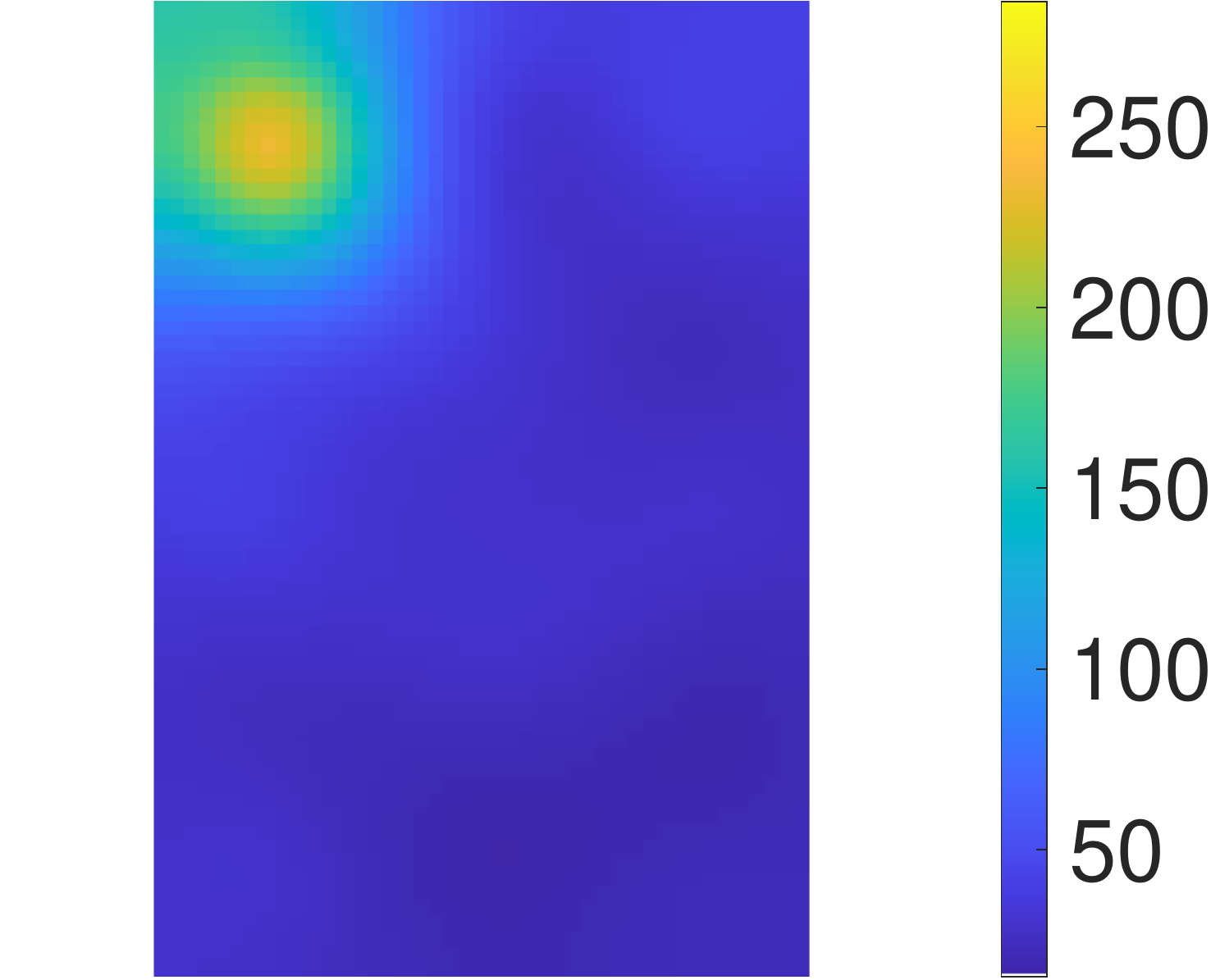}}

%----------------- 
% TV
\put(290,0){\includegraphics[width=90pt]{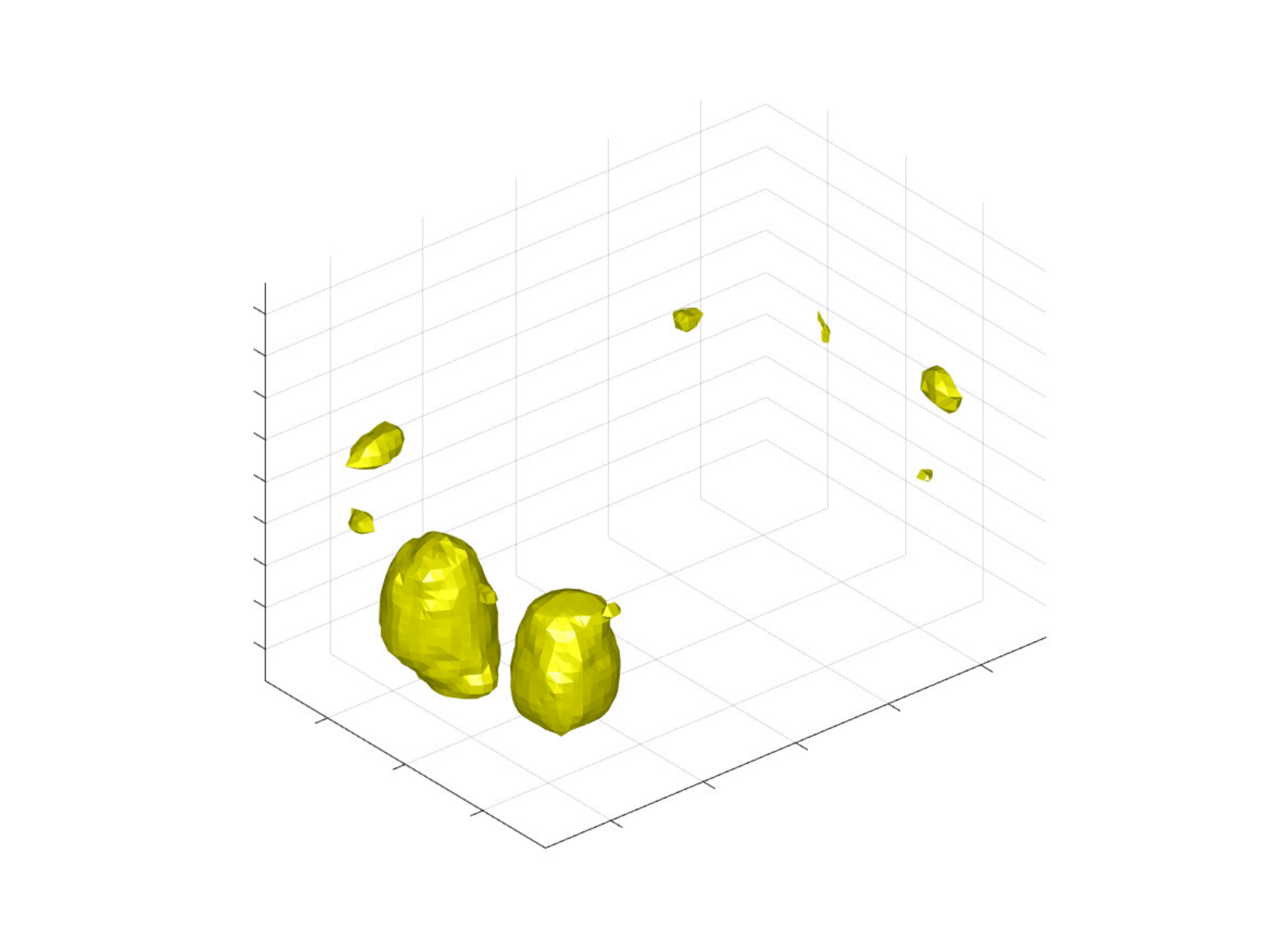}}
\put(350,0){\includegraphics[width=100pt]{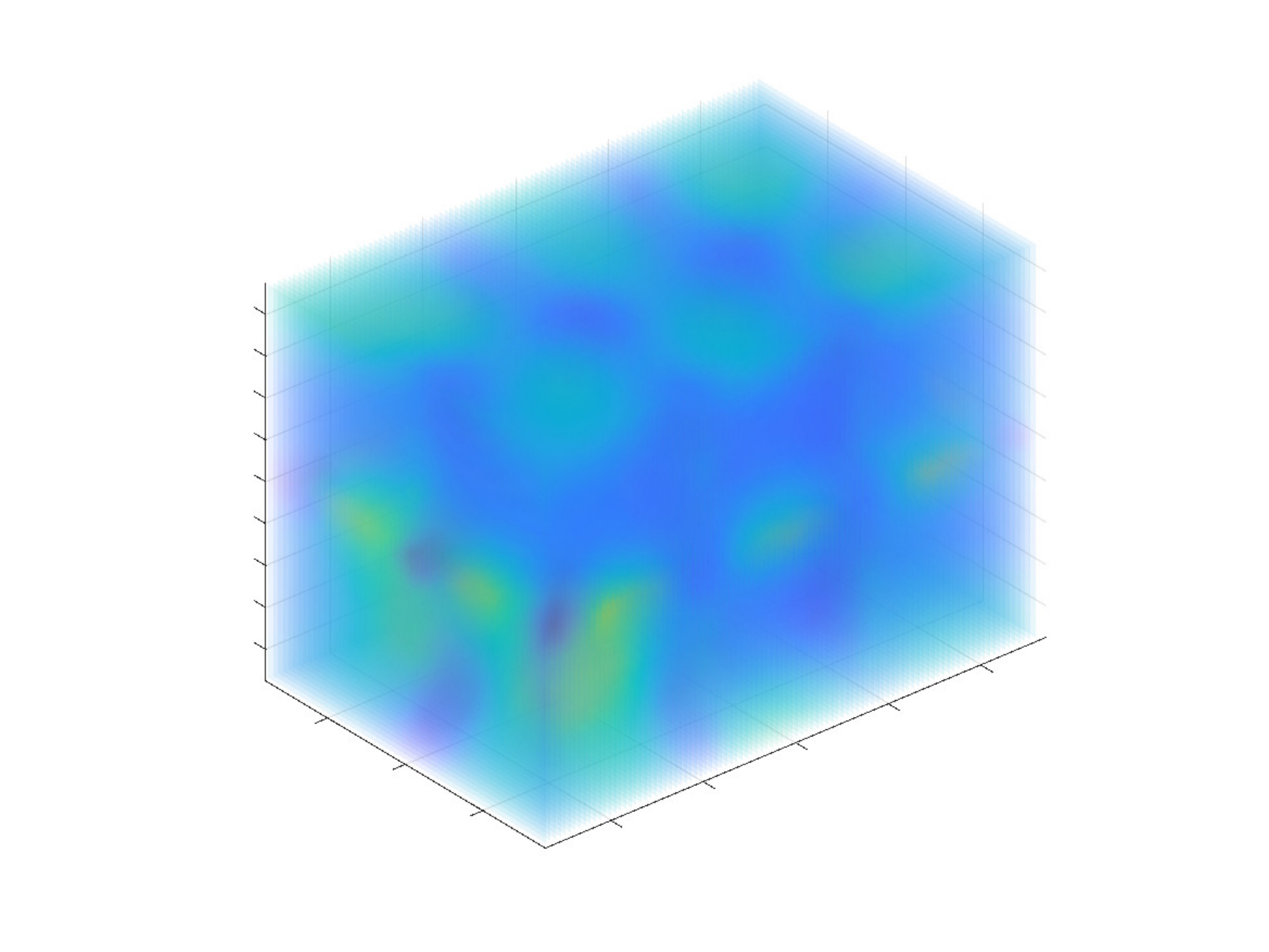}}
\put(230,0){\includegraphics[width=70pt]{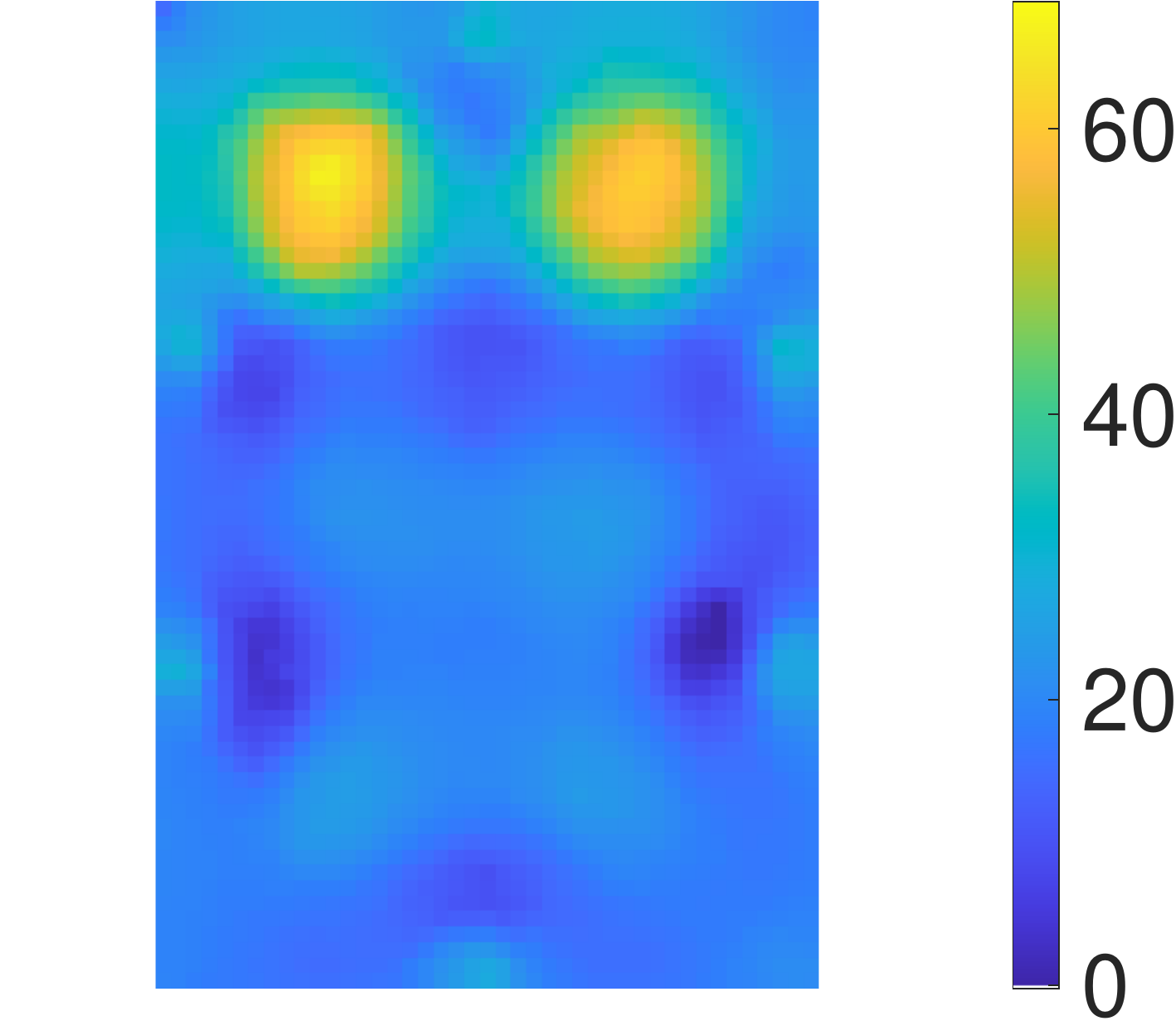}}

\put(60,0){\includegraphics[width=90pt]{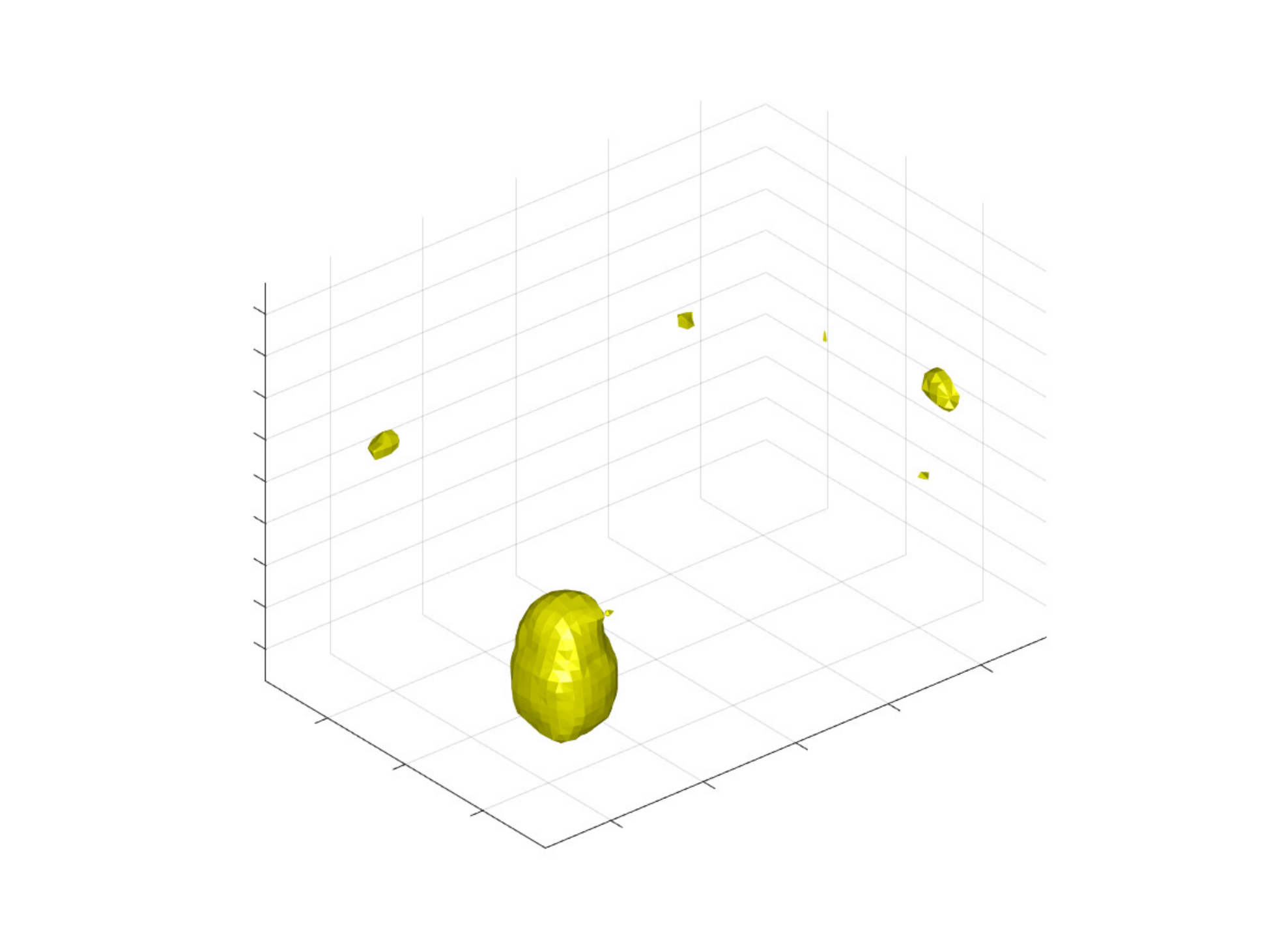}}
\put(120,0){\includegraphics[width=100pt]{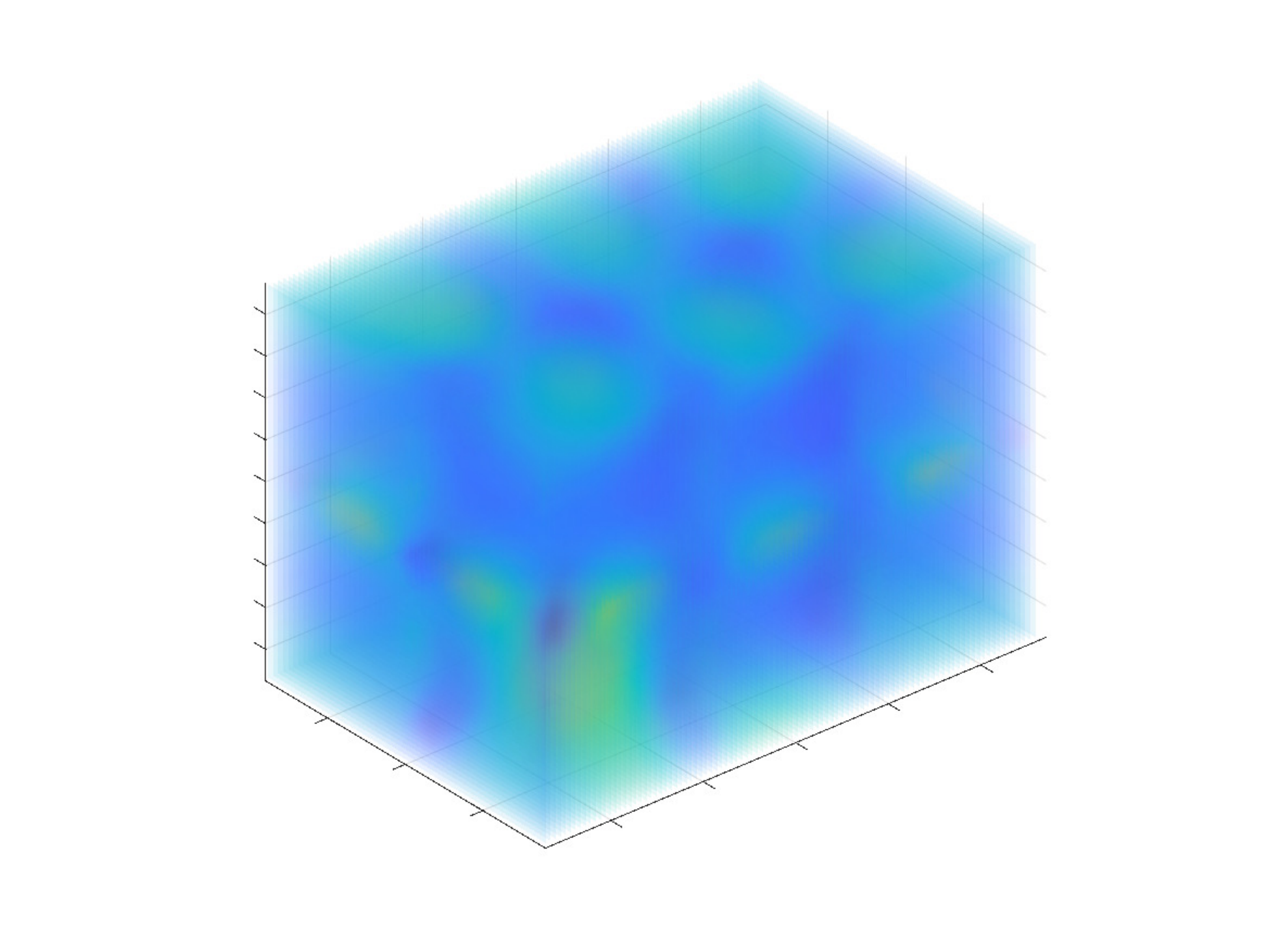}}
\put(0,0){\includegraphics[width=70pt]{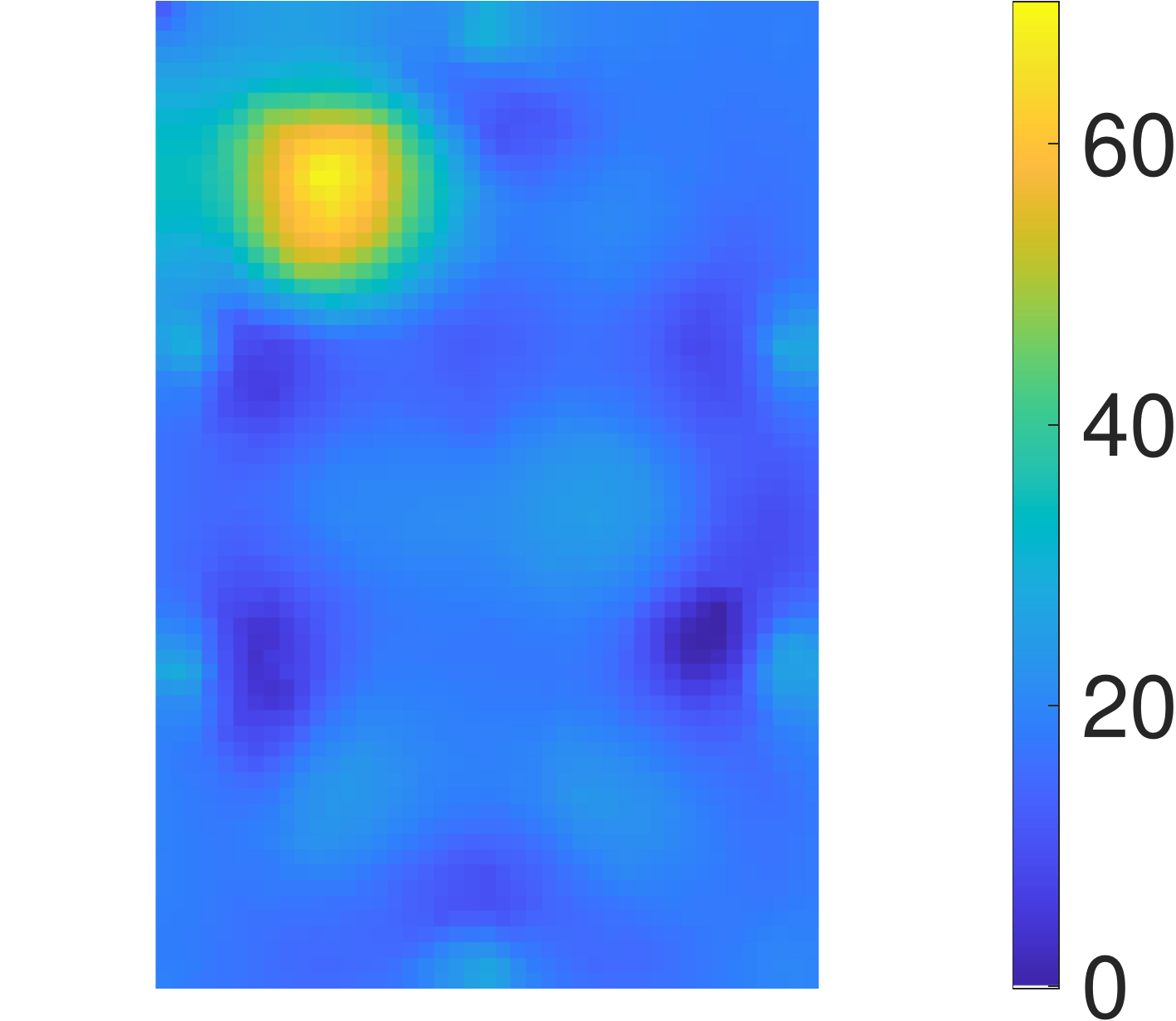}}

\put(3,0){\line(0,1){400}}
%\put(145,0){\line(0,1){400}}
\put(220,0){\line(0,1){400}}
%\put(445,0){\line(0,1){400}}
%----------------- 
% Top labels:
\put(20,385){{{\sc\footnotesize Slice}}}
\put(75,385){{{\sc\footnotesize Isosurface}}}
\put(160,385){{{\sc\footnotesize 3D}}}
\put(250,385){{{\sc\footnotesize Slice}}}
\put(310,385){{{\sc\footnotesize Isosurface}}}
\put(390,385){{{\sc\footnotesize 3D}}}

\put(68,400){{\underline{\sc One Target}}}
\put(298,400){{\underline{\sc Two Targets}}}
%\put(320,395){{\underline{\sc 20x35x25 modeling}}}

%----------------- 
% Side labels:
\put(-10,330){\rotatebox{90}{\sc Truth}}
\put(-10,250){\rotatebox{90}{\sc Cald}}
\put(-10,180){\rotatebox{90}{\sc $\texp$}}
\put(-10,105){\rotatebox{90}{\sc $\tzero$}}
\put(-10,25){\rotatebox{90}{\sc TV}}
%----------------- 

\put(445,0){\line(0,1){400}}

\end{picture}
\caption{\label{fig:abs_182719} {\bf Absolute image} reconstructions comparing the CGO methods to the regularized method with {\it moderately incorrect domain modeling}, using a box of size 18cm x 27cm x 19cm.  Note the truth targets had a measured conductivity of approx 290 mS/m.}
\end{figure}

\begin{figure}[ht]
\centering
\begin{picture}(450,410)
\linethickness{.3mm}
%----------------- 
%----------------- 
% Truth
\put(300,320){\includegraphics[width=60pt]{truth_2targ_iso.png}}
%\put(350,320){\includegraphics[width=100pt]{cald_true_1_targ_ABS_3D-eps-converted-to.pdf}}
\put(240,320){\includegraphics[angle=-90,origin=c,width=42pt]{TwoTargs_TopView.jpeg}}

\put(70,320){\includegraphics[width=60pt]{truth_1targ_iso.png}}
%\put(120,240){\includegraphics[width=100pt]{cald_true_2_targ_ABS_3D-eps-converted-to.pdf}}
\put(10,320){\includegraphics[angle=-90,origin=c,width=42pt]{OneTarg_TopView.jpeg}}
% %----------------- 
% cald
\put(290,240){\includegraphics[width=90pt]{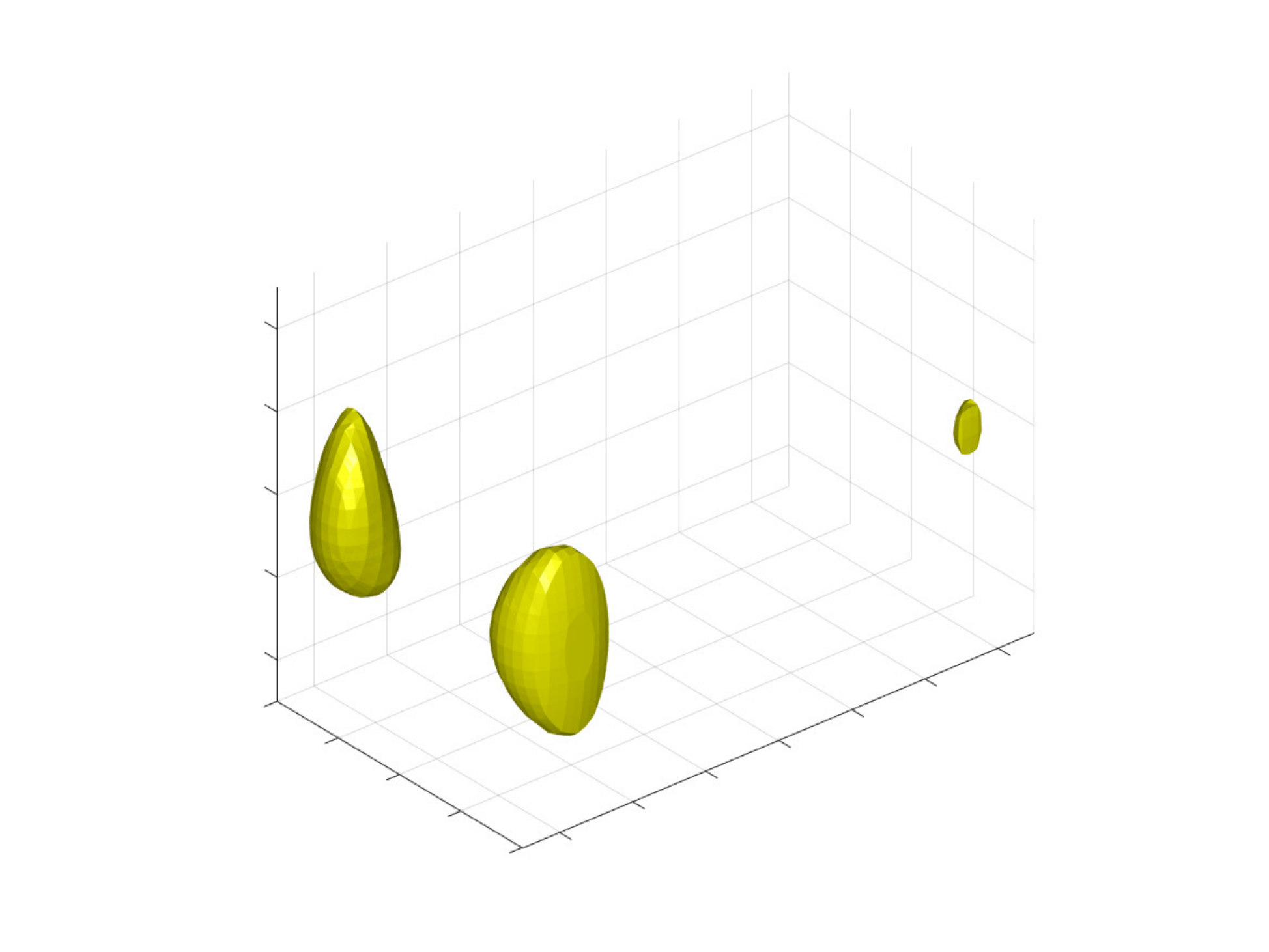}}
\put(350,240){\includegraphics[width=100pt]{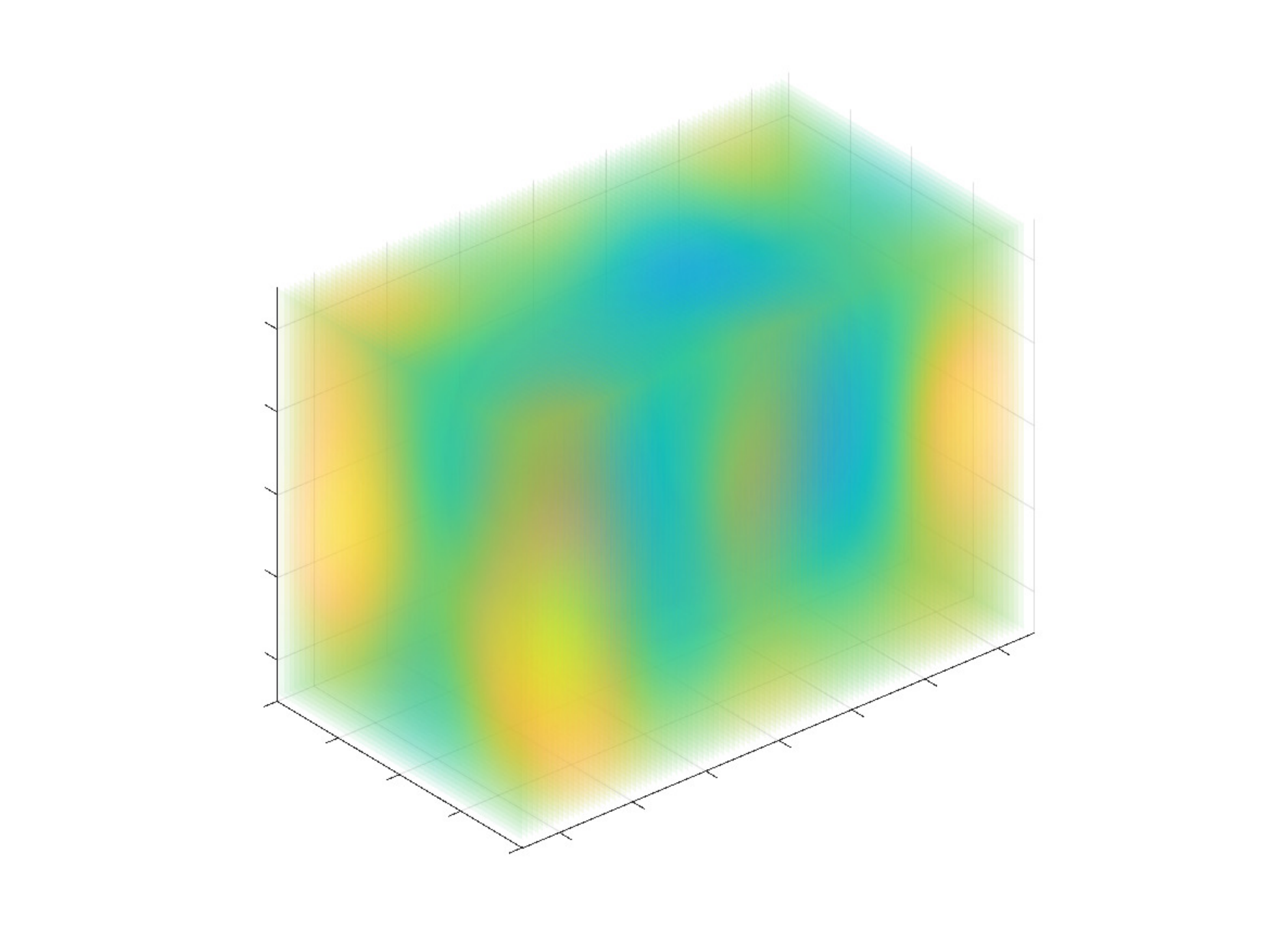}}
\put(230,240){\includegraphics[width=70pt]{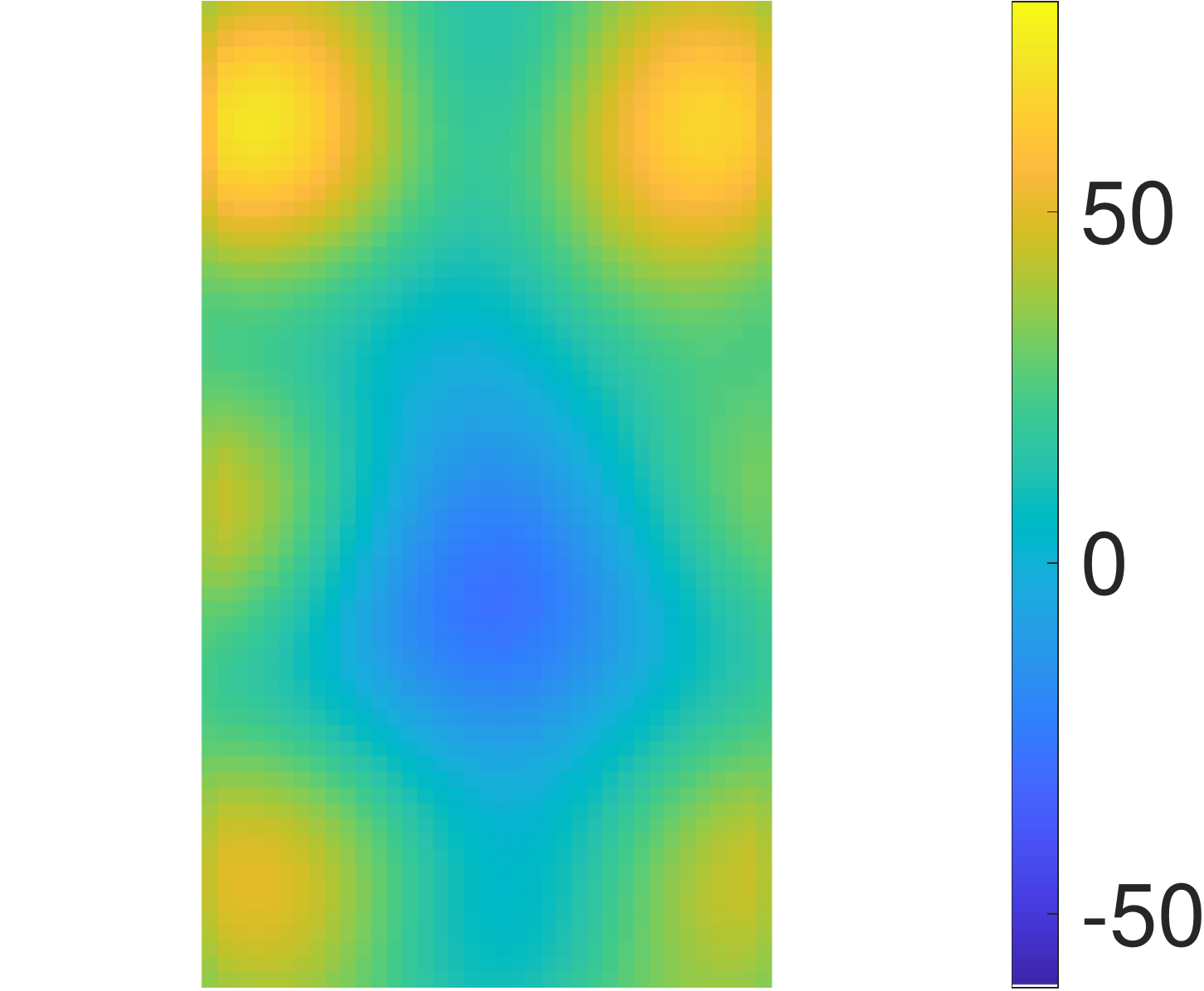}}

\put(60,240){\includegraphics[width=90pt]{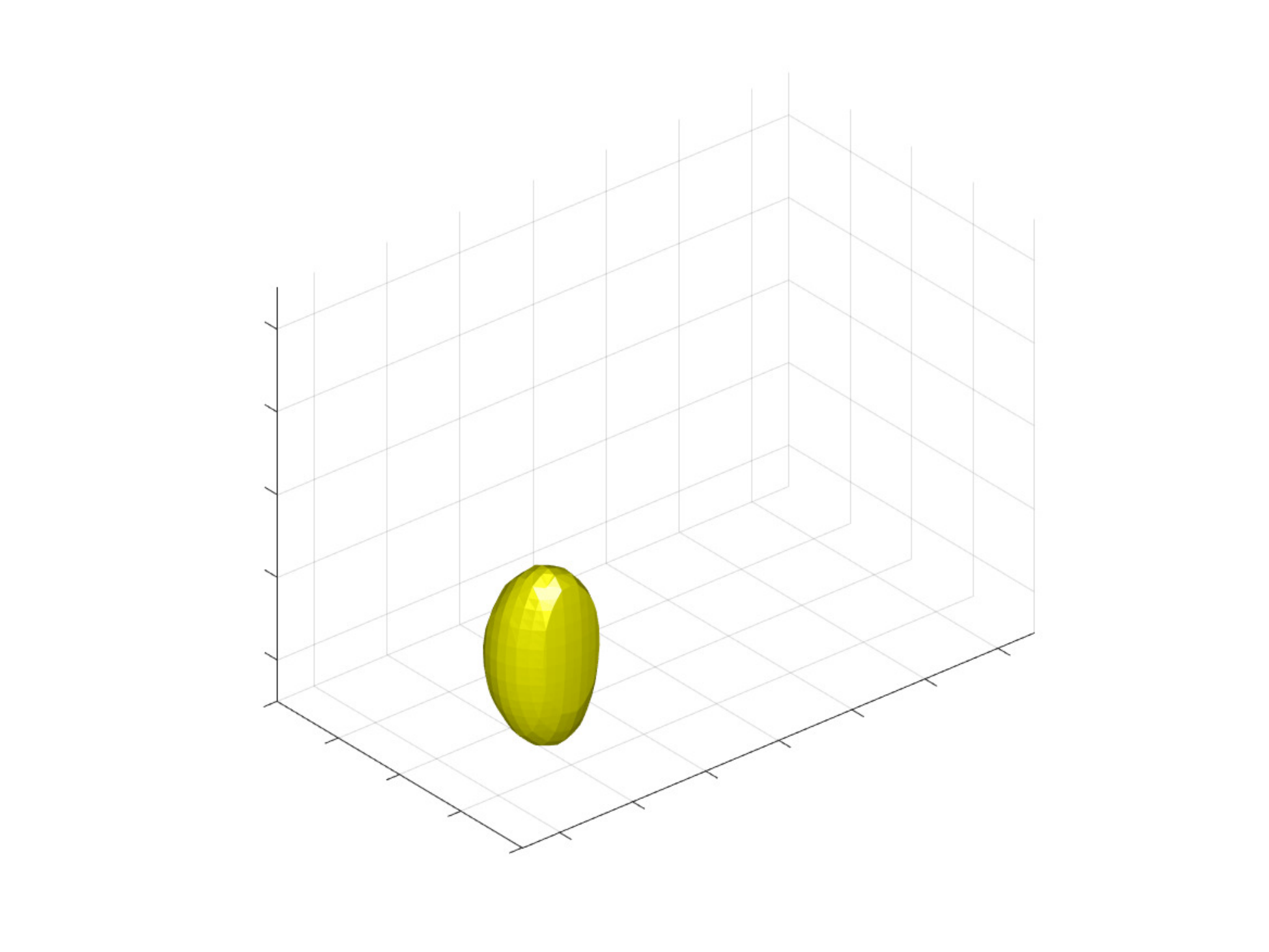}}
\put(120,240){\includegraphics[width=100pt]{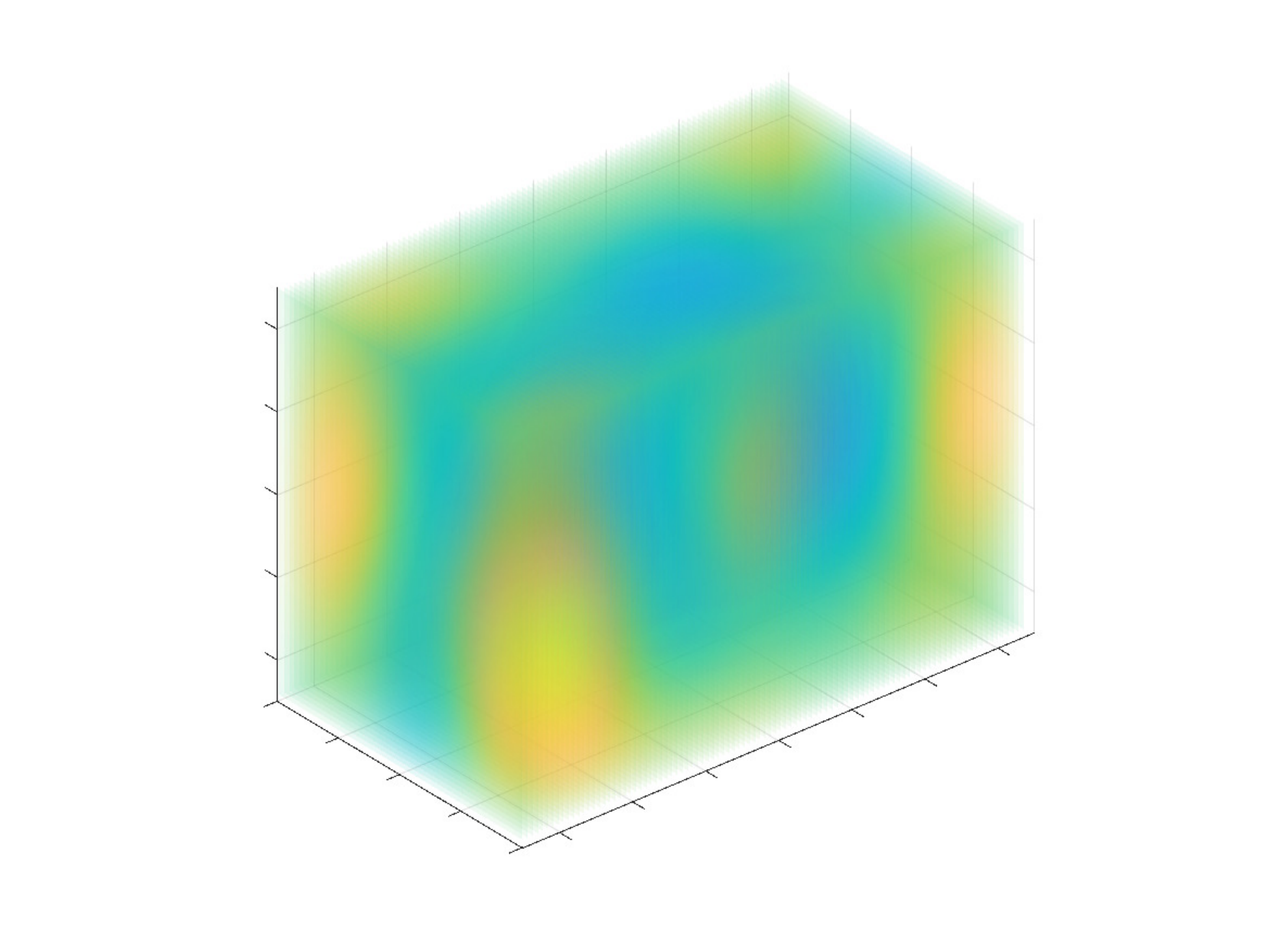}}
\put(0,240){\includegraphics[width=70pt]{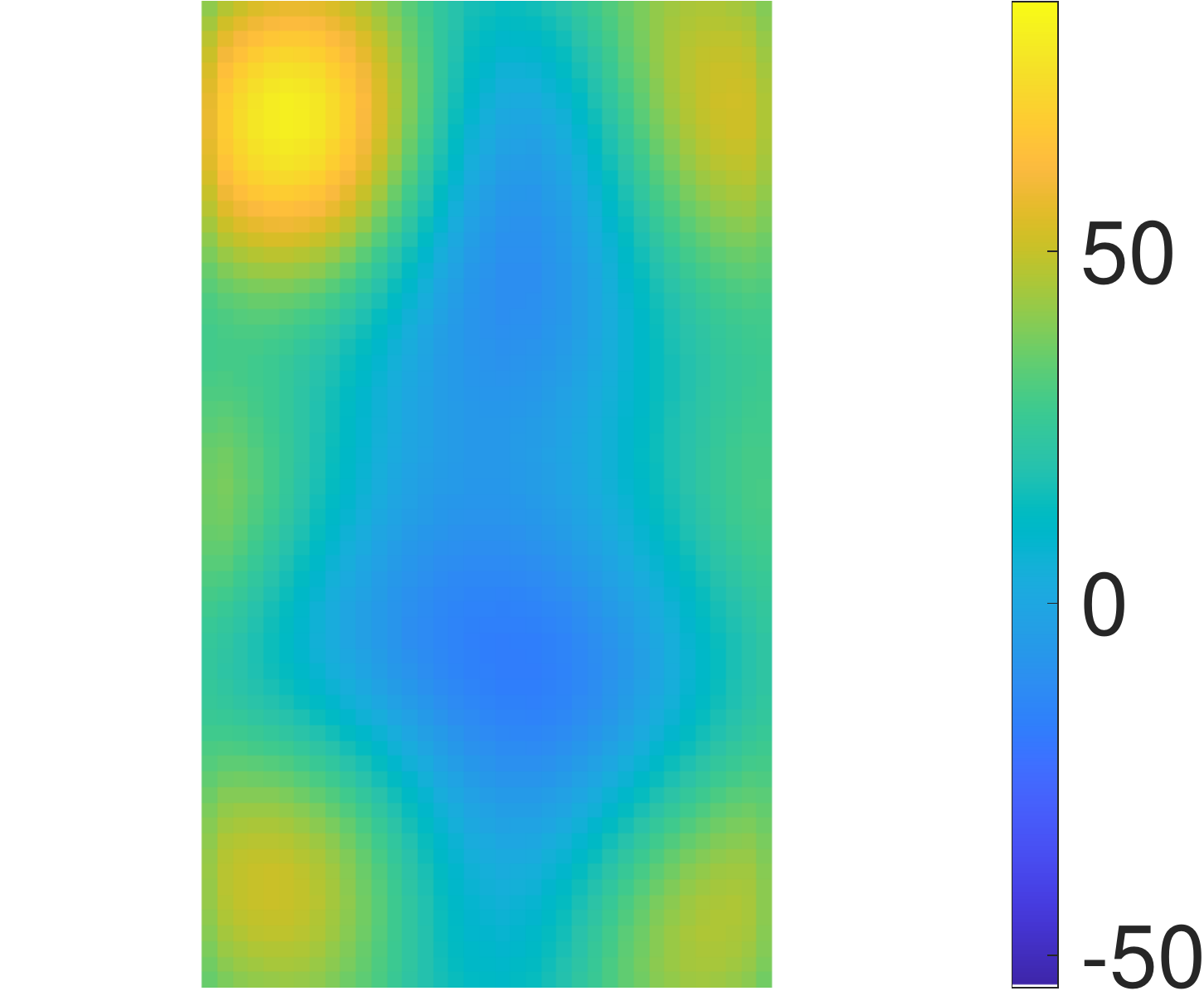}}
% %----------------- 
% texp
\put(290,160){\includegraphics[width=90pt]{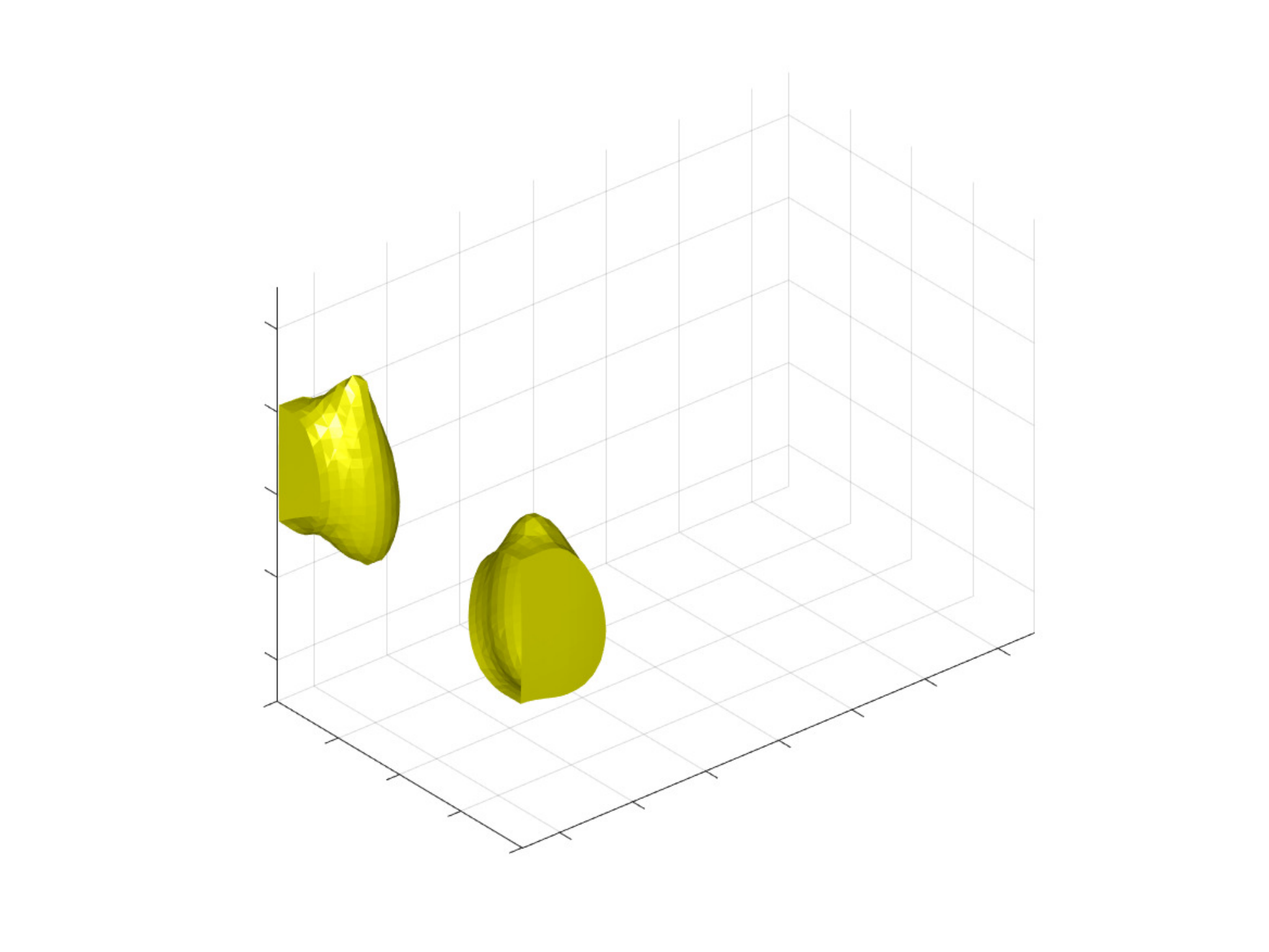}}
\put(350,160){\includegraphics[width=100pt]{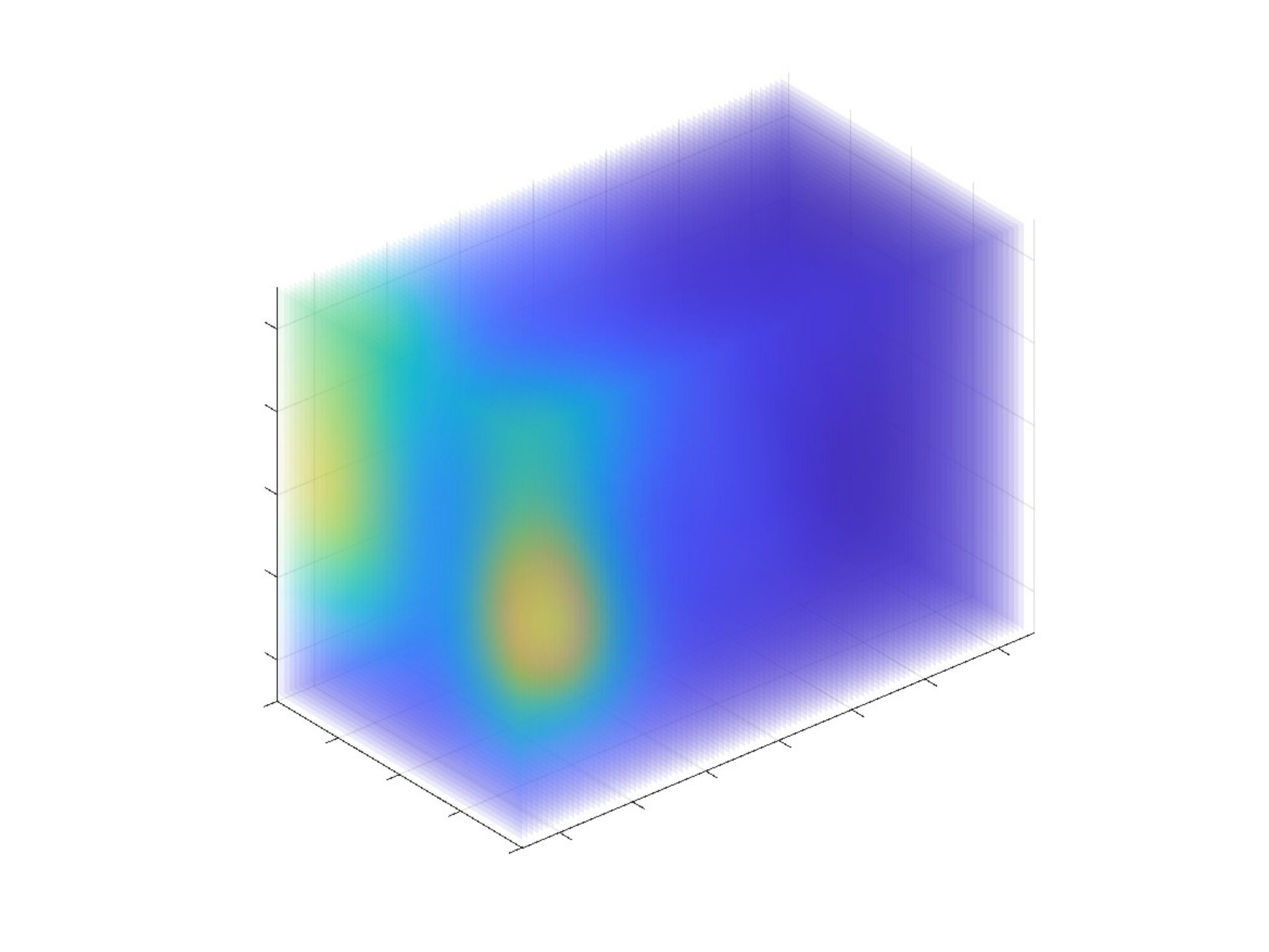}}
\put(230,160){\includegraphics[width=70pt]{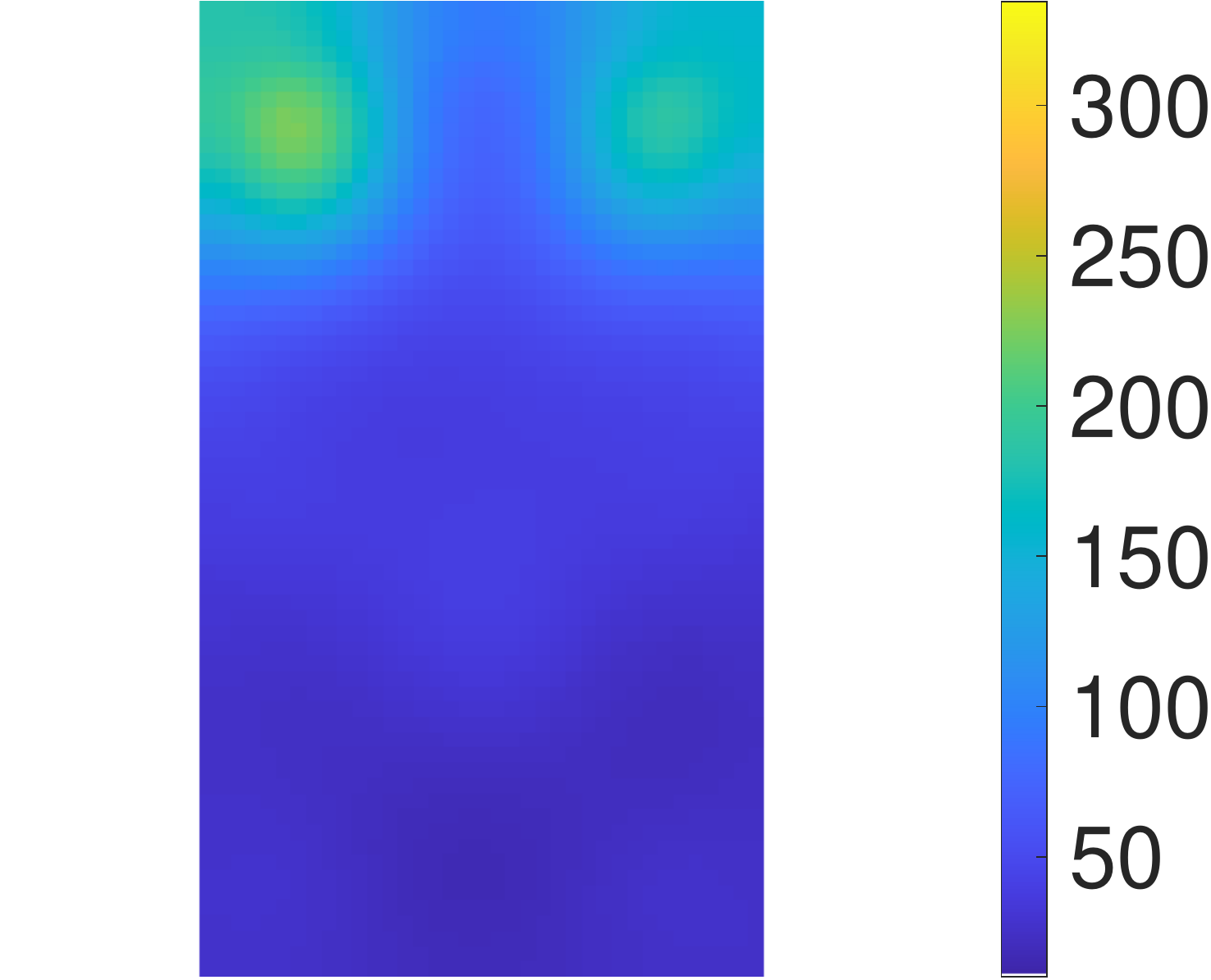}}

\put(60,160){\includegraphics[width=90pt]{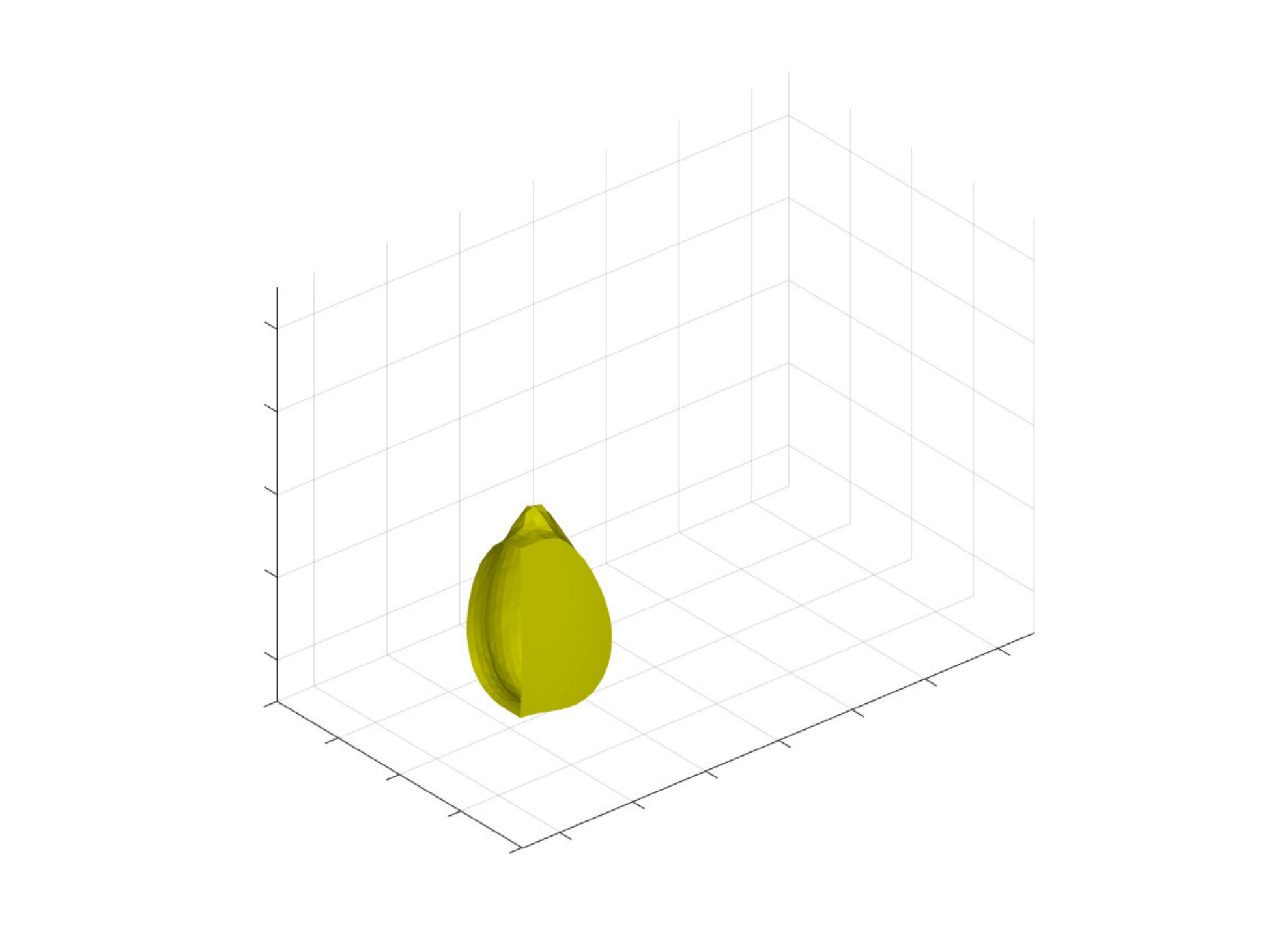}}
\put(120,160){\includegraphics[width=100pt]{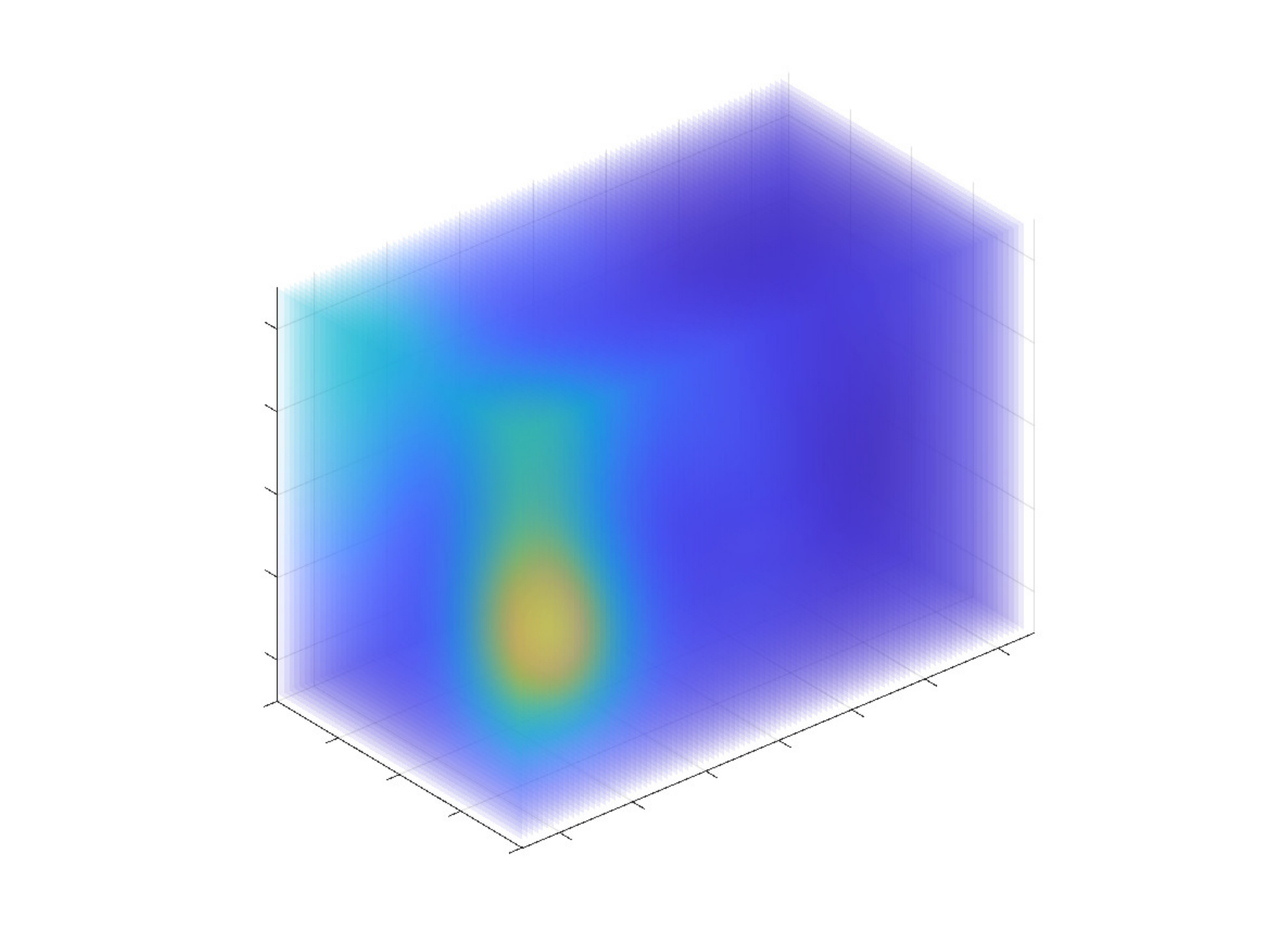}}
\put(0,160){\includegraphics[width=70pt]{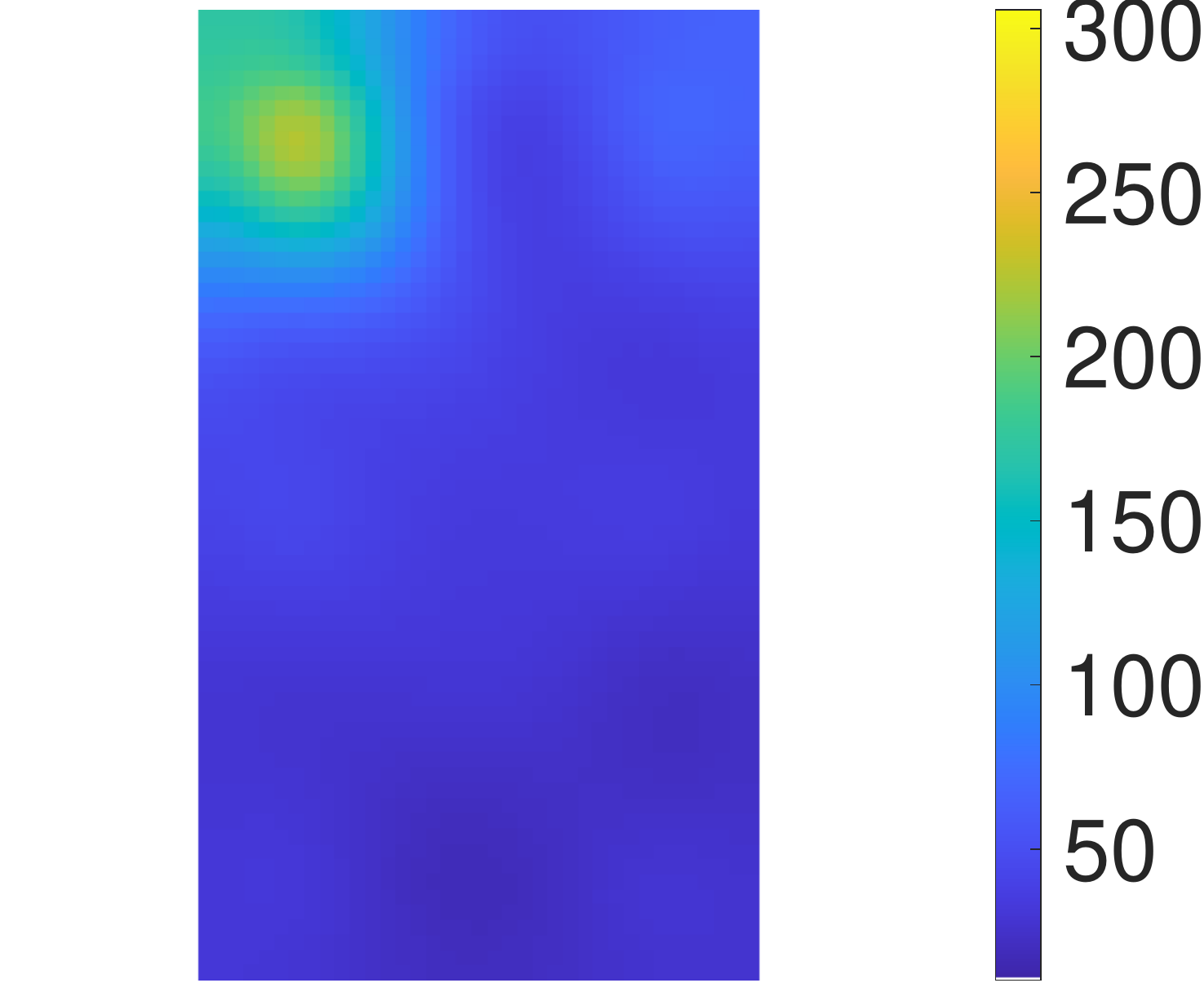}}

% %----------------- 
% t0
\put(290,80){\includegraphics[width=90pt]{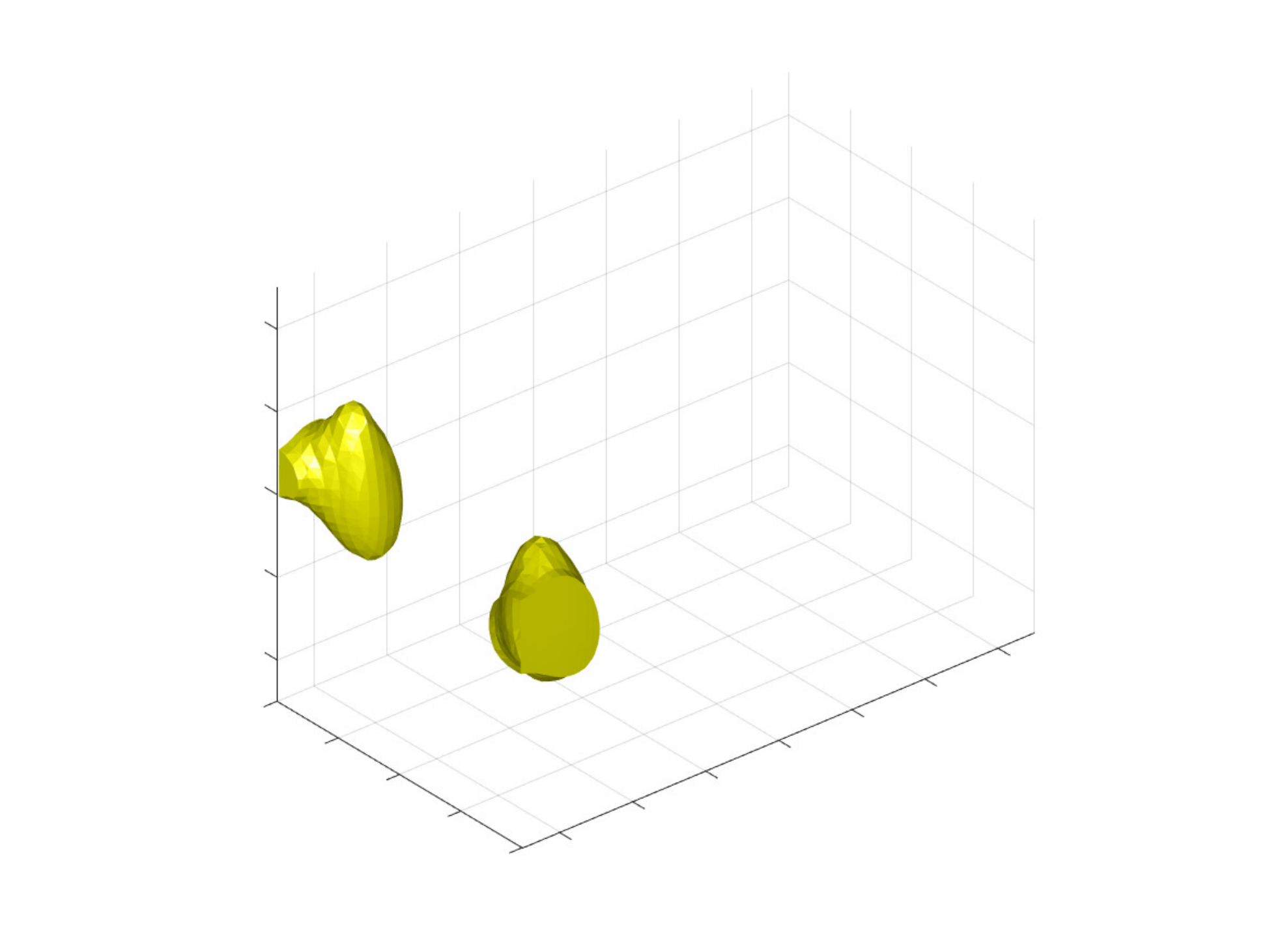}}
\put(350,80){\includegraphics[width=100pt]{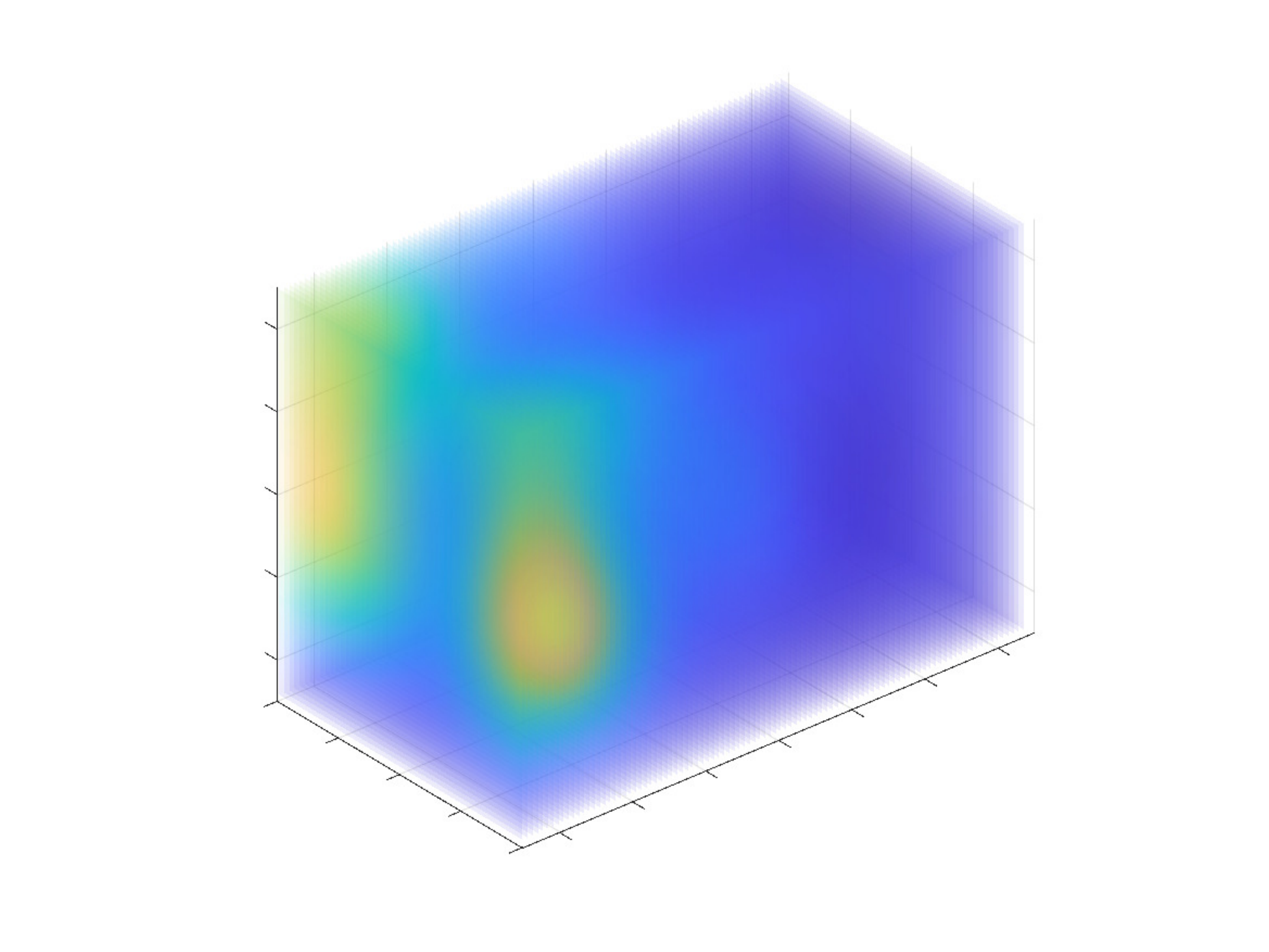}}
\put(230,80){\includegraphics[width=70pt]{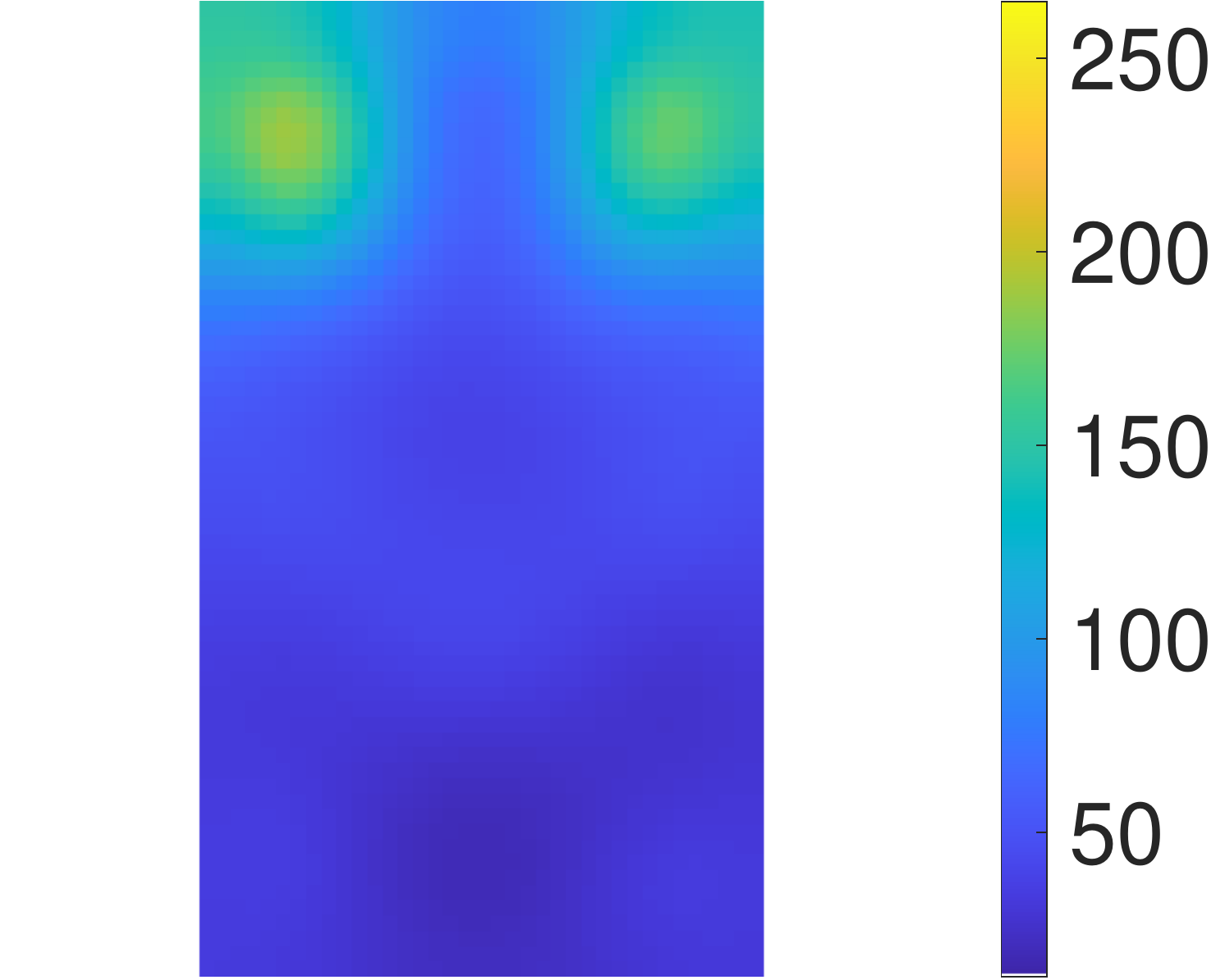}}

\put(60,80){\includegraphics[width=90pt]{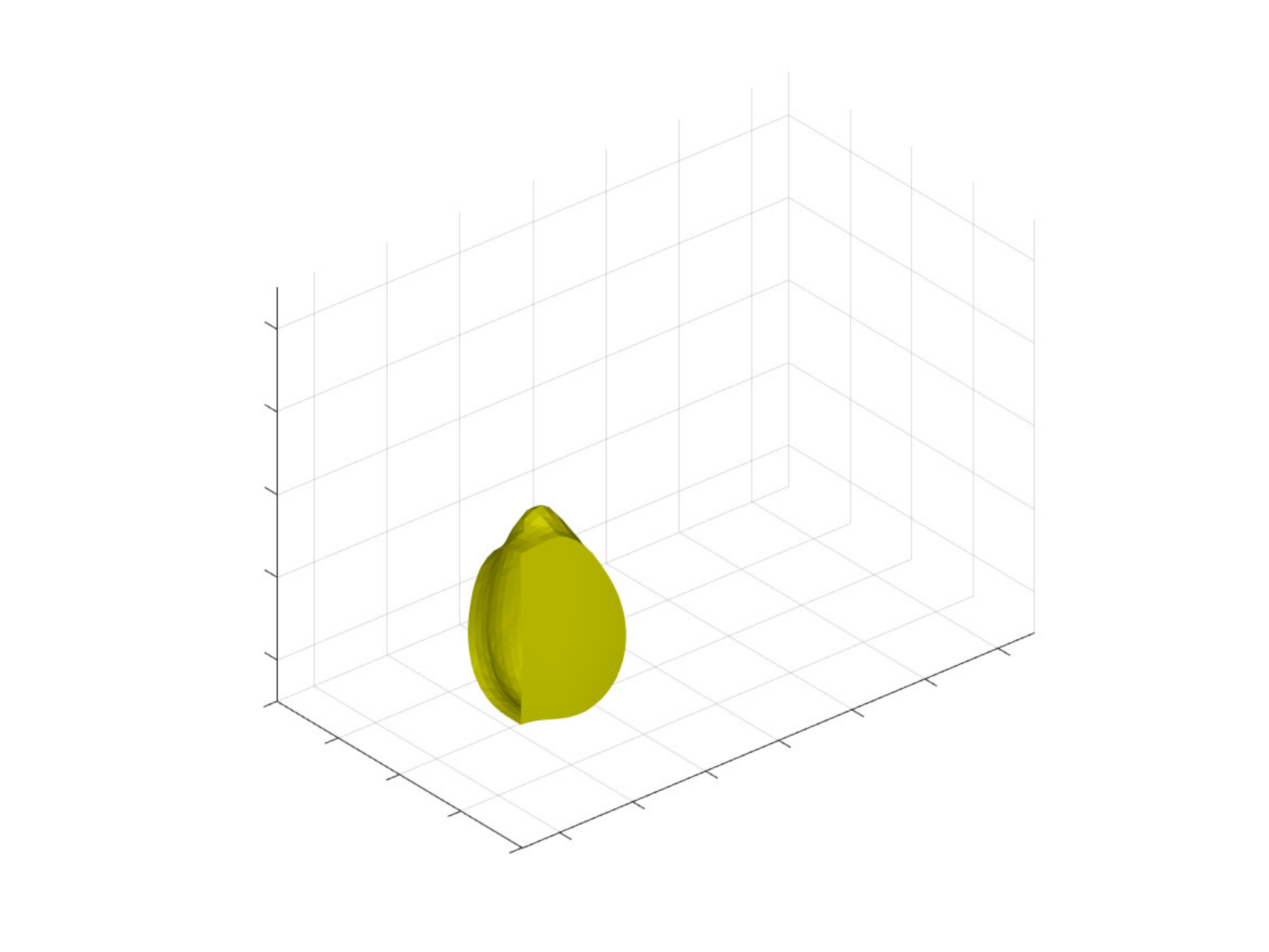}}
\put(120,80){\includegraphics[width=100pt]{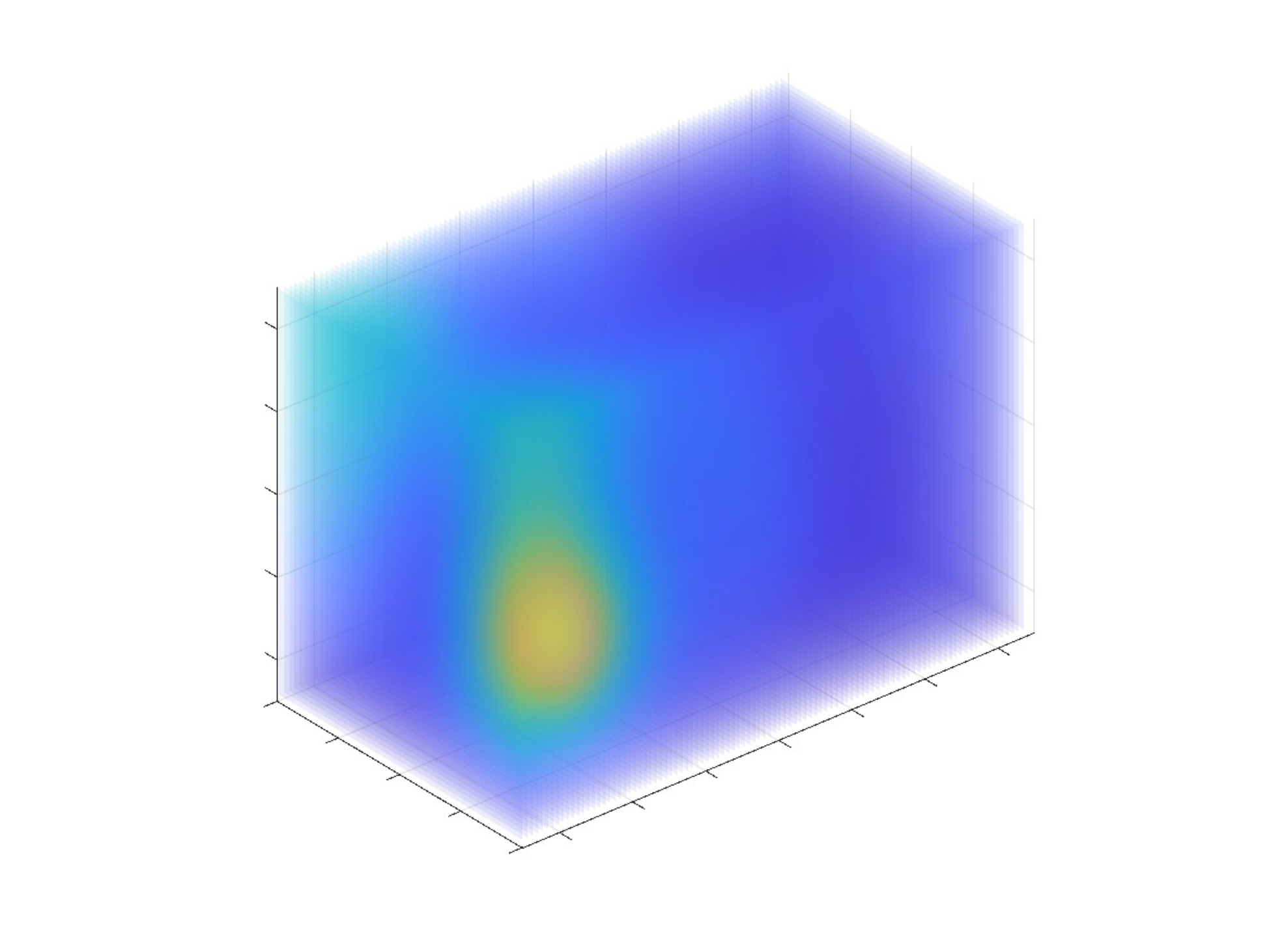}}
\put(0,80){\includegraphics[width=70pt]{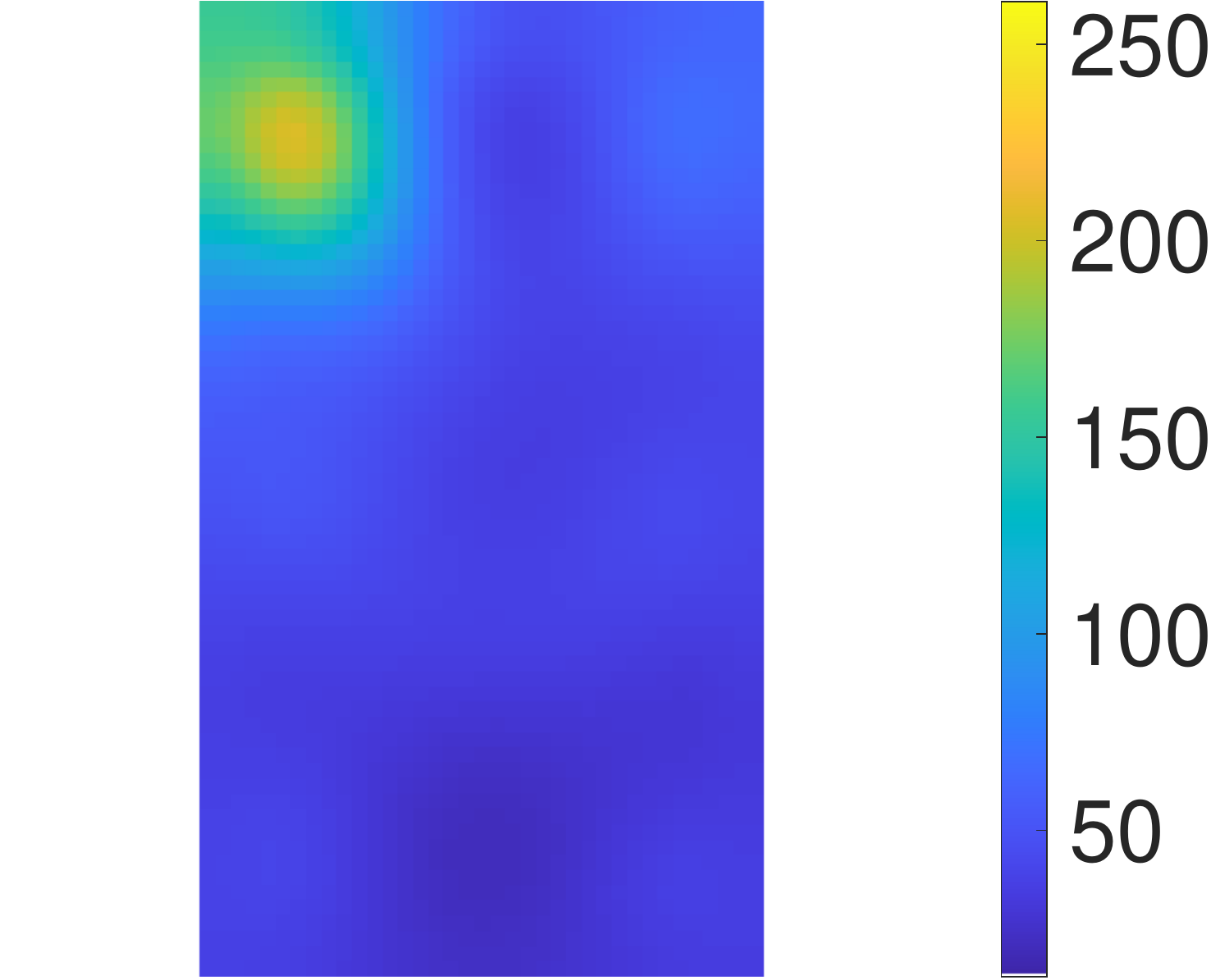}}

%----------------- 
% TV
\put(290,0){\includegraphics[width=90pt]{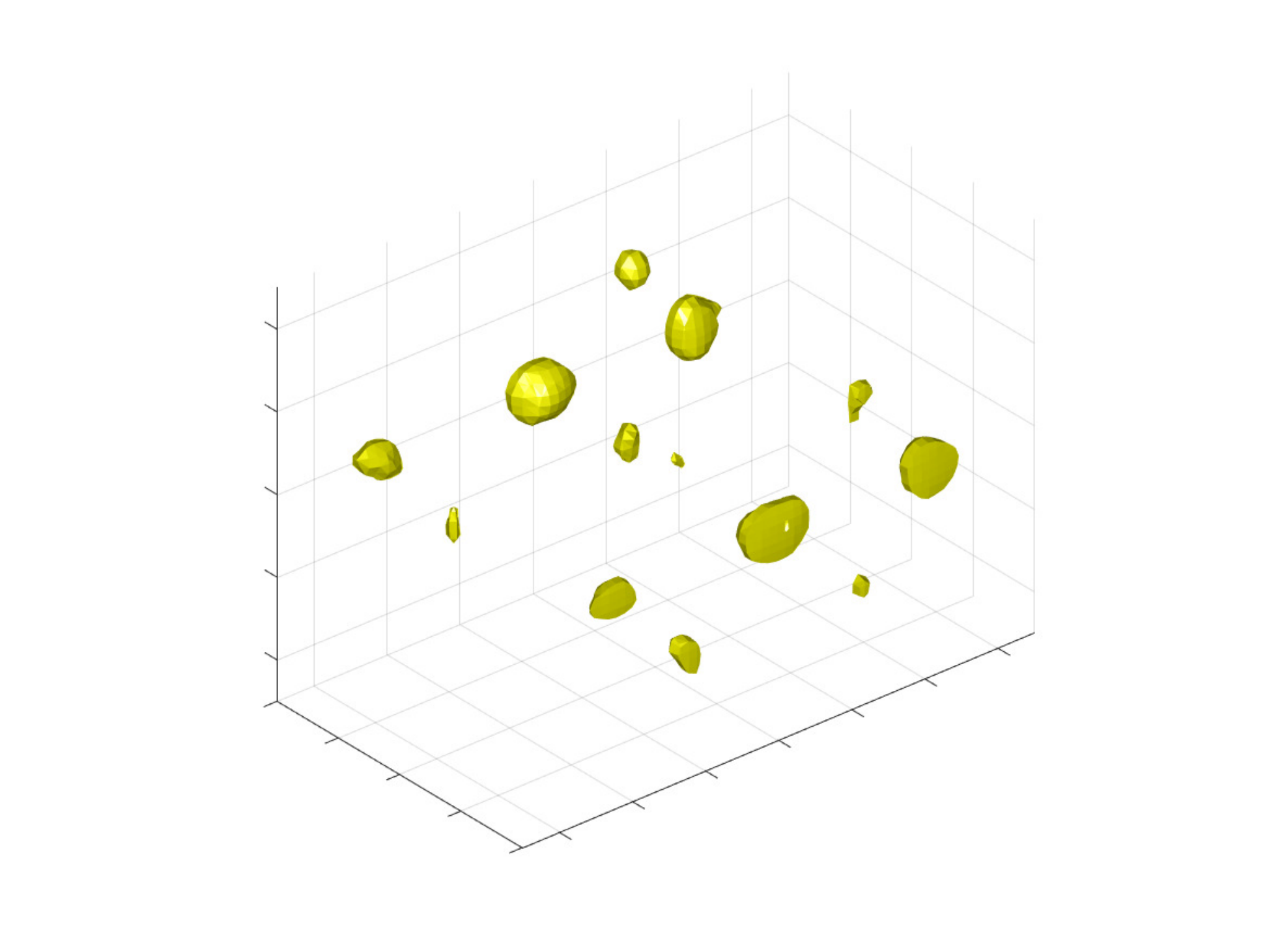}}
\put(350,0){\includegraphics[width=100pt]{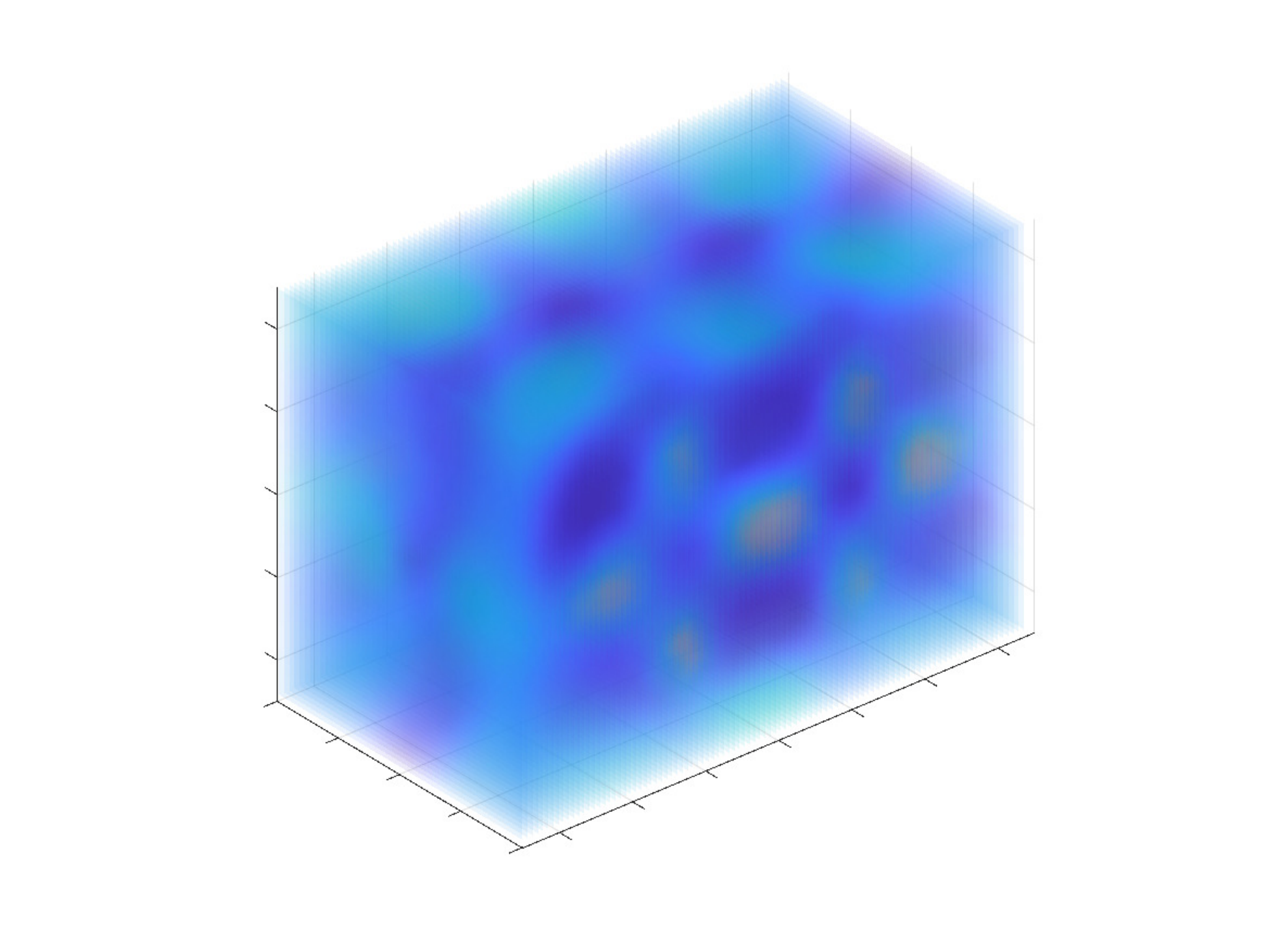}}
\put(230,0){\includegraphics[width=70pt]{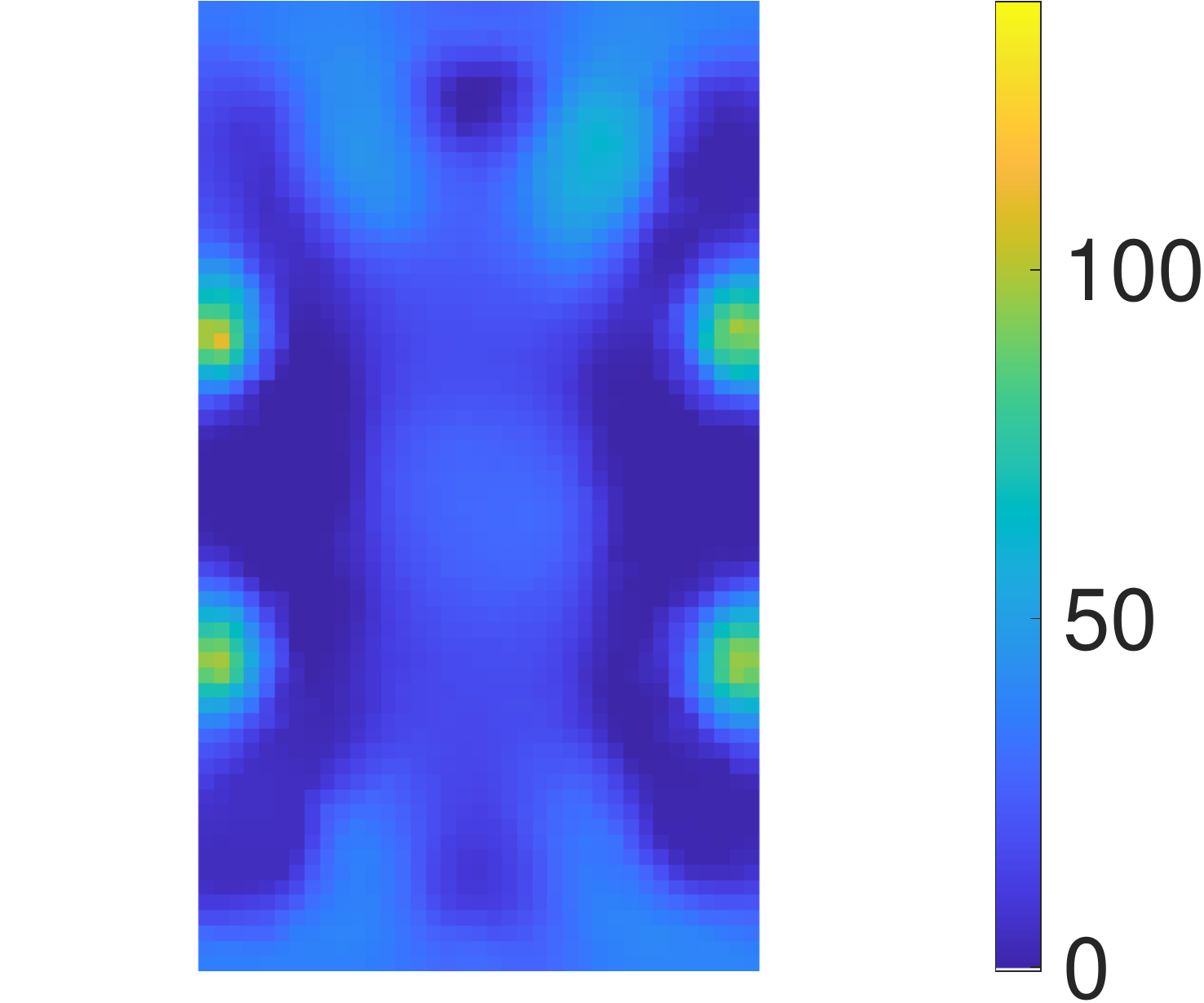}}

\put(60,0){\includegraphics[width=90pt]{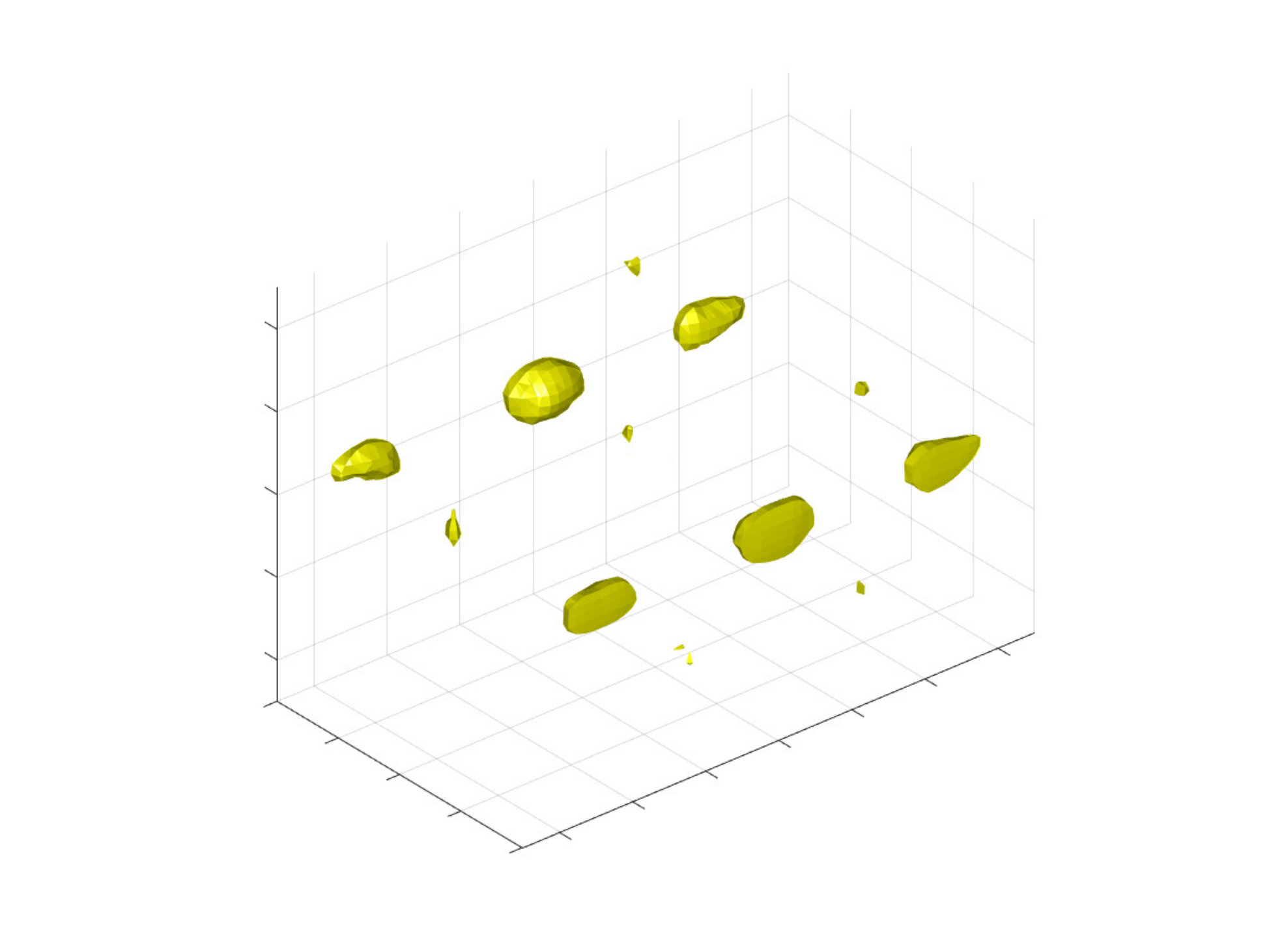}}
\put(120,0){\includegraphics[width=100pt]{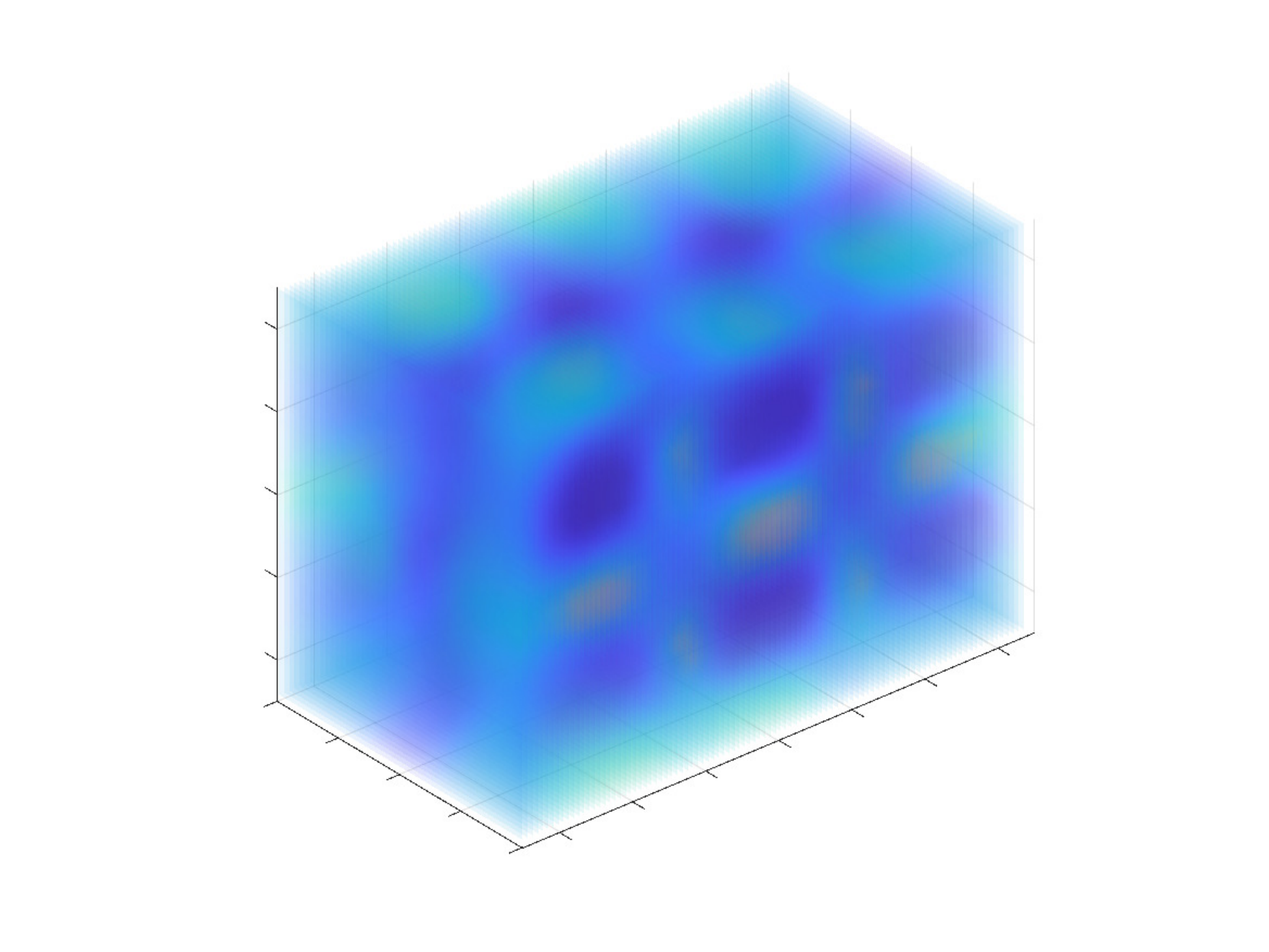}}
\put(0,0){\includegraphics[width=70pt]{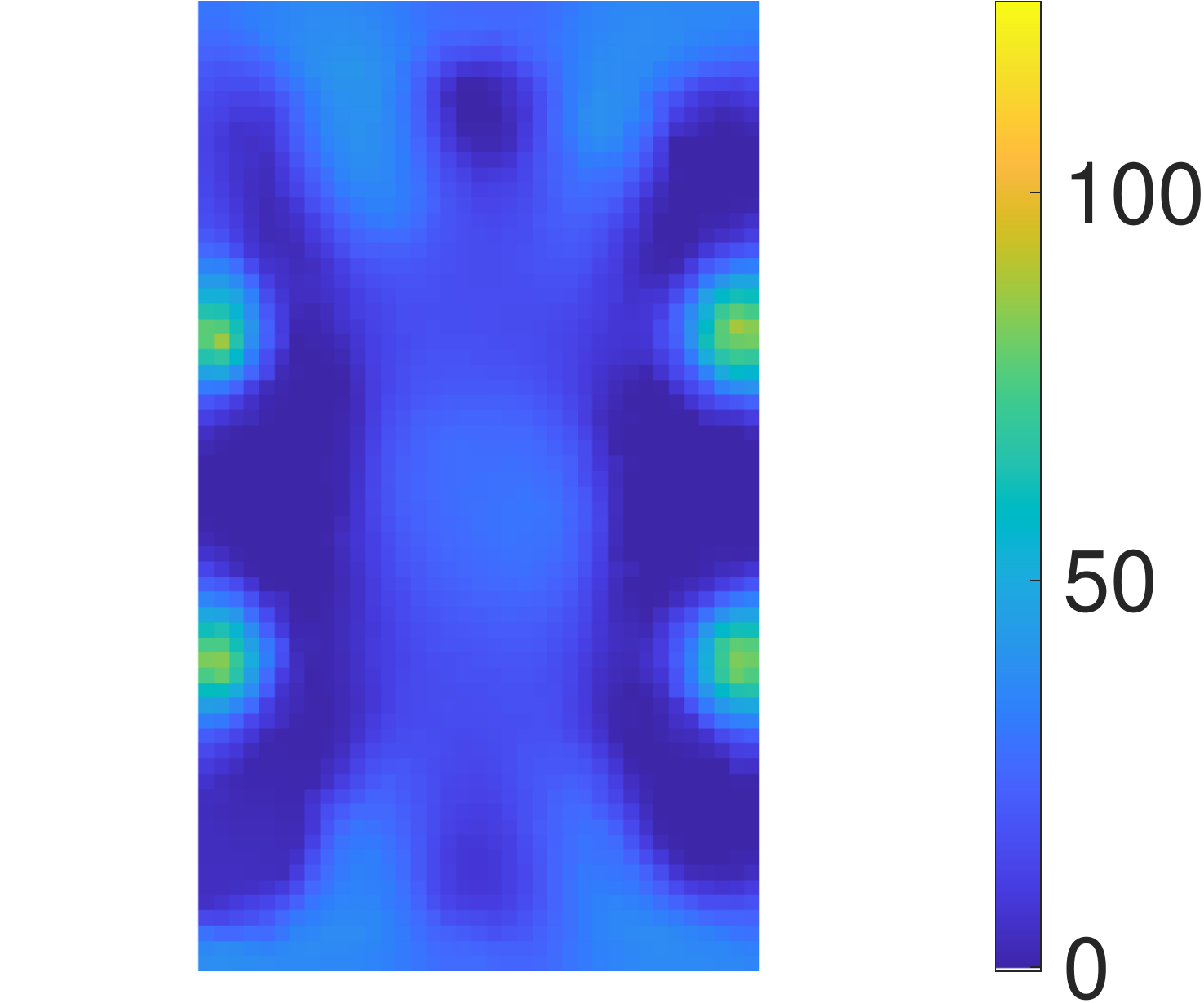}}

\put(3,0){\line(0,1){400}}
%\put(145,0){\line(0,1){400}}
\put(220,0){\line(0,1){400}}
%\put(445,0){\line(0,1){400}}
%----------------- 
% Top labels:
\put(20,385){{{\sc\footnotesize Slice}}}
\put(75,385){{{\sc\footnotesize Isosurface}}}
\put(160,385){{{\sc\footnotesize 3D}}}
\put(250,385){{{\sc\footnotesize Slice}}}
\put(310,385){{{\sc\footnotesize Isosurface}}}
\put(390,385){{{\sc\footnotesize 3D}}}

\put(68,400){{\underline{\sc One Target}}}
\put(298,400){{\underline{\sc Two Targets}}}
%\put(320,395){{\underline{\sc 20x35x25 modeling}}}

%----------------- 
% Side labels:
\put(-10,330){\rotatebox{90}{\sc Truth}}
\put(-10,250){\rotatebox{90}{\sc Cald}}
\put(-10,180){\rotatebox{90}{\sc $\texp$}}
\put(-10,105){\rotatebox{90}{\sc $\tzero$}}
\put(-10,25){\rotatebox{90}{\sc TV}}
%----------------- 

\put(445,0){\line(0,1){400}}

\end{picture}
\caption{\label{fig:abs_203525} {\bf Absolute image} reconstructions comparing the CGO methods to the regularized method with {\it largely incorrect domain modeling}, using a box of size 20cm x 35cm x 25cm.  Note the truth targets had a measured conductivity of approx 290 mS/m.}
\end{figure}

\begin{figure}[ht]
\centering
\begin{picture}(450,420)
\linethickness{.3mm}
%----------------- 
% truth
\put(70,320){\includegraphics[width=60pt]{truth_2targ_iso.png}}
\put(8,330){\includegraphics[angle=-90,origin=c,width=42pt]{TwoTargs_TopView.jpeg}}

\put(220,320){\includegraphics[width=60pt]{truth_2targ_iso.png}}
\put(158,330){\includegraphics[angle=-90,origin=c,width=42pt]{TwoTargs_TopView.jpeg}}

\put(370,320){\includegraphics[width=60pt]{truth_2targ_iso.png}}
\put(308,330){\includegraphics[angle=-90,origin=c,width=42pt]{TwoTargs_TopView.jpeg}}
%----------------- 
% cald
\put(360,240){\includegraphics[width=90pt]{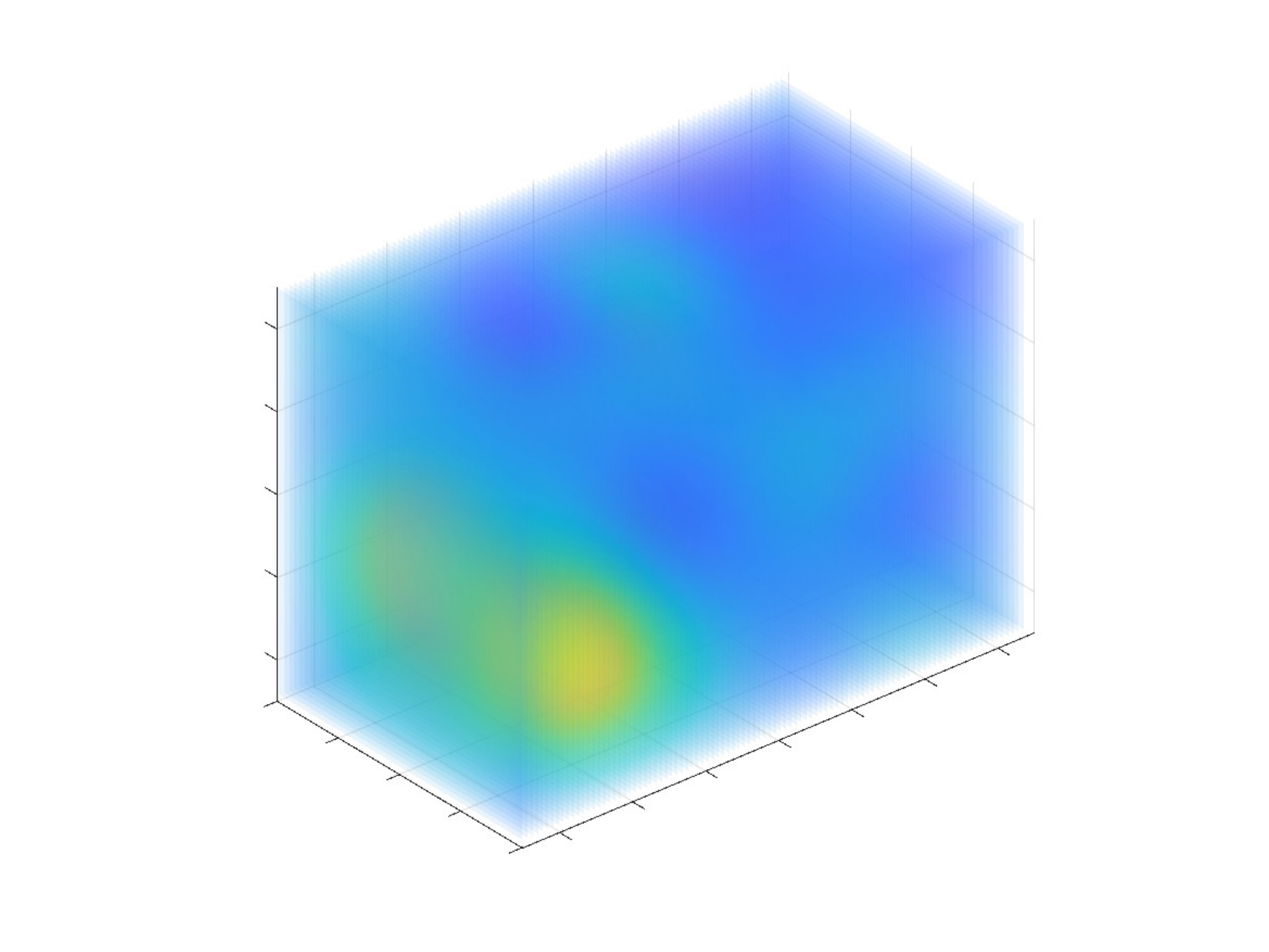}}
\put(300,240){\includegraphics[width=70pt]{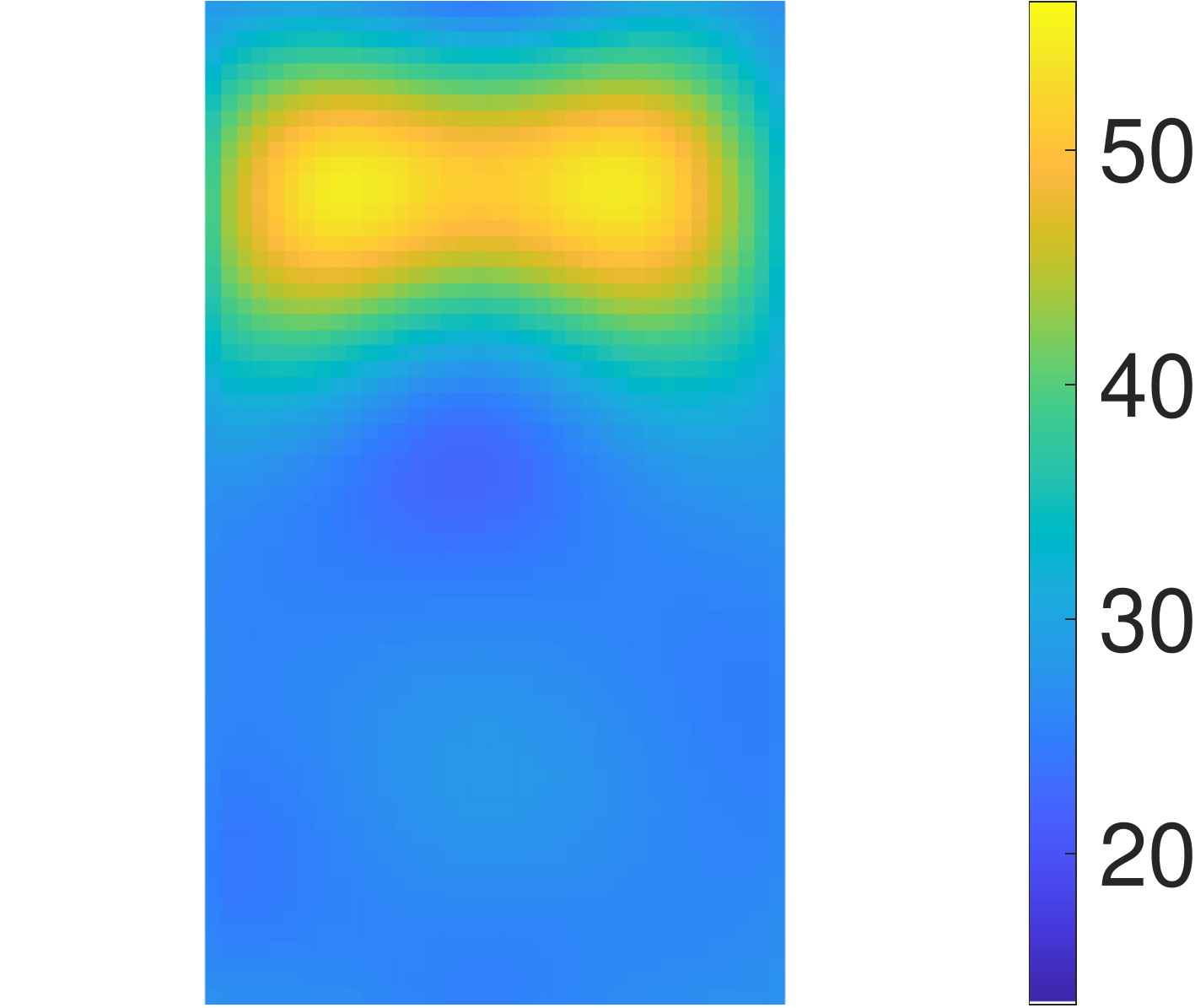}}

\put(210,240){\includegraphics[width=90pt]{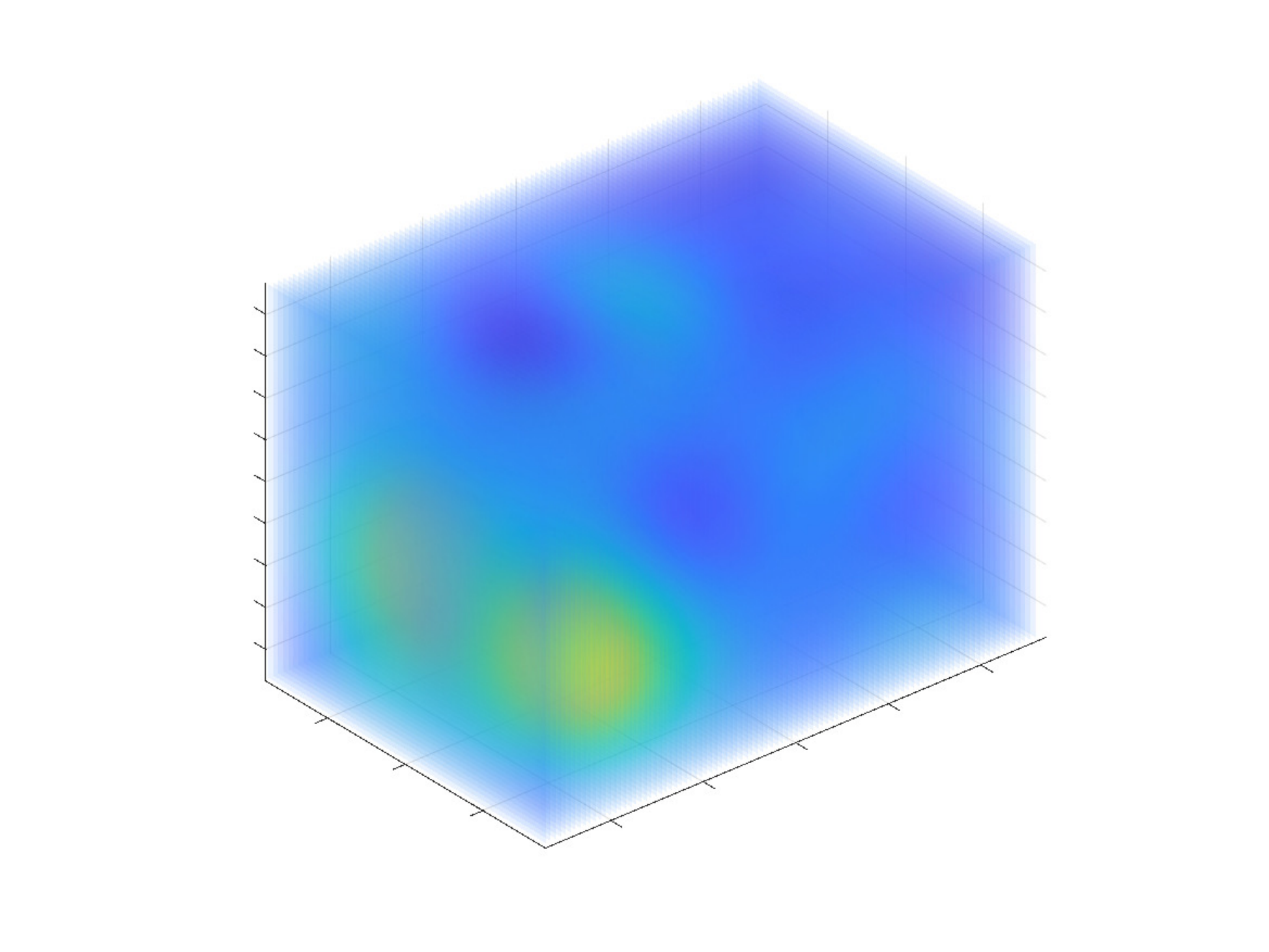}}
\put(150,240){\includegraphics[width=70pt]{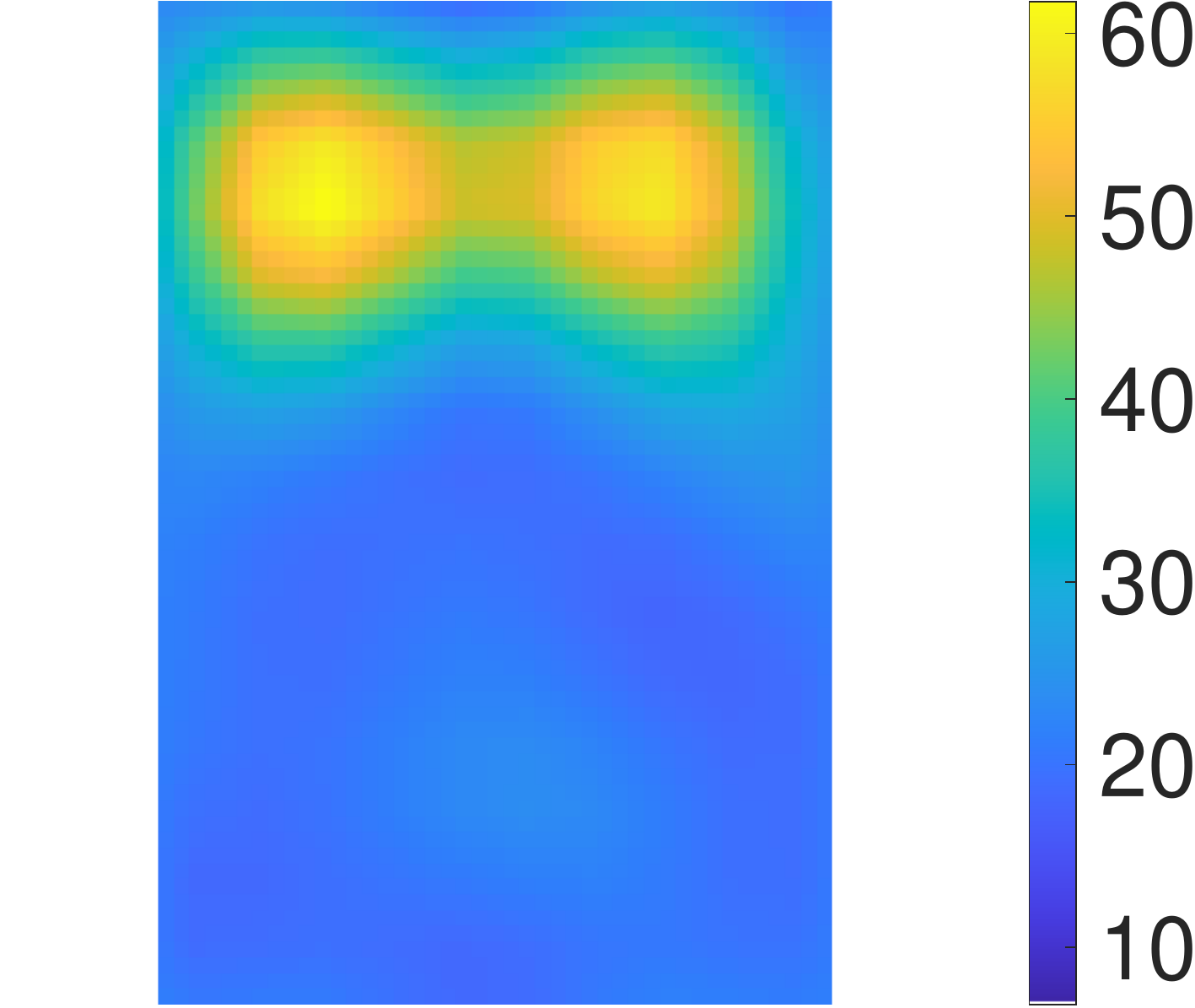}}

\put(60,240){\includegraphics[width=90pt]{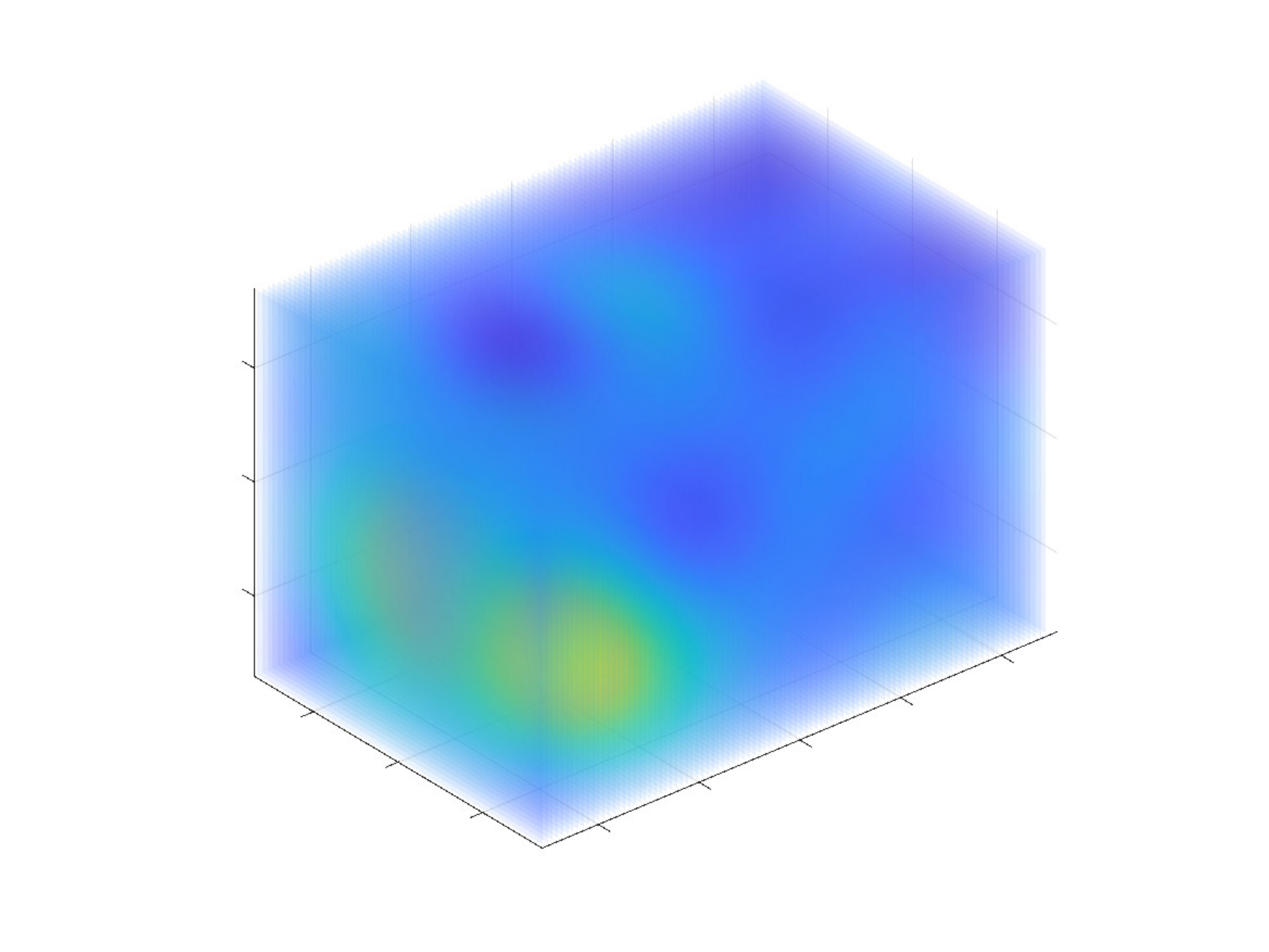}}
\put(0,240){\includegraphics[width=70pt]{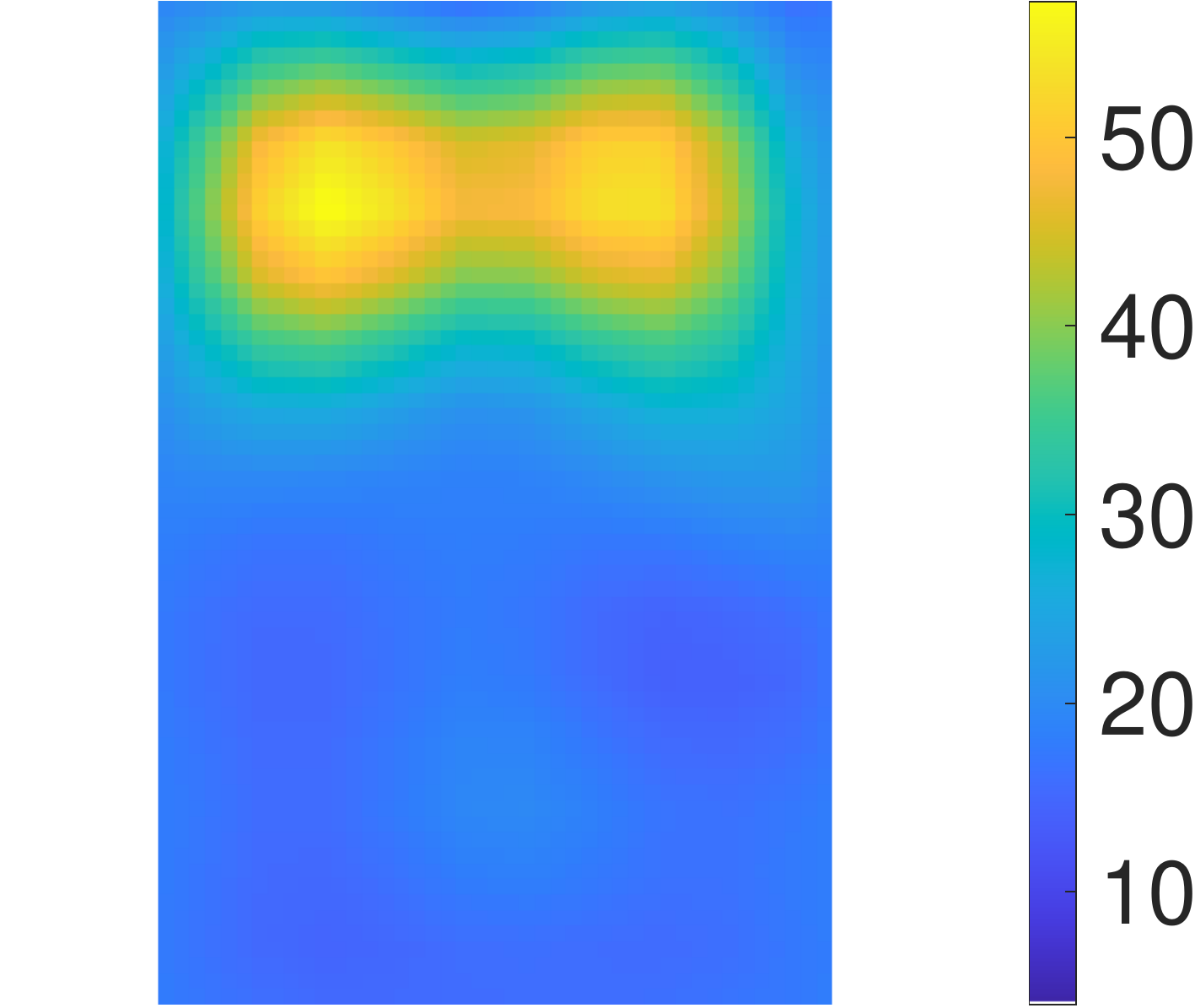}}
%----------------- 
% texp
\put(360,160){\includegraphics[width=90pt]{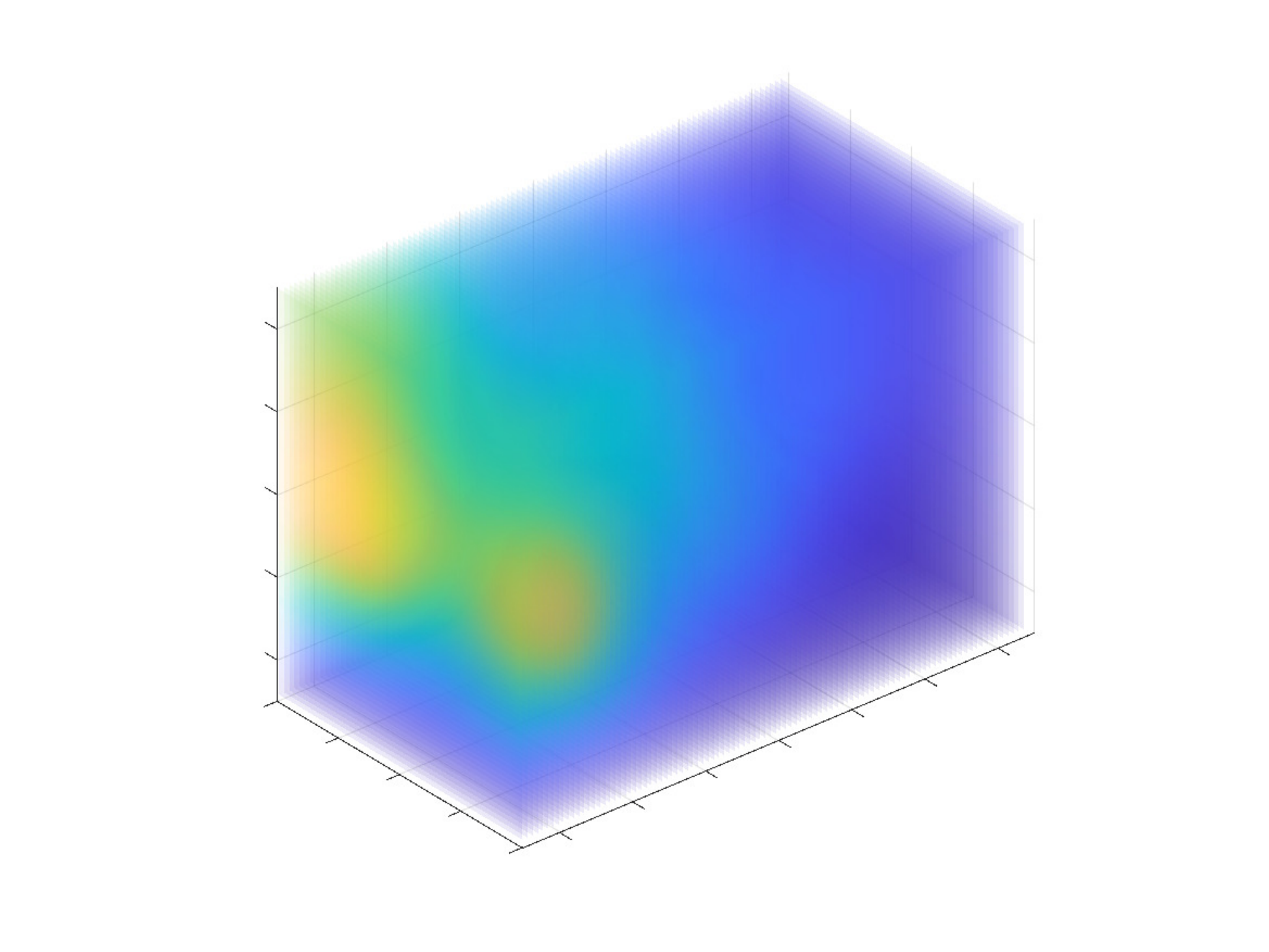}}
\put(300,160){\includegraphics[width=70pt]{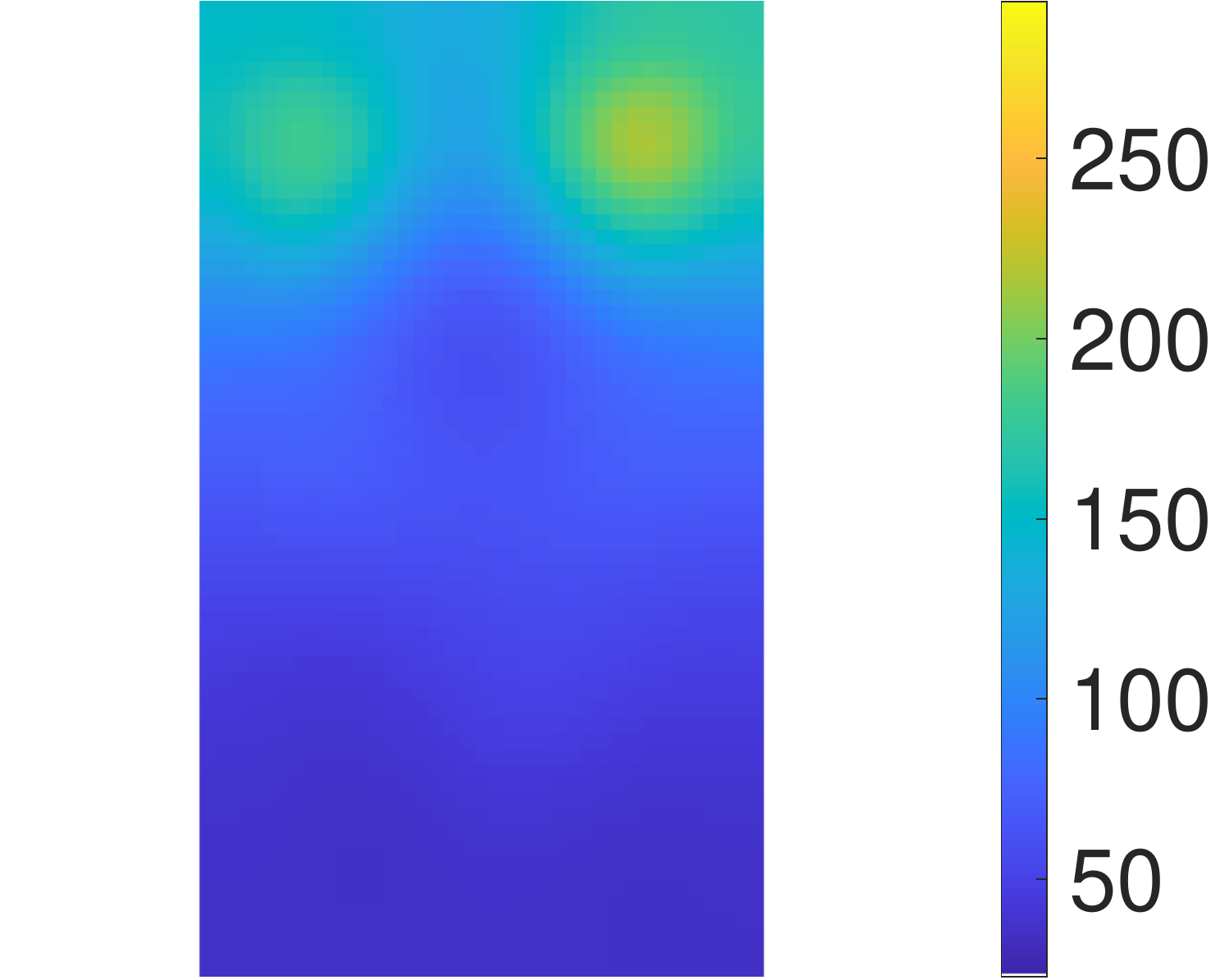}}

\put(210,160){\includegraphics[width=90pt]{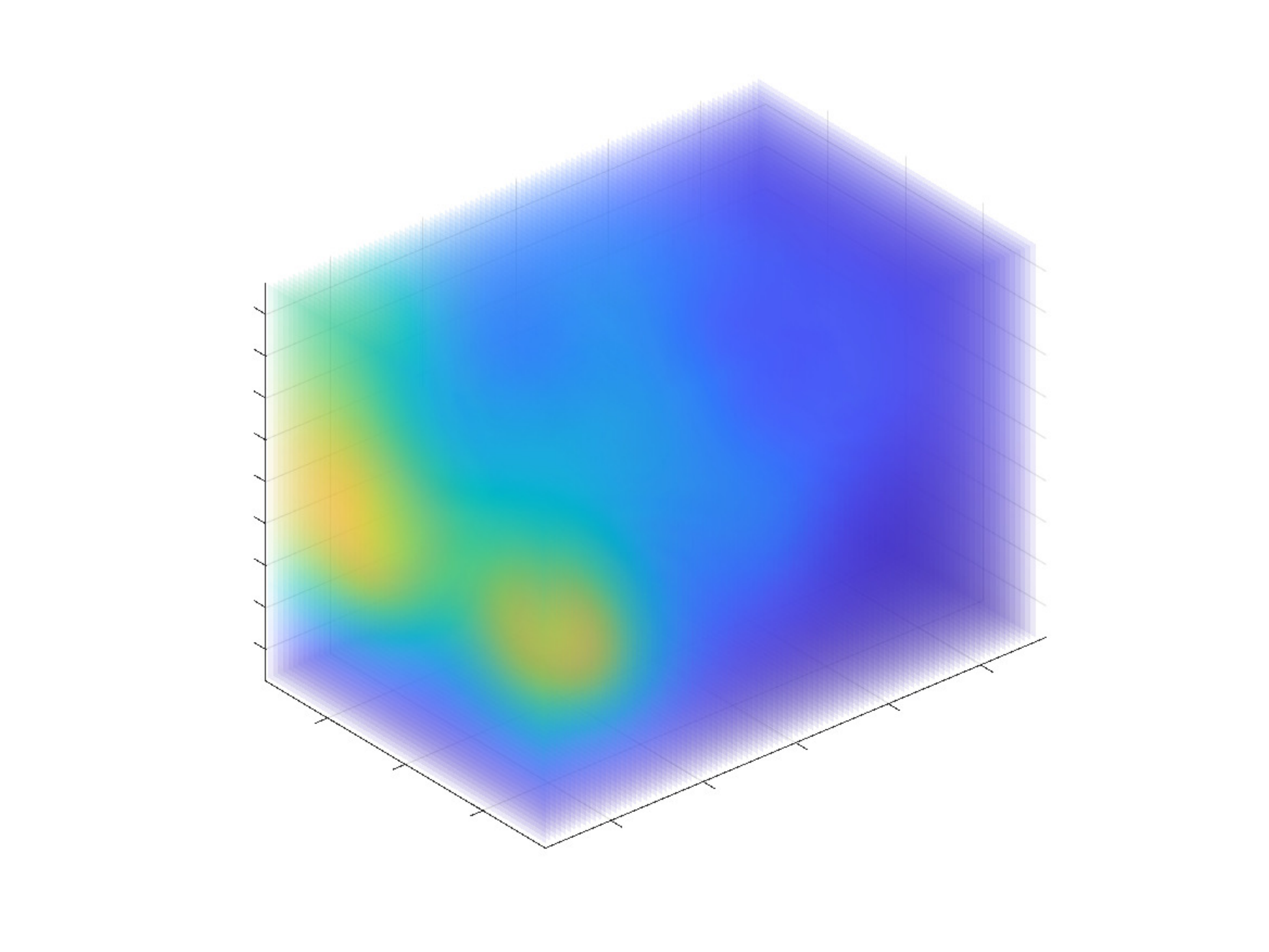}}
\put(150,160){\includegraphics[width=70pt]{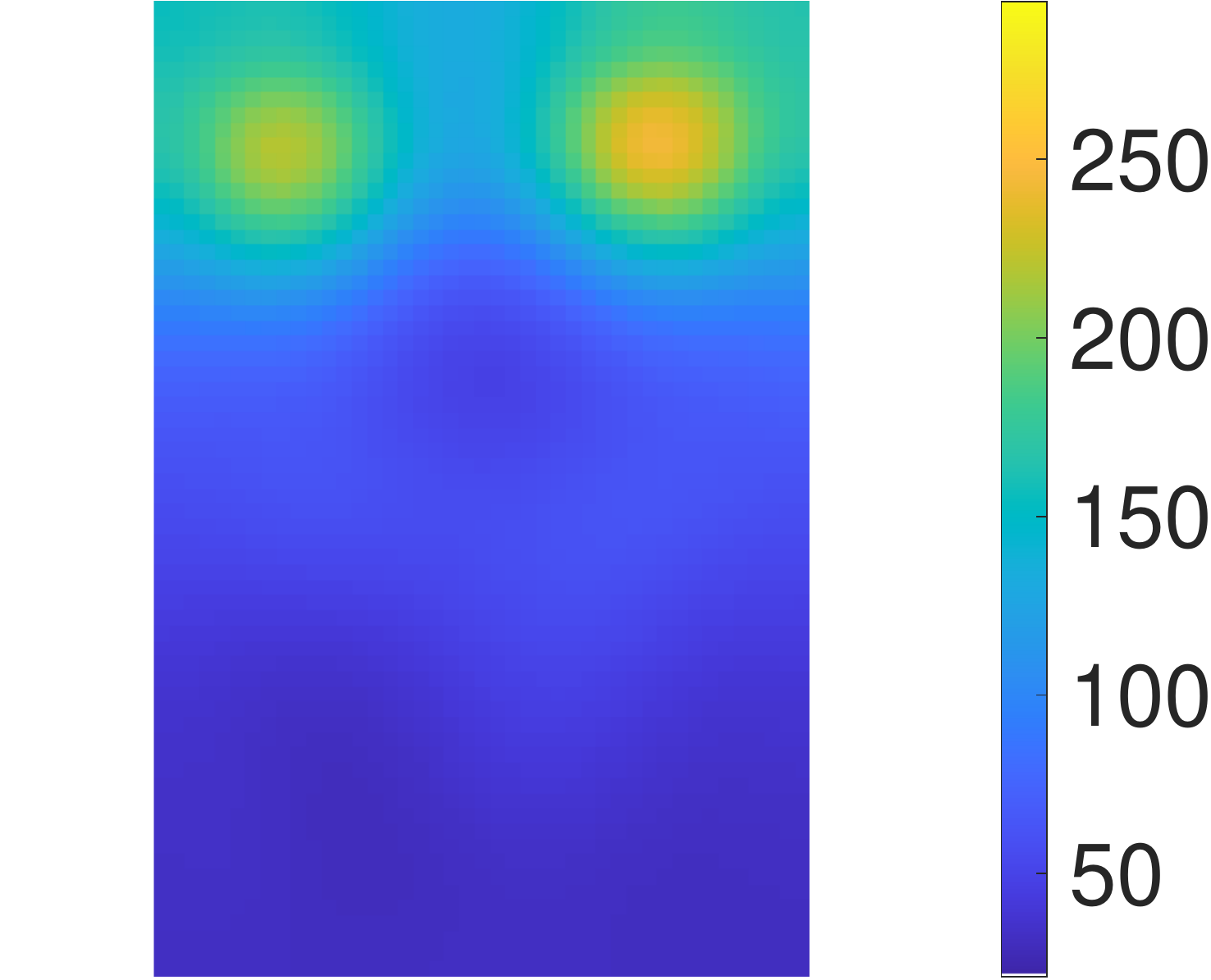}}

\put(60,160){\includegraphics[width=90pt]{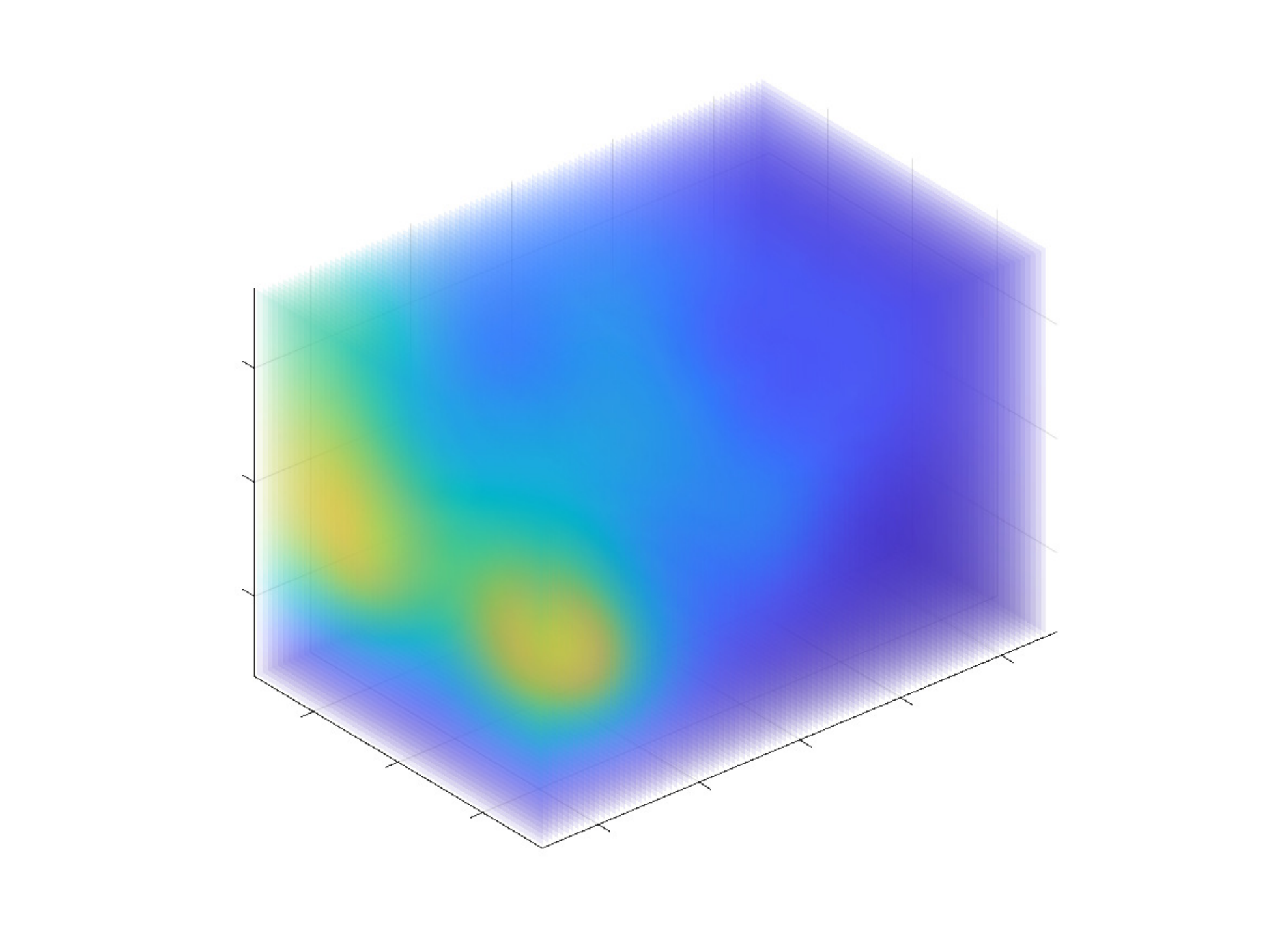}}
\put(0,160){\includegraphics[width=70pt]{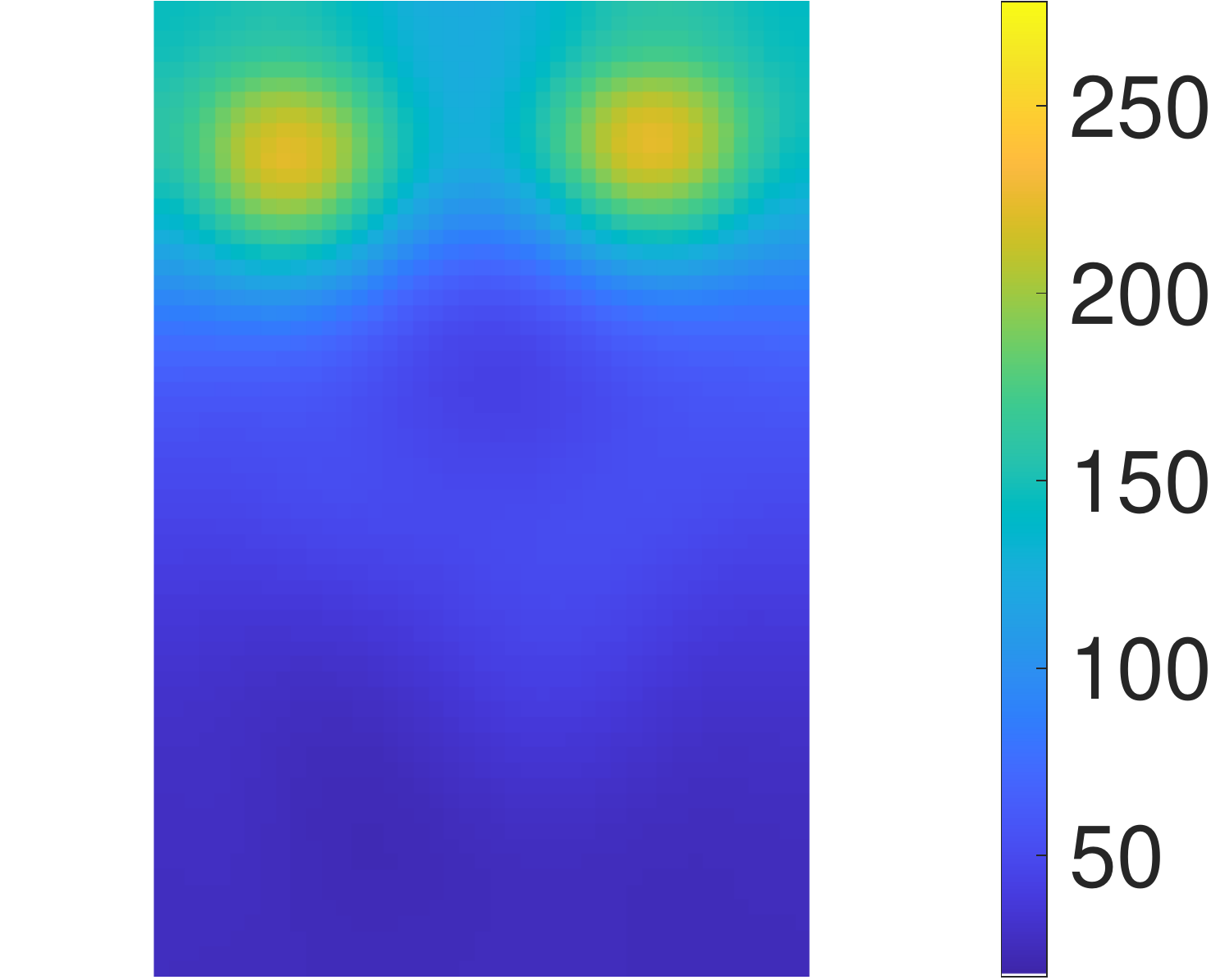}}

%----------------- 
% t0
\put(360,80){\includegraphics[width=90pt]{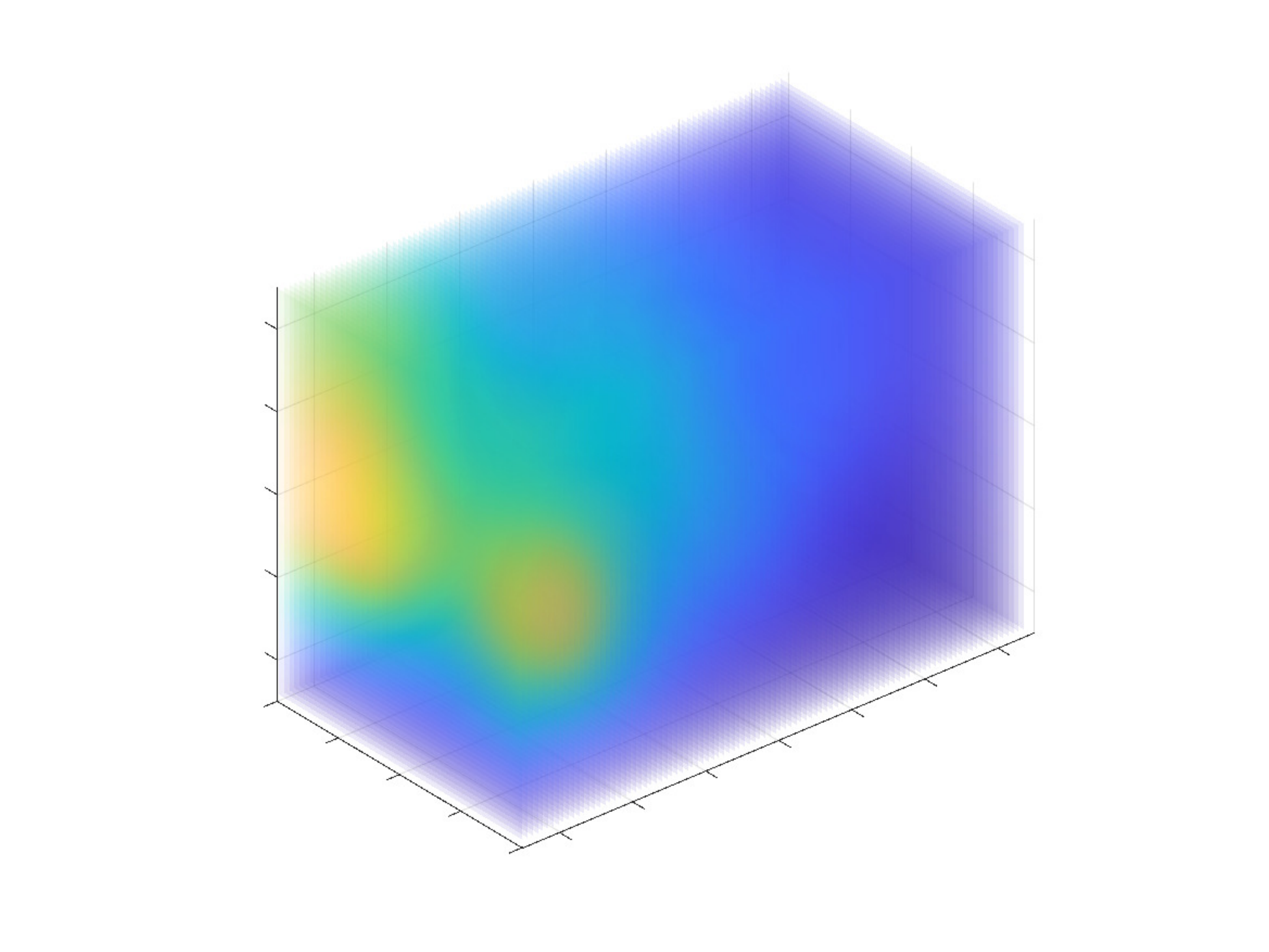}}
\put(300,80){\includegraphics[width=70pt]{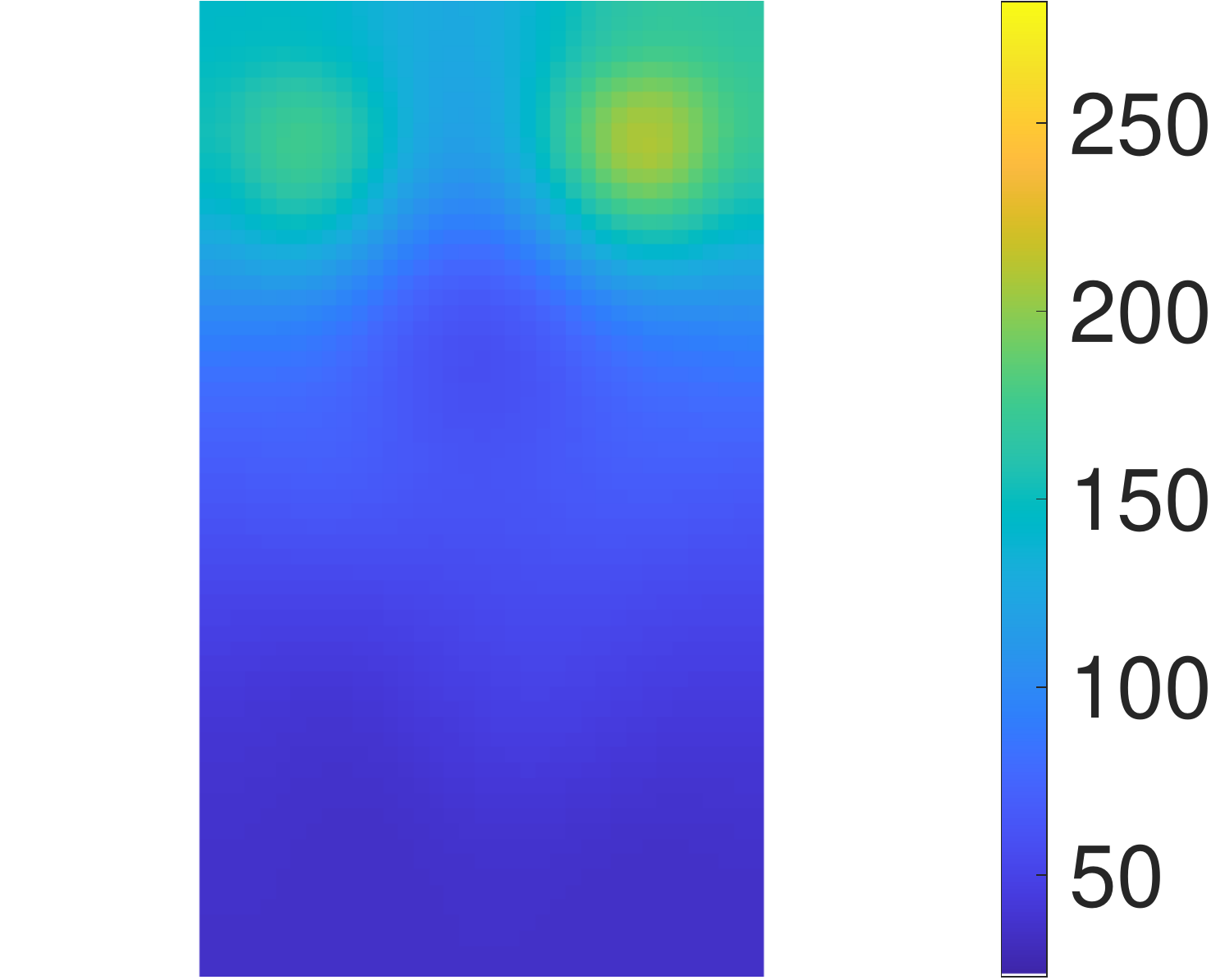}}

\put(210,80){\includegraphics[width=90pt]{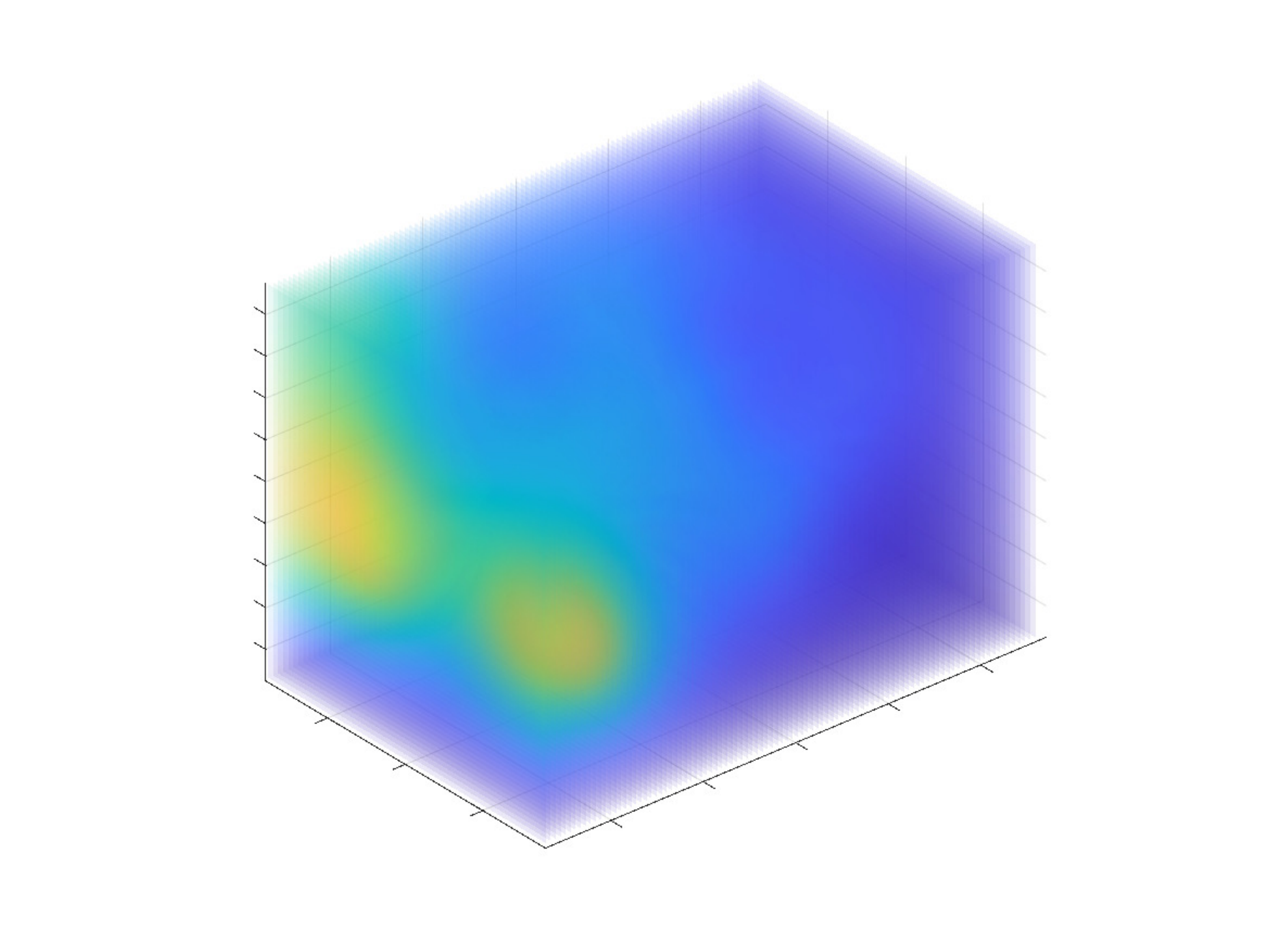}}
\put(150,80){\includegraphics[width=70pt]{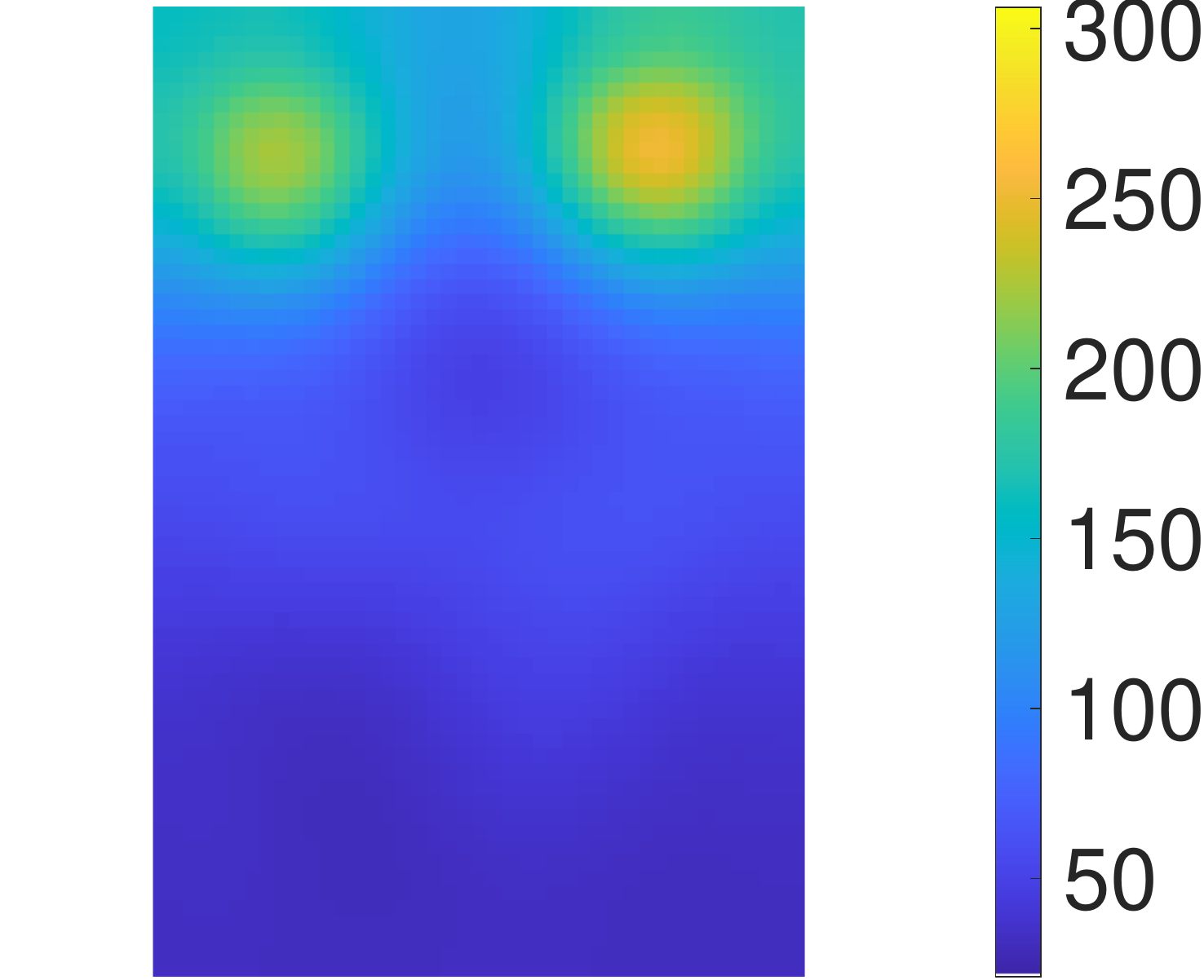}}

\put(60,80){\includegraphics[width=90pt]{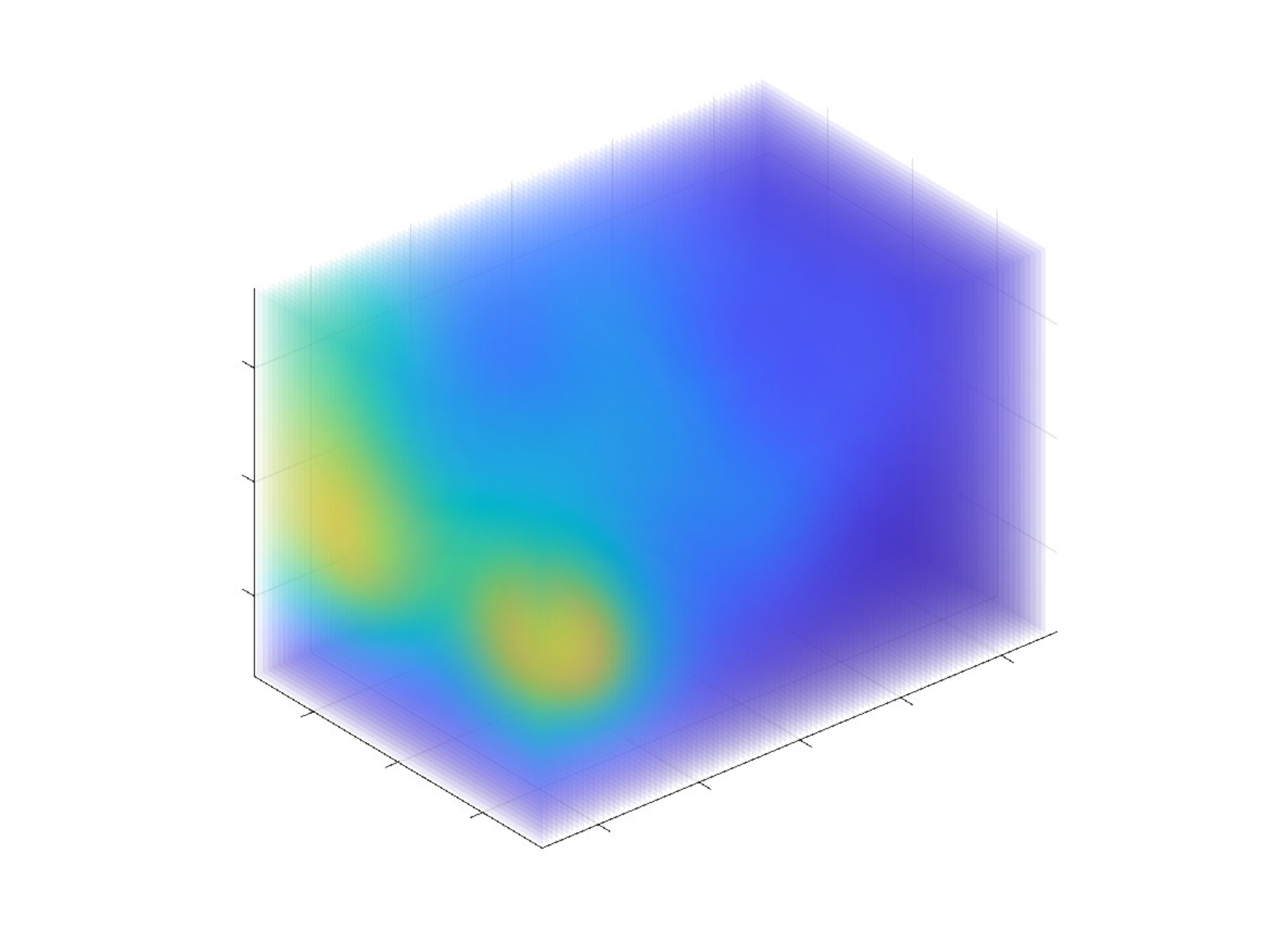}}
\put(0,80){\includegraphics[width=70pt]{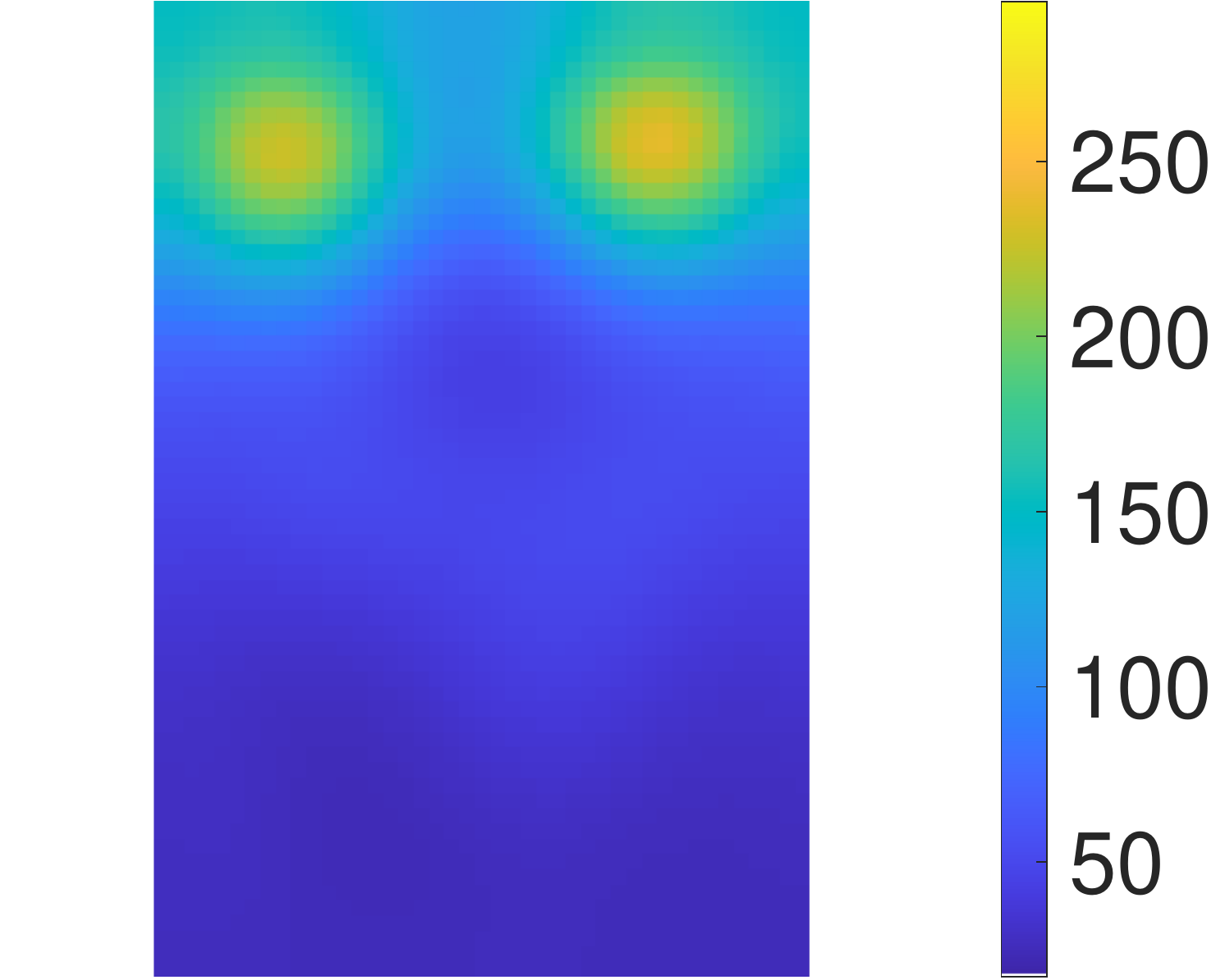}}

%----------------- 
% TV
\put(360,0){\includegraphics[width=90pt]{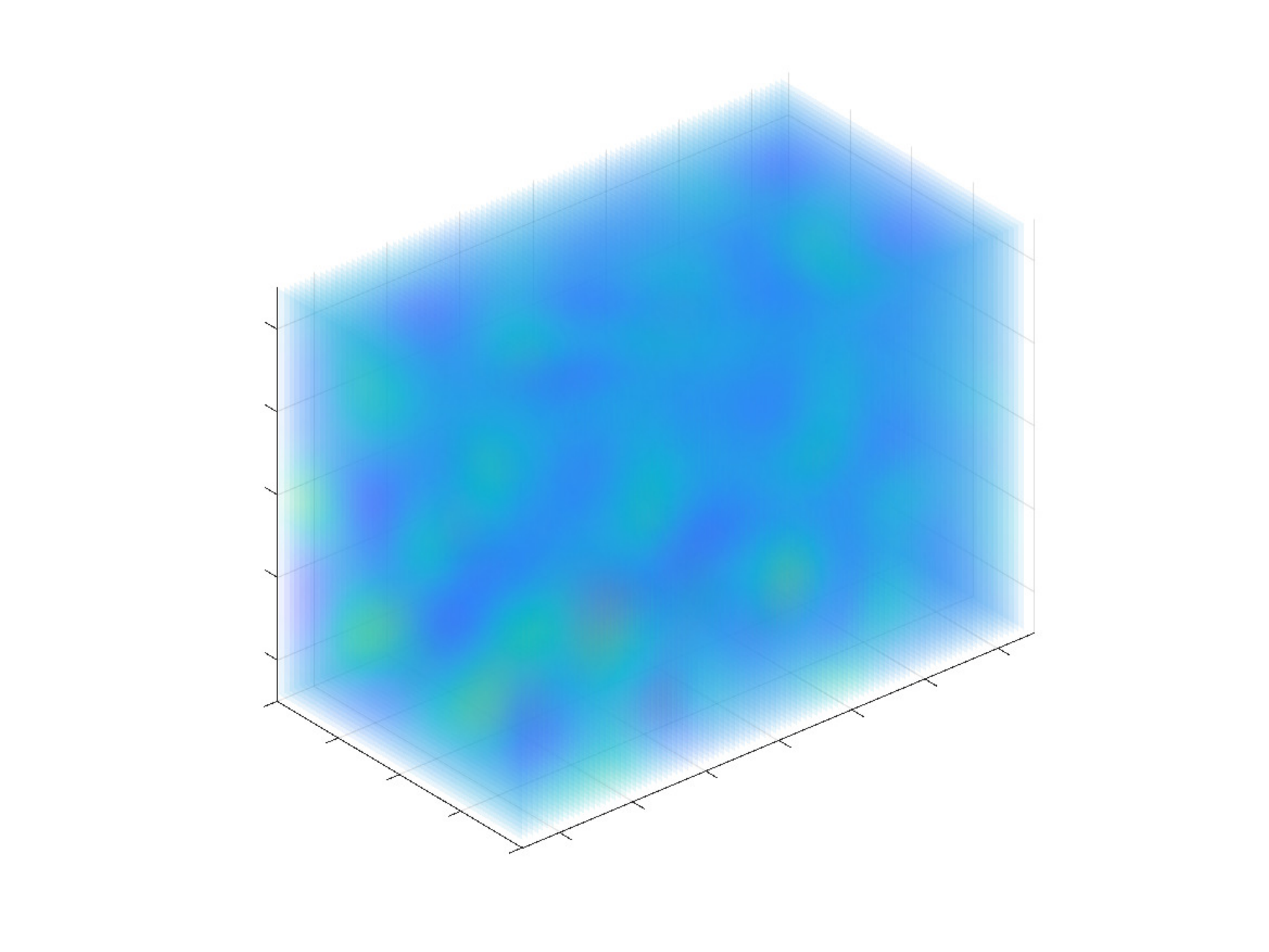}}
\put(300,0){\includegraphics[width=70pt]{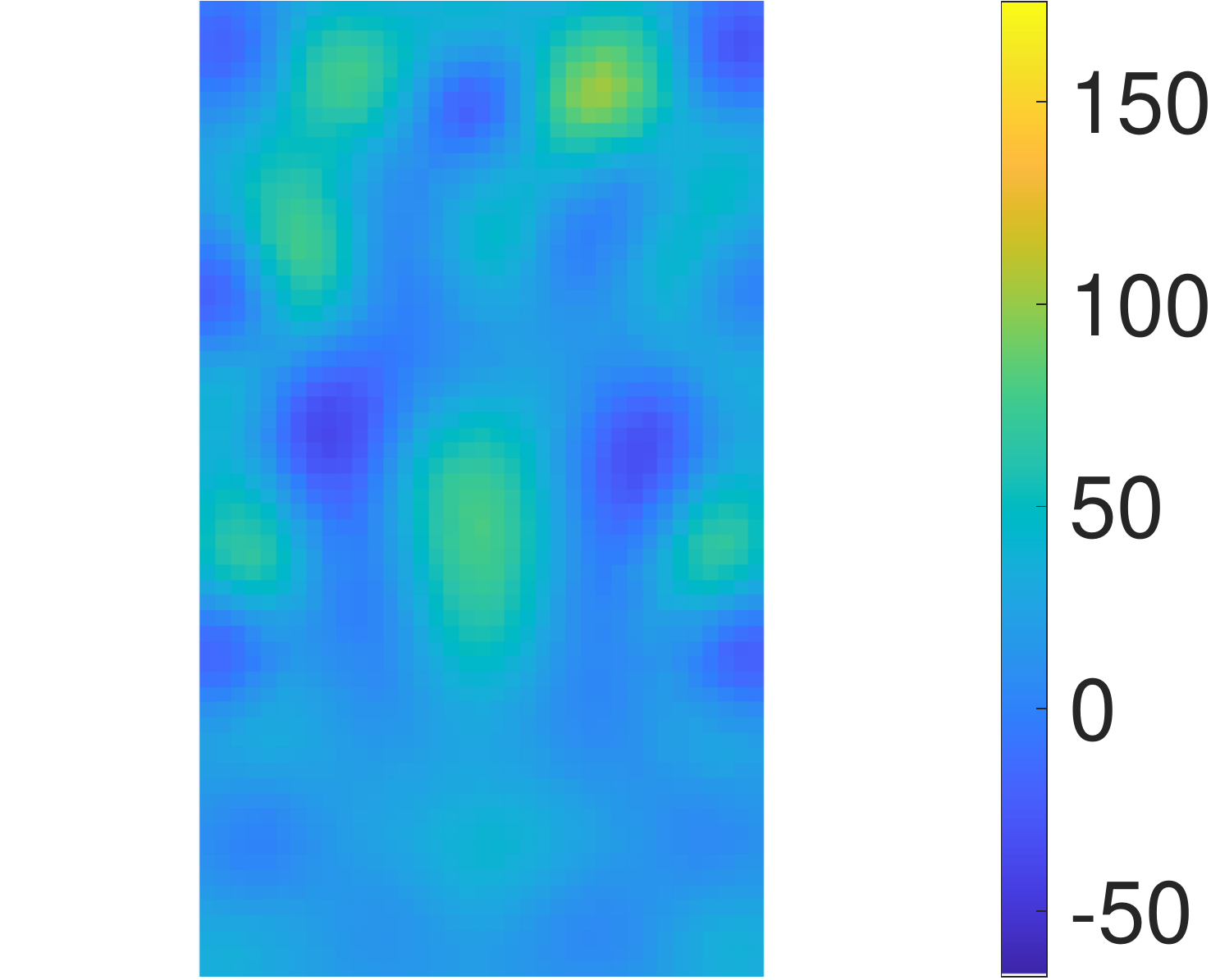}}

\put(210,0){\includegraphics[width=90pt]{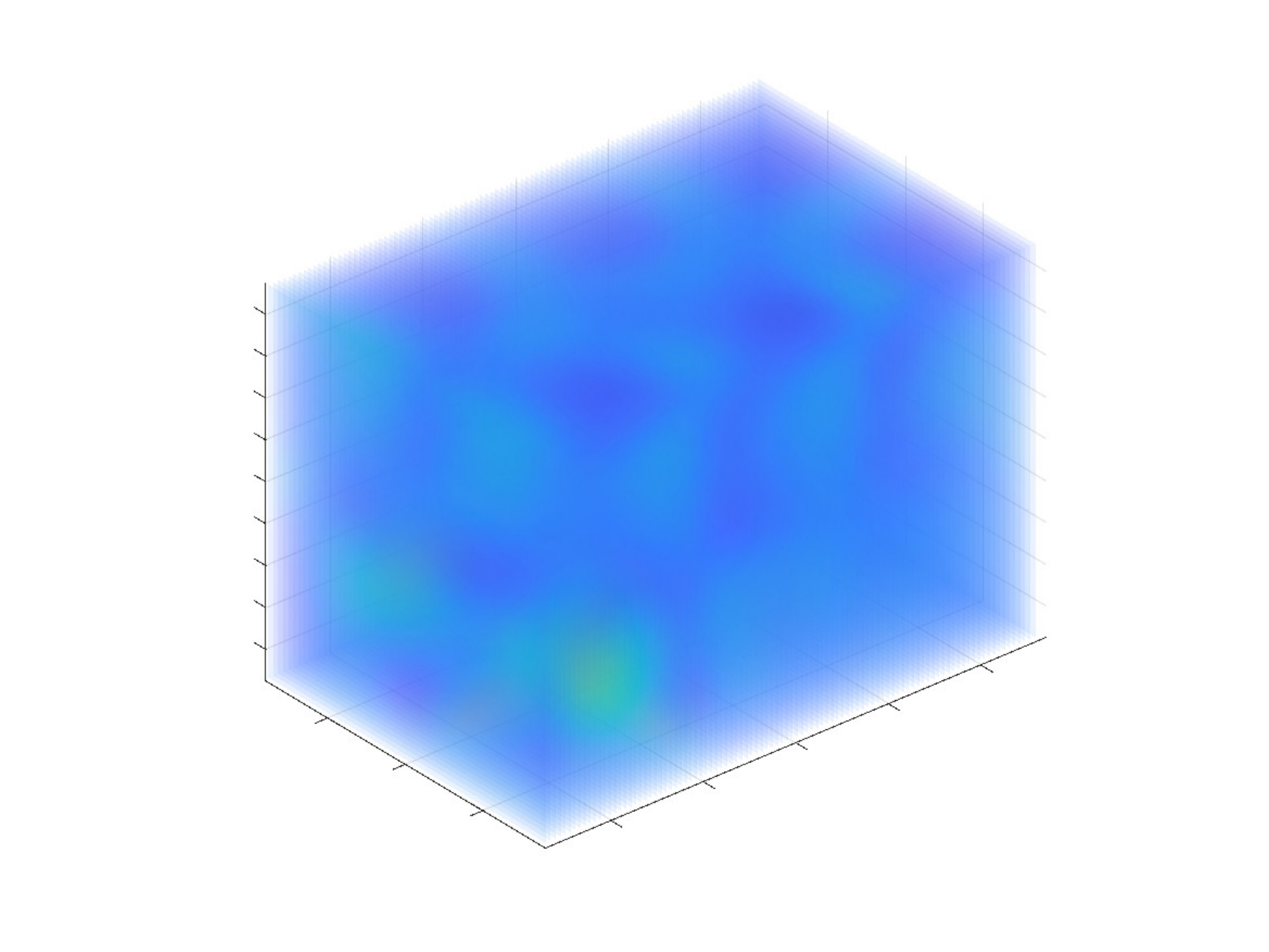}}
\put(150,0){\includegraphics[width=70pt]{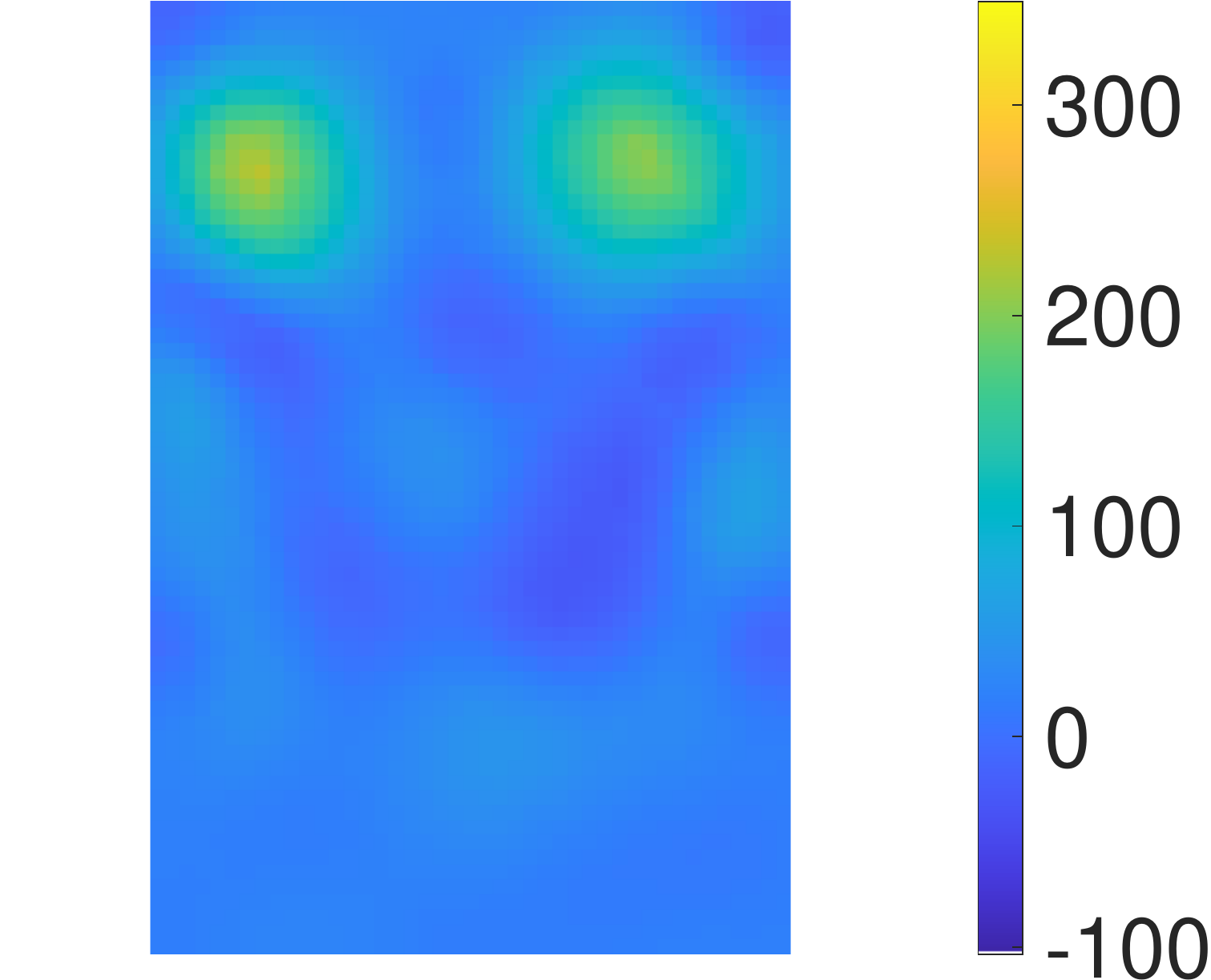}}

\put(60,0){\includegraphics[width=90pt]{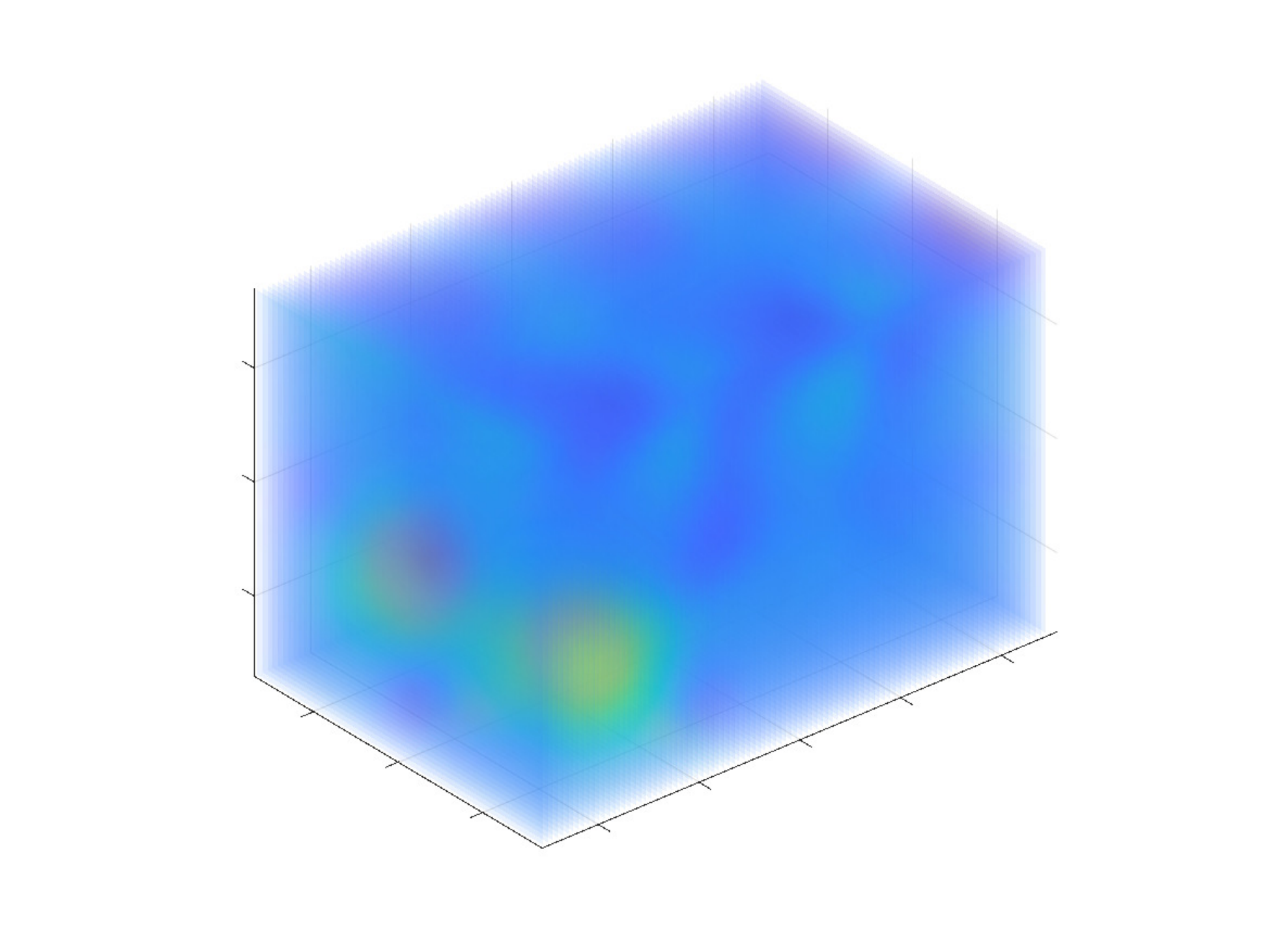}}
\put(0,0){\includegraphics[width=70pt]{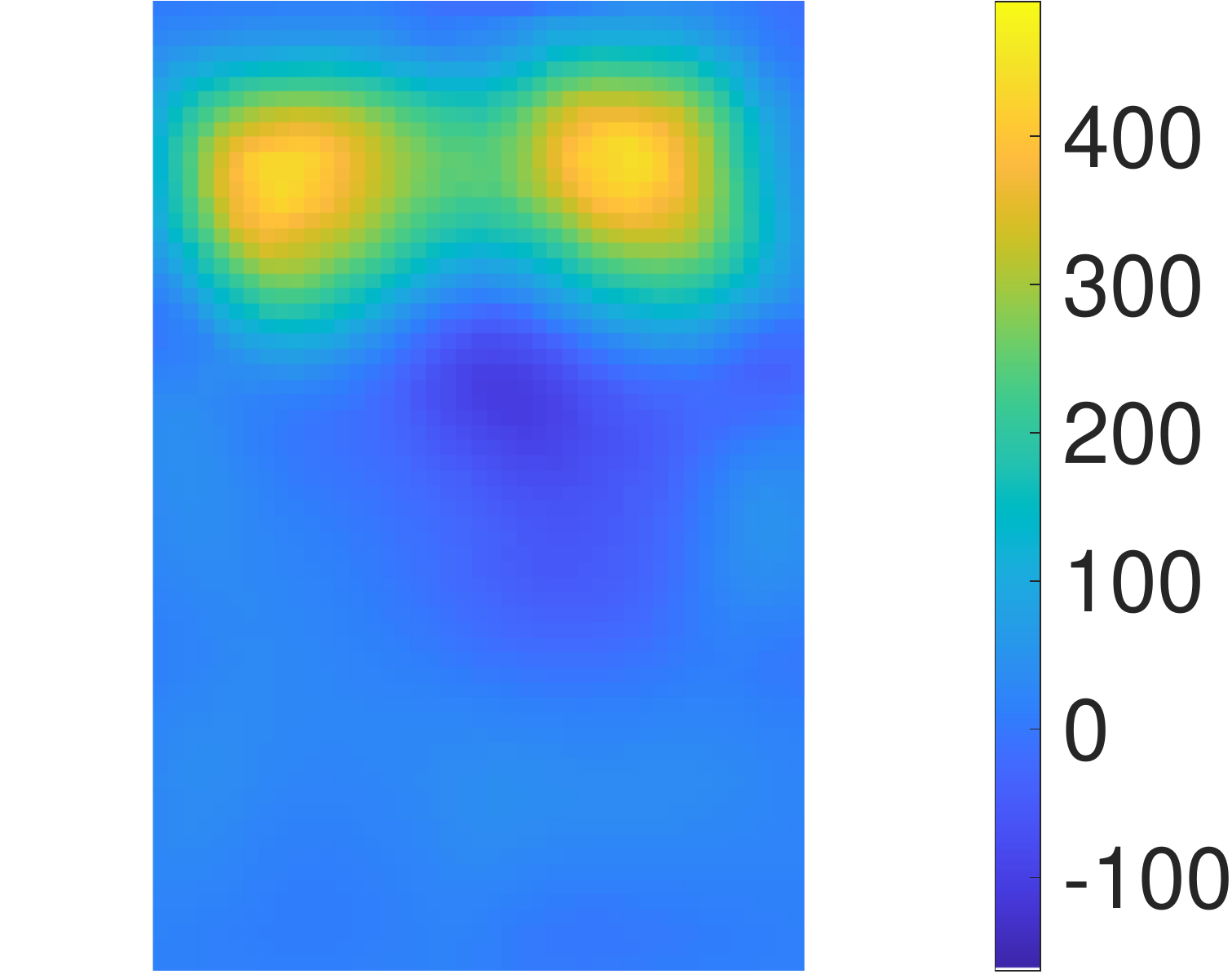}}

\put(3,0){\line(0,1){400}}
\put(145,0){\line(0,1){400}}
\put(295,0){\line(0,1){400}}
\put(445,0){\line(0,1){400}}
%----------------- 
% Top labels:
\put(20,395){{\underline{\sc Correct modeling}}}
\put(170,395){{\underline{\sc 18x27x19 modeling}}}
\put(320,395){{\underline{\sc 20x35x25 modeling}}}
%\put(170,395){\sc }
%\put(270,395){\sc 3D}
%----------------- 
% Side labels:
\put(-10,330){\rotatebox{90}{\sc Truth}}
\put(-10,250){\rotatebox{90}{\sc Cald}}
\put(-10,180){\rotatebox{90}{\sc $\texp$}}
\put(-10,105){\rotatebox{90}{\sc $\tzero$}}
\put(-10,15){\rotatebox{90}{\sc Linear}}
%----------------- 

\put(445,0){\line(0,1){400}}

\end{picture}
\caption{\label{fig:2targ_DIFF} {\bf Difference image} reconstructions comparing the CGO methods to a typical linear method.  Slices and 3D renderings of the conductivity are shown for the correct domain modeling, and increasing levels of error in domain modeling.  Note the truth targets had a measured conductivity difference from the background of approx 266 mS/m.}
\end{figure}

Figure~\ref{fig:abs_182719} shows the absolute EIT reconstructions for the moderately mismodelled domain whereas Figure~\ref{fig:abs_203525} shows the absolute EIT reconstructions for the significantly mismodelled domain.  Difference images are shown in Figure~\ref{fig:2targ_DIFF}.  Metrics of maximum conductivity per target and \trev{scaled} localization error are shown in Tables~\ref{tab:1Targ_errs} and \ref{tab:2Targ_errs} for the absolute images and Table~\ref{tab:2Targ_errs_diff} for the difference images. Bolded entries correspond to the best metric value in each row.

% ---------------------------------------------------------------------------------
%%%%%%%%%%%%%%%%%%%%%%%%%%%%%%%%%%%%
% Correct modeling case - new format with BOTH 1 and 2 targets here.
%%%%%%%%%%%%%%%%%%%%%%%%%%%%%%%%%%%%
%%%%%%%%%%%%%%%%%%%%%%%%%%%%%%%%%%%%
% ---------------------------------------------------------------------------------

% ---------------------------------------------------------------------------------
%%%%%%%%%%%%%%%%%%%%%%%%%%%%%%%%%%%%
% 18x27x19 modeling case - new format with BOTH 1 and 2 targets here.
%%%%%%%%%%%%%%%%%%%%%%%%%%%%%%%%%%%%
%%%%%%%%%%%%%%%%%%%%%%%%%%%%%%%%%%%%
% ---------------------------------------------------------------------------------

% ---------------------------------------------------------------------------------
%%%%%%%%%%%%%%%%%%%%%%%%%%%%%%%%%%%%
% 20x35x25 modeling case - new format with BOTH 1 and 2 targets here.
%%%%%%%%%%%%%%%%%%%%%%%%%%%%%%%%%%%%
%%%%%%%%%%%%%%%%%%%%%%%%%%%%%%%%%%%%
% ---------------------------------------------------------------------------------

% \subsection{Difference Images}
% ---------------------------------------------------------------------------------
%%%%%%%%%%%%%%%%%%%%%%%%%%%%%%%%%%%%
% 2 Target - New format comparing TRUE, 18x27x19, and 20x35x25
%%%%%%%%%%%%%%%%%%%%%%%%%%%%%%%%%%%%
%%%%%%%%%%%%%%%%%%%%%%%%%%%%%%%%%%%%
% ---------------------------------------------------------------------------------

 %\clearpage

%--------------------------------------------------------------------
\section{Discussion}\label{sec:discussion}
%--------------------------------------------------------------------
Beginning with the correct domain modeling scenario, each method, Calder\'on, $\texp$, $\tzero$, and TV produced reconstructions that clearly show the target(s) with good localization. The TV images are sharpest, as expected.  The $\texp$ and $\tzero$ methods achieved the best contrast and best approximated the estimated experimental conductivity.  The Calder\'on and TV methods achieved lower contrast but better \trev{scaled} localization error.  Notably, as shown in Figure~\ref{fig:2targ_DIFF}, in terms of artefacts, the absolute CGO reconstructions are as clean as their corresponding difference images.  The CGO methods required simulated voltages data for the same experimental setup but with a conductivity $\sigma\equiv1$ S/m.  Even though this data was not tuned to the EIT machine (noise and contact impedances) the methods produced high quality absolute reconstructions. 

Moving into incorrect domain modeling, we see all methods handled the moderate domain mismodeling quite well and the CGO methods handled the severe domain mismodeling as well.  No significant boundary artefacts are seen in the CGO reconstructions for the 18cm x 27cm x 19cm and only moderate boundary artefacts in the TV reconstructions.  In the more severe mis-modeling case, assuming the data is coming from a much larger box (20cm x 35cm x 25cm) instead of the true domain (17cm x 25.5cm x 17cm), the $\texp$ and $\tzero$ methods still produce good reconstructions showing the correct number of targets and correct region of the tank.  The localization is worst in the $x_3$ direction showing the targets slightly elongated.  While artefacts in the Calder\'on 20x35x25 reconstructions are present, the method still does detect the true objects as the most conductive.  The TV algorithm failed to identify any of the targets for both cases with the 20x35x25 box.
\trev{The artefacts in the TV reconstructions in these cases with incorrect domain model could be expected as  
the regularized non-linear least squares minimization-based absolute imaging EIT algorithms are known to be highly sensitive to modeling errors such as errors in modelling of the domain shape, see e.g. \cite{Nissinen10compensation,kolehmainen2008electrical,Lehikoinen07approximation}.}

We can see the effect of the underlying assumption for Calder\'on's method that the conductivity is a small perturbation from a constant.  In the physical experiments, the conductor(s) had roughly twelve times the conductivity of the background saline. This is seen in the underestimation of maximum conductivity from Calder\'on's method.  However, unsurprisingly, this method does a good job of localizing the targets, in all cases except for the absolute image from the significant domain modeling error.  The artefacts seen in the Calder\'on reconstructions are consistent with what has been observed for absolute imaging with Calder\'on in 2D, such as Gibbs phenomenon for all absolute reconstructions and higher reconstructed conductivity at electrode locations when the domains were incorrectly modelled.  Across all modeling scenarios, the regularized TV method struggled with the high contrast targets, giving maximum target conductivities well below the measured values. The same effect was seen in simulated scenarios where, however, maximum conductivities of lower contrast targets were reliably recovered. This is at least partially explained by saturation of contrast distinquishability of the measurements, i.e. the effect of increasing the conductivity of the target inclusion on the measurements gradually diminishes and until certain level the TV regularized method can no longer discern between high and even higher conductivities.

We note that the electrodes used in this experiment were very large and the structure of the domain, a box with corners, may exacerbate some of the modeling and/or hardware challenges.  Nevertheless, the study provides informative results on the feasibility of absolute EIT reconstruction in 3D.

The difference images from Calder\'on are able to handle the stronger mismodeling of the domain, as are the $\texp$ and $\tzero$ methods.  The strong mismodeling proved too severe for the linear difference imaging reference method, which did not manage to identify the targets, see Figure~\ref{fig:2targ_DIFF}.

\begin{figure}[h]
    \centering
    \begin{picture}(320,150)
    \put(0,10){\includegraphics[height=120pt]{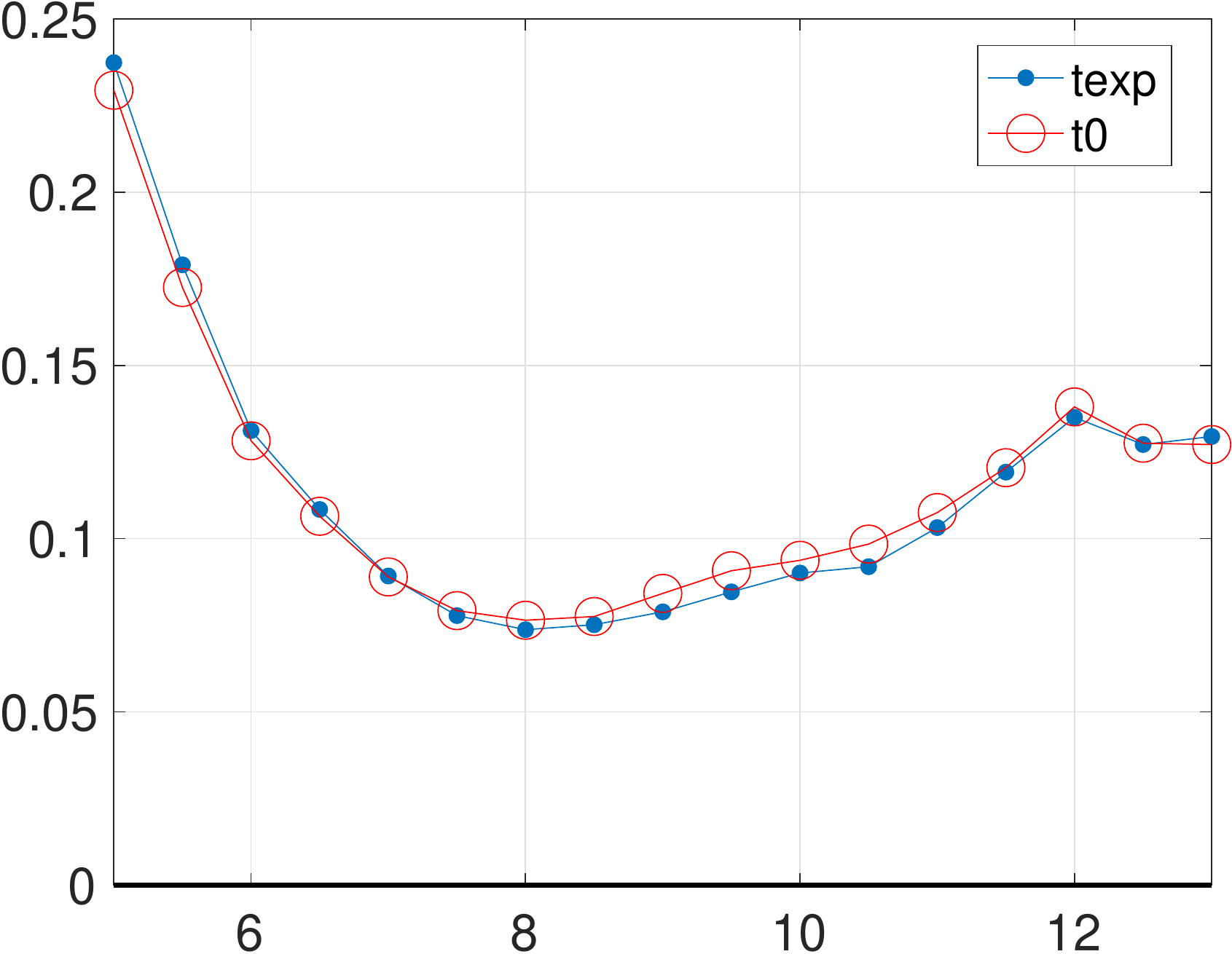}}
    \put(180,10){\includegraphics[height=120pt]{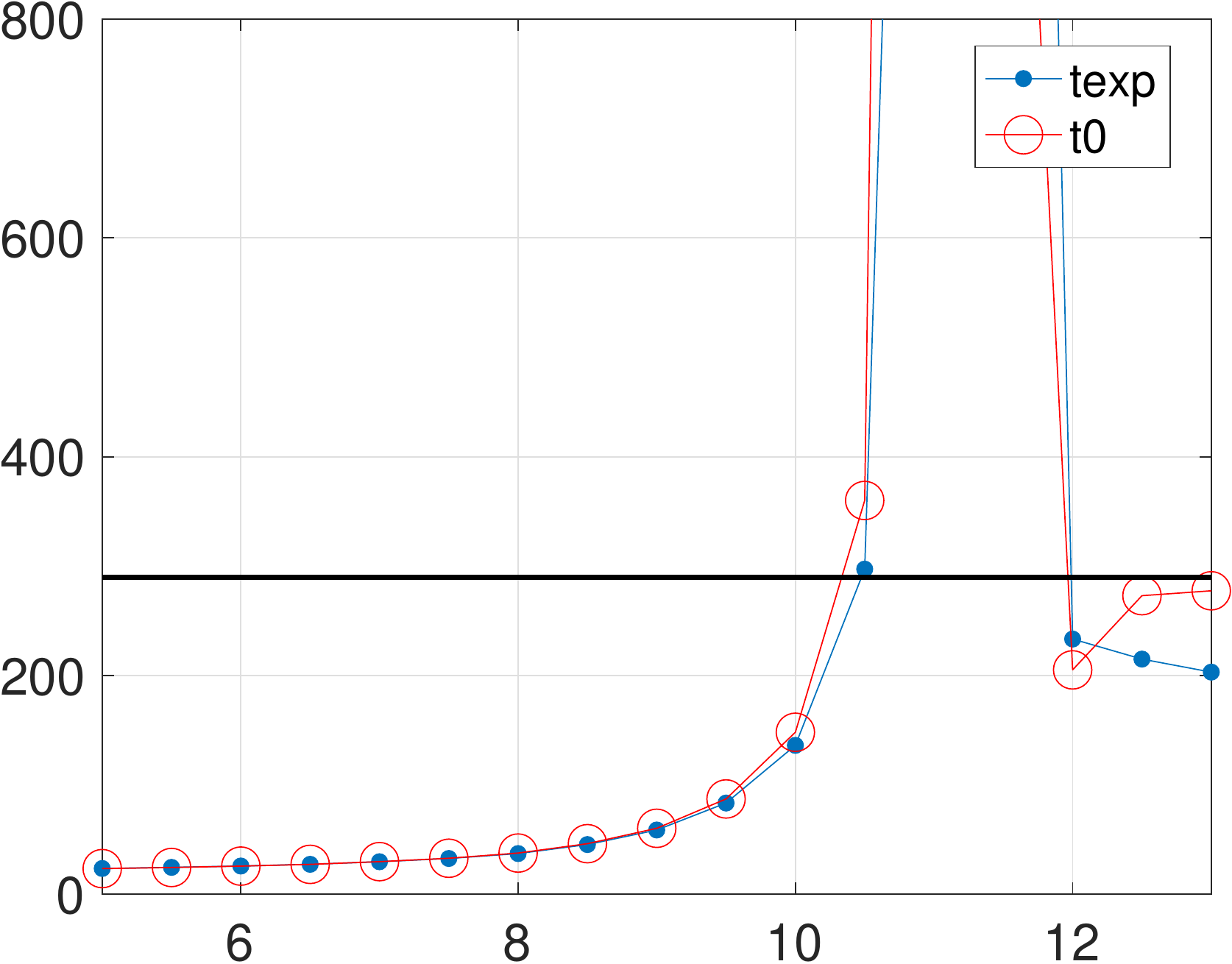}}
    
    \put(15,135){\footnotesize \sc \trev{Scaled Localization Error}}
    \put(200,135){\footnotesize \sc Max Target Cond. (mS/m)}
    
    \put(65,0){\scriptsize $T_\xi$}
    \put(250,0){\scriptsize $T_\xi$}
    
    % \put(-10,50){\scriptsize $M_\xi$}
    
    \end{picture}
    \caption{Comparison of the effect of the truncation value $T_\xi$ of the scattering radius in the $\texp$ and $\tzero$ CGO methods for \trev{Scaled} Localization Error (left) and Maximum value of the recovered target.  Max conductivity values (right) for $T_\xi=11,\;11.5$, off the plot, spiked into the 1500-5000~mS/m range. }
    \label{fig:Mxi_sweep_texp_t0}
\end{figure}

While the $\texp$ and $\tzero$ CGO methods did reliably recover the contrast and approximate location of the targets across examples studied here, they do appear more sensitive than their 2D \trev{D-bar based} counterparts in regards to the regularization parameter, $T_\xi$, used in the truncation of the nonlinear scattering data.  Figure~\ref{fig:Mxi_sweep_texp_t0} displays the effect that $T_\xi$ plays on the \trev{scaled} localization error and maximum recovered conductivity value for the single target, correct domain modeling case. Recall that a secondary nonuniform truncation is also enforced where scattering data with magnitudes exceeding 20 for the real or imaginary parts are set to zero.  Adjusting that value will also have an effect on the reconstruction.  The value of 20 was chosen in this work for its overall reliability across examples. The contrast appears more sensitive than the localization error.  As in \cite{Delbary2014}, the minimum $\zeta$ parameterization of $\xi$ was used, as the scattering data is a function of both $\zeta\in\C^3$ and $\xi\in\R^3$ in 3D instead of just $k\in\C$ as in 2D \trev{D-bar based methods}.  Alternative parameterizations and a more detailed study of the effect of the regularization parameters, while interesting, are outside the scope of this work.

%--------------------------------------------------------------------
%--------------------------------------------------------------------
%\subsection{Computational Cost}\label{sec:disc_cost}
%--------------------------------------------------------------------
In terms of speed, when running reconstructions on a MacBook Pro with a 2.3 GHz Dual-Core Intel\textsuperscript{\textregistered} Core i5 processor, Calder\'on reconstructions on a $16\times16\times16$ $x$-grid which are interpolated to a $64\times64\times64$ $x$-grid take 1 to 2 seconds without optimizing for parallelization. This increases to 6-8 seconds per reconstruction when the initial $x$-grid is $32\times32\times32$. When running on a PC with a AMD EPYC 7702P 64-Core Processor 2.00 GHz, the reconstruction times are 0.6-0.7 seconds and 4 seconds, respectively, again without optimizing for parallelization.  On an 2015 iMac with a 4 GHz Quad-Core Intel\textsuperscript{\textregistered} Core i7 processor, the $\texp$ and $\tzero$ methods require 3-4 seconds/recon using a $21\by21\by21$ $x$-grid for the potential $q(x)$ or 6-8 seconds/recon when using a $41\by41\by41$ $x$-grid for $q(x)$.  The timings are non-optimized with the highest computational cost coming from computing the inverse Fourier transform and solving the boundary value problem using FEM.  The regularized TV, and linear difference, reconstructions averaged 2-3 hrs, and 3 minutes, respectively, when computed on a server with 256GB of RAM and two 10~core Intel\textsuperscript{\textregistered} Xeon\textsuperscript{\textregistered} CPU E5-2630 v4 @2.20GHz processors. \trev{We remark that the rather long computation times of the TV regularized non-linear least squares approach are caused by the 3D problem with a computationally rather challenging geometry as the sufficient accuracy of the CEM  forward model \eqref{eq:CEM} necessitates significant mesh refinement near the electrodes, leading to the large number of degrees of freedom (approx 250,000 nodes) in the FEM based forward model. The FEM model needs to be solved multiple times in the line search at each iteration of the Gauss-Newton method and with the mesh used each forward solution takes approximately 80s computation time.}

%--------------------------------------------------------------------
\section{Conclusions}\label{sec:conclusion}
%--------------------------------------------------------------------
In this work, we presented the first 3D absolute EIT reconstructions from CGO-based methods on experimental 3D tank data, and compared them to the current standard, a total variation regularized non-linear least squares approach.  We demonstrated that, with correct domain modeling, quality 3D absolute reconstructions can be obtained by all of the methods, comparable to the quality seen in linear difference imaging. All methods, Calder\'on, $\texp$, $\tzero$, and TV reasonably handled the moderate domain modeling error within little noticeable change in localization error and target contrast.  For the large modeling error case, the $\texp$ and $\tzero$ methods correctly identified the targets with high contrast, additional artefacts were introduced into the Calder\'on reconstruction, and the error proved too significant for the TV method.  The computational cost of the CGO reconstruction is trivial compared to TV (non-optimized, less than 1 sec/recon for Calder\'on, approximately 5~sec/recon for $\texp$ and $\tzero$, compared to 2-3 hours per reconstruction for TV).  

%--------------------------------------------------------------------
\section*{Acknowledgments}
Research reported in this paper was supported by the National Institute of Biomedical Imaging and Bioengineering of the National Institutes of Health under award numbers R21EB028064 (SH and JN) 1R01EB026710-01A1 (GS, DI, JN, ORS, and the development of the ACT5 device). The content is solely the responsibility of the authors and does not necessarily represent the official views of the National Institutes of Health.  JT and VK were supported by the Academy of Finland (Project 336791, Finnish Centre of Excellence in Inverse modeling and Imaging), the Jane and Aatos Erkko Foundation and Neurocenter Finland.  

%--------------------------------------------------------------------

%--------------------------------------------------------
% BIBLIOGRAPHY
%--------------------------------------------------------
%\bibliographystyle{amsalpha}
\bibliographystyle{plain}
\small
\bibliography{biblographyRefs.bib}
%--------------------------------------------------------

\end{document}